\documentclass{article}
\usepackage{amssymb}
\usepackage{amsmath}
\usepackage{enumitem}
\usepackage{bbm}
\usepackage{bookmark}
\usepackage{tikz-cd}
\usepackage{hyperref}

\usepackage[titletoc]{appendix}
\usepackage{amsthm}
\usepackage{pdflscape}
\usepackage{fancyhdr}
\usepackage{parskip}
\usepackage{adjustbox}
\usepackage{amsfonts}
\usepackage{caption}
\usepackage{mathrsfs}
\usepackage{afterpage}
\usepackage{float}
\usepackage{multirow}
\usepackage{graphicx}
\usepackage{color}
\usepackage[margin=1.0in]{geometry}
\newcommand*{\SHom}{\mathscr{H}\kern -.5pt om}
\newcommand*{\SExt}{\mathscr{E}\kern -.5pt xt}

\begin{document}
\title{The Getzler-Gauss-Manin connection and Kontsevich-Soibelman operations on the periodic cyclic homology}
\date{}
\author{Zihong Chen}
\maketitle
\theoremstyle{definition}
\newtheorem{mydef}{Definition}[section]
\numberwithin{mydef}{section}
\newtheorem{rmk}[mydef]{Remark}
\newtheorem{conj}[mydef]{Conjecture}
\theoremstyle{plain}
\newtheorem{question}[mydef]{Question}
\newtheorem{cor}[mydef]{Corollary}
\newtheorem{ex}[mydef]{Example}
\newtheorem{lemma}[mydef]{Lemma}
\newtheorem{thm}[mydef]{Theorem}
\newtheorem{prop}[mydef]{Proposition}
\begin{abstract}
    We study equivariant operations on the periodic cyclic homology of dg algebras that arise from the chain level action of the two-colored Kontsevich-Soibelman operad. The first main result is that these operations are covariantly constant with respect to the Getzler-Gauss-Manin connection on the periodic cyclic homology of a family of dg algebras. Then, using classical computations of Cohen \cite{Coh}, we explicitly compute a set of generators for these operations under composition, and show that these generators are closely related to the $p$-fold equivariant cap products previously studied by the author \cite{Che2} in relation to equivariant Gromov-Witten theory with mod $p$ coefficients. The main technical novelty is a re-formulation of the Kontsevich-Soibelman operad in terms of a two-colored version of the cacti operad, and a proof that it is \emph{equivariantly} quasi-equivalent to the two-colored operad of little disks on a disk/cylinder.  
\end{abstract}

\def\acts{\curvearrowright}
\tableofcontents

\numberwithin{equation}{section}
\renewcommand{\theHequation}{\thesection.\arabic{equation}}

\section{Introduction}
This paper is concerned with equivariant operations on the Hochschild homology and cohomology of a dg algebra, specifically those coming from the chain level action of the $E_2$-operad and its generalizations. Special attention will be paid to the case where the ground ring has positive characteristic. 
\subsection{Context}
In the early 1970s, the language of operads was introduced \cite{BV}\cite{May} as a package 
to describe algebraic operations with a certain topological `shape', such as the action of little $n$-disks operads on iterated loop spaces. \par\indent
The primary development relevant for us is Deligne's Hochschild cohomology conjecture, which postulates that the Hochschild cochains of an associative ring admit an action by the little $2$-disks operad (which from now on, will be called \emph{the little disks operad}), and various proofs have appeared in the works of Kontsevich-Soibelman \cite{KS2}, Tamarkin \cite{Tam}, McClure-Smith \cite{MS} and Voronov \cite{Vo2}. As a consequence, the equivariant homology of the $E_2$-operad, which was fully computed in \cite{Coh}, gives rise to an abundance of operations on Hochschild cohomology. While some of these operations, such as the cup product and Lie bracket (which form part of the Gerstenhaber structure \cite{Ger}), were classically known; the solution of Deligne's conjecture provided a new set of operations coming from the mod $p$ equivariant homology of the $E_2$-operad, called the \emph{Dyer-Lashof-Cohen operations}, which was explicitly computed by Tourtchine \cite{Tou}. \par\indent
The current paper aims to classify and compute equivariant operations on the \emph{periodic cyclic homology} of a dg algebra that parallels the Dyer-Lashof-Cohen operations in the following context. In \cite{KS1}, Kontsevich and Soibelman generalized Deligne's conjecture by constructing a two-colored operad $\{\mathbf{KS}(l,0),\mathbf{KS}(k,1)\}$ (where the suboperad $\{\mathbf{KS}(l,0)\}$ is quasi-equivalent to the little disks operad) that acts on the pair $(CC^*(\mathcal{A}),CC_*(\mathcal{A}))$ consisting of the Hochschild cochains and chains of a dg (or $A_{\infty}$) algebra $\mathcal{A}$; a quasi-equivalent model for this two-colored operad is sometimes called the \emph{operad of homotopy calculus}, cf. \cite{DTT1}. In particular, there is an action map
\begin{equation}\label{eq:the action map of KS}
\mathbf{KS}(k,1)\longrightarrow \mathrm{Hom}(CC^*(\mathcal{A})^{\otimes k}\otimes CC_*(\mathcal{A}),CC_*(\mathcal{A}))
\end{equation}
which is furthermore equivariant \eqref{eq:ksn1 action again} with respect to certain $S^1$ and $\Sigma_k$ actions on both sides (this will be explained in detail in Section 4.2). Taking (Tate) equivariant homology with respect to these actions, one obtains a map\footnote{$\mathbf{KS}(k,1)^{t S^1}_{h\Sigma_k}$ is our abuse of notation for $(\mathbf{KS}(k,1)_{h\Sigma_k})^{tS^1}$.}
\begin{equation}\label{eq:equiv KS action in intro}
H^*(\mathbf{KS}(k,1)^{t S^1}_{h\Sigma_k})\rightarrow \mathrm{Hom}(H^*_{\Sigma_k}(CC^*(\mathcal{A})^{\otimes k})\otimes HH_*^{tS^1}(\mathcal{A}), HH^{tS^1}_*(\mathcal{A})),
\end{equation}
where $HH^{tS^1}_*(\mathcal{A})=HH^{per}_*(\mathcal{A})$ is the periodic cyclic homology of $\mathcal{A}$. Therefore, any choice of equivariant homology classes $[\alpha]\in H^*(\mathbf{KS}(k,1)^{t S^1}_{h\Sigma_k})$ and $[\phi]\in H^*_{\Sigma_k}(CC^*(\mathcal{A})^{\otimes k})$ gives rise to an endomorphism of $HH_*^{per}(\mathcal{A})$. These endomorphisms (across all $k\geq 0$) generate a subalgebra of $\mathrm{End}_{H^*_{S^1}(pt)}(HH_*^{per}(\mathcal{A}))$ which we will call the \emph{Kontsevich-Soibelman operations}. \par\indent
A natural question is: 
\begin{center}
Can we describe all the Kontsevich-Soibelman operations for a dg algebra $\mathcal{A}$?     
\end{center}
One of the goals of this paper is to answer this question in the affirmative. On the other hand, the exact description turns out to very much depend on the coefficient ring over which $\mathcal{A}$ is defined, and is significantly more interesting when the ground ring has characteristic $p>0$. In the latter case, our main result Theorem \ref{thm:main theorem} identifies a simple set of generators for these operations, analogous to the classical case of Dyer-Lashof-Cohen operations. As we will see, the study of Kontsevich-Soibelman operations on $HH^{per}_*(\mathcal{A})\simeq CC_*(\mathcal{A})^{tS^1}$ can be illuminated by first replacing $(-)^{tS^1}$ by $(-)^{tC_p}$, using the well known relation (cf. Appendix A.3)
\begin{equation}\label{eq:HH^per compare with HH^{C_p,per}}
HH^{C_p,per}_*(\mathcal{A}) \cong HH^{per}_*(\mathcal{A}) \oplus  HH^{per}_*(\mathcal{A}) \theta,
\end{equation}
where $\theta$ is a formal variable of degree $1$. More precisely, we describe a set of generators for $H^*(\mathbf{KS}(k,1)^{t S^1}_{h\Sigma_k})$ under \emph{operadic composition} and explicitly compute their associated Kontsevich-Soibelman operations, which are most easily expressed using relation \eqref{eq:HH^per compare with HH^{C_p,per}} and explicit chain models for $HH^{C_p,per}_*(\mathcal{A})$. \par\indent
Along the way, we also study the interaction of these endomorphisms with other well-known structures on the periodic cyclic homology, such as the Getzler-Gauss-Manin connection \cite{Get}.
\subsection{Main results}
Let $\mathcal{A}$ be a dg algebra (or d$(\mathbb{Z}/2)$g algebra) over a ground ring $R$. In \cite{Che2}, the author introduced an explicit chain complex computing the Tate fixed points of the sub $C_p\subset S^1$-action on $CC_*(\mathcal{A})$, which is called the \emph{periodic $C_p$-equivariant Hochschild complex} and denoted $CC^{C_p,per}_*(\mathcal{A})$; cf. Section 2.3 for the definition. Moreover, motivated by the study of equivariant genus $0$ Gromov-Witten theory, the author defined an operation 
\begin{equation}\label{eq:prod^C_p}
\prod\nolimits^{C_p}: H^*_{C_p}(CC^*(\mathcal{A})^{\otimes p})\otimes HH_*^{C_p,per}(\mathcal{A})\longrightarrow HH_*^{C_p,per}(\mathcal{A}),
\end{equation}
called the \emph{$p$-fold equivariant cap product} (cf. Definition \autoref{thm:Cp equivariant cap product}), which is a $C_p$-equivariant analogue of the usual cap product action $\cap: HH^*(\mathcal{A})\otimes HH_*(\mathcal{A})\rightarrow HH_*(\mathcal{A})$. Incidentally, $\prod^{C_p}$ will play a crucial role in the classification of Kontsevich-Soibelman operations. Let $\prod^{C_p}_1$ (resp. $\prod^{C_p}_2$) denote the result of pre-composing $\prod^{C_p}$ with the restriction $HH^{per}_*(\mathcal{A})\rightarrow HH^{C_p,per}_*(\mathcal{A})$ and post-composing with the map in the left (resp. right) direction in \eqref{eq:projections to HH^{per} factors of HH^{C_p,per}} (the projections $pr_1,pr_2$ are taken with respect to the decomposition \eqref{eq:HH^per compare with HH^{C_p,per}})
\begin{equation}\label{eq:projections to HH^{per} factors of HH^{C_p,per}}
 HH^{per}_*(\mathcal{A})   \xleftarrow{pr_1} HH^{C_p,per}_*(\mathcal{A})\xrightarrow{pr_2} HH^{per}_*(\mathcal{A})\theta\xrightarrow{\textrm{divide by}\;\theta} HH^{per}_*(\mathcal{A}). 
\end{equation}

Our main result is the following theorem:\\
\begin{thm}\label{thm:main theorem}
\begin{enumerate}[label=\roman*)]
\item When $\mathbb{Q}\subset R$, there are no nontrivial Kontsevich-Soibelman operations on $HH^{per}_*(\mathcal{A})$. 
    \item When $\mathbb{F}_p\subset R$, let $[\phi]\in H^*_{\Sigma_p}(CC^*(\mathcal{A})^{\otimes p})$ and let $[\phi]_{C_p}$ denote its restriction to $H^*_{C_p}(CC^*(\mathcal{A})^{\otimes p})$. Then \begin{equation}
        \prod\nolimits_1^{C_p}([\phi]_{C_p},-)\quad \textrm{and} \quad\prod\nolimits_2^{C_p}([\phi]_{C_p},-)
    \end{equation} are Kontsevich-Soibelman operations.
    \item When $\mathbb{F}_p\subset R$, and assuming that the natural map
\begin{equation}
(CC^*(\mathcal{A})^{\otimes k_1})^{h\Sigma_{k_1}}\otimes\cdots\otimes (CC^*(\mathcal{A})^{\otimes k_l})^{h\Sigma_{k_l}}\rightarrow (CC^*(\mathcal{A})^{\otimes k_1+\cdots+k_l})^{h(\Sigma_{k_1}\times \cdots\times\Sigma_{k_l})}     
\end{equation}
is a quasi-isomorphism for all $k_1,\cdots,k_l,l\geq 1$ and that the left hand side above satisfies the K\"{u}nneth isomorphism, then the algebra of Kontsevich-Soibelman operations on $HH^{per}_*(\mathcal{A})$ is generated by endomorphisms of the form $\prod^{C_p}_1([\phi]_{C_p},-)$ and $\prod^{C_p}_2([\phi]_{C_p},-)$.\\
\end{enumerate}
\end{thm}
It is a standard fact that when $X$ (resp. $Y$) is a bounded chain complex with finite group action $G$ (resp. $H$), the natural map $X^{hG}\otimes Y^{hH}\rightarrow (X\otimes Y)^{h(G\times H)}$ is a quasi-isomorphism. In particular,\\
\begin{lemma}\label{thm:scenario where main thm apply}
The assumptions in Theorem \ref{thm:main theorem} iii) are satisfied for a $\mathbb{Z}$-graded dg algebra $\mathcal{A}$ over a field $\mathbf{k}\supset \mathbb{F}_p$ such that $HH^*(\mathcal{A})$ lies in bounded degrees (for instance, when $\mathcal{A}$ is smooth proper and $\mathbb{Z}$-graded).   \qed
\end{lemma}
The proof of Theorem \autoref{thm:main theorem} will take up most of the paper and be divided into two parts.
\begin{itemize}
    \item Identifying the two types of abstract generators for Kontsevich-Soibelman operations under compositions. This is obtained by computing the cohomology of the $S^1$-Tate fixed points of a convenient topological model (cf. Section 4.3) for the Kontsevich-Soibelman operad and finding suitable generators under operadic compositions, cf. Theorem \ref{thm:classification of untwisted KS operations}. In fact, from the topological point of view it is more natural to start from certain $C_p$-equivariant classes, cf. Theorem \autoref{thm:classification of untwisted KSp operations}, and then transfer this information via comparisons of $S^1$ and $C_p$-fixed points, cf. Corollary \ref{thm:restriction of epsilon classes} and Corollary \ref{thm:formula for action of epsilon classes}.
    \item Showing that the prior identified abstract generators indeed agree with (factors of) the explicitly defined $p$-fold equivariant cap product on $HH_*^{C_p,per}(\mathcal{A})$ using the appropriate chain model. This is the more technical part of the paper, and is done by introducing a good combinatorial model (called \emph{cyclic cacti with spines}, cf. Section 4.2) for the Kontsevich-Soibelman operad that is simultaneously fine enough to encode the various operadic compositions and $S^1$ or $C_p$-actions at once and discrete enough to see (variants of) the `cyclic bar complex' model for Hochschild homology; these requirements will be explained in more detail in Section 1.3. The proof that this combinatorial model agrees with the topological one, and the proof that certain abstract Kontsevich-Soibelman operations agree with $C_p$-equivariant cap products, will be completed in Proposition \ref{thm:e0 operation is p fold equivariant cap product}.  
\end{itemize}
On the other hand, the assumptions of Theorem \ref{thm:main theorem} iii) are rarely satisfied for a d$(\mathbb{Z}/2)$g algebra, which may be a case of interest (e.g. in mirror symmetry). In this case, one can still obtain a weaker decomposition by restricting to a smaller class of operations: fix a cocycle $\varphi\in CC^{even}(\mathcal{A})$, a \emph{$\varphi$-Kontsevich-Soibelman operation} is an endomorphism of $HH^{per}_*(\mathcal{A})$ obtained by plugging $[\varphi^{\otimes k}]\in H^*_{\Sigma_k}(CC^*(\mathcal{A})^{\otimes k})$ and some class $[\alpha]\in H^*(\mathbf{KS}(k,1)^{t S^1}_{h\Sigma_k})$ into \eqref{eq:equiv KS action in intro}.\\
\begin{thm}\label{thm:main theorem without homological bounds assumption}
Let $\varphi\in CC^{even}(\mathcal{A})$ be a cocycle. Then any $\varphi$-Kontsevich-Soibelman operation on $\mathcal{A}$ can be written as a linear combination of endomorphisms of the form
\begin{equation}
  \bigcap\nolimits^{C_p}_2([\varphi],-)\circ  \cdots \circ\bigcap\nolimits^{C_p}_2([\varphi],-)\circ \prod\nolimits^{C_p}_1([\phi]_{C_p},-)
\end{equation}
for some class $[\phi]\in H^*_{\Sigma_p}(CC^*(\mathcal{A})^{\otimes p})$\footnote{$[\phi]$ may be explicitly constructed out of $[\varphi]$ and the equivariant homology class $[\alpha]$ giving rise to the $\mathrm{KS}_{\varphi}$-operation in question, cf. Theorem \ref{thm:classification of untwisted KSp operations}. However, $\prod\nolimits^{C_p}_1([\phi]_{C_p},-)$ itself may not be a $\varphi$-Kontsevich-Soibelman operation.}, where $\bigcap^{C_p}$ denotes the \emph{$C_p$-equivariant cap product} (cf. Definition \ref{thm:Cp equivariant cap product}) and $(-)_1$ (resp. $(-)_2$) are as in Theorem \ref{thm:main theorem}. 
\end{thm}

Another result we prove concerns Kontsevich-Soibelman operations in the relative setting, i.e. when $\mathcal{A}$ is a dg algebra over $R$, together with an extension $\mathbf{k}\subset R$ where $\mathbf{k}$ is a field. Intuitively, we think of $HH^{per}_*(\mathcal{A}/R)$ as a vector bundle over $\mathrm{Spec}\,R$, whose `total space' is $HH^{per}_*(\mathcal{A}/\mathbf{k})$. \par\indent
In this setting, Getzler \cite{Get} defined a non-commutative analogue of the Gauss-Manin connection on the cohomology of a family of complex varieties, the \emph{Getzler-Gauss-Manin connection}
\begin{equation}
\nabla^{GGM}: HH^{per}_*(\mathcal{A}/R)\rightarrow HH^{per}_*(\mathcal{A}/R)\otimes \Omega^1_{R/\mathbf{k}}.      
\end{equation}
This forms part of a (partially conjectural) variation of semi-infinite Hodge structure on the periodic cyclic homology \cite{KKP} and plays an important role in mirror symmetry \cite{GPS}. In this setting, we say that \emph{a Kontsevich-Soibelman operation is $\mathbf{k}$-linear} if it is induced from a Kontsevich-Soibelman class in the image of  $H^*((\mathbf{KS}(n,1)_{\mathbf{k}})^{tS^1}_{h\Sigma_n})\rightarrow H^*((\mathbf{KS}(n,1)_{R})^{tS^1}_{h\Sigma_n})$ and an equivariant Hochschild cohomology class in the image of $H^*_{\Sigma_n}(CC^*(\mathcal{A}/R)^{\otimes_\mathbf{k} n})\rightarrow H^*_{\Sigma_n}(CC^*(\mathcal{A}/R)^{\otimes_R n})$. \\
\begin{thm}(Automatic covariant constancy) \label{thm:covariant constancy of KS operation}
Any $\mathbf{k}$-linear Kontsevich-Soibelman operation commutes with $\nabla^{GGM}$. 
\end{thm}
Perhaps surprisingly, the proof of Theorem \autoref{thm:covariant constancy of KS operation} (cf. Theorem \autoref{thm:automatic covariant constancy}) does not use any explicit relations in $\mathbf{KS}$, nor does it use the classification of Theorem \autoref{thm:main theorem}. Instead, it follows from the reformulation \cite{Kal2}\cite{PVV} of the Getzler-Gauss-Manin connection as a Grothendieck connection, and the observation that any Kontsevich-Soibelman operation (over $\mathbf{k}$) automatically `extends to the first order neighborhood'. Hence we call Theorem \autoref{thm:covariant constancy of KS operation} `automatic covariant constancy'. 
\subsection{Methods of proof}
\textbf{Equivariant homology of $\mathbf{KS}$}. As was explained in the previous subsection, the first step towards proving Theorem \autoref{thm:main theorem} is a computation of the (Tate) equivariant homology of $\mathbf{KS}(k,1)$. The idea is to reduce this to studying the equivariant homology of certain configuration spaces, at which point it follows from an explicit computation in algebraic topology aided by the classical work of Cohen \cite{Coh} on the homology of $E_n$-spaces; the latter computation is done in Section 5. \par\indent
On the other hand, the key to the reduction step is a set of \emph{equivariant} weak equivalences
\begin{equation}\label{eq:weak equivalence between KS and Cyl}
\mathbf{KS}(l,0)\simeq \mathrm{Cyl}(l,0)\quad,\quad\mathbf{KS}(k,1)\simeq \mathrm{Cyl}(k,1),    
\end{equation}
where $\mathrm{Cyl}(l,0)$ is the configuration space of $l$ disks on a disk, and $\mathrm{Cyl}(k,1)$ is the configuration space of $k$ disks on a cylinder, together with one marked point moving along each boundary component (modulo an overall rotation). Moreover, \eqref{eq:weak equivalence between KS and Cyl} should be a weak equivalence of \emph{operads}. \par\indent
\textbf{Cactus and its cyclic version}. The weak equivalence of operads in \eqref{eq:weak equivalence between KS and Cyl} is well-known and proved in \cite[Section 11]{KS1}, \cite[Section 7]{KS2} and \cite[Section 6]{Wil} using an explicit model for $\mathbf{KS}$ based on the `minimal operads' of Kontsevich-Soibelman \cite[Section 5]{KS2}, which are cellular operads whose cells correspond to certain trees or forests on a cylinder that have tautological actions on the Hochschild (co)chains of an $A_{\infty}$-category. \par\indent
However, the aforementioned proofs are not sufficient for our purpose of studying equivariant homology operations as they do not keep track of the appropriate equivariance conditions. To elaborate, what we need is a weak equivalence of operads as in \eqref{eq:weak equivalence between KS and Cyl} such that: 
\begin{itemize}
    \item The weak equivalence $\mathbf{KS}(k,1)\simeq \mathrm{Cyl}(k,1)$ is $\Sigma_k\times S^1\times S^1$-equivariant, where on the right hand side, $\Sigma_k$ permutes the labels of the disks and $S^1\times S^1$ rotates the two marked points on the boundaries of the cylinder. 
    \item The action map \eqref{eq:the action map of KS} is also $\Sigma_k\times S^1\times S^1$-equivariant, where on the right hand side, $\Sigma_k$ permutes the factors of $CC^*(\mathcal{A})^{\otimes k}$ and the two circles act naturally on the two copies of $CC_*(\mathcal{A})$ (e.g. via the Connes operator).
    \item Moreover, to obtain the explicit formulae of Theorem \autoref{thm:main theorem} in terms of the $p$-fold equivariant cap product, the proclaimed $S^1\times S^1$-action on $\mathbf{KS}(k,1)$ must be fine enough to `see' the subgroup $C_p\subset S^1$. For instance, the algebra of cellular chains of $S^1$ with one $0$-cell and one $1$-cell would not serve as a good model for this circle action.
\end{itemize}
Based on these requirements, we turn to a different, but quasi-equivalent formulation of the two-colored operad $\mathbf{KS}$ using the operad of `cactus' instead of the `minimal operad'. \emph{Cactus} was first introduced combinatorially by \cite{MS} in their solution of Deligne's conjecture and studied extensively in \cite{Vo1} \cite{Kau} \cite{Sal1}. Roughly speaking, a cactus is a connected union of embedded circles (called its \emph{lobes}) in $\mathbb{C}$ which configure like, well, \emph{a cactus} (see Figure \ref{fig:cactus}). It has two important features relevant for us. First, the space of cacti is a realization of a multi-semisimplicial set (which we extend to a cosimplicial multi-simplicial set, cf. Definition \autoref{thm:Cacti with spines as a cosimplicial multisimplicial set}); second, the space of cacti is a model for the configuration space on $\mathbb{C}$. The proof of the latter was obtained in \cite{Sal1}, and the idea is to associate to each configuration $(z_1,\cdots,z_k)\in \mathrm{Conf}_k(\mathbb{C})$ and positive weights $a_1,\cdots,a_k>0, \sum_i a_i=1$, a vector field
\begin{equation}
E(z)=\sum_{i=1}^ka_i\frac{(z-z_i)}{|z-z_i|^2}.
\end{equation}
Its critical graph (i.e. the union of flow lines connecting zeros and poles of $E$) gives the dual graph of the desired cactus; see Section 6.1 for a detailed review.\par\indent
Our main construction is a generalization of cactus to the cyclic setting, done in a simplicial framework. To be more precise, we construct two types of functors
\begin{equation}\label{eq:cactus and cyclic cactus with spines}
\widehat{\mathfrak{Cact}^l}:\Delta\times (\Delta^{op})^l\rightarrow\mathrm{Sets} \quad ,\quad \widehat{\mathfrak{Cact}^k_{\circlearrowright}}:\Lambda^{op}\times (\Delta^{op})^k\times \Lambda\rightarrow\mathrm{Sets},
\end{equation}
called \emph{cacti with spines} and \emph{cyclic cacti with spines}, where $\Delta$ denotes the simplex category and $\Lambda$ denotes Connes' cyclic category \cite{Con}, which is a small combinatorial category that models circle action out of finite cyclic group actions. We define appropriate structure maps making \eqref{eq:cactus and cyclic cactus with spines} into a two-colored operad in an appropriate simplicial sense, and show that it acts \emph{tautologically} on the pair $(CC^*(\mathcal{A}),CC_*(\mathcal{A}))$, each viewed with its (co)simplicial or (co)cyclic structure. \par\indent
On the other hand, it is well-known \cite{Jon} that given a functor $\Lambda\rightarrow \mathrm{Sets}$ (resp. $\Lambda^{op}\rightarrow \mathrm{Sets}$), its realization (resp. totalization) is equipped with a canonical $S^1$-action. Moreover, one can recover the sub $C_p\subset S^1$ action combinatorially via a map of small categories ${}_p\Lambda\rightarrow \Lambda$, where ${}_p\Lambda$ is the \emph{finite $p$-cyclic category}, cf. Section 2.4 b). This observation is crucial for obtaining the explicit formulae in Theorem \autoref{thm:main theorem}.\par\indent
Finally, we mimic the method of \cite{Sal1} and show that the totalization/realization of \eqref{eq:cactus and cyclic cactus with spines} gives a two-colored topological operad that is quasi-equivalent to $\mathrm{Cyl}$, where the $S^1\times S^1$-equivariance is a built-in feature of the construction. In fact, we go through an intermediate step and first show it is equivariantly quasi-equivalent to the two-colored Fulton-MacPherson operad $\{FM(l),FM_{\circlearrowright}(k)\}$ \eqref{eq:FM_circ_space}, which is further quasi-equivalent to $\mathrm{Cyl}$. 
\subsection{Relation to other directions}
We explain two directions which have largely motivated this project, even though neither will occur substantial in the body of this paper. \par\indent
\textbf{Equivariant Gromov-Witten theory with mod $p$ coefficients}. Initially, this paper grew out of an attempt to study equivariant operations in symplectic Gromov-Witten theory via Fukaya categories, cf. \cite{Sei2}\cite{Che2} for some relevant results. For technical simplicity, we focus on the case of a closed monotone (i.e. $c_1(TX)=[\omega]$) symplectic manifold $(X,\omega)$.  \par\indent
Let $QH^*(X;R)$ denote the quantum cohomology of $X$ with coefficient ring $R$. When $R=\mathbf{k}$ is a field of characteristic $p$, Fukaya \cite{Fuk} and Wilkins \cite{Wilk} introduced a Frobenius $p$-linear multiplicative action called the \emph{quantum Steenrod operations}
\begin{equation}
Q\Sigma: QH^*(X;\mathbf{k})\rightarrow \mathrm{End}_{\mathbf{k}[[t,\theta]]}(H^*(X;\mathbf{k})[[t,\theta]]), \quad \mathbf{k}[[t,\theta]]=H^*_{C_p}(pt;\mathbf{k}) 
\end{equation}
built out of $C_p$-equivariant genus $0$ Gromov-Witten invariants with mod $p$ coefficients. The quantum Steenrod operations are a quantization of the classical Steenrod operations: the `classical part' (i.e. contributions from constant $J$-holomorphic curves) of $Q\Sigma_b$ is cup product with the total Steenrod power\footnote{Conventionally Steenrod powers are defined for cohomology classes with $\mathbb{F}_p$-coefficients, and here we simply extend them over $\mathbf{k}$ by Frobenius $p$-linearity.}
(cf. \cite[(1d)]{Sei2} for the sign convention)
\begin{equation}\label{eq:classical Steenrod powers}
St(b):=\begin{cases}
 \sum Sq^i(b)t^{\frac{|b|-i}{2}}\;,\quad\textrm{when $p=2$}  \\
 \pm(\frac{p-1}{2}!)^{|b|}\cdot \sum (-1)^i P^i(b)t^{\frac{(p-1)(|b|-2i)}{2}}\pm  \sum (-1)^i\beta P^i(b) t^{\frac{(p-1)(|b|-2i)-2}{2}}\theta\;,   \quad  \textrm{when $p>2$}
\end{cases}   
\end{equation}
where $Sq^i$ is the $i$-th Steenrod square, $P^i$ is the $i$-th Steenrod power (for $p$ odd) and $\beta$ is the Bockstein (for $p=2$, $\beta$ agrees with $Sq^1$). Recall that $t$ and $\theta$ are the degree $2$ and $1$ generators of $H^*_{C_p}(pt;\mathbf{k})$; in particular, $\theta^2=0$ if $p$ is odd and $\theta^2=t$ (or $\theta=t^{1/2}$) if $p=2$.\par\indent
Despite their relatively recent appearance in symplectic topology, quantum Steenrod operations have already seen an abundance of links and applications to Hamiltonian dynamics \cite{Shel}\cite{CGG}; arithmetic mirror symmetry \cite{Sei2}; representation theory \cite{Lee}\cite{BL1}; and the study of quantum connection \cite{SW}\cite{Lee}\cite{BL2}\cite{Che3}. \par\indent
In \cite{Che2}, the author attempted an explanation of the theoretical importance of quantum Steenrod operations by interpreting them as equivariant operations on the Hochschild homology of the Fukaya category of $X$:\\
\begin{thm}\label{thm:rough statement that QSigma is Cp equivariant cap}(\cite{Che2})
Under suitable assumptions on the Fukaya category $\mathrm{Fuk}(X;\mathbf{k})$, the quantum Steenrod operations $Q\Sigma: QH^*(X;\mathbf{k})\acts H^*(X;\mathbf{k})[[t,\theta]]$ can be recovered from 
$\prod^{C_p}$ \eqref{eq:prod^C_p} applied to $\mathcal{A}=\mathrm{Fuk}(X;\mathbf{k})$. 
\end{thm}
Perhaps on a philosophical level, the current paper further affirms this point by proving that (cf. Theorem \autoref{thm:main theorem}) factors of $\prod^{C_p}$ generate all the `natural' endomorphisms of $HH^{per}_*(\mathcal{A})$ from the operadic point of view.  \\
\begin{rmk}
Given Theorem \autoref{thm:main theorem} and Theorem \autoref{thm:rough statement that QSigma is Cp equivariant cap}, it is natural to wonder if there is a direct proof (provided that it is true at all) that the $Q\Sigma_b$'s generate all endomorphisms of $H^*(X;\mathbf{k})((t,\theta))$ coming from the mod $p$ equivariant homology of the genus $0$ modular operad $\overline{\mathcal{M}}_{0,-}$ on a monotone symplectic manifold $X$. To the author's knowledge there are at least two challenges.
\begin{itemize}
    \item The technical/analytic setup needed to assign Floer operations to all elements in the equivariant homology of $\overline{\mathcal{M}}_{0,n}$.
    \item We (the author) do not know of a computation of the equivariant homology of $\overline{\mathcal{M}}_{0,n}$ (for all $n$) or an analogue of Lemma \autoref{thm:generation of equiv cohomology of Cyl} in this context. 
\end{itemize}
\end{rmk}
\textbf{$p$-curvature in non-commutative Hodge theory}. $p$-curvature is a central invariant of a differential equation defined over characteristic $p$. Its definition rests on the following simple observation: if $X$ be a (formal) scheme over a field $\mathbf{k}$ of characteristic $p>0$, and $D: \mathcal{O}_X\rightarrow \mathcal{O}_X$ be a (continuous) $\mathbf{k}$-linear derivation on $X$, then its $p$-th iteration 
\begin{equation}
D^p:=\overbrace{D\circ D\circ\cdots\circ D}^{p\;times}: \mathcal{O}_X\rightarrow \mathcal{O}_X  
\end{equation}
remains a derivation (an easy application of the Leibniz rule), sometimes called the \emph{$p$-th power of $D$}. \par\indent
In the above setting, further fix $(M,\nabla)$ an $\mathcal{O}_X$-module with connection. The \emph{$p$-curvature} of $\nabla$ along $D$ is defined to be
\begin{equation}
F^{\nabla}_D:=\nabla^p_{D}-\nabla_{D^p}: M\rightarrow M.    
\end{equation}
We refer the readers to \cite[Section 5]{Ka} for its key properties and a more in-depth introduction.  \par\indent
One of the most important features of $p$-curvature is that it gives a criterion (known as \emph{Cartier descent}, cf. \cite[Theorem 5.1]{Ka}) for when a differential equation in characteristic $p$ has a full set of algebraic solutions. Namely, this is the case if and only if its $p$-curvature vanishes. This starkly contrasts the case of characteristic $0$, where no such simple criterion for the existence of algebraic solutions to a general differential equation (with singularity) is known. A well-known (and still wildly open) conjecture of Grothendieck and Katz predicts that the existence of algebraic solutions to a differential equation in characteristic $0$ is encoded in the $p$-curvature of its mod $p$ reduction across all primes $p$. \par\indent
The differential equation of primary concern in non-commutative Hodge theory is the Getzler-Gauss-Manin connection (cf. Section 3 for more details) on the relative periodic cyclic homology of a dg algebra $\mathcal{A}$ over a base $R$. When $\mathcal{A}$ is smooth proper and $\mathbb{Z}$-graded over a base $R$ of characteristic $p$, the $p$-curvature of its Getzler-Gauss-Manin connection was extensively studied in \cite{PVV}, which led to interesting results such as a non-commutative analogue of the local monodromy theorem. \par\indent
The relevance of this to the current paper is the structural similarities between the $p$-curvature of $\nabla^{GGM}$ and the $p$-fold equivariant cap product. 
\begin{enumerate}
    \item Their `non-equivariant parts' agree. In Section 3.2 and 3.3, we transfer the Getzler-Gauss-Manin connection in characteristic $p$ to $HH^{C_p,per}_*(\mathcal{A})$ via \eqref{eq:HH^per compare with HH^{C_p,per}} and then explicitly write down a chain level formula \eqref{eq:p GGM}. It follows from the expression \eqref{eq:p GGM} that 
\begin{equation}
(t\nabla^{GGM,p})^p|_{(t,\theta)=0}={}_p\prod(\kappa^{\otimes p},-)=\prod\nolimits^{C_p}(\kappa^{\otimes p},-)|_{(t,\theta)=0},
\end{equation}
where $\kappa$ denotes the Kodaira-Spencer class of the family $\mathcal{A}/R$ and ${}_p\prod$ is the (non-equivariant) $p$-fold cap product \eqref{eq:p fold cap product}. 
\item Both $\prod^{C_p}$ and the $F^{GGM}$ are covariantly constant with respect to $\nabla^{GGM}$: the former because of $\prod^{C_p}$ being a Kontsevich-Soibelman operation and Theorem \ref{thm:automatic covariant constancy}, and the latter because $\nabla^{GGM}$ is flat. 
\end{enumerate}
Given these, it is natural to ask further:\\
\begin{question}\label{thm:p curv and equivariant cap product}
Does $\prod^{C_p}(\kappa^{\otimes p},-)$ agree with the $p$-curvature of $\nabla^{GGM}$ in general? 
\end{question}
A positive answer to this question is significant at various levels. First, it gives the $p$-curvature of Getzler-Gauss-Manin connection, which is a fundamental invariant of a family of dg algebras, a clean topological interpretation in terms of the Kontsevich-Soibelman operad. More importantly, such topological interpretation sheds light on properties of the $p$-curvature that are otherwise elusive from explicit formulae, such as its nilpotence\footnote{Roughly speaking, the idea is: while the Getzler-Gauss-Manin connection is most naturally viewed as part of an $L_{\infty}$-action of Hochschild cochains on the negative cyclic chains \cite[Section 11.4]{KS1} \cite[Section 3]{DTT2}, its \emph{multiplicative} properties are only seen after passing to the bigger $\mathbf{KS}$-operad (aka the homotopy calculus operad), where both $A_{\infty}$ and $L_{\infty}$-module structures live.}. In a followup paper \cite{Che4}, we will affirm Question \ref{thm:p curv and equivariant cap product} for certain families of differential $\mathbb{Z}/2$-graded algebras and study its consequences for classical questions in non-commutative Hodge theory.\\
\begin{rmk}
\begin{itemize}
\item In the context of symplectic topology (roughly speaking, setting $\mathcal{A}$ to be the Fukaya category of a closed symplectic manifold), an analogue of Question \ref{thm:p curv and equivariant cap product} was asked by Lee \cite{Lee} and affirmed in several cases \cite{Lee}\cite{BLP}\cite{PS}.
    \item When $\mathcal{A}$ is smooth proper and moreover can be lifted to the sphere spectrum (e.g. when $\mathcal{A}$ has a $\mathbb{Z}$-grading and certain homological bounds on $\mathcal{A}$ and $HH^*(\mathcal{A})$ are satisfied, cf. \cite[Theorem 9.1]{Rez}), a recent work of Rezchikov \cite{Rez} gives a positive answer to Question \ref{thm:p curv and equivariant cap product}. An earlier work of \cite{PVV} also addresses this (in slightly different formulation) in the $\mathbb{Z}$-graded setting when similar homological bounds are imposed. Both works, explicitly or implicitly, use structures beyond just the Kontsevich-Soibelman operad (e.g. cyclotomic structure, and other ideas from stable homotopy theory). 
\end{itemize}
\end{rmk}

\subsection*{Organization}
The organization of this paper is as follows. In Section 2, we review the definition of Hochschild (co)homology of a dg algebra. In particular, we will study them from the perspective of (co)simplicial and (co)cyclic objects. In Section 3, we review the construction of the Getzler-Gauss-Manin connection, as well as Kaledin's reformulation of it. While this is unrelated to the other parts of the paper, we also give an explicit formula of a Getzler-Gauss-Manin connection on the $C_p$-equivariant Hochschild complex, which may be of interest. In Section 4, we study cacti and their variants and use them to define a version of the Kontsevich-Soibelman operad. We then define Kontsevich-Soibelman operations on the periodic cyclic homology and prove Theorem \autoref{thm:covariant constancy of KS operation}. In Section 5, we compute the equivariant homology of the Kontsevich-Soibelman operad building on Cohen's work on configuration spaces \cite{Coh}, and identify a set of generators under operadic composition. In Section 6, we prove the equivalence between (cyclic) cacti and certain configuration spaces, and derive explicit formulae for the generators of Kontsevich-Soibelman operations, thus concluding the proof of Theorem \autoref{thm:main theorem}.  

\subsection*{Acknowledgements}
We would like to thank Paul Seidel, Dmitry Kaledin, Oscar Randal-Williams, Tommaso Rossi, Paolo Salvatore, John Christian Ottem and Mark Grant for enlightening discussions or correspondences at various stages of this project. We especially thank Semon Rezchikov and Paul Seidel for pointing out an error in a prior version of this manuscript. The author is partially supported by a Title A Fellowship from Trinity College, Cambridge and a Postdoctoral Fellowship from the Herchel-Smith Fund.

\section{Hochschild invariants of a dg category}
Let $\mathcal{A}$ be a differential ($\mathbb{Z}$ or $\mathbb{Z}/2$-) graded algebra over a ground ring $R$, i.e. an $R$-linear chain complex $(\mathcal{A},d)$ with an $R$-bilinear associative product $\cdot$ such that
\begin{equation}
d(x\cdot y)=dx\cdot y+(-1)^{|x|}x\cdot dy. \label{eq:leibniz}    
\end{equation}
For simplicity, we will assume that $\mathcal{A}$ contains a strict unit $1$. \par\indent
The goal of this section is to review the definition of Hochschild (co)homology, as well as suitable equivariant versions, of $\mathcal{A}$ in terms of explicit chain models. 
\subsection{Hochschild cohomology}
As a graded $R$-module, the \emph{Hochschild cochain complex} of $\mathcal{A}$ is
\begin{equation}
CC^*(\mathcal{A}):=\prod_{k\geq 0}\mathrm{Hom}_{R}^*(\mathcal{A}^{\otimes_R k}[1],\mathcal{A}).\label{eq:CC^*}
\end{equation}
There is a chain level binary operation on \eqref{eq:CC^*}, called the \emph{circle product}, given by
\begin{equation}
(\varphi\circ\phi)(x_1,\cdots,x_n):=\sum_{j\leq k}(-1)^{|\phi|\cdot \sum_{i=k}^n\|x_i\|} \varphi(x_1,\cdots,x_{j-1},\phi(x_j,\cdots,x_{k-1}),x_k,\cdots,x_n),  
\end{equation}
where $\|x\|:=|x|-1$ denotes the reduced degree. This gives rise to a degree $-1$ chain level Lie bracket by
\begin{equation}
[\varphi,\phi]:=\varphi\circ \phi-(-1)^{\|\varphi\|\|\phi\|}\phi\circ\varphi.
\end{equation}
Define $\mu_{\mathcal{A}}\in CC^{2}(\mathcal{A})$ to be the Hochschild cochain given by
\begin{equation}
\mu_{\mathcal{A}}^0=0,\;\mu^1_{\mathcal{A}}(x)=dx,\;\mu^2_{\mathcal{A}}(x,y)=xy,\;\mu^k_{\mathcal{A}}=0\;\;\textrm{for}\;\;k\geq 3.    
\end{equation}
Then associativity and the Leibniz rule \eqref{eq:leibniz} implies that $\mu_{\mathcal{A}}\circ \mu_{\mathcal{A}}=0$, and in particular, $[\mu_{\mathcal{A}},\mu_{\mathcal{A}}]=0$. One then defines the differential on $CC^*(\mathcal{A})$ to be
\begin{equation}
[\mu_{\mathcal{A}},-]: CC^*(\mathcal{A})\rightarrow CC^{*+1}(\mathcal{A}),
\end{equation}
which squares to zero as a consequence of the Jacobi identity. We denote the cohomology of this complex by $HH^*(\mathcal{A})$. \par\indent
There is a product structure on the Hochschild cochain complex $CC^*(\mathcal{A})$ called the \emph{cup product} (or \emph{Yoneda product}). On the chain level, it is given by
\begin{equation}
\phi\cup\psi(x_1,\cdots,x_k):=\sum_{i=1}^{k+1}(-1)^{|\phi|\cdot\sum_{j=i}^k\|x_j\|}\phi(x_1,\cdots,x_{i-1})\psi(x_{i},\cdots,x_k).
\end{equation}
The cup product descends to cohomology and defines an associative (and in fact graded commutative) algebra structure on $HH^*(\mathcal{A})$, for which $1\in\mathcal{A}$ (viewed as a Hochschild cocycle of length $0$) is the unit.

\subsection{Hochschild homology and cyclic homology}
The \emph{Hochschild chain complex} of $\mathcal{A}$ 
\begin{equation}
CC_*(\mathcal{A}):=\bigoplus_{n\geq 0}\mathcal{A}\otimes_R\mathcal{A}[1]^{\otimes_R n},
\end{equation}
with differential given by
$$ b(x_0\otimes x_1\otimes\cdots\otimes x_n)=\sum_{i=0}^n (-1)^{\sum_{j=i+1}^n\|x_j\|} x_0\otimes x_1\otimes\cdots \otimes dx_i\otimes\cdots\otimes x_n$$
\begin{equation}
+\sum_{i=0}^{n-1} (-1)^{\sum_{j=i+2}^n\|x_j\|}\mathbf{x}\otimes x_1\otimes \cdots\otimes x_{i}x_{i+1}\otimes\cdots\otimes x_n+(-1)^{\sum_{j=1}^n\|x_j\|+\|x_n\|\cdot(|x_0|+\sum_{i=1}^n\|x_i\|)} x_n x_0\otimes x_1\otimes\cdots\otimes x_{n-1}.\label{eq:b}
\end{equation}
The cohomology of this complex is called the \emph{Hochschild homology} of $\mathcal{A}$, denoted $HH_*(\mathcal{A})$, and the above chain complex is also called the cyclic bar model.  \par\indent
A key feature of Hochschild homology is that it admits a chain level $S^1$-action. There are various ways to interpret what this means, which we further explore in later parts of the paper. For the moment, we recall one such interpretation via the \emph{Connes operator}.\par\indent
For the purpose of describing the Connes operator, it would be more convenient to use a quasi-isomorphic chain model for Hochschild homology, namely the \emph{normalized Hochschild complex}. It is defined as
\begin{equation}
\overline{CC}_*(\mathcal{A}):=CC_*(\mathcal{A})/N,     
\end{equation}
where $N$ denotes the subcomplex generated by elements of the form $\mathbf{x}\otimes x_1\otimes\cdots\otimes 1\otimes\cdots\otimes x_n$ ($1$ sits in the $i$-th entry, $i>0$). The differential on $\overline{CC}_*(\mathcal{A})$ is inherited from $CC_*(\mathcal{A})$ under the quotient, which we still denote by $b$. It is a standard exercise using the simplicial description of the Hochschild chain complex that the quotient map $CC_*(\mathcal{A})\rightarrow \overline{CC}_*(\mathcal{A})$ is a quasi-isomorphism. \par\indent
The Connes operator is the degree $-1$ operator on $\overline{CC}_*(\mathcal{A})$ given by
\begin{equation}
B(x_0\otimes x_1\otimes\cdots\otimes x_n):=\sum_{i=0}^n (-1)^{(\sum_{i}^n(|x_i|-1))(|x_0|+\sum_{1}^{i-1}(|x_i|-1))} 1\otimes x_i\otimes x_{i+1}\otimes\cdots\otimes x_n\otimes x_0\otimes\cdots\otimes x_{i-1}.\label{eq:B}
\end{equation}
$B$ satisfies 
\begin{equation}
b\circ B+B\circ b=0,\;\;\;B^2=0.    
\end{equation}
In particular, it defines an action of the dg algebra $R[\epsilon]/\epsilon^2\simeq C_*(S^1)$ (where $|\epsilon|=-1$) on $\overline{CC}_*(\mathcal{A})$.\par\indent 
The \emph{negative cyclic chain complex} is the following explicit complex computing the homotopy fixed point of this $S^1$-action: 
\begin{equation}
CC^{-}_*(\mathcal{A}):=(\overline{CC}_*(\mathcal{A})[[t]],b+tB)\cong \mathrm{RHom}_{R[\epsilon]/\epsilon^2}(R,\overline{CC}_*(\mathcal{A})),\;|t|=2
\end{equation}
where we think of $t$ as the generator of $H^2_{S^1}(\mathrm{pt})$. The homology of this complex is called the \emph{negative cyclic homology} of $\mathcal{A}$ and denoted $HH^-_*(\mathcal{A})$. \par\indent
Finally, the \emph{periodic cyclic chain complex} is defined as 
\begin{equation}
CC^{per}_*(\mathcal{A}):=(\overline{CC}_*(\mathcal{A})((t)),b+tB),
\end{equation}
whose homology is called the \emph{periodic cyclic homology} of $\mathcal{A}$ and denoted $HH^{per}_*(\mathcal{A})$.
\subsection{The $C_N$-equivariant Hochschild complex}
Let $N$ be any positive integer. As seen in the previous subsection, the Hochschild chain complex of a dg algebra is equipped with an $S^1$-action. It is natural to ask the following question: consider the inclusion $C_N\subset S^1$ of the finite cyclic group of order $N$ as the $N$-th roots of unity, how can we describe induced $C_N$-action? To answer this question, we consider the following variant of the Hochschild chain complex, cf. \cite{Che1}, \cite{Che2} for a more detailed study.\\
\begin{mydef}
The \emph{$N$-fold Hochschild chain complex} ${}_NCC_*(\mathcal{A})$ is defined as
\begin{equation}
{}_NCC_*(\mathcal{A}):=\bigoplus_{k_1,\cdots,k_N\geq 0} \mathcal{A}\otimes_R\mathcal{A}[1]^{\otimes_R k_1}\otimes_R\cdots\otimes_R\mathcal{A}\otimes_R\mathcal{A}[1]^{\otimes_R k_N},
\end{equation}
\end{mydef}
The differential, denoted $b^N$, is given by a sum $b^N=\sum_{i=1}^N b^N_i$, where
$$
b^N_i(\mathbf{x}^1\otimes x_1^1\otimes \cdots\otimes x^1_{k_1}\otimes\mathbf{x}^2\otimes x_1^2\otimes \cdots\otimes x^2_{k_2}\otimes \cdots\otimes \mathbf{x}^N\otimes x_1^N\otimes \cdots\otimes x^N_{k_N}):=
$$
$$\sum_{l=0}^{k_i}\pm\mathbf{x}^1\otimes x_1^1\otimes \cdots\otimes x^1_{k_1}\otimes\cdots\otimes\mathbf{x}^i\otimes x^i_1 \otimes\cdots\otimes dx^i_l\otimes \cdots\otimes x^i_{k_i}\otimes \cdots\otimes \mathbf{x}^N\otimes x_1^N\otimes \cdots\otimes x^N_{k_N} +$$
\begin{equation}
\sum_{l=0}^{k_i}\pm\mathbf{x}^1\otimes x_1^1\otimes \cdots\otimes x^1_{k_1}\otimes\cdots\otimes\mathbf{x}^i\otimes x^i_1 \otimes\cdots\otimes x^i_{l-1}x^i_l\otimes \cdots\otimes x^i_{k_i}\otimes \cdots\otimes \mathbf{x}^N\otimes x_1^N\otimes \cdots\otimes x^N_{k_N},    
\end{equation}
with the usual Koszul signs and where $x^i_{-1}$ is interpreted as $x^{(i-1)}_{k_{i-1}}$ if $k_{i-1}>0$ and $0$ if $k_{i-1}=0$. \par\indent
There is a quasi-isomorphism for $N\geq 2$, given by 
\begin{equation}
\epsilon^0_{N-1,N}: {}_{N}CC_*(\mathcal{A})\rightarrow {}_{N-1}CC_*(\mathcal{A})
\end{equation}
defined by
$$
\epsilon^0_{N-1,N}(\mathbf{x}^1\otimes x_1^1\otimes \cdots\otimes x^1_{k_1}\otimes\mathbf{x}^2\otimes x_1^2\otimes \cdots\otimes x^2_{k_2}\otimes \cdots\otimes \mathbf{x}^N\otimes x_1^N\otimes \cdots\otimes x^N_{k_N})
$$
\begin{equation}
:=\begin{cases}
\mathbf{x}^1\mathbf{x}^2\otimes x_1^2\otimes \cdots\otimes x^2_{k_2}\otimes \cdots\otimes \mathbf{x}^N\otimes x_1^N\otimes \cdots\otimes x^N_{k_N}\quad\mathrm{if}\;\;k_1=0\\
0\quad\mathrm{if}\;\;k_1>0
\end{cases}
\end{equation}
where in the first case, $\mathbf{x}^1_0\mathbf{x}^2_0$ becomes a new distinguished entry. It can be easily checked that $\epsilon^0_{N-1,N}$ is a quasi-isomorphism for all $N\geq 2$. Taking the composition of $\epsilon^0_{k-1,k}$ for $k=2,3,\cdots,N$, one obtains a quasi-isomorphism
\begin{equation}
\Phi^0_N:\epsilon^0_{1,2}\circ\cdot\circ\epsilon^0_{N-1,N}: {}_NCC_*(\mathcal{A})\rightarrow CC_*(\mathcal{A}).\label{eq:Phi^0_N}
\end{equation}
There is a chain level $\mathbb{Z}/N$-action on ${}_NCC_*(\mathcal{A})$, where the generator $\tau\in\mathbb{Z}/N$ acts by cyclically permuting the $N$ `blocks':
$$\tau: \mathbf{x}^1\otimes x_1^1\otimes \cdots\otimes x^1_{k_1}\otimes\mathbf{x}^2\otimes x_1^2\otimes \cdots\otimes x^2_{k_2}\otimes \cdots\otimes \mathbf{x}^N\otimes x_1^N\otimes \cdots\otimes x^N_{k_N}\mapsto $$
\begin{equation}\label{eq:tau}
(-1)^{\dag}\mathbf{x}^N\otimes x_1^N\otimes \cdots\otimes x^N_{k_N}\otimes \mathbf{x}^1\otimes x_1^1\otimes \cdots\otimes x^1_{k_1}\otimes\cdots\otimes \mathbf{x}^{N-1}\otimes x_1^{N-1}\otimes \cdots\otimes x^{N-1}_{k_{N-1}},    
\end{equation}
where
\begin{equation}
\dag=\big(|\mathbf{x}^N|+\sum_{i=1}^{k_N} \|x^N_i\|\big)\cdot\big(\sum_{j=1}^{N-1} |\mathbf{x}^j|+\sum_{j=1}^{N-1}\sum_{i=1}^{k_j}\|x^j_{i}\|\big)    
\end{equation}
is the Koszul sign. One easily verifies that $\tau\circ d_{{}_NCC}=d_{{}_NCC}\circ\tau$ and $\tau^{N}=1$. There is a precise sense in which ${}_NCC_
*(\mathcal{A})$ equipped with the above $C_N$ action is induced from the $S^1$-action on $CC_*(\mathcal{A})$, under the quasi-isomorphism of the underlying chain complexes \eqref{eq:Phi^0_N}, cf. \cite[Section 5]{Che1} for more details.   \par\indent
As in the cyclic case, there is an explicit model that represents the homotopy fixed point of the above $C_N$-action on ${}_NCC_*(\mathcal{A})$, or equivalently,
\begin{equation}
\mathrm{Rhom}_{R[C_N]}(R,{}_NCC_*(\mathcal{A})).    
\end{equation}
Recall that if $X$ is a cohomologically graded complex with a $C_N$-action, there is a canonical complex computing $\mathrm{Rhom}_{R[C_N]}(k,X)$ given by $X[[t,\theta]]$, where $|t|=2,|\theta|=1,\theta^2=0$. 
Let $\tau\in C_N$ be the standard generator, then the differential is given by
\begin{equation}
\begin{cases}
d_{eq}(x)=dx+(-1)^{|x|}(\tau x-x)\theta,\\
d_{eq}(x\theta)=dx\,\theta+(-1)^{|x|}(x+\tau x+\cdots+\tau^{N-1}x)t.    \label{eq:eq diff} 
\end{cases}  
\end{equation} 
and extended $t$-linearly.\par\indent
The \emph{negative $C_N$-equivariant Hochschild complex} of $\mathcal{A}$ is defined to be   
\begin{equation}CC^{C_N,-}_*(\mathcal{A}):={}_NCC_*(\mathcal{A})[[t,\theta]]\label{eq:C_N negative}
\end{equation} 
equipped with the differential \eqref{eq:eq diff}. The \emph{periodic $C_N$-equivariant Hochschild complex} is obtained from \eqref{eq:C_N negative} by inverting $t$, and is denoted $CC_*^{C_N,per}(\mathcal{A})$. We denote their respective homologies by $HH^{C_N,-}_*$ and $HH^{C_N,per}_*$. 

\subsection{A simplicial approach to Hochschild invariants}
For the purpose of comparing the various equivariant (e.g. $S^1$ and $C_p$) versions of the Hochschild chain complex and certain operations on them, it is often convenient to use an alternative definition via simplicial or cyclic objects. In this subsection, we recall the relevant constructions.\par\indent
\textbf{a) Cyclic homology via Connes' cyclic category $\Lambda$}. First, recall that the \emph{simplex category} $\Delta$ is the category with objects $[n], n\geq 0$, and morphisms from $[n]$ to $[m]$ given by order preserving maps from $\{0\leq 1\leq 2\leq \cdots\leq n\}$ to $\{0\leq 1\leq 2\leq \cdots\leq m\}$. There are distinguished morphisms $\sigma_i:[n]\rightarrow [n+1], 0\leq i\leq n+1$, $\epsilon_i: [n]\rightarrow [n-1], 0\leq i\leq n-1$ defined by 
\begin{equation}
\sigma_i(j)=\begin{cases} j,\;\mathrm{if}\;j< i\\
j+1, \;\mathrm{if}\;j\geq i
\end{cases}   \quad\quad,\quad\quad \epsilon_i(j)=\begin{cases} j,\;\mathrm{if}\;j\leq i\\
j-1, \;\mathrm{if}\;j> i
\end{cases}\quad,
\end{equation}
which generate $\Delta$. The \emph{semi-simplex category} $\overrightarrow{\Delta}\subset \Delta$ is the subcategory on the same set of objects but whose morphisms are injective order preserving maps; in particular, it is generated by the $\sigma_i$'s. \par\indent
The cyclic category $\Lambda$ can be formally obtained from $\Delta$ by adjoining the cyclic permutations of $\{0\leq 1\leq 2\leq\cdots\leq n\}$ as automorphisms of $[n]$. Instead, we give a geometric description of $\Lambda$ following \cite{Con}, \cite[E.6.1.2]{Lod}:
\begin{itemize}
    \item Objects of $\Lambda$ are labeled by natural numbers $[n], n\geq 0$. One thinks of $[n]$ as a configuration of $n+1$ marked points $z_0<\cdots<z_n$ on a circle.
    \item Morphisms from $[n]$ to $[m]$ are homotopy classes of degree $1$ nondecreasing maps from $S^1$ to $S^1$ that send marked points to marked points.
\end{itemize}
Combinatorially, a morphism from $[n]$ to $[m]$ is uniquely determined by 
\begin{itemize}
    \item a cyclic reordering of $[n]$ given by $[n]_{\sigma}=\{\sigma_0<\sigma_1<\cdots<\sigma_n\}$, where $\sigma\in \frac{\mathbb{Z}}{n+1}\subset S_{n+1}$ and
    \item  an (ordered) partition of $[n]_{\sigma}$ into $m+1$ subsets, i.e. (possibly empty) order subsets $f_0,\cdots,f_m\subset [n]_{\sigma}$ such that $[n]_{\sigma}=f_0\star f_1\star\cdots\star f_m$, where
$\star$ denotes the join of partially ordered sets.
\end{itemize}
To make the identification, given a homotopy class of $f:S^1\rightarrow S^1$, let $f_i$ be the set of (indices of) marked points that are sent to $z_i$ by $f$. By abuse of notation, we also write the ordered set $f_i$ as $f^{-1}(i)$.
\begin{figure}[H]
 \centering
 \includegraphics[width=0.8\textwidth]{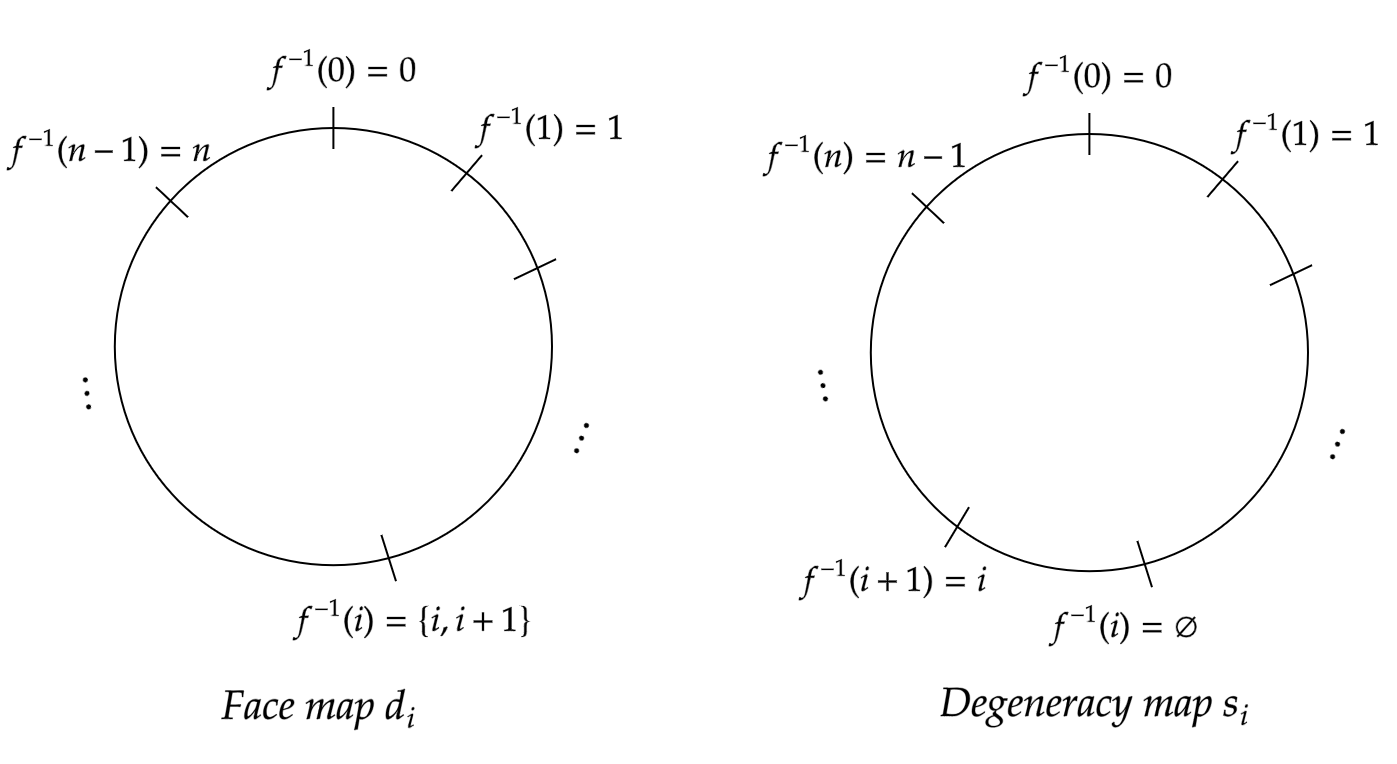}
 \caption{}
 \label{fig:figure1}
\end{figure}
Schematically, we represent a morphism $f: [n]\mapsto [m]$ by drawing a circle with $m+1$ ordered marked points (thought of as the object $[m]$), and write the elements of $f^{-1}(i)$ in order next to the $i$-th marked point. There is an inclusion of $i:\Delta^{op}\subset\Lambda$ as the subcategory whose morphisms are (homotopy classes of) maps that sends $z_0$ to $z_0$, and a further inclusion $\overrightarrow{\Delta}^{op}\subset \Delta^{op}\subset \Lambda$ as the subcategory whose morphisms are (homotopy classes of) maps which send $z_0$ to $z_0$ and are surjective on marked points. The face maps $d_i=\sigma_i^{op}: [n]\rightarrow [n-1], 0\leq i\leq n$ and degeneracy maps $s_i=\epsilon_i^{op}:[n-1]\rightarrow [n], 0\leq i\leq n-1$ are presented in Figure \ref{fig:figure1}. The face and degeneracy maps, together with the morphisms $\tau_n:[n]\rightarrow [n]$ determined by cyclically permuting the $n+1$ marked points on $S^1$ (i.e. $\tau_n(i)=i+1$), generate $\Lambda$. \par\indent
Let $(\mathcal{A},d,\cdot)$ be a dg algebra, and let $\mathrm{Mod}_R$ denote the dg category of chain complexes over $R$. \\
\begin{mydef}
The \emph{Hochschild chain functor} of $\mathcal{A}$ is the functor
\begin{equation}
\mathcal{A}_{\sharp}: \Lambda\rightarrow \mathrm{Mod}_R \label{eq:HH functor}
\end{equation}
given by
\begin{itemize}
    \item On objects, $\mathcal{A}_{\sharp}([n]):=(\mathcal{A},d)^{\otimes n+1}$ as chain complex.
    \item On morphisms, $\mathcal{A}_{\sharp}$ sends $f\in \Lambda([n],[m])$ to the map $\mathcal{A}^{\otimes n+1}\rightarrow \mathcal{A}^{\otimes m+1}$ given by
\begin{equation}
a_0\otimes a_1\otimes \cdots\otimes a_n\mapsto \prod_{j\in f^{-1}(0)}a_j\otimes \prod_{j\in f^{-1}(1)}a_j\otimes\cdots\otimes \prod_{j\in f^{-1}(m)}a_j,\label{eq:HH functor on morphism}
\end{equation}
where each product is taken with respect to the induced ordering on the elements of $f^{-1}(j)$, and if $f^{-1}(j)=\emptyset$, it is understood to be the unit $1\in \mathcal{A}$. 
\end{itemize}
\end{mydef}
By inspecting the formulae of the Hochschild and cyclic differentials \eqref{eq:b} \eqref{eq:B}, it is clear that one can associate to any functor $E:\Lambda\rightarrow \mathrm{Mod}_R$ chain complexes (cf. \cite[Definition 3.5]{Che1})
\begin{equation}
CC_*(E),\quad CC_*^{-}(E),\quad CC_*^{per}(E),  
\end{equation}
which will be called its Hochschild, negative cyclic and periodic cyclic complex, such that when applied to $E=\mathcal{A}_{\sharp}$ recover the usual notions. \par\indent
By \cite[Proposition 1.1]{Hoy}, there exists a homotopy exact square of $\infty$ categories
\begin{equation}
\begin{tikzcd}[row sep=1.2cm, column sep=0.8cm]
\Delta^{op}\arrow[r,"i"]\arrow[d]&\Lambda\arrow[d]\\
*\arrow[r]& BS^1\label{eq:Lambda homotopy exact sq}
\end{tikzcd},
\end{equation}
where $S^1$ denotes the circle group. Given $X:\Lambda\rightarrow \mathrm{Mod}_R$ a cyclic chain complex, let $\int_{S^1}X\in\mathrm{Fun}(BS^1,\mathrm{Mod}_R)$ denote its left Kan extension along the map $\Lambda\rightarrow BS^1$ (which is an infinity groupoid completion). \par\indent
The following result from \cite{Hoy} allows one to recover various versions of cyclic homology as abstract $\infty$ categorical (co)-limits.\\
\begin{thm}\cite[Theorem 2.3]{Hoy}\label{thm:Hoyois theorem}
Let $E:\Lambda\rightarrow \mathrm{Mod}_R$ be a cyclic chain complex. Under the identification $\gamma: R[\epsilon]/\epsilon^2\simeq C_{-*}(S^1;R)$ and the induced equivalence of categories $\gamma^*:\mathrm{Fun}(BS^1,\mathrm{Mod}_R)\simeq \mathrm{Mod}_{R[\epsilon]/\epsilon^2}$,
there is an equivalence of mixed complexes
\begin{equation}
CC_*(E)\simeq \gamma^*\int_{S^1}E.\label{eq:Connes operator via left kan extension}
\end{equation}\qed
\end{thm}
It follows from Theorem \autoref{thm:Hoyois theorem} that there are quasi-isomorphisms of chain complexes
\begin{equation}
CC_*^-(E)\simeq (\int_{S^1}E)^{hS^1}\qquad\mathrm{and}\qquad CC^{per}_*(E)\simeq (\int_{S^1}E)^{tS^1},\label{eq:CC(A) as abstract categorical (co)limits}
\end{equation}
where $hS^1$ and $tS^1$ denotes the homotopy and Tate fixed points, respectively. \par\indent
\textbf{b) The finite $p$-cyclic category ${}_p\Lambda$}. Let $p$ be a prime number. The \emph{finite $p$-cyclic category} ${}_p\Lambda$, introduced in \cite{Che1}, has
\begin{itemize}
    \item Objects are $p$-tuples of natural numbers $[k_1,\cdots,k_p]$, thought of as a configuration of $k_1+\cdots+k_p+p$ increasing marked points on the circle, with the $0, k_1+1, k_1+k_2+2,\cdots,k_1+\cdots+k_{p-1}+p-1$-th points marked as distinguished.
    \item Morphisms from $[k_1,\cdots,k_p]$ to $[k_1',\cdots,k_p']$ are homotopy classes of degree $1$ nondecreasing maps $f: S^1\rightarrow S^1$ that sends marked points to marked points, and is further required to send the distinguished marked points bijectively to distinguished marked points. In particular, by the nondecreasing condition, $f$ must act as a $C_p$ cyclic permutation on the distinguished points. 
\end{itemize}
There is an obvious functor $j:{}_p\Lambda\rightarrow \Lambda$ which on objects sends $[k_1,\cdots,k_p]$ to $[k_1+\cdots+k_p+p-1]$, where the $p$ distinguished marked points become ordinary marked points. Moreover, $i_p:(\Delta^{op})^p\subset {}_p\Lambda$ sits as the subcategory whose morphisms fix each distinguished marked point. Let $\mathfrak{j}: C_p\subset S^1$ denote the standard inclusion as $p$-th roots of unity. These functors fit into homotopy exact squares \cite[Lemma 3.18, Lemma 5.2]{Che1}
\begin{equation}
\begin{tikzcd}[row sep=1.2cm, column sep=0.8cm]\label{eq:two homotopy exact sq}
(\Delta^{op})^p\arrow[r,"{i_p}"]\arrow[d]&{}_p\Lambda\arrow[d]\arrow[r,"j"]&\Lambda\arrow[d]\\
*\arrow[r]& B C_p\arrow[r,"{\mathfrak{j}}"]&BS^1
\end{tikzcd}.
\end{equation}
The following result is a counterpart of Theorem \autoref{thm:Hoyois theorem} in the finite cyclic context.\\
\begin{thm}\label{thm:Cp version of Hoyois}\cite[Theorem 3.20]{Che1}
Let $E: {}_p\Lambda\rightarrow \mathrm{Mod}_R$ be a finite $p$-cyclic chain complex. Under the identification $\mathrm{Fun}(BC_p,\mathrm{Mod}_R)\simeq \mathrm{Mod}_{R[C_p]}$, there is an equivalence of $C_p$-chain complexes
\begin{equation}
{}_pCC_*(E)\simeq \int_{C_p}E,     
\end{equation}
where $\int_{C_p}$ denotes the left Kan extension along ${}_p\Lambda\rightarrow BC_p$. 
\end{thm}
It follows from Theorem \autoref{thm:Cp version of Hoyois} that
 there are quasi-isomorphisms of chain complexes
\begin{equation}
CC_*^{C_p,-}(E)\simeq (\int_{C_p}E)^{hC_p}\qquad\mathrm{and}\qquad CC^{C_p,per}_*(E)\simeq (\int_{C_p}E)^{tC_p}.\label{eq:pCC(A) as abstract categorical (co)limits}
\end{equation}
One consequence of reinterpreting cyclic homology and $C_p$-equivariant Hochschild homology in terms of abstract $\infty$ categorical (co)limits is the following comparison result, which will be used extensively in later parts of this paper. \\
\begin{prop}\label{thm:pGysin for CC}
Let $R$ be a ring of characteristic $p$. Let $E: \Lambda\rightarrow \mathrm{Mod}_{R}$ be a cyclic chain complex. Then there is a quasi-isomorphism of $R$-chain complexes, natural in $E$, 
\begin{equation}
\phi_p: CC^{per}_*(E)\otimes_{R((t))}R((t,\theta))\simeq CC^{C_p,per}_*(j^*E)\label{eq:S^1 C_p comparison}.
\end{equation}
\end{prop}
\emph{Proof}. Since the right square in \eqref{eq:two homotopy exact sq} is homotopy exact, there is a natural quasi-isomorphism
\begin{equation}
\mathfrak{j}^*\int_{S^1}E\simeq \int_{C_p}j^*E.      
\end{equation}
Therefore, by \cite[Lemma IV.4.12]{NS} applied to $\int_{S^1}E$, there is a natural quasi-isomorphism
\begin{equation}
 (\int_{S^1}E)^{tS^1}\otimes_{R((t))}R((t,\theta))\simeq (\int_{C_p}j^*E)^{tC_p}. \label{eq:p Gysin}  
\end{equation}
The proof follows from combining \eqref{eq:p Gysin}, Theorem \autoref{thm:Hoyois theorem} and \autoref{thm:Cp version of Hoyois}.\qed\par\indent
\textbf{c) Hochschild cohomology as a cosimplicial object}. We recall the standard construction of the Hochschild cochain complex as the realization of a cosimplicial chain complex, cf. \cite{Lod}. The \emph{Hochschild cochain functor} of a dg algebra $\mathcal{A}$ is the cosimplicial chain complex
\begin{equation}
\mathcal{A}^{\sharp}: \Delta\rightarrow \mathrm{Mod}_R  \label{eq:HH cochain functor}   
\end{equation}
given by
\begin{itemize}
    \item On objects, $\mathcal{A}^{\sharp}([n]):=\mathrm{Hom}_R(\mathcal{A}^{\otimes n},\mathcal{A})$ as chain complex.
    \item On morphisms, if $\partial_i\in\Delta([n],[n+1]), 0\leq i\leq n+1$ is the map in that skips $i$, then
    \begin{equation}\mathcal{A}^{\sharp}(\partial_i)(\phi)(a_1,\cdots,a_{n+1}):=
     \begin{cases}
      a_1\phi(a_2,\cdots,a_{n+1}), \quad\qquad\qquad i=0\\
      \phi(a_1,\cdots,a_{i}a_{i+1},\cdots,a_{n+1}),\quad 1\leq i\leq n\\
      \phi(a_1,\cdots,a_n)a_{n+1},\quad\qquad\qquad i=n+1
     \end{cases} ;
    \end{equation}
    if $\sigma_i\in\Delta([n],[n-1]), 0\leq i\leq n-1$ is the map in that doubles $i$, then
    \begin{equation}
      \mathcal{A}^{\sharp}(\sigma_i)(\phi)(a_1,\cdots,a_{n-1}):=\phi(a_1,\cdots,1,\cdots,a_{n-1}),
    \end{equation}
    where the unit $1$ is inserted in the $i$-th position. 
\end{itemize}
It is evident from the definition that the realization of $\mathcal{A}^\sharp$ recovers the Hochschild cochain complex from \eqref{eq:CC^*}.\par\indent
\textbf{d) Topological realization of (co)cyclic spaces.} a)-c) can be viewed as specific instances of realization of (co)cyclic and (co)simplicial objects in the category of chain complexes. Here we discuss the topological analogue. \par\indent
Let $\Delta^n=\{\sum_{i=0}^n t_i=1,t_i\geq 0\}\subset \mathbb{R}^{n+1}$ be the standard $n$-simplex. The assignment $[n]\mapsto \Delta^n$ defines a cosimplicial space $\Delta^{\bullet}$ (since linear maps between simplices are uniquely determined by their values on vertices). \par\indent
Let $X: \Delta^{op}\rightarrow \mathrm{Spaces}$ be a simplicial space, then its \emph{realization} is defined as
\begin{equation}\label{eq:realization}
|X|:=\coprod_{n\geq 0}X_n\times \Delta^n/(f^*(x),y)\sim (x,f_*(y)), \;\;f\in\Delta([n],[m]).    
\end{equation}
Let $Y: \Delta\rightarrow \mathrm{Spaces}$ be a cosimplicial space, then its \emph{totalization} is defined as the space of cosimplicial maps
\begin{equation}\label{eq:totalization}
\mathrm{Tot}(Y):=\mathrm{Hom}_{\Delta}(\Delta^{\bullet},Y)\subset \prod_{n\geq 0}\mathrm{Map}(\Delta^n,Y^n)
\end{equation}
with the subspace topology.\par\indent
We recall the following well known result, cf. for instance \cite{Jon}. \\
\begin{lemma}\label{thm:circle actions of realization of (co)cyclic sets}
i) If $X: \Lambda\rightarrow\mathrm{Spaces}$ is a cyclic space, then there is a canonical $S^1$-action on $|X|$.\par\indent
ii)  If $X: \Lambda^{op}\rightarrow\mathrm{Spaces}$ is a cocyclic space, then there is a canonical $S^1$-action on $\mathrm{Tot}(Y)$.
\end{lemma}
Before giving the proof, we remark that in the rest of the paper, $S^1$ will be given the clockwise orientation, i.e. we identify $\mathbb{R}/\mathbb{Z}$ with $\{|z|=1\}\subset \mathbb{C}^*$ via $\theta\mapsto e^{-2\pi i\theta}$.\par\indent
\emph{Proof}. The key is the cocyclic space $\Lambda^{\bullet}$, which as a cosimplicial space, is $S^1\times \Delta^{\bullet}$. The cocyclic structure is defined such that for $\tau_n\in\mathrm{Aut}_{\Lambda}([n],[n])=\mathbb{Z}/(n+1)\mathbb{Z}$ the standard generator,
\begin{equation}
\tau_n^*(z,t_0,\cdots,t_n):=(ze^{-2\pi i t_n},t_1,\cdots,t_{n},t_0).    
\end{equation}
One observes that the obvious level-wise $S^1$-action on $\Lambda^{\bullet}$ is compatible with its cocyclic structure maps. Therefore, there are natural $S^1$-actions on
\begin{equation}
|X|=\Delta^{\bullet}\otimes_{\Delta^{op}}X=\Lambda^{\bullet}\otimes_{\Lambda}X    
\end{equation}
and
\begin{equation}
\mathrm{Tot}(Y)=\mathrm{Hom}_{\Delta}(\Delta^{\bullet},Y)=\mathrm{Hom}_{\Lambda^{op}}(\Lambda^{\bullet},Y).    
\end{equation}\qed\par\indent
Let $\mathrm{Mon}(I,\partial I)$ denote the space of weakly monotone self-maps of $I=[0,1]$ fixing the endpoints. There is a homeomorphism $\mathrm{Mon}(I,\partial I)\cong\mathrm{Tot}(\Delta^{\bullet})$ given by
\begin{equation}\label{eq:monotonic map of interval as totalization}
 f\mapsto \big((0\leq x_1\leq \cdots\leq x_k\leq 1)\mapsto (0\leq f(x_1)\leq \cdots\leq f(x_k)\leq 1)\big)_{k\geq 0}, 
\end{equation}
where we identify $\Delta^k$ with $\{0\leq x_1\leq \cdots\leq x_k\leq 1\}\subset [0,1]^k$ via $x_1=t_0,x_2=t_0+t_1,\cdots,x_k=t_0+\cdots+t_{k-1}$. Let $\mathrm{Mon}(S^1)$ denote the weakly \emph{monotonically increasing degree $1$} self-maps of $S^1=I/\partial I=\mathbb{R}/\mathbb{Z}$. There is a chain of homeomorphisms
\begin{equation}\label{eq:monotonic map of circle as totalization}
\mathrm{Mon}(S^1)\cong S^1\times \mathrm{Mon}(I,\partial I)\cong S^1\times \mathrm{Tot(\Delta^{\bullet})}\cong \mathrm{Tot}(\Lambda^{\bullet}).  
\end{equation}
By Lemma \autoref{thm:circle actions of realization of (co)cyclic sets} ii), there is a natural $S^1$-action on $\mathrm{Tot}(\Lambda^{\bullet})$. The following statement is an easy exercise from the definition.\\
\begin{lemma}\label{thm:circle action on Mon(S^1)}
Under the homeomorphism of \eqref{eq:monotonic map of circle as totalization}, the $S^1$-action on $\mathrm{Tot}(\Lambda^{\bullet})$ corresponds to the $S^1$-action on $\mathrm{Mon}(S^1)$ given by 
\begin{equation}
\theta\cdot f:=(x\mapsto f(x+\theta)),\quad\theta\in S^1, f\in\mathrm{Mon}(S^1).
\end{equation}\qed
\end{lemma}

\section{The Getzler-Gauss-Manin connection}
In this section, we work in the following relative setting: fix $R$ a ring over a field $\mathbf{k}$ (the primary example would be the ring of functions on a smooth affine irreducible curve). Let $(\mathcal{A}, d,\cdot)$ be a $\mathbb{Z}$ or $\mathbb{Z}/2$-graded, smooth and proper dg-algebra over $R$. We remark that all the formulas in this section work more generally for an $A_{\infty}$-algebra $(\mathcal{A}, \{\mu^d\}_{d\geq 1})$. \par\indent
Given a Hochschild cocycle $D\in CC^*(\mathcal{A})$, Getzler defined a `contraction' operator
\begin{equation}e_D:CC_*(\mathcal{A})\rightarrow CC_{*+|D|}(\mathcal{A}),\label{eq:e_D}
\end{equation}
a `Lie derivative' operator 
\begin{equation}\mathcal{L}_D: CC_*(\mathcal{A})\rightarrow CC_{*+|D|-1}(\mathcal{A}),\label{eq:L_D}
\end{equation}
and an auxiliary operator
\begin{equation}
E_D: CC_*(\mathcal{A})\rightarrow CC_{*+|D|-2}(\mathcal{A})  \label{eq:E_D}  
\end{equation}
which satisfy 
\begin{equation}
[e_D,b]=0, \;\;[E_D,B]=0   \label{eq:e_D,E_D basic}
\end{equation}
and an analogue of the Cartan homotopy formula (up to homotopy): 
\begin{equation}
[e_D, B]=\mathcal{L}_D-[E_D,b].     \label{eq:cartan}
\end{equation}
We now recall Getzler's construction of the Gauss-Manin connection in \cite{Get}. First, after replacing $\mathcal{A}$ by a semi-free resolution, we fix an $R$-basis of $\mathcal{A}=\bigoplus_i \mathcal{A}_i$ as a graded $R$-module. Let 
\begin{equation}\nabla': \mathcal{A}\rightarrow \mathcal{A}\otimes_R\Omega^1_R\label{eq:nabla'}
\end{equation}
be the trivial connection on $\mathcal{A}$ (as a graded $R$-module) with respect to this basis. Extend $\nabla'$ to a connection $\bigoplus_n \mathcal{A}^{\otimes n+1}\rightarrow \bigoplus_n \mathcal{A}^{\otimes n+1}\otimes_R\Omega^1_R$ by the Leibniz rule. Define
 \begin{equation}
  \kappa:=[\nabla',\mu]: \bigoplus_n\mathcal{A}^{\otimes n}\rightarrow \mathcal{A}\otimes_R\Omega^1_R.\label{eq:kappa}
 \end{equation}
 Equivalently, \eqref{eq:kappa} is obtained by expanding $\mu$ as a matrix with respect to the chosen $R$-basis, and differentiating entry-wise. In particular, \eqref{eq:kappa} gives rise to a Hochschild cocycle and its cohomology class $[\kappa]\in HH^2(\mathcal{A})\otimes_R\Omega^1_R$, which is independent of the choice of an $R$-basis, is called the \emph{Kodaira-Spencer class} of $\mathcal{A}$. \par\indent
 $\nabla'$ induces a connection on the negative cyclic chain complex $CC^{-}_*(\mathcal{A})$ as a graded $R$-module, which however does not commute with the differential $b+tB$. The failure to commute is measured by
 \begin{equation}
 [\nabla', b]=\mathcal{L}_{\kappa},\;\;[\nabla',B]=0.   \label{eq:nabla' property} 
 \end{equation}
\eqref{eq:e_D,E_D basic}, \eqref{eq:cartan} and \eqref{eq:nabla' property} imply that the chain level connection on $CC^{-}_*(\mathcal{A})$
 \begin{equation}
 \nabla^{GGM}:=\nabla'-\frac{1}{t}(e_{\kappa}+tE_{\kappa})    \label{eq:GGM}
 \end{equation}
 descends to the negative cyclic homology. \\
 \begin{mydef}
The connection $HH^{-}_*(\mathcal{A})\rightarrow HH^{-}_*(\mathcal{A})\otimes_R\Omega^1_R$ induced by \eqref{eq:GGM} is called the \emph{Getzler-Gauss-Manin connection} of $\mathcal{A}$. We also refer to the induced (by inverting $t$) connection on $HH^{per}_*(\mathcal{A})$ as the Getzler-Gauss-Manin connection.  
 \end{mydef}
 
\subsection{Kaledin's definition}
There is a formula-free definition of the Getzler-Gauss-Manin connection on the periodic cyclic homology due to Kaledin \cite{Kal2}; we follow the exposition in \cite[section 3.2]{PVV}. As a preliminary, we recall the definition of a \emph{Grothendieck connection}.\par\indent
Let $X=\mathrm{Spec}\,R$ and $X^{[2]}$ denote the first order infinitesimal neighborhood of the diagonal $\Delta\subset X\times X$. Let $p_1, p_2: X^{[2]}\rightarrow X$ denote the two coordinate projections. Let $E$ be a vector bundle on $X$, then a Grothendieck connection on $E$ is an $R$-linear splitting $s: E\rightarrow {p_1}_*p_2^*E$ of the Atiyah sequence
\begin{equation}
0\rightarrow E\otimes_{\mathcal{O}_X} \Omega^1_X\rightarrow {p_1}_*p_2^*E\xrightarrow{\pi} E\rightarrow 0.   
\end{equation}
It is a classical fact that Grothendieck's definition of a connection is equivalent to the usual definition. Explicitly, given a connection $\nabla: E\rightarrow E\otimes_{\mathcal{O}_X}\Omega^1_X$, we obtain a section of $\pi$ via
\begin{equation}
s:=j^1+\nabla,\label{eq:first jet}
\end{equation}
where $j^1: E\rightarrow {p_1}_*p_2^*E$ is the $\mathbf{k}$-linear section of $\pi$ known as the \emph{first jet}: writing ${p_1}_*p_2^*E= E\otimes_R(R\otimes_{\mathbf{k}} R/I^2)$, where $I\subset R\otimes_{\mathbf{k}} R$ is the ideal corresponding to the diagonal embedding, then $j^1(e):=e\otimes 1\otimes 1$. Conversely, given a Grothendieck connection $s$, one obtains a connection in the usual sense via $\nabla:=s-j^1$.  \par\indent
Consider the two-term filtration on ${p_1}_*p_2^*\mathcal{A}$ 
\begin{equation}
0\rightarrow \mathcal{A}\otimes_{\mathcal{O}_X} \Omega^1_X\rightarrow {p_1}_*p_2^*\mathcal{A}\rightarrow \mathcal{A}\rightarrow 0.\label{eq:first jet sequence}
\end{equation}
\eqref{eq:first jet sequence} induces a filtration $I^i({p_1}_*p_2^*\mathcal{A})_{\sharp}$ on the cyclic object 
\begin{equation}
({p_1}_*p_2^*\mathcal{A})_{\sharp}: \Lambda\rightarrow \mathrm{Ch}_{R}
\end{equation}
associated to ${p_1}_*p_2^*\mathcal{A}$, cf. \cite{Kal2}. Recall that object-wise, $({p_1}_*p_2^*\mathcal{A})_{\sharp}([n]):=({p_1}_*p_2^*\mathcal{A})^{\otimes n+1}$, and $I$ is just the tensor product of the two-term filtration \eqref{eq:first jet sequence}. This further induces a filtration on the periodic cyclic complex, which by an abuse of notation we also denote $I^i$. \par\indent
Consider the following diagram (where $I^i$ is shorthand for $I^iCC^{per}_*({p_1}_*p_2^*\mathcal{A}):=CC^{per}_*(I^i{p_1}_*p_2^*\mathcal{A})$) 
\begin{equation}
\begin{tikzcd}[row sep=1.2cm, column sep=0.8cm]
I^1/I^2\arrow[r]& I^0/I^2\arrow[r,"\pi"] \arrow[d,"m"]&I^0/I^1=CC^{per}_*(\mathcal{A})\\
& {p_1}_*p_2^*CC^{per}_*(\mathcal{A})  &
\end{tikzcd}\label{eq:GGM diagram}
\end{equation}
The chain complex $I^1/I^2$ is contractible by \cite[Lemma 3.1]{PVV}, which in particular implies that $\pi$ is a quasi-isomorphism. Kaledin's definition of the Gauss-Manin connection on $CC^{per}_*(\mathcal{A})$ is $\nabla^{Kal}:=m\pi^{-1}$ (in the derived category). \\
\begin{prop}\label{thm:Kaledin connection agrees with GGM}
Kaledin's connection agrees with Getzler's connection as a morphism $CC^{per}_*(\mathcal{A})\rightarrow {p_1}_*p_2^*CC^{per}_*(\mathcal{A})$ in the derived category.
\end{prop}
\emph{Proof}. For completeness, we review the proof from \cite[Proposition 3.2]{PVV}. \par\indent
Under \eqref{eq:first jet}, we identify the connection $\nabla'$ from \eqref{eq:nabla'} with an $R$-linear section of ${p_1}_*p_2^*\mathcal{A}\rightarrow \mathcal{A}$, denoted
\begin{equation}
\varphi: \mathcal{A}\rightarrow {p_1}_*p^*_2\mathcal{A}. \label{eq:first jet of nabla'}
\end{equation}
By taking tensor products, \eqref{eq:first jet of nabla'} induces an $R$-linear section (as a map of graded $R$-modules) of the projection $\pi$ in \eqref{eq:GGM diagram}, which we also denote as $\varphi: CC^{per}_*(\mathcal{A})\rightarrow I^0/I^2$. Then, 
\begin{equation}
[\varphi, b]=\tilde{\mathcal{L}}_{\kappa}\;, \;\;[\varphi, B]=0,  \label{eq:varphi basic property}
\end{equation}
where 
\begin{equation}
\tilde{\mathcal{L}}_{\kappa}: CC^{per}_*(\mathcal{A})\rightarrow I^1/I^2\simeq \bigoplus_{n\geq 0}\bigoplus_{i=0}^n\mathcal{A}\otimes\cdots\mathcal{A}\otimes \overbrace{(\mathcal{A}\otimes_{\mathcal{O}_X}\Omega^1_X)}^{i\mathrm{th\;position}}\otimes\mathcal{A}\otimes\cdots\otimes \mathcal{A}    
\end{equation} is the Lie derivative along the Kodaira-Spencer element (more precisely, it is a lift of $\mathcal{L}_{\kappa}$ in \eqref{eq:L_D} along the natural projection $I^1/I^2\rightarrow CC_*^{per}(\mathcal{A})\otimes_{\mathcal{O}_X}\Omega^1_X$, i.e. $m\circ\tilde{\mathcal{L}}_{\kappa}=\mathcal{L}_{\kappa}$ where $m$ is as in \eqref{eq:GGM diagram}). Similarly, one can define $\tilde{e}_{\kappa}, \tilde{E}_{\kappa}$ that lifts $e_{\kappa}, E_{\kappa}$. \par\indent
By \eqref{eq:varphi basic property} and the Cartan homotopy formula, $\varphi-\frac{1}{t}(\tilde{e}_{\kappa}+t\tilde{E}_{\kappa}): CC^{per}_*(\mathcal{A})\rightarrow I^0/I^2$ is a chain map which splits $\pi$, which immediately implies that in the derived category,
\begin{equation}
\nabla^{GGM}=m\circ (\varphi-\frac{1}{t}(\tilde{e}_{\kappa}+t\tilde{E}_{\kappa}))=m\circ \pi^{-1}=\nabla^{Kal}. 
\end{equation}\qed

\subsection{The $C_p$-Getzler-Gauss-Manin connection}
For the sake of completeness, we record an explicit formula for a version of the Getzler-Gauss-Manin connection on the negative (and periodic) $C_p$-equivariant Hochschild complex of an $A_{\infty}$-algebra; this will not appear in the rest of the paper. \par\indent
Let $D\in CC^*(\mathcal{A})$ be a Hochschild cochain. One can define the chain level `$1/p$ of a contraction' operator, of degree $|D|$, $e^p_D: {}_pCC_*(\mathcal{A})\rightarrow {}_pCC_*(\mathcal{A})$ by
$$e^p_D(\mathbf{x}^1\otimes x_1^1\otimes \cdots\otimes x^1_{k_1}\otimes\mathbf{x}^2\otimes x_1^2\otimes \cdots\otimes x^2_{k_2}\otimes \cdots\otimes \mathbf{x}^p\otimes x_1^p\otimes \cdots\otimes x^p_{k_p}):=$$
\begin{equation}
\pm \sum_{0\leq i_1<i_2<i_3\leq k_p, 1\leq l\leq k_1} \mu(x^p_{i_1+1},\cdots,x^p_{i_2},D(x^p_{i_2+1},\cdots,x^p_{i_3}),x^p_{i_3+1},\cdots,x^p_{k_p},\mathbf{x^1},x^1_1,\cdots,x^1_{l})\otimes x^1_{l+1}\otimes\cdots \otimes\mathbf{x}^2\otimes\cdots\otimes\mathbf{x}^p\otimes x^p_1\otimes\cdots x^p_{i_1},   \label{eq:e^p_D}
\end{equation}
with the usual Koszul signs. If $D$ is a Hochschild cocycle, it is easy to check that $e^p_D$ commutes with the $p$-fold Hochschild differential $b^p$. However, it does not commute with the action of $\tau\in C_p$. The required nullhomotopy is given by `$1/p$ of a Lie derivative' $E^p_D$, of degree $|D|-1$. Explicitly, 
$$E^p_D(\mathbf{x}^1\otimes x_1^1\otimes \cdots\otimes x^1_{k_1}\otimes\mathbf{x}^2\otimes x_1^2\otimes \cdots\otimes x^2_{k_2}\otimes \cdots\otimes \mathbf{x}^p\otimes x_1^p\otimes \cdots\otimes x^p_{k_p}):=$$
$$\sum_{1\leq i_1<i_2\leq k_1} \pm\mathbf{x}^1\otimes x_1^1\otimes\cdots\otimes x_{i_1}\otimes D(x^1_{i_1+1},\cdots,x^1_{i_2})\otimes x^1_{i_2+1}\otimes \cdots\otimes x^1_{k_1}\otimes\mathbf{x}^2\otimes x_1^2\otimes \cdots\otimes x^2_{k_2}\otimes \cdots\otimes \mathbf{x}^p\otimes x_1^p\otimes \cdots\otimes x^p_{k_p}+$$
\begin{equation}
\sum_{1\leq i_1\leq k_p, 1\leq i_2\leq k_1}\pm D(x^p_{i_1+1},\cdots,\mathbf{x}^1,\cdots,x^1_{i_2})\otimes x^1_{i_2+1}\otimes \cdots\otimes x^1_{k_1}\otimes\mathbf{x}^2\otimes x_1^2\otimes \cdots\otimes x^2_{k_2}\otimes \cdots\otimes \mathbf{x}^p\otimes x_1^p\otimes \cdots\otimes x^p_{i_1}.   \label{eq:E^p_D}
\end{equation}
$E^p_D$ and $e^p_D$ satisfy a $C_p$-analogue of the Cartan homotopy formula: up to sign, these are (where $\delta$ denotes the differential on $CC^*(\mathcal{A})$):
\begin{equation}
[e^p_D, b^p]\pm e^p_{\delta D}=0.  \label{eq:e^p_D,E_D^p basic property}  
\end{equation}
\begin{equation}
[E^p_D,b^p]\pm [e^p_D\tau^{-1},\tau]\pm E^p_{\delta D}=0.\label{eq:p Cartan}
\end{equation}
We define the operation $\widetilde{E^p_D\theta}$ on ${}_pCC_*(\mathcal{A})((t,\theta))$ to be the $t$-linear extension of
\begin{equation}
\begin{cases}
 x\mapsto E^p_Dx\theta\\
 x\theta\mapsto \sum_{0\leq i<j\leq p-1} \tau^j E^p_D\tau^{i-j} xt.
\end{cases}    
\end{equation}
Then, there is an additional relation 
\begin{equation}
[(\tau-1)\tilde{\theta},\widetilde{E^p_D\theta}]=\sum_{j=0}^{p-1}\tau^{j}E^p_D\tau^{-j} x t=:\mathcal{L}^p_Dx t.   \label{eq:additional relation}
\end{equation}
We are now ready to define the $\mathbb{Z}/p$-version of the Getzler-Gauss-Manin connection. 
$\nabla'$ induces a connection on $CC^{C_p,-}_*(\mathcal{A})$ as a graded $R$-module:
\begin{equation}
\nabla': CC^{C_p,-}_*(\mathcal{A})\rightarrow   CC^{C_p,-}_*(\mathcal{A})\otimes_R\Omega^1_R,
\end{equation}
which however is incompatible with the differential. More precisely, one has 
\begin{equation}
 [\nabla',b^p]=\mathcal{L}^p_{\kappa}\;,\;\;[\nabla',\tau]=0. \label{eq:p nabla' property}  
\end{equation}
Therefore, by \eqref{eq:e^p_D,E_D^p basic property}, \eqref{eq:p Cartan}, \eqref{eq:additional relation}, and \eqref{eq:p nabla' property}, we deduce that the chain level connection on $CC^{C_p,-}(\mathcal{A})$ given by
\begin{equation}
\nabla^{GGM,p}:=\nabla'+\frac{1}{t}(e^p_{\kappa}\tau^{-1}+\widetilde{E^p_{\kappa}\theta})   \label{eq:p GGM}
\end{equation}
commutes with the differential $b^p+(\tau-1)\tilde{\theta}$.\\
\begin{mydef}
The connection on $HH^{C_p,-}_*(\mathcal{A})$ (or $HH^{C_p,per}_*(\mathcal{A})$) induced by \eqref{eq:p GGM} is called the \emph{$C_p$-Getzler-Gauss-Manin connection}.   
\end{mydef}
At the level of cohomology, $\nabla^{GGM,p}$ is independent of the choice of $\nabla'$, and moreover a map of $A_{\infty}$-algebras induces a map on $HH^{C_p,per}_*$ that is compatible with $\nabla^{GGM,p}$. Both of those facts can either be checked by explicit but tedious formulae, or follow from an invariant definition of the connection given below.\par\indent
We imitate Section 3.1 and give a `Kaledin's re-formulation' of the $C_p$-Getzler-Gauss-Manin connection, using the finite $p$-cyclic category introduced in Section 2.4 b). \par\indent
Given any functor $E: {}_p\Lambda\rightarrow \mathrm{Mod}_{R}$, one can form 
an analogue of the negative (resp. periodic) $C_p$-equivariant Hochschild complex $CC^{C_p,-}_*(E)$ (resp. $CC^{C_p,per}_*(E)$), cf. \cite[Definition 3.17]{Che1}, such that when $E=j^*\mathcal{A}_{\sharp}$ for $\mathcal{A}$ a $\mathbb{Z}$ or $\mathbb{Z}/2$-graded dg algebra, it recovers \eqref{eq:C_N negative}. \par\indent
Now, we consider the two-step filtration \eqref{eq:first jet sequence} and its induced filtration on the cyclic object $({p_1}_*p_2^*\mathcal{A})_{\sharp}$. By an abuse of notation, we denote $I^ij^*({p_1}_*p_2^*\mathcal{A})_{\sharp}:=j^*I^i({p_1}_*p_2^*\mathcal{A})_{\sharp}$ the induced filtration on the finite $p$-cyclic object $j^*({p_1}_*p_2^*\mathcal{A})_{\sharp}$, and also use $I^i$ to denote the induced filtration on the associated periodic $C_p$-equivariant Hochschild complex.
There is an analogue of diagram \eqref{eq:GGM diagram} in this case (where now $I^i$ stands for $I^iCC^{C_p,per}_*({p_1}_*p_2^*\mathcal{A}_{\sharp})$):
\begin{equation}
\begin{tikzcd}[row sep=1.2cm, column sep=0.8cm]
I^1/I^2\arrow[r]& I^0/I^2\arrow[r,"\pi_p"] \arrow[d,"m_p"]&I^0/I^1=CC^{C_p,per}_*(\mathcal{A})\\
& {p_1}_*p_2^*CC^{C_p,per}_*(\mathcal{A})  &
\end{tikzcd}\label{eq:p GGM diagram}
\end{equation}
In \eqref{eq:p GGM diagram}, 
\begin{equation}
I^1CC^{C_p,per}_*({p_1}_*p_2^*\mathcal{A}_{\sharp})/I^2CC^{C_p,per}_*({p_1}_*p_2^*\mathcal{A}_{\sharp})=CC^{C_p,per}_*(j^*(I^1({p_1}_*p_2^*\mathcal{A}_{\sharp})/I^2({p_1}_*p_2^*\mathcal{A}_{\sharp})))
\end{equation}
is contractible by Proposition \autoref{thm:pGysin for CC} and the fact that $CC^{per}_*(I^1({p_1}_*p_2^*\mathcal{A}_{\sharp})/I^2({p_1}_*p_2^*\mathcal{A}_{\sharp}))$ is contractible, cf. \cite[Lemma 3.1]{PVV}. Thus, $\pi_p$ in \eqref{eq:p GGM diagram} is a quasi-isomorphism. We define $\nabla^{Kal, p}:=m_p\pi_p^{-1}$, which gives a Grothendieck connection on $CC^{tC_p}_*(\mathcal{A})$ in the derived category. \\
\begin{prop}
$\nabla^{Kal,p}$ agrees with $\nabla^{GGM, p}$ 
as a morphism $CC^{tC_p}_*(\mathcal{A})\rightarrow {p_1}_*p_2^*CC^{tC_p}_*(\mathcal{A})$ in the derived category.   
\end{prop}
\emph{Proof}. The proof is exactly parallel to that of Proposition \autoref{thm:Kaledin connection agrees with GGM}. First, extending the connection $\nabla'$ to tensor products via the Leibniz rule, then under the identification \eqref{eq:first jet} we obtain an $R$-linear section of $\pi_p$, denoted $\varphi: CC^{C_p,per}_*(\mathcal{A})\rightarrow I^0/I^2$. Similar to \eqref{eq:varphi basic property}, we have
\begin{equation}
[\varphi, b^p]=\tilde{\mathcal{L}}^p_{\kappa}\;,\;\;[\varphi, \tau]=0,\label{eq:p varphi basic property}
\end{equation}
where $\tilde{\mathcal{L}}^p_{\kappa}$ is the canonical lift of $\mathcal{L}^p_{\kappa}$ along $I^1/I^2\rightarrow CC^{C_p,per}_*(\mathcal{A})\otimes_{\mathcal{O}_X}\Omega^1_X$. There are also canonical lifts $\tilde{E}^p_{\kappa}, \tilde{e}^p_{\kappa}$ of $E^p_{\kappa}, e^p_{\kappa}$, respectively. The $C_p$-analogue of Cartan homotopy formula (cf. \eqref{eq:p Cartan}) then implies that $\varphi-\frac{1}{t}(\tilde{e}^p_{\kappa}+\widetilde{\tilde{E}^p_{\kappa}\theta})$ is a chain map that splits $\pi_p$. In particular, this implies that $\nabla^{GGM,p}=m_p\circ (\varphi-\frac{1}{t}(\tilde{e}^p_{\kappa}+\widetilde{\tilde{E}^p_{\kappa}\theta}))=m_p\circ \pi^{-1}_p=\nabla^{Kal,p}$.\qed

\subsection{Comparison of $\nabla^{GGM,p}$ and $\nabla^{GGM}$}
In this section we show that the comparison map of Proposition \autoref{thm:pGysin for CC} intertwines the Getzler-Gauss-Manin connection $\nabla^{GGM}$ with its $C_p$-variant $\nabla^{GGM,p}$ at the level of homology. \par\indent  
Indeed, apply the quasi-isomorphism $\phi_p$ of \eqref{eq:S^1 C_p comparison} to the cyclic object ${p_1}_*p_2^*\mathcal{A}_{\sharp}$, as well as the various associated graded (e.g. $I^0/I^2, I^0/I^1, I^1/I^2)$, we obtain an object-wise quasi-isomorphism
\begin{equation}
\phi_p: \textrm{diagram \eqref{eq:GGM diagram}}\otimes_{R((t))}R((t,\theta))\simeq \textrm{diagram \eqref{eq:p GGM diagram}}.\label{eq:Gysin diagram map}
\end{equation}
In fact, \eqref{eq:Gysin diagram map} is a map of diagrams (i.e. all the squares commute). To see this, for instance, applying naturality of $\phi_p$ to the map of cyclic chain complexes 
\begin{equation}
I^0({p_1}_*p_2^*\mathcal{A})_{\sharp}/I^2({p_1}_*p_2^*\mathcal{A})_{\sharp}\rightarrow   I^0({p_1}_*p_2^*\mathcal{A})_{\sharp}/I^1({p_1}_*p_2^*\mathcal{A})_{\sharp}\cong \mathcal{A}_{\sharp}
\end{equation}
shows that $\phi_p$ intertwines $\pi$ with $\pi_p$. Similarly, applying naturality of $\phi_p$ to the map of cyclic chain complexes
\begin{equation}
I^0({p_1}_*p_2^*\mathcal{A})_{\sharp}/I^2({p_1}_*p_2^*\mathcal{A})_{\sharp}\rightarrow {p_1}_*p_2^*\mathcal{A}_{\sharp} 
\end{equation}
induced by (where $[d]\in \Lambda, a_i\in \mathcal{A}$ and $r_i\in R_{[2]}:=R\otimes R/I^2$ is the first order neighborhood of the diagonal in $R\otimes R$)
\begin{equation}
[(a_0\otimes r_0)\otimes \cdots (a_d\otimes r_d)]\mapsto (a_0\otimes\cdots\otimes a_d)\otimes (r_0\cdots r_d),     
\end{equation}
one obtains that $\phi_p$ intertwines $m$ with $m_p$.\par\indent
Since $\nabla^{Kal}$ (resp. $\nabla^{Kal,p}$) is equal to $m\circ \pi^{-1}$ (resp. $m_p\circ \pi_p^{-1}$), cf. diagram \eqref{eq:GGM diagram} (resp. \eqref{eq:p GGM diagram}), the existence of the map of diagrams \eqref{eq:Gysin diagram map} immediately implies that\\
\begin{cor}
For a dg algebra $\mathcal{A}$ over $R$ over a characteristic $p$ field $\mathbf{k}$, the quasi-isomorphism $\phi_p: CC^{per}_*(\mathcal{A})\otimes_{R((t))}R((t,\theta))\rightarrow CC^{C_p,per}_*(\mathcal{A})$ intertwines the connections $\nabla^{Kal}\otimes_{R((t))}R((t,\theta))$ with $\nabla^{Kal,p}$. \qed
\end{cor}
Combining the previous results, we deduce that\\
\begin{cor}
For a dg algebra $\mathcal{A}$ over $R$ over a characteristic $p$ field $\mathbf{k}$, the quasi-isomorphism $\phi_p$ intertwines $\nabla^{GGM}\otimes_{R((t))}R((t,\theta))$ with $\nabla^{GGM,p}$ in the derived category. \qed
\end{cor}

\section{Kontsevich-Soibelman operations on the periodic cyclic homology}
The purpose of this section is to introduce a two-colored operad $\{\mathbf{KS}(l,0),\mathbf{KS}(k,1)\}_{l\geq 1,k
\geq 0}$ valued in chain complexes constructed by Kontsevich-Soibelman \cite{KS1}, which parametrizes `universal operations' on the pair $(CC^*(\mathcal{A}),CC_*(\mathcal{A}))$ of Hochschild (co)chain complexes of a dg algebra $\mathcal{A}$. In Section 4.1 and 4.2, we give explicit combinatorial descriptions of $\mathbf{KS}(n,0),\mathbf{KS}(n,1)$, respectively. However, instead of using the `minimal operad' in \cite{KS1}\cite{KS2}, we use a different combinatorial model known as \emph{cacti} (and its variations) first introduced combinatorially by \cite{MS} in their solution of Deligne's conjecture and subsequently developed in more geometric terms in \cite{Vo1} \cite{Kau} \cite{Sal1} \cite{RS}. In Section 4.3, we review the two colored operad of configuration spaces of disks on a disk/cylinder, which is a topological counterpart of $\mathbf{KS}$. Finally in Section 4.4, we study endomorphisms of the periodic cyclic homology of a dg algebra induced from the action of $\mathbf{KS}$, and in particular, show that they are covariantly constant with respect to the Getzler-Gauss-Manin connection.

\subsection{Cacti and their actions on Hochschild cochains}
Recall that an \emph{$n$-fold semi-simplicial (resp. semi-cosimplicial) object in a category $C$} is a covariant (resp. contravariant) functor from $(\overrightarrow{\Delta}^{op})^n$ (cf. Section 2.4(a)) to $C$. The \emph{geometric realization} of an $n$-fold semi-simplicial set $X$ is defined to be
\begin{equation}\label{eq:geometric realization of n fold semisimplicial set}
|X|:=\coprod_{a_1,\cdots,a_n\geq 0} X_{a_1,\cdots,a_n}\times \Delta^{a_1}\times\cdots\Delta^{a_n}/\sim,  
\end{equation}
where $\sim$ is the equivalence relation given by 
\begin{equation}
(f^*(x),y)\sim (x,f_*(y)),\;\mathrm{where}\;f\in \overrightarrow{\Delta}([a_1,\cdots,a_n],[b_1,\cdots,b_n]), x\in X_{b_1,\cdots,b_n}, y\in\prod_{i=1}^n\Delta^{a_i}.
\end{equation}
\begin{mydef}(\cite[Section 4.2]{MS},\cite[Definition 4.2]{Sal1})\label{thm:cacti}
$\mathfrak{Cact}^k=\mathfrak{Cact}^k_{\bullet,\cdots,\bullet}$ is the $k$-fold semi-simplicial set defined as follows. For $[m_1,\cdots,m_k]\in (\overrightarrow{\Delta}^{op})^k$, $\mathfrak{Cact}^k_{m_1,\cdots,m_k}$ is the set of surjective maps $f: \{1,\cdots,m+k\}\rightarrow \{1,\cdots,k\}$ with $m=\sum_{i=1}^k m_i$ such that
\begin{enumerate}[label=\arabic*)]
    \item $f^{-1}(j)$ has cardinality $m_j+1$ for each $1\leq j\leq k$.
    \item $f(a)\neq f(a+1)$ for $1\leq a<m+k$.
    \item There are no values $1\leq a<b<c<d\leq m+k$ such that $f(a)=f(c)\neq f(b)=f(d).$
\end{enumerate}
We will equivalently represent $f$ by the sequence $f(1)f(2)\cdots f(m+k)$. Fix $1\leq j\leq k$, and denote $f^{-1}(j)=\{a_0<a_1<\cdots<a_{m_j}\}$. Then the $i$-th ($0\leq i\leq m_j$) face map of the $j$-th component is defined by
\begin{equation}
d^{(j)}_i(f):=f(1)\cdots\widehat{f(a_i)}\cdots f(m+k): \{1,\cdots,m+k-1\}\rightarrow \{1,\cdots,k\}.
\end{equation}
The geometric realization of this $k$-fold semi-simplicial set $\mathrm{Cact}^k:=|\mathfrak{Cact}^k|$ is called the \emph{space of Cacti with $k$ lobes}. \par\indent
Note that $\mathrm{Cact}^k$ has a natural cellular structure induced from the (multi-)semi-simplicial structure of $\mathfrak{Cact}^k$. Then $C^{cell}_*(\mathrm{Cact}^k;R)$ agrees with the realization of $R\langle\mathfrak{Cact}^k\rangle$ in the category of $R$-chain complexes. \\
\end{mydef}
\begin{rmk} \label{thm:geometric interpretation of cacti}
$\mathrm{Cact}^k$ is called the space of cacti for the following reason. An element $(f,t^1,\cdots,t^k)\in \mathfrak{Cact}^k_{m_1,\cdots,m_k}\times \Delta^{m_1}\times\cdots\Delta^{m_k}$ may be viewed as a partition of $[0,k]$ into $m+k$ closed sub-intervals $[0,t_1], [t_1,t_1+t_2],\cdots, [\sum_{i=1}^{m+k-1}t_i,\sum_{i=1}^{m+k}t_i=k]$ such that if we denote $f^{-1}(j)=\{a_0<a_1<\cdots<a_{m_j}\}$ then $(t_{a_0},\cdots,t_{a_{m_j}})=t^j$, for each $1\leq j\leq k$. The $i$-th interval $[\sum_{l=1}^{i-1}t_l,\sum_{l=1}^it_l]$ is said to have `color' $f(i)\in\{1,\cdots,k\}$. Note the sum of lengths of intervals of each color is $1$. By identifying the endpoints of $[0,k]$, we view this as a partition of $S^1$. We define an equivalence relation on $S^1=[0,k]/0\sim k$ where $z\sim z'$ if $z$ and $z'$ are the boundary points of the same connected component of $S^1\backslash\mathrm{int}(I_j)$ for some $1\leq j\leq k$, where $I_j$ is the union of all intervals of color $j$. Condition 2) of Definition \autoref{thm:cacti} implies that the quotient under the prior equivalence relation is a connected union of $k$ circles of length $1$. This quotient is called the \emph{cactus associated with $(f,t^1,\cdots,t^k)$} (we will also abuse terminology and refer to $(f,t^1,\cdots,t^k)$ itself as a cactus); the image of $I_j$ is called \emph{the $j$-th lobe} of the cactus; the image of $0$ is \emph{the basepoint of the cactus} and the image of an endpoint of an interval $[\sum_{l=1}^{i-1}t_l,\sum_{l=1}^it_l]$ is called a \emph{marked point on the cactus}. See Figure \ref{fig:cactus} for an illustration.  
\begin{figure}[H]
 \centering
 \includegraphics[width=0.9\textwidth]{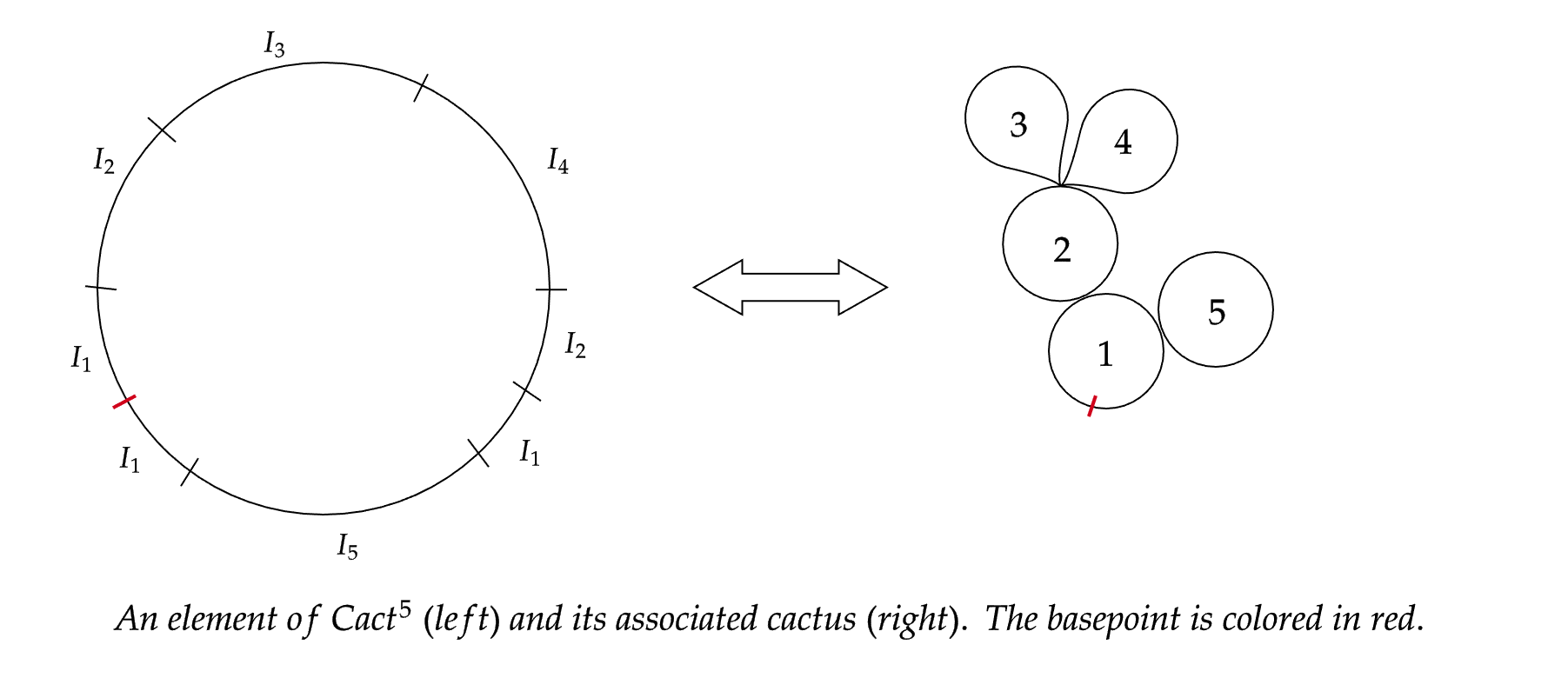}
 \caption{}
 \label{fig:cactus}
\end{figure}
\end{rmk}
Using the notations from Remark \autoref{thm:geometric interpretation of cacti}, \\
\begin{mydef}\label{thm:cactus map}
Let $x=(f,t^1,\cdots,t^k)\in \mathrm{Cact}^k$ be a cactus. The \emph{cactus map} is a piecewise (linear) oriented isometry
\begin{equation}\label{eq:cactus map}
c_x=(c_x^1,\cdots,c_x^k): [0,k]\rightarrow [0,1]^k   
\end{equation}
defined by
\begin{itemize}
    \item $c_x(0)=(0,\cdots,0)$.
    \item If $t\in [\sum_{l=1}^{i-1}t_l,\sum_{l=1}^{i}t_l]$, then for $1\leq j\leq k$, $c^j_x(t)=\begin{cases} c^j_x(\sum_{l=1}^{i-1}t_l)\,,\,\mathrm{if}\;j\neq f(i)\\
    c_x^j(\sum_{l=1}^{i-1}t_l)+(t-\sum_{l=1}^{i-1}t_l)\,,\,\mathrm{if}\;j=f(i)     
    \end{cases}$.
\end{itemize}
\end{mydef}
Intuitively, imagine a point going around the cactus clockwise with unit speed starting from its basepoint. Then as the point moves on the $j$-th lobe, the $j$-th coordinate of $c_x$ increases with unit speed, while all other coordinates remain constant. It is clear from this description that a cactus $x$ is uniquely determined by its cactus map $c_x$. \par\indent
\textbf{The operadic structure on cacti}. For positive integers $k,n_1,\cdots,n_k$, there are operadic structure maps
\begin{equation}\label{eq:operadic structure on cacti}
\theta_{k,n_1,\cdots,n_k}: \mathrm{Cact}^{n_1}\times\cdots\times\mathrm{Cact}^{n_k}\times \mathrm{Cact}^k\rightarrow \mathrm{Cact}^{n},\;\mathrm{where}\;n=\sum_{i=1}^kn_i.    
\end{equation}
In words, the operadic composition of $(x_1,\cdots,x_k,x)$ is obtained from $x$ by inserting the cactus $x_j$ into its $j$-th lobe (via a piecewise linear identification) in a way such that the basepoint of $x_j$ coincides with the \emph{local basepoint} of the $j$-th lobe of $x$. The local basepoint of a lobe is given by the basepoint if the lobe contains the basepoint; else it is the intersection of that lobe with the connected component of the closure of its complement that contains the basepoint. The precise definition of $\theta$ can be formulated in terms of cactus maps.\par\indent
Namely, there exists a unique piecewise linear map $\alpha: [0,k]\rightarrow [0,n]$ and a piecewise oriented isometry $c: [0,n]\rightarrow\prod_{j=1}^k[0,n_j]$ making the following diagram commute
\begin{center}
\begin{tikzcd}[row sep=1.2cm, column sep=0.8cm]
[0,k]\arrow[r,"{c_x}"]\arrow[d,"{\alpha}"]& {[0,1]}^k\arrow[d,"{(t_1,\cdots,t_k)\mapsto(n_1t_1,\cdots,n_kt_k)}"] \\
{[0,n]} \arrow[r,"c"] &\prod_{j=1}^k[0,n_j]
\end{tikzcd}\label{eq:cacti operadic composition diagram}
\end{center}\par
\begin{figure}[H]
 \centering
 \includegraphics[width=0.9\textwidth]{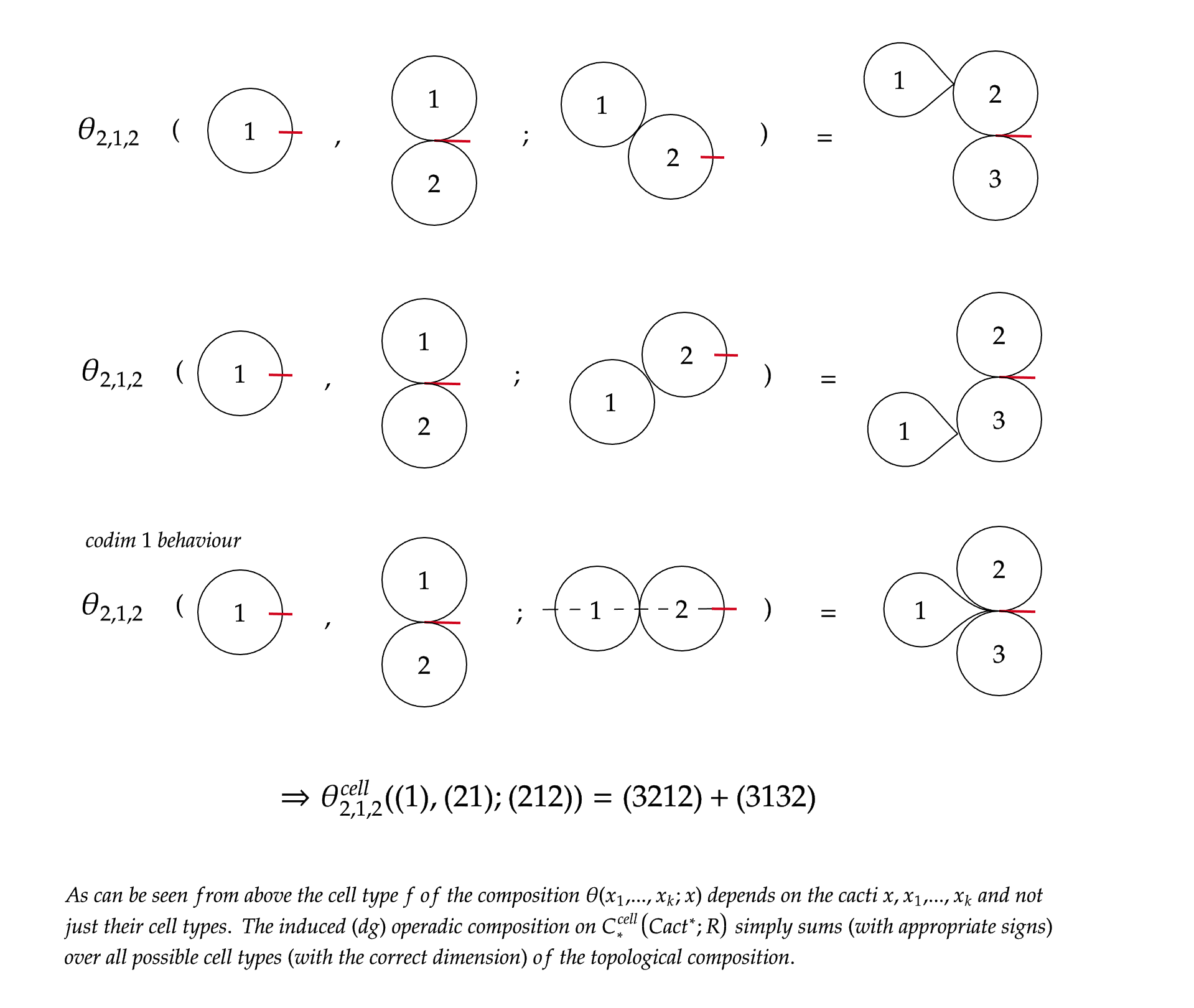}
 \caption{}
 \label{fig:cacti_composition}
\end{figure}

\indent
\begin{mydef}\label{thm:operadic structure on cacti}
$\theta(x_1,\cdots,x_k;x)$ is the cactus whose cactus map is
\begin{equation}\label{eq:cactus map of operadic composition}
(\prod_{j=1}^kc_{x_j})\circ c.    
\end{equation}
\end{mydef}
It turns out that $\theta$ is only associative up to homotopy.\\
\begin{thm}\cite[Theorem 4.10, Corollary 4.11]{Sal1}\label{thm:associativity of cacti operad}
$\theta$ defines an operadic structure up to homotopy on $\{\mathrm{Cact}^k\}_{k\geq 1}$. Furthermore, it induces a (strictly associative) dg operadic structure on $\{C^{cell}_*(\mathrm{Cact}^k;R)\}_{k\geq 1}$. 
\end{thm}
See Figure \ref{fig:cacti_composition} for an illustration of the operadic composition. \par\indent
There are two ways to turn cacti into a strictly associative operad. The first is to give `weights' to the lobes of a cactus and suitably multiply these weights when composing, cf. \cite[Section 4]{MS}). The spaces of weights $\mathcal{P}=\{\mathcal{P}(n)\}_{n\geq 1}$ themselves form a topological operad where $\mathcal{P}(n):=\{a_1+\cdots+a_n=1,a_i>0\}$ and
\begin{equation}\label{eq:operadic composition of weights}
(b_1,\cdots,b_m)\circ_i(a_1,\cdots,a_n):=(a_1,\cdots,a_{i-1},a_ib_1,\cdots,a_ib_m,a_{i+1},\cdots,a_n).    
\end{equation}
\begin{mydef}\label{thm:weighted cacti}
The \emph{space of weighted cacti with $k$ lobes} is $\widetilde{\mathrm{Cact}^k}:=\mathrm{Cact}^k\times \mathcal{P}(k)$.    
\end{mydef}
One thinks of $(x,\mathbf{a}=(a_1,\cdots,a_k))\in \widetilde{\mathrm{Cact}^k}$ as rescaling the cactus $x$ so that the length of its $i$-th lobe is $a_i$, and its total length is $1$. By rescaling Definition \autoref{thm:cactus map}, one can associate to a weighted cactus $(x,\mathbf{a})$ its cactus map (which uniquely determines it) of the form
\begin{equation}\label{eq:weighted cactus map}
c_{x,\mathbf{a}}:[0,1]\rightarrow \prod_{i=1}^k[0,a_i].  
\end{equation}
In terms of cactus maps, the operadic composition of weighted cacti \begin{equation}
\widetilde{\theta_j}:\widetilde{\mathrm{Cact}^l}\times \widetilde{\mathrm{Cact}^k}\rightarrow\widetilde{\mathrm{Cact}^{k+l-1}}    
\end{equation}is given by
\begin{equation}\label{eq:operadic composition of weighted cacti}
\widetilde{\theta_j}(c_{y,\mathbf{b}},   c_{x,\mathbf{a}}):=[0,1]\xrightarrow{c_{(x,\mathbf{a})}}\prod_{i=1}^k[0,a_i]\xrightarrow{\mathbbm{1}\times\cdots \times \overbrace{a_jc_{y,\mathbf{b}}a_j^{-1}}^{\textrm{jth entry}}\times\cdots\mathbbm{1}} [0,a_1]\times\cdots\times [0,a_{j-1}]\times\prod_{s=1}^l[0,a_jb_s]\times[0,a_{j+1}]\times\cdots\times[0,a_k],
\end{equation}
whose associativity is evident. Intuitively, $\circ_j$ rescales the first cactus by $a_j$ and then inserts it into the $j$-th lobe of the second cactus (in a way that matches up the (local) basepoints). 
\textbf{Cacti with spines and their operadic structure.} The second way to strictify the operadic composition of cacti is to thicken $\mathfrak{Cact}^k$ to a suitable cosimplicial $k$-fold simplicial set, called `cacti with spines', which form a strict operad in a suitable simplicial sense. Its topological realization, which turns out to be quasi-equivalent to the space of cacti (without spines), is then naturally a strict topological operad. This thickening has another advantage that it acts \emph{tautologically} on the Hochschild cochains of a dg algebra.  \\
\begin{rmk}\label{thm:cacti with spine terminology caveat}
Our definition below of `cacti with $n$ spines and $k$ lobes' recovers McClure-Smith's `Formula of type $(k,n)$' \cite[Section 5.1]{MS}. \emph{However, we immediately caveat that this notion of cacti with spines is different from that in \cite{Vo1} or \cite{Kau}; the latter is sometimes also called framed cacti \cite[Remark 4.7 (1)]{RS}.}\\
\end{rmk}
\begin{mydef}\label{thm:Cacti with spines as a cosimplicial multisimplicial set}
We define a cosimplicial k-fold simplicial set $\widehat{\mathfrak{Cact}^k}=\widehat{\mathfrak{Cact}^k_\bullet}_{;\bullet,\cdots,\bullet}: \Delta\times (\Delta^{op})^k\rightarrow \mathrm{Sets}$ as follows. Given non-negative integers $n,m_1,\cdots,m_k$, $\widehat{\mathfrak{Cact}^k_n}_{;m_1,\cdots,m_k}$ is the set of pairs $(f,S)$ consisting of 
\begin{itemize}
    \item A surjective map $f: \{1,\cdots,m+k\}\rightarrow \{1,\cdots,k\}$ with $m=\sum_{i=1}^k m_i$ satisfying conditions 1) and 3) of Definition \autoref{thm:cacti}.
    \item A map $S:\{0,1,\cdots,m+k\}\rightarrow \mathbb{Z}_{\geq 0}$ satisfying $\sum_{i=0}^{m+k}S(i)=n$.
\end{itemize}
We represent $f, S$ as $f(1)\cdots f(m+k), S(0)S(1)\cdots S(m+k)$, respectively. The cosimplicial and $k$-fold simplicial structure maps are defined as follows. Fix $1\leq j\leq k$, and denote $f^{-1}(j)=\{a_0<a_1<\cdots<a_{m_j}\}$. Then the $i$-th ($0\leq i\leq m_j$) face map of the $j$-th component is defined by
\begin{equation}\label{eq:face maps of thickened cacti}
\begin{cases}
d^{(j)}_i(f):= f(1)\cdots \widehat{f(a_i)}\cdots f(m+k)\\
d^{(j)}_i(S):= S(0)\cdots (S(a_{i}-1)+S(a_i))\cdots S(m+k)
\end{cases} 
\end{equation}
and the $i$-th ($0\leq i\leq m_j$) degeneracy map of the $j$-th component is defined by
\begin{equation}\label{eq:degeneracy maps of thickened cacti}
\begin{cases}
s^{(j)}_i(f):= f(1)\cdots f(a_i)f(a_i)\cdots f(m+k)\\
s^{(j)}_i(S):= S(0)\cdots 0S(a_i)\cdots S(m+k)
\end{cases} 
\end{equation}
For $0\leq i\leq n+1$, let $l_i:=\min\{q:S(0)+\cdots+S(q)\geq i\}$ if $i\leq n$ and $l_i:=m+k$ if $i=n+1$. The $i$-th ($0\leq i\leq n+1$) coface map is defined by
\begin{equation}\label{eq:coface maps of thickened cacti}
\begin{cases}
\sigma_i(f):=f\\
\sigma_i(S):=S(0)\cdots (S(l_i)+1)\cdots S(m+k)
\end{cases}  
\end{equation}
and the $i$-th ($0\leq i\leq n-1$) codegeneracy map is defined by
\begin{equation}\label{eq:codegeneracy maps of thickened cacti}
\begin{cases}
\epsilon_i(f):=f\\
\epsilon_i(S):=S(0)\cdots (S(l_{i+1})-1)\cdots S(m+k).
\end{cases}  
\end{equation}
Note the definition of $l_i$ guarantees that $S(l_i)>0$ if $i>0$, so \eqref{eq:codegeneracy maps of thickened cacti} makes sense. \par\indent
The multi-simplicial realization $\widehat{\mathrm{Cact}^k_n}:=|\widehat{\mathfrak{Cact}^k_n}_{;\bullet,\cdots,\bullet}|$ is called \emph{the space of cacti with $n$ spines and $k$ lobes}. It has a natural cellular structure with cells labeled by pairs $(f,S)$. We then denote the totalization of the associated cosimplicial space by $\widehat{\mathrm{Cact}^k}:=\mathrm{Tot}([n]\mapsto |\widehat{\mathrm{Cact}^k_n}|)$. 
\end{mydef}
\begin{figure}[H]
 \centering
 \includegraphics[width=0.9\textwidth]{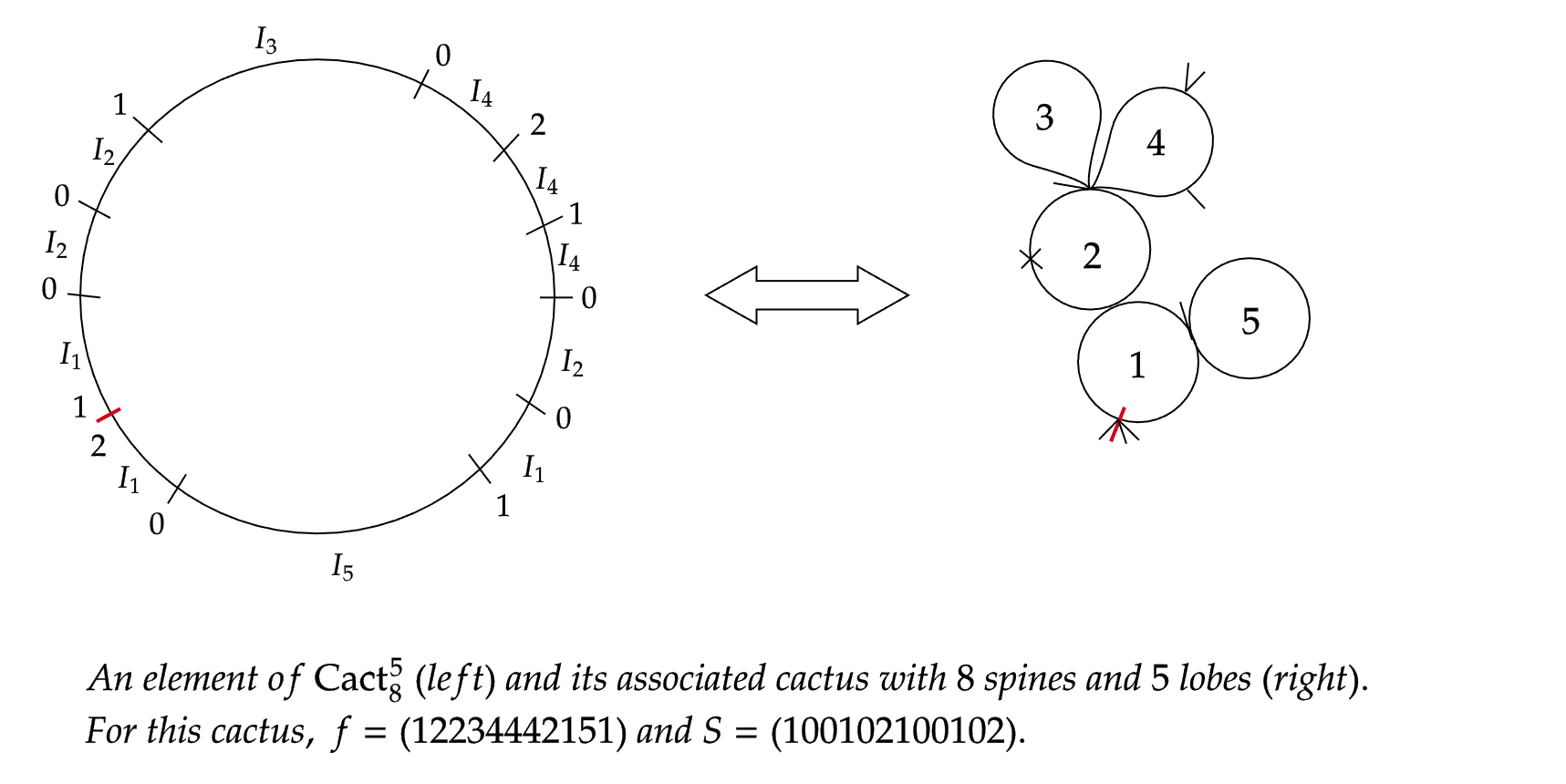}
 \caption{}
 \label{fig:cactus_with_spine}
\end{figure}
\begin{rmk}\label{thm:geometric interpretation of cacti with spines}
We explain the terminology `cactus with spines'. Following the procedure of Remark \autoref{thm:geometric interpretation of cacti}, to each element $((f,S),t^1,\cdots,t^k)\in \widehat{\mathfrak{Cact}^k_n}_{;m_1,\cdots,m_k}\times \Delta^{m_1}\times\cdots\times \Delta^{m_k}$ one associates a cactus with $k$ lobes (recall this is a quotient of $S^1$ homeomorphic to a connected union of circles), which is now decorated with some extra data. 
Namely, the map $S$ specifies a way to attach `spines' to the cactus, so that $S(i)$ spines are attached at the $i$-th marked point (when $S(i)=0$, the empty spine is attached, which is denoted by a cross if no lobe is based at $i$), see Figure \ref{fig:cactus_with_spine}.\par\indent
As a convention, we will draw a cactus such as the right hand side of Figure \ref{fig:cactus_with_spine} to represent its corresponding cell $(f,S)$. Schematically, the (co)face and (co)degeneracy maps of $\widehat{\mathfrak{Cact}^k}$ may be interpreted as in Figure \ref{fig:cactus_with_spine_co_face_and_co_degeneracy}. 
\begin{figure}[H]
 \centering
 \includegraphics[width=1.0\textwidth]{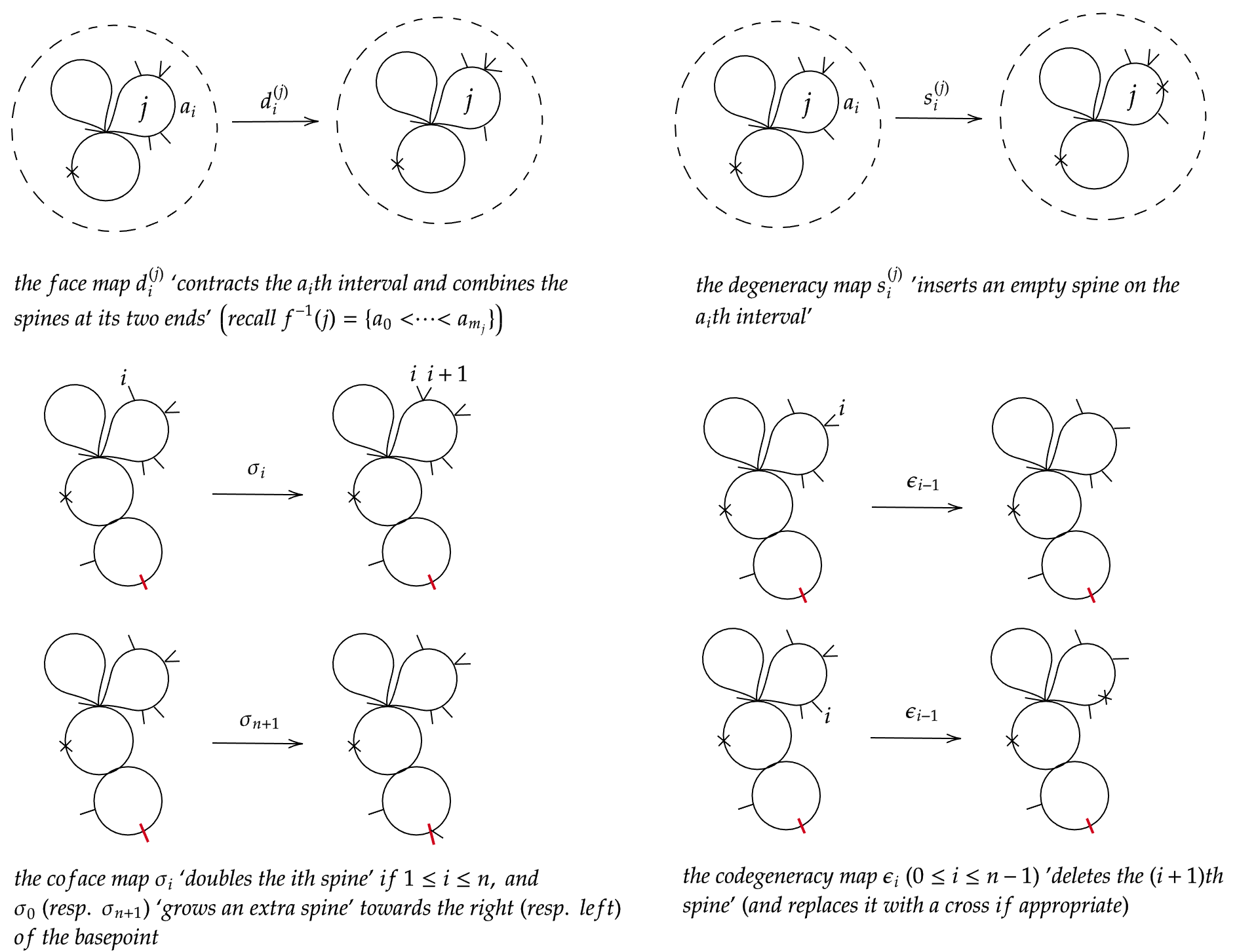}
 \caption{}
 \label{fig:cactus_with_spine_co_face_and_co_degeneracy}
\end{figure}
\end{rmk}
The sequence of cosimplicial multi-simplicial sets $\{\widehat{\mathfrak{Cact}^k}\}_{k\geq 1}$ has a (strictly associative) operadic structure which we now describe.\par\indent
First we recall some basic facts about operads in the cosimplicial multi-simplicial context. Given a sequence $\{X^k\}_{k\geq 1}$, where $X^k:=X^k_{\bullet;\bullet,\cdots,\bullet}: \Delta\times(\Delta^{op})^k\rightarrow \mathrm{Sets}$, the appropriate notion of operadic structure maps is a collection of natural transformations
\begin{equation}\label{eq:operadic structure for cosimplicial multisimplcial sets}
\gamma_{k,n_1,\cdots,n_k}: (X^{n_1}\times\cdots\times X^{n_k})\times_{(\Delta^{op})^k}X^k\rightarrow X^{n},\quad (n=\sum_{i=1}^k n_i)
\end{equation}
where both sides of \eqref{eq:operadic structure for cosimplicial multisimplcial sets} are viewed as functors $\Delta\times (\Delta^{op})^n\rightarrow \mathrm{Sets}$, that satisfies the usual operadic associativity (and unitality, if applicable) properties. We explain the notation in \eqref{eq:operadic structure for cosimplicial multisimplcial sets}: $(\Delta^{op})^k$ acts on the left of $X^k$ through its $k$-fold simplicial structure, and on the right of $X^{n_1}\times\cdots\times X^{n_k}$ through the cosimplicial structures on $X^{n_i}$ ($1\leq i\leq k$); $-\times_{(\Delta^{op})^k}-$ denotes the coend (explicitly, given $X\,(\mathrm{resp.}\;Y): \Delta\,(\mathrm{resp.}\;\Delta^{op})\rightarrow \mathrm{Sets}$, $X\times_{\Delta^{op}}Y=\coprod_n X^n\times Y_n/(x,f_*(y))\sim (f^*(x),y)$ for $f\in \Delta([m],[n])$), and in particular the left hand side of \eqref{eq:operadic structure for cosimplicial multisimplcial sets} has a residual (left) $\Delta\times(\Delta^{op})^n$-action.\par\indent
Given an operadic structure $\gamma$ on $\{X^k\}_{k\geq 1}$ as in \eqref{eq:operadic structure for cosimplicial multisimplcial sets}, by taking totalization-realization in the category of topological spaces (resp. $R$-chain complexes) one obtains \\
\begin{lemma}\label{thm:induced operadic structure on totalization-realization}
$\gamma$ induces the structure of a topological (resp. dg) operad on the sequence
\begin{equation}\label{eq:totalization realization operad}
\{\mathrm{Tot}([n]\mapsto |X^k_{n;\bullet,\cdots,\bullet}|)\}_{k\geq 1}  \qquad (\mathrm{resp.}\;  \{\mathrm{Tot}_R([n]\mapsto |X^k_{n;\bullet,\cdots,\bullet}|_R)\}_{k\geq 1}.
\end{equation}\qed
\end{lemma}
In \eqref{eq:totalization realization operad}, $\mathrm{Tot}_R$ and $|\cdot|_R$ denote the totalization and realization in the category of $R$-chain complexes. Similar to their topological counterparts, they can be defined as abstract right and left Kan extensions, or be computed via canonical chain models (cf. Section 2.4 a)-c)). \par\indent
We now define the operadic structure maps 
\begin{equation}\label{eq:operadic structure on cacti with spine}
\widehat{\theta}_{k,n_1,\cdots,n_k}: (\widehat{\mathfrak{Cact}^{n_1}}\times\cdots\times \widehat{\mathfrak{Cact}^{n_k}})\times_{(\Delta^{op})^k} \widehat{\mathfrak{Cact}^k}\rightarrow\widehat{\mathfrak{Cact}^n}\quad(n=\sum_{i=1}^k n_i). 
\end{equation}
Via the usual operadic yoga, this is equivalent to specifying natural transformations
\begin{equation}\label{eq:operadic structure on cacti with spine version 2}
\widehat{\theta}_{j}:=\widehat{\theta}(\mathbbm{1},\cdots,\mathbbm{1},\overbrace{-}^{j\mathrm{th}\;\mathrm{entry}},\mathbbm{1},\cdots,\mathbbm{1};-):\widehat{\mathfrak{Cact}^l}\times_{j,\Delta^{op}}\widehat{\mathfrak{Cact}^k}\rightarrow \widehat{\mathfrak{Cact}^{k+l-1}}\quad(1\leq j\leq k)
\end{equation}
between functors from $\Delta\times(\Delta^{op})^{k+l-1}$ to $\mathrm{Sets}$, where $\times_{j,\Delta^{op}}$ means that $\Delta^{op}$ acts on the left of $\widehat{\mathfrak{Cact}^k}$ via the $j$-th component of its $k$-fold simplicial structure. The unit $\mathbbm{1}$ represents the cactus with $1$ lobe.\\
\begin{mydef}\label{thm:pictorial composition of cacti with spines}
Given $x_1\in \widehat{\mathfrak{Cact}^{l}}_{m_j;m_1',\cdots,m'_{l}}$ and $x_2\in \widehat{\mathfrak{Cact}^k_n}_{;m_1,\cdots,m_k}$, their composition $\widehat{\theta}_j(x_1,x_2)\in \widehat{\mathfrak{Cact}}^{\,k+l-1}_{n;m_1,\cdots,m_{j-1},m_1',\cdots,m_l',m_j+1,\cdots,m_k}$ is the cactus with spines obtained by inserting $x_1$ into the $j$-th lobe of $x_2$ so that 1) their (local) basepoints match and 2) the spines of $x_1$ precisely match up the marked points on the $j$-th lobe of $x_2$, see Figure \ref{fig:composition_of_cacti_with_spines}. Note that $\widehat{\theta}_j(x_1,x_2)$ has the same number of spines as $x_2$.
\begin{figure}[H]
 \centering
 \includegraphics[width=1.0\textwidth]{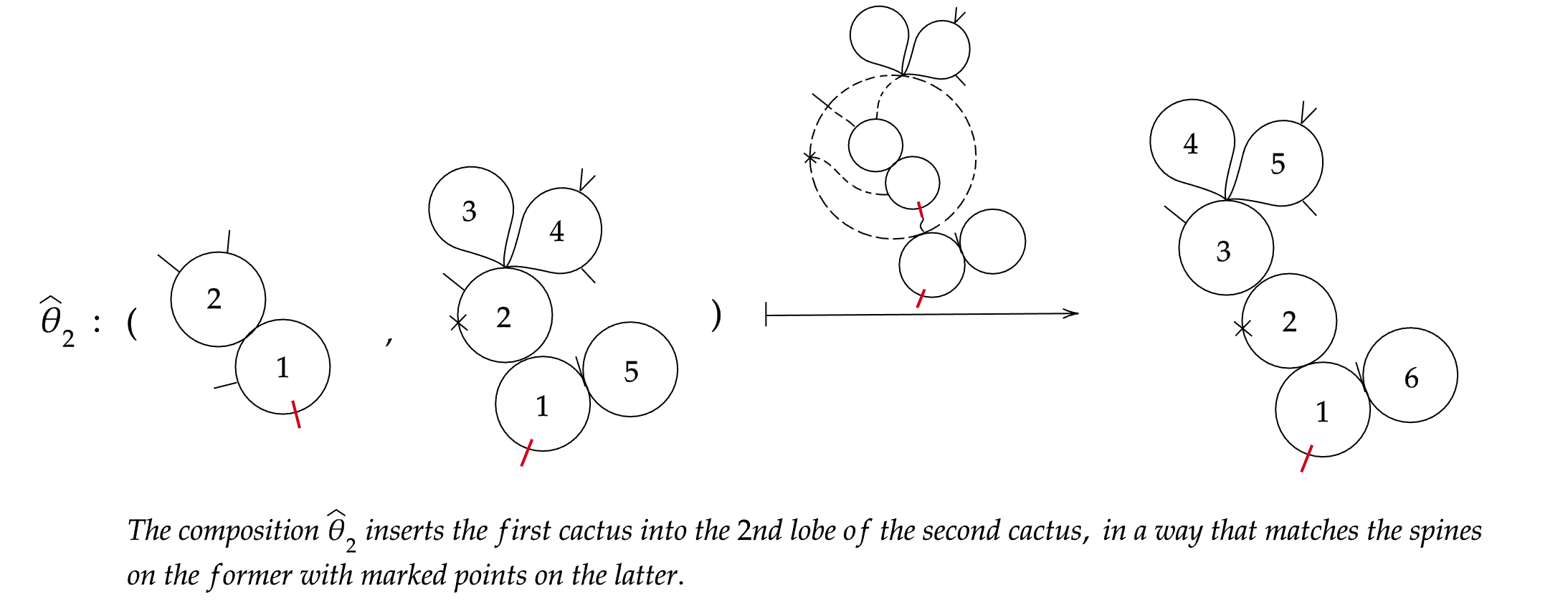}
 \caption{}
 \label{fig:composition_of_cacti_with_spines}
\end{figure}
\end{mydef}
From Figure \ref{fig:cactus_with_spine_co_face_and_co_degeneracy} and \ref{fig:composition_of_cacti_with_spines}, it is evident that $\widehat{\theta}_j$ \eqref{eq:operadic structure on cacti with spine version 2} is well defined (i.e. factors through the coend $-\times_{j,\Delta}-$) defines a map of cosimplicial $k+l-1$-fold simplicial sets. As a result of Lemma \autoref{thm:induced operadic structure on totalization-realization}, one obtains the structure of a topological (resp. dg) operad on the sequence
\begin{equation}\label{eq:totalization realization operad for cacti}
\{\widehat{\mathrm{Cact}^k}=\mathrm{Tot}([n]\mapsto|\widehat{\mathfrak{Cact}^k_n}_{;\bullet,\cdots,\bullet}|)\}_{k\geq 1} \qquad (\mathrm{resp.}\;\{\widehat{\mathrm{Cact}^k_R}=\mathrm{Tot}_R([n]\mapsto|\widehat{\mathfrak{Cact}^k_n}_{;\bullet,\cdots,\bullet}|_R)\}_{k\geq 1}).
\end{equation}
\begin{mydef}\label{thm:definition of KS(n,0)_R}
Define $\{\mathbf{KS}(k,0)\}_{k\geq 1}$ (resp. $\{\mathbf{KS}(n,0)_R\}_{k\geq 1}$) to be the topological (resp. dg) operad of \eqref{eq:totalization realization operad for cacti}.  \\ 
\end{mydef}
\begin{prop}\label{thm:KS equals weighted cacti}
There is a quasi-equivalence of topological operads 
\begin{equation}\label{eq:KS equals weighted cacti}
\{\mathbf{KS}(k,0)\}_{k\geq 1}\simeq \{\widetilde{\mathrm{Cact}^k}\}_{k\geq 1}.
\end{equation}
\end{prop}
Proposition \autoref{thm:KS equals weighted cacti} follows by embedding both sides of \eqref{eq:KS equals weighted cacti} into the coendomorphism operad of $S^1$, which we now recall.\\
\begin{mydef}\label{thm:coendomorphism operad}
Let $X$ be a topological space. The \emph{two-colored coendomorphism operad of $X$} is the two-colored operad $\{\mathrm{CoEnd}_X(n,0),\mathrm{CoEnd}_X(m,1)\}_{n\geq 1,m\geq 0}$ defined by 
\begin{equation}
\mathrm{CoEnd}_X(n,0):=\mathrm{Map}(X,X^{\times n})\quad,\quad\mathrm{CoEnd}_X(m,1):=\mathrm{Map}(X,X\times X^{\times m});   
\end{equation}
let $f=(f_1,\cdots,f_k)\in\mathrm{CoEnd}_X(k,0), g=(g_1,\cdots,g_l)\in\mathrm{CoEnd}_X(l,0), F=(F_0,F_1,\cdots,F_m)\in\mathrm{CoEnd}_X(m,1), G=(G_0,G_1,\cdots,G_n)\in\mathrm{CoEnd}_X(n,1)$, then the operadic structure maps are given by
\begin{equation}
g\circ_j f:= (f_1,\cdots,f_{j-1},g_1\circ f_j,\cdots,g_l\circ f_j,f_{j+1},\cdots,f_k),  \quad 1\leq j\leq k,
\end{equation}
\begin{equation}
g\circ_j^0 F:= (F_0,F_1,\cdots,F_{j-1},g_1\circ F_j,\cdots,g_l\circ F_j,F_{j+1},\cdots,F_m),\quad 1\leq j\leq m 
\end{equation}
and
\begin{equation}
G\circ^1 F:= (F_0\circ G_0,F_1\circ G_0,\cdots,F_{m}\circ G_0,G_1,\cdots,G_n).
\end{equation}
The \emph{coendomorphism operad of $X$} is the suboperad $\{\mathrm{CoEnd}_X(n,0)\}_{n\geq 1}$.
\end{mydef}
\emph{Proof of Proposition \autoref{thm:KS equals weighted cacti}.}  The simple but key observation is that there is a homeomorphism $g=(g_1,g_2):\widehat{\mathrm{Cact}^k_n}\xrightarrow{\cong}\mathrm{Cact}^k\times \Delta^n$. Namely, $g_1$ records the underlying cactus without spine, and $g_2$ records $1/k$ of the lengths of the intervals $(0,y_1),(y_1,y_2),\cdots,(y_n,k)$, where $y_i$ denotes the position of the $i$-th spine. Moreover, these homeomorphisms assemble into an isomorphism of cosimplicial spaces $\widehat{\mathrm{Cact}^k_{\bullet}}\cong  \mathrm{Cact}^k\times \Delta^{\bullet}$. Therefore,
\begin{equation}\label{eq:cacti with spine as cacti times Mon(I,partial I)}
\mathbf{KS}(k,0)=\mathrm{Tot}(\widehat{\mathrm{Cact}^k_{\bullet}})\cong  \mathrm{Cact}^k\times \mathrm{Tot}(\Delta^{\bullet})\cong \mathrm{Cact}^k\times \mathrm{Mon}(I,\partial I).   
\end{equation}
Given $x\in \mathrm{Cact}^k$, we view its cactus map $c_x$ as an element of $\mathrm{CoEnd}_{S^1}(k,0)=\mathrm{Map}(S^1,(S^1)^{\times k})$ by rescaling each interval in \eqref{eq:cactus map} to have length $1$ and identifying its endpoints. Given $(x,f)\in \mathrm{Cact}^k\times \mathrm{Mon}(I,\partial I)$, define its \emph{corrected cactus map} by $c_{x,f}:=c_x\circ f\in \mathrm{CoEnd}_{S^1}(k,0)$. By a straightforward computation, the map
\begin{equation}
(x,f)\mapsto c_{x,f}   
\end{equation}
induces an embedding of operads $\{\mathbf{KS}(k,0)\}_{k\geq 1}\hookrightarrow\{\mathrm{CoEnd}_{S^1}(k,0)\}_{k\geq 1}$. Moreover, by \cite[Proposition 4.5]{Sal2}, its image is given by the space of \emph{based} maps $g=(g_1,\cdots,g_k):S^1\rightarrow (S^1)^{\times k}$ such that
\begin{enumerate}[label=\arabic*)]
    \item Each $g_i$ is weakly monotonically increasing of degree $1$;
    \item there is a partition of $S^1$ into closed intervals
such that on each interval $g_j$ is constant for all except one special index $j=i$;
 \item the clockwise sequence of special indices does not contain a subsequence of the form $\{i, j, i, j\}$ with $i\neq j$.
\end{enumerate}
On the other hand, sending $(x,\mathbf{a})\in \widetilde{\mathrm{Cact}^k}$ to its (rescaled) cactus map $c_{x,\mathbf{a}}$ gives an embedding of operads $\{\widetilde{\mathrm{Cact}^k}\}_{k\geq 1}\hookrightarrow\{\mathrm{CoEnd}_{S^1}(k,0)\}_{k\geq 1}$ whose image is the space of maps $g:S^1\rightarrow (S^1)^{\times k}$ satisfying Conditions 1) and 3) above together with 
\begin{enumerate}[label=\arabic*')]
 \setcounter{enumi}{1}
    \item there exists a multi-set $\{a_1,\cdots,a_k\}$ with $a_j>0,\sum_{j=1}^k a_j=1$ and a partition of $S^1$ into closed intervals
such that on each interval $g_j$ is constant for all except one special index $j=i$ \emph{and} $g_i$ is linear of slope $1/a_i$.  
\end{enumerate}
Clearly 2') is stronger than 2), and thus there is an embedding of operads $\{\widetilde{\mathrm{Cact}^k}\}_{k\geq 1}\hookrightarrow\{\mathbf{KS}(k,0)\}_{k\geq 1}$. Moreover, for each $k\geq1$ the embedding fits into a commutative diagram 
\begin{equation}
\begin{tikzcd}[row sep=1.2cm, column sep=0.8cm]
\widetilde{\mathrm{Cact}^k}\arrow[hookrightarrow,rr]\arrow[dr]& &\mathbf{KS}(k,0)\arrow[dl]\\
 & \mathrm{Cact}^k  &
\end{tikzcd}
\end{equation}
Since both downward arrows are homotopy equivalences, the result follows.\qed\par\indent
\textbf{Action on Hochschild cochain}. Let $\mathcal{A}$ be a dg algebra over $R$ and $\mathcal{A}^{\sharp}$ be the cosimplicial chain complex from \eqref{eq:HH cochain functor} whose totalization is the Hochschild cochain complex of $\mathcal{A}$. \\
\begin{mydef}\label{thm: action of KS(n,0)_R on Hochschild cochains}
Fix $x=(f,S)\in \widehat{\mathfrak{Cact}^k_n}_{;m_1,\cdots,m_k}$, $\varphi_j\in\mathcal{A}^{\sharp}([m_j])=\mathrm{Hom}_R(\mathcal{A}^{\otimes m_j},\mathcal{A}), 1\leq j\leq k$, we define 
\begin{equation}
\varphi:=\mathrm{Act}_{x}(\varphi_1,\cdots,\varphi_k) \in \mathrm{Hom}_R(\mathcal{A}^{\otimes n},\mathcal{A})   
\end{equation}
to be the following multi-linear map. Given $a_1\otimes\cdots\otimes a_n\in\mathcal{A}^{\otimes n}$, attach $a_i$ to the $i$-th spine of $x$ (and attach the unit to an empty spine), and insert $\varphi_j$ into the $j$-th lobe of $x$. The cactus with spines $x$ specifies a way of successively applying $\varphi_j$ or multiplication in $\mathcal{A}$ to the tuple $(a_1,\cdots,a_n)$, and $\varphi(a_1,\cdots,a_n)$ is defined to be the output at the basepoint adjusted by the Koszul sign. See Figure \ref{fig:action_of_cacti_with_spines_on_HH_cochain} for an illustration.
\begin{figure}[H]
 \centering
 \includegraphics[width=1.0\textwidth]{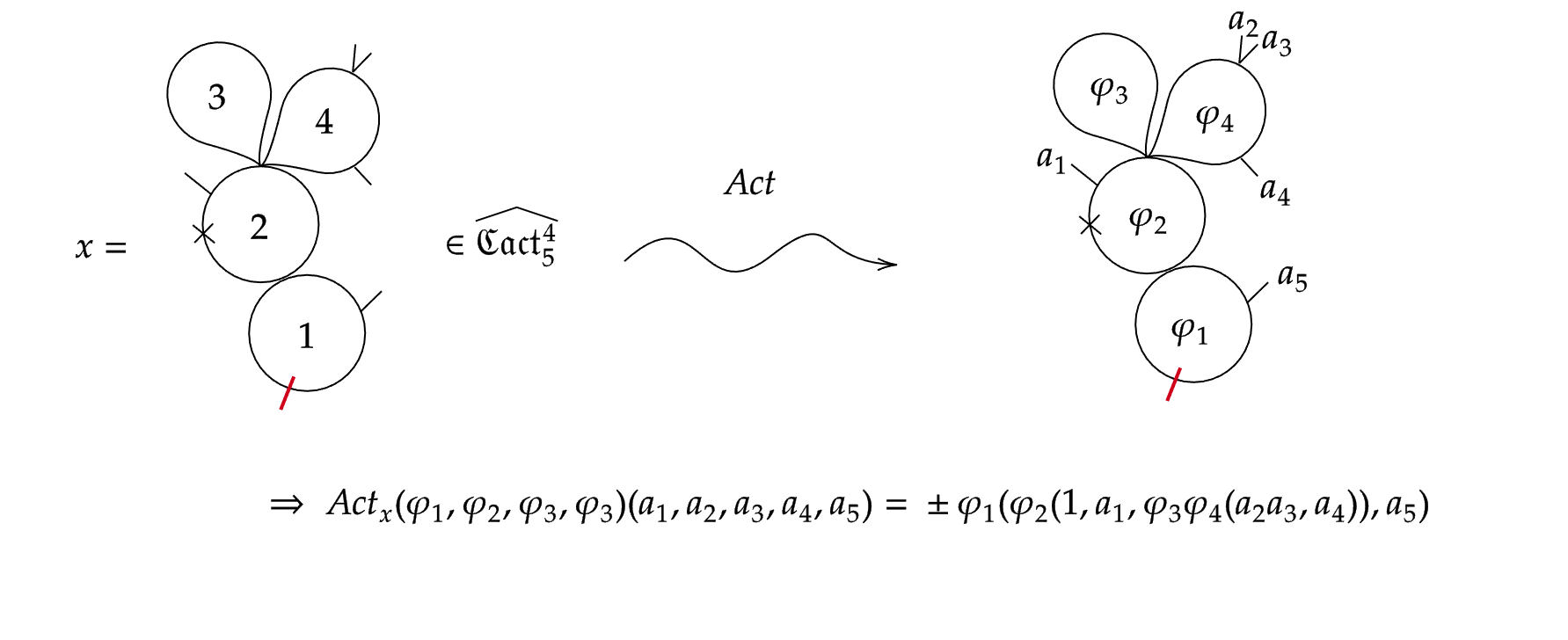}
 \caption{}
 \label{fig:action_of_cacti_with_spines_on_HH_cochain}
\end{figure}
\end{mydef}
The following lemma is immediate from Definition \autoref{thm: action of KS(n,0)_R on Hochschild cochains}. \\
\begin{lemma}
$\mathrm{Act}: \widehat{\mathfrak{Cact}^k}\rightarrow\mathrm{Hom}_R((\mathcal{A}^{\sharp})^{\otimes k},\mathcal{A}^{\sharp})$ defines a map of cosimplicial $k$-fold simplicial sets. Moreover, it intertwines the operadic structure $\widehat{\theta}$ with the operadic structure on $\{\mathrm{Hom}_R((\mathcal{A}^{\sharp})^{\otimes k},\mathcal{A}^{\sharp})\}_{k\geq 1}$ given by insertion of multi-linear maps. \qed  
\end{lemma}
Taking totalization-realization in the category of $R$-chain complexes gives rise to a map of dg operads
\begin{equation}\label{eq:KS(n,0)_R action on Hocschild cochains, after taking totalization}    
\mathrm{Act}: (\{\mathbf{KS}(n,0)_R\}_{n\geq 1}\rightarrow \{\mathrm{Hom}_R(CC^*(\mathcal{A})^{\otimes n},CC^*(\mathcal{A}))\}_{n\geq 1}.
\end{equation}
In view of the quasi-equivalence of Lemma \autoref{thm:KS equals weighted cacti}, this gives an action of the singular chains on the space of weighted cacti on $CC_*(\mathcal{A})$. It was shown in e.g. \cite{MS} that $\{\widetilde{\mathrm{Cact}^k}\}_{k\geq 1}$ is weakly equivalent to the operad of little disks, from which they were able to conclude Deligne's conjecture. 

\subsection{Cyclic cacti and $\mathbf{KS}(n,1)$}
In Section 4.1, we studied the suboperad $\{\mathbf{KS}(n,0)_R\}_{n\geq 1}$ of Kontsevich-Soibelman's two colored operad $\{\mathbf{KS}(n,0)_R,\mathbf{KS}(m,1)_R\}_{n\geq 1,m\geq 0}$ and its action on $CC^*(\mathcal{A})$. In this subsection, we use a variation of cacti to model the other part $\mathbf{KS}(n,1)_R$ (and its topological version $\mathbf{KS}(n,1)$), and construct the full action of the two colored operad on the pair $(CC^*(\mathcal{A}),CC_*(\mathcal{A}))$. In comparison with \cite{KS1}'s original definition via `the minimal operad', our cactus model makes it easier to study circle actions (and its finite cyclic subgroup actions), which will play an important role in the classification of equivariant Kontsevich-Soibelman operations in Section 5. \par\indent
For $k\geq 0$, consider the space of cacti with $k+1$ lobes $\mathrm{Cact}^{k+1}$, with the lobes labeled $0,1,\cdots,k$. As in Remark \autoref{thm:geometric interpretation of cacti}, each cactus is uniquely determined by a partition of $S^1\cong \mathbb{R}/(k+1)\mathbb{Z}$ into closed $1$-dimensional submanifolds $I_0,\cdots,I_k$ (recall $I_j$ is the union of intervals with color $j$).\\
\begin{mydef}\label{thm:space of cyclic cacti}
The \emph{space of cyclic cacti with $k$ lobes} $\mathrm{Cact}^k_{\circlearrowright}$ is the space $\mathrm{Cact}^{k+1}\times S^1$. For $k\geq 1$, we equip it with the $S^1\times S^1$-action given by the product of the $S^1$-action on $\mathrm{Cact}^{k+1}$ that clockwisely rotates the basepoint along the periphery of the cactus, i.e.  
\begin{equation}
t\cdot(I_0,\cdots,I_k)= (I_0+t,\cdots,I_k+t) \quad (t\in S^1\cong \mathbb{R}/(k+1)\mathbb{Z})
\end{equation}
with the regular action of $S^1$. 
Pictorially, we think of the second $S^1$ factor of $\mathrm{Cact}^{k}_{\circlearrowright}$ as an extra marked point moving along the $0$-th lobe; this requires a choice of a reference point on the $0$-th lobe, which we take to be the intersection of the $0$-th lobe with the connected component of its complement closure that contains the $1$-st lobe. We also call the $0$-th lobe the \emph{base lobe}. The basepoint of the underlying (non-cyclic) cactus $x$ will be called the \emph{input basepoint}; the extra marked point on the base lobe will be called the \emph{output basepoint}; the \emph{adjusted local basepoint} of the $j$-th lobe ($1\leq j\leq k$) is its intersection with the connected component of its complement closure that contains the base lobe, and if $j=0$ it is the output basepoint. See Figure \ref{fig:cyclic_cacti} for an example of a cyclic cactus.\par\indent
For $k=0$, we think of $S^1=\mathrm{Cact}^0_{\circlearrowright}$ as measuring the clockwise angle from the input basepoint to the output basepoint; it is equipped with the $S^1\times S^1$-action via the quotient $S^1\times S^1/\mathrm{diag}\cong S^1$.  
\begin{figure}[H]
 \centering
 \includegraphics[width=0.9\textwidth]{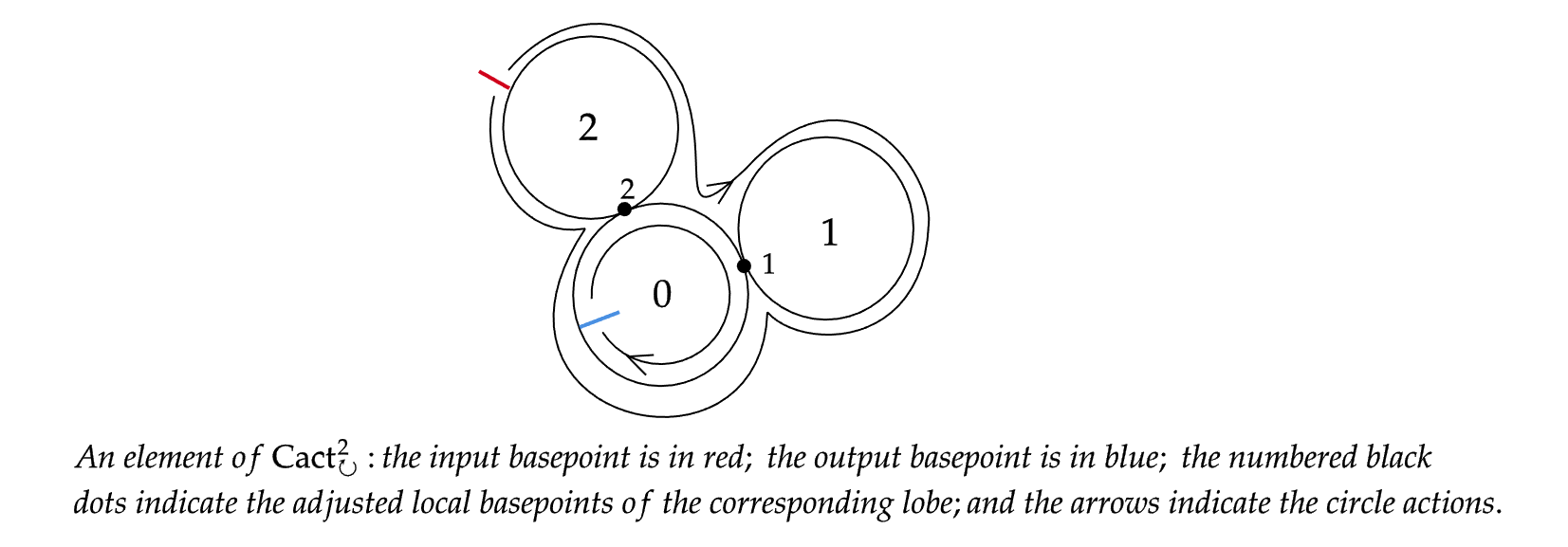}
 \caption{}
 \label{fig:cyclic_cacti}
\end{figure}
\end{mydef}
The `cactus map' of a cyclic cactus $(x,\theta)\in \mathrm{Cact}^k_{\circlearrowright}$ is defined to be the pair $(c_x,\theta)$, where $c_x$ is as in Definition \autoref{thm:cactus map}.\par\indent
\textbf{The operadic structure on cyclic cacti}. We define structure maps
\begin{equation}\label{eq:cyclic cacti operadic structure type 0}
\theta^{\circlearrowright,0}_j: \mathrm{Cact}^l\times \mathrm{Cact}^k_{\circlearrowright}\rightarrow \mathrm{Cact}^{k+l-1}_{\circlearrowright}, \quad1\leq j\leq k    
\end{equation}
and 
\begin{equation}\label{eq:cyclic cacti operadic structure type 1}
\theta^{\circlearrowright,1}: \mathrm{Cact}^l_{\circlearrowright}\times \mathrm{Cact}^k_{\circlearrowright}\rightarrow \mathrm{Cact}^{k+l}_{\circlearrowright} 
\end{equation}
such that when combined with $\theta$ from \eqref{eq:operadic structure on cacti} makes $\{\mathrm{Cact}^l,\mathrm{Cact}^k_{\circlearrowright}\}_{l,k\geq 1}$ into a two colored topological operad up to homotopy. \\
\begin{mydef}\label{thm:operadic structure on cyclic cacti}
\begin{itemize}
    \item $\theta^{\circlearrowright,0}_j$ inserts the first cactus into the $j$-th lobe of the second cyclic cactus (upon applying a piecewise linear oriented isometry) so that the basepoint of the former matches the adjusted local basepoint of the latter. 
    \item $\theta^{\circlearrowright,1}$ inserts the second cyclic cactus into the base lobe of the first cyclic cactus (upon applying a piecewise linear oriented isometry) so that the input basepoint of the second cyclic cactus matches the output basepoint of the first cyclic cactus. The base lobe of the second cyclic cactus will become the base lobe of the composition. 
\end{itemize}
See Figure \ref{fig:composition_of_cyclic_cacti} for an illustration. 
\begin{figure}[H]
 \centering
 \includegraphics[width=0.9\textwidth]{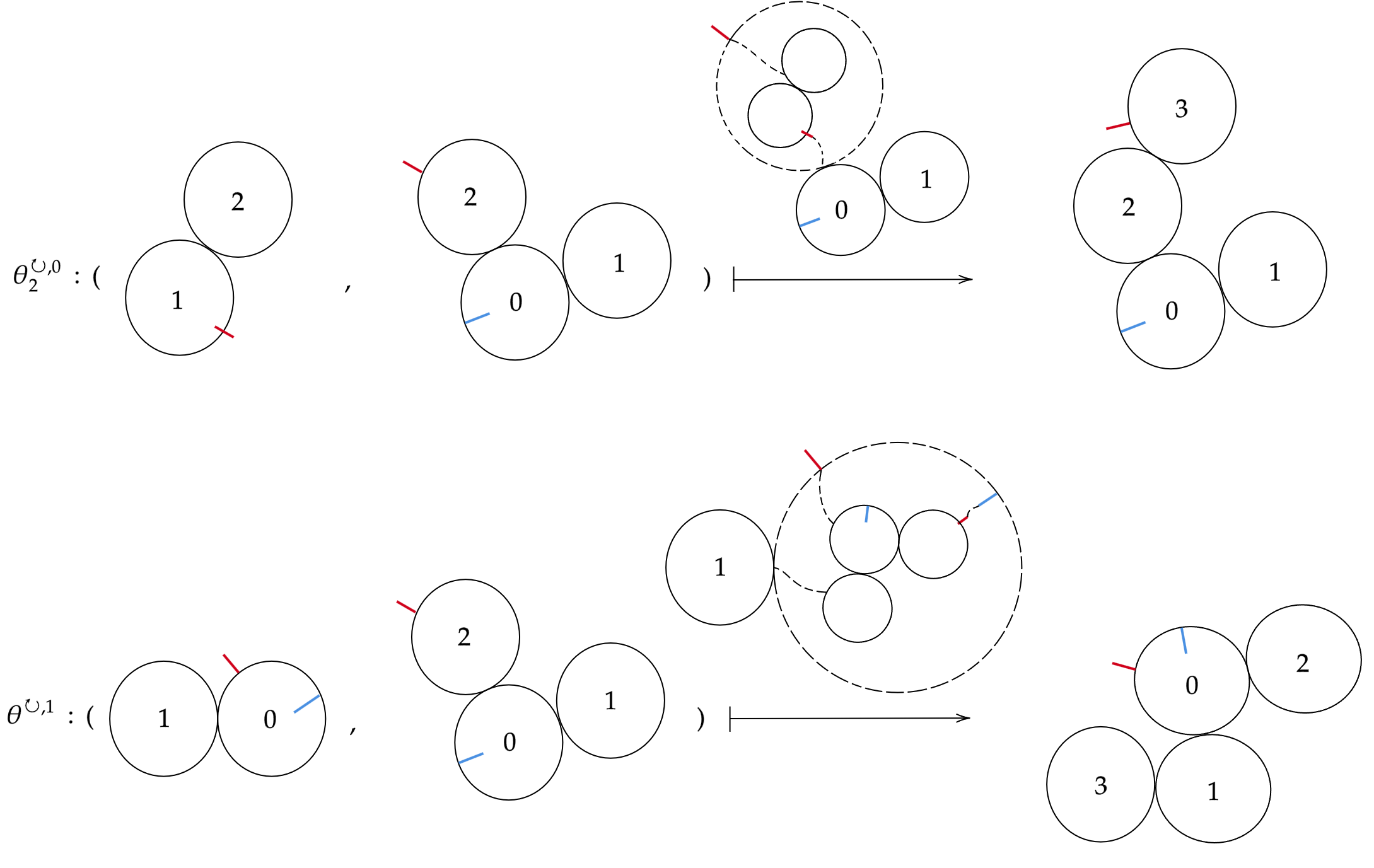}
 \caption{}
 \label{fig:composition_of_cyclic_cacti}
\end{figure}
\end{mydef}
We record the following lemma in analogy with Theorem \autoref{thm:associativity of cacti operad} and omit the proof. \\
\begin{lemma}\label{thm: associativity and equivariance of cyclic cacti operad}
\begin{itemize}
    \item Up to rescaling the $S^1$'s, $\theta^{\circlearrowright,0}_j$ is $S^1\times S^1$-equivariant; $\theta^{\circlearrowright,1}$ factors through 
    \begin{equation}\label{eq:quotient of theta circle 1}
    \theta^{\circlearrowright,1} :\mathrm{Cact}^l_{\circlearrowright}\times_{S^1} \mathrm{Cact}^k_{\circlearrowright}\rightarrow   \mathrm{Cact}^{k+l}_{\circlearrowright},
    \end{equation}
    where $S^1$ rotates the output basepoint of $\mathrm{Cact}^l_{\circlearrowright}$ and the input basepoint of $\mathrm{Cact}^k_{\circlearrowright}$. Moreover, \eqref{eq:quotient of theta circle 1} is also $S^1\times S^1$-equivariant where the left hand side is equipped with the obvious residual $S^1\times S^1$-action.  
    \item The structure maps $\theta, \theta^{\circlearrowright,0},\theta^{\circlearrowright,1}$ turns $\{\mathrm{Cact}^l,\mathrm{Cact}^k_{\circlearrowright}\}_{l,k\geq 1}$ into a two-colored topological operad up to homotopy.
\end{itemize}\qed
\end{lemma}
Similar as before, there are two ways to modify cyclic cacti to strictify the (equivariant) two-colored topological operad up to homotopy in Lemma \autoref{thm: associativity and equivariance of cyclic cacti operad}. The first is by assigning weights to the lobes. \\
\begin{mydef}\label{thm:weighted cyclic cacti}
The \emph{space of weighted cyclic cacti with $k$ lobes} is $\widetilde{\mathrm{Cact}_{\circlearrowright}^k}:=\mathrm{Cact}_{\circlearrowright}^k\times \mathcal{P}(k+1)$. 
\end{mydef}
Again, we think of $(x,\theta,\mathbf{a}=(a_0,a_1,\cdots,a_k))\in \widetilde{\mathrm{Cact}_{\circlearrowright}^k}=\mathrm{Cact}^{k+1}\times S^1\times \mathcal{P}(k+1)$ as a rescaling of $(x,\theta)$ so that the $i$-th lobe of $x$ has length $a_i$. A weighted cyclic cacti $(x,\theta,\mathbf{a})$ is uniquely determined by its \emph{cactus map}, which is the pair $(c_{x,\mathbf{a}},\theta)$, where $c_{x,\mathbf{a}}$ is the cactus map of the weighted cactus $(x,\mathbf{a})$, cf. \eqref{eq:weighted cactus map}. Let $\theta_j\in S^1, 0\leq j\leq k$ denote the clockwise angle from the local basepoint to the adjusted local basepoint (cf. Definition \autoref{thm:space of cyclic cacti}) of the $j$-th lobe. 
\par\indent
In terms of cactus maps, the operadic compositions involving weighted cyclic cacti are
\begin{equation}
\widetilde{\theta_j^{\circlearrowright,0}}:\widetilde{\mathrm{Cact}^l}\times \widetilde{\mathrm{Cact}^k_{\circlearrowright}}\rightarrow\widetilde{\mathrm{Cact}^{k+l-1}_{\circlearrowright}} ,\quad 1\leq j\leq k   
\end{equation}\
defined by $\widetilde{\theta_j^{\circlearrowright,0}}(c_{y,\mathbf{b}},   (c_{x,\mathbf{a}},\theta)):=$
\begin{equation}\label{eq:operadic composition of weighted cyclic cacti 0}
\big([0,1]\xrightarrow{c_{(x,\mathbf{a})}}\prod_{i=0}^k[0,a_i]\xrightarrow{\mathbbm{1}\times\cdots \times \overbrace{a_j(T_{-\theta_j}c_{y,\mathbf{b}})a_j^{-1}}^{\textrm{jth entry}}\times\cdots\mathbbm{1}} [0,a_0]\times\cdots\times [0,a_{j-1}]\times\prod_{s=1}^l[0,a_jb_s]\times[0,a_{j+1}]\times\cdots\times[0,a_k],\theta\big),
\end{equation}
where $(T_{\theta}c)(t):=c(t+\theta)-c(\theta)$ (the expression is considered cyclically);
\begin{equation}
\widetilde{\theta^{\circlearrowright,1}}:\widetilde{\mathrm{Cact}^l_{\circlearrowright}}\times \widetilde{\mathrm{Cact}^k_{\circlearrowright}}\rightarrow\widetilde{\mathrm{Cact}^{k+l}_{\circlearrowright}}   
\end{equation}\
defined by $\widetilde{\theta^{\circlearrowright,1}}((c_{x',\mathbf{a}'},\theta'),   (c_{x,\mathbf{a}},\theta)):=$
\begin{equation}\label{eq:operadic composition of weighted cyclic cacti 1}
\big([0,1]\xrightarrow{c_{(x',\mathbf{a}')}}\prod_{i=0}^l[0,a_i']\xrightarrow{a_0'(T_{-\theta'_0}c_{x,\mathbf{a}})(a_0')^{-1}\times\mathbbm{1}\times\cdots \times\mathbbm{1}} \prod_{j=0}^k[0,a_0'a_j]\times[0,a_1']\times\cdots\times[0,a'_l],\theta\big).
\end{equation}
We remark that $\theta'$ might not agree with $\theta'_0$, but their difference $\theta'-\theta_0'$ only depends on the underlying non-cyclic cactus $x'$. The following statement is evident from the definition.\\
\begin{lemma}\label{thm:weighted cyclic cacti form two colored operad}
The maps $\{\widetilde{\theta},\widetilde{\theta^{\circlearrowright,0}},\widetilde{\theta^{\circlearrowright,1}}\}$ define the structure of a two-colored topological operad on $\{\widetilde{\mathrm{Cact}^l},\widetilde{\mathrm{Cact}^k_{\circlearrowright}}\}_{l\geq 1,k\geq 0}$. Moreover, $\widetilde{\theta^{\circlearrowright,1}}$ factors through $\widetilde{\mathrm{Cact}^l_{\circlearrowright}}\times_{S^1} \widetilde{\mathrm{Cact}^k_{\circlearrowright}}$ and $\widetilde{\theta^{\circlearrowright,0}}, \widetilde{\theta^{\circlearrowright,1}}$ are $S^1\times S^1$-equivariant.  \qed
\end{lemma}
\textbf{Cyclic cacti with spines and their operadic structure}. There is an obvious way to combine Definition \autoref{thm:Cacti with spines as a cosimplicial multisimplicial set} with Definition \autoref{thm:space of cyclic cacti}.\\
\begin{mydef}\label{thm:cyclic cacti with spines as a (co)cyclic multisimplicial object}
We define a cocyclic cyclic k-fold simplicial set $\widehat{\mathfrak{Cact}^k_{\circlearrowright}}=(\widehat{\mathfrak{Cact}^k_\circlearrowright})_{\bullet;\bullet,\cdots,\bullet;\bullet}: \Lambda^{op}\times (\Delta^{op})^k\times \Lambda\rightarrow \mathrm{Sets}$ as follows. Given non-negative integers $n,m_0,m_1,\cdots,m_k$, $(\widehat{\mathfrak{Cact}^k_\circlearrowright})_{n;m_1,\cdots,m_k;m_0}$ is the set of pairs $(f,S,\iota)$ consisting of 
\begin{itemize}
    \item A surjective map $f: \{1,\cdots,m+k+1\}\rightarrow \{0,1,\cdots,k\}$ with $m=\sum_{i=0}^k m_i$ satisfying conditions 1) and 3) of Definition \autoref{thm:cacti}.
    \item A map $S:\{0,1,\cdots,m+k+1\}\rightarrow \mathbb{Z}_{\geq 0}$ satisfying $\sum_{i=0}^{m+k+1}S(i)=n$.
    \item A distinguished element $\iota\in \{0,1,\cdots,m_0\}$.
\end{itemize}    
The \emph{space of cyclic cacti with $k$ lobes and $n$ spines} is the $(k+1)$-fold simplicial realization $\widehat{\mathrm{Cact}^k_{\circlearrowright,n}}:=|\widehat{(\mathfrak{Cact}^k_{\circlearrowright}})_{n;\bullet,\cdots,\bullet;\bullet}|$. Denote $\widehat{\mathrm{Cact}^k_{\circlearrowright}}:=\mathrm{Tot}([n]\mapsto\widehat{\mathrm{Cact}^k_{\circlearrowright,n}})$; by Lemma \autoref{thm:circle actions of realization of (co)cyclic sets}, $\widehat{\mathrm{Cact}^k_{\circlearrowright}}$ has an induced $S^1\times S^1$-action. 
\begin{figure}[H]
 \centering
 \includegraphics[width=0.9\textwidth]{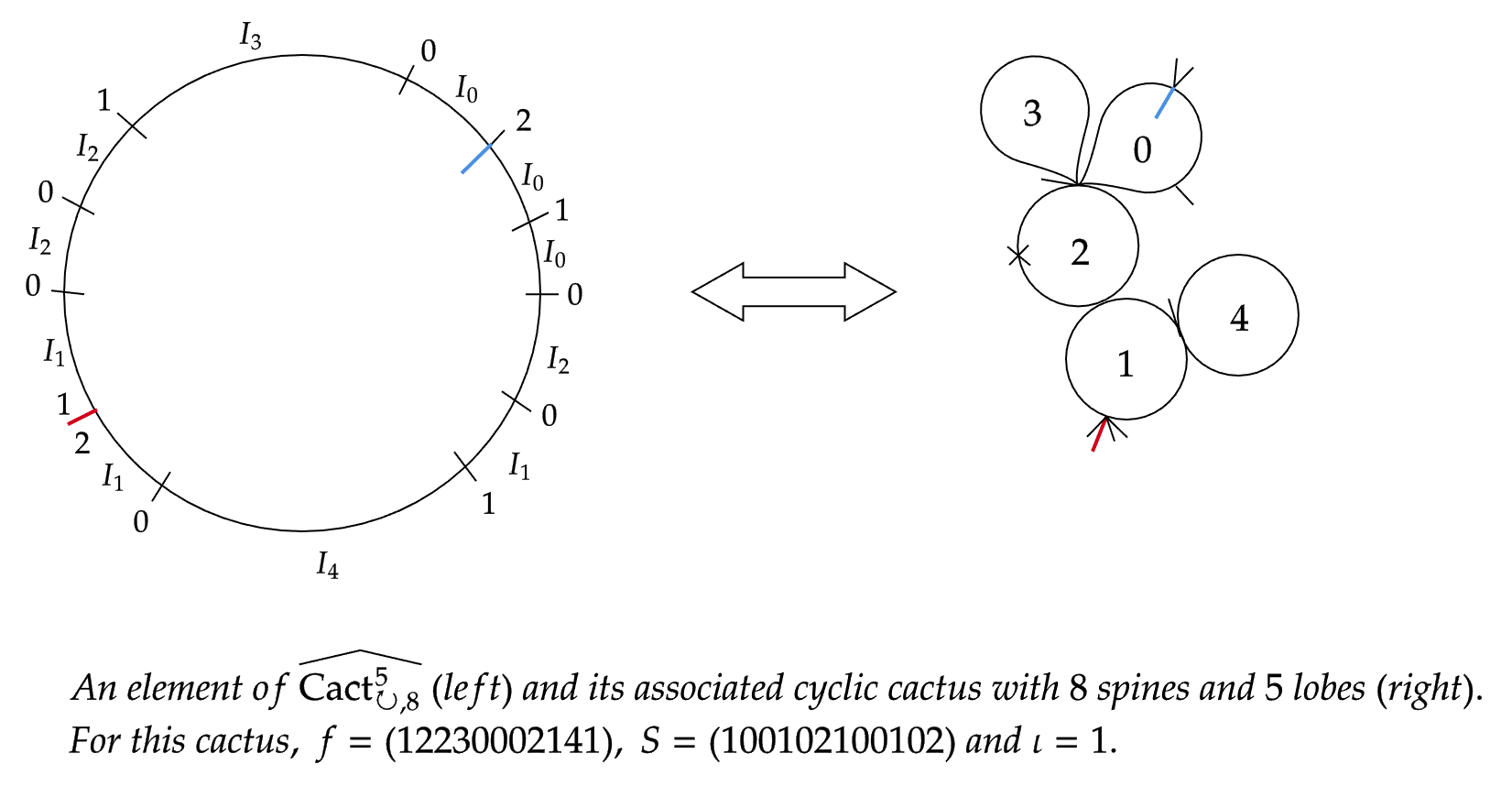}
 \caption{}
 \label{fig:cyclic_cactus_with_spine}
\end{figure}
We think of the input basepoint as `the $0$-th spine' and if $f^{-1}(0)=\{a_0<\cdots<a_{m_0}\}$ then the marked point $a_{\iota}$ corresponds to the output basepoint, see Figure \ref{fig:cyclic_cactus_with_spine}.
\begin{figure}[H]
 \centering
 \includegraphics[width=0.9\textwidth]{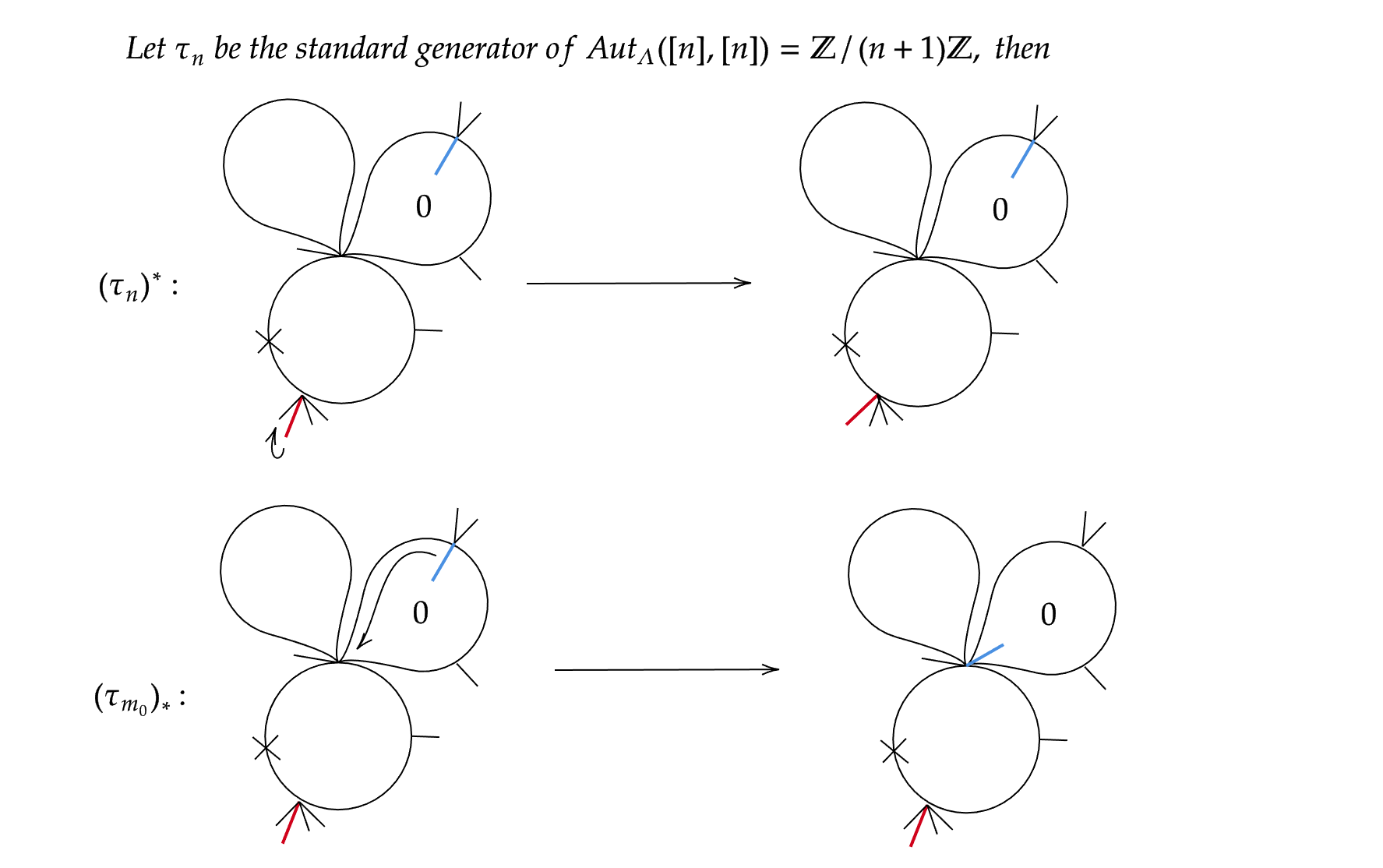}
 \caption{}
 \label{fig:co_cyclic_maps_of_cyclic_cacti_with_spines}
\end{figure}
\end{mydef}
The (co)face and (co)degeneracy maps are defined the same way as in Definition \autoref{thm:Cacti with spines as a cosimplicial multisimplicial set} (also see Figure \ref{fig:cactus_with_spine_co_face_and_co_degeneracy}) with the caveat that on the $j$-th lobe one needs to cyclically reorder $f^{-1}(j)=\{a^j_0<\cdots<a^j_{m_j}\}$ as $\{a^j_{\eta_j}<a^j_{\eta_j+1}<\cdots\}$ where $a^j_{\eta_j}$ is the adjusted local basepoint on the $j$-th lobe (e.g. the $0$-th face map `contracts' the interval between the adjusted local basepoint and the next marked point). The cyclic (resp. cocyclic) structure is induced by cyclically permuting the spines (resp. the marked points on the base lobe), see Figure \ref{fig:co_cyclic_maps_of_cyclic_cacti_with_spines}.\par\indent
Next we define operadic compositions for cyclic cacti with spines. In analogy with \eqref{eq:operadic structure on cacti with spine version 2}, but now taking into account the (co)cyclic structures, there are two types of composition maps:
\begin{equation}\label{eq:cyclic cacti with spine operadic composition 1}
\widehat{\theta^{\circlearrowright,0}_j}:\widehat{\mathfrak{Cact}^l}\times_{j,\Delta^{op}}\widehat{\mathfrak{Cact}^k_{\circlearrowright}}\rightarrow   \widehat{\mathfrak{Cact}^{k+l-1}_{\circlearrowright}},\quad(1\leq j\leq k)
\end{equation}
which is a map of $\Lambda^{op}\times (\Delta^{op})^{k+l-1}\times\Lambda$-sets, and
\begin{equation}\label{eq:cyclic cacti with spine operadic composition 1}
\widehat{\theta^{\circlearrowright,1}}:\widehat{\mathfrak{Cact}^l_{\circlearrowright}}\times_{\Lambda^{op}}\widehat{\mathfrak{Cact}^k_{\circlearrowright}}\rightarrow   \widehat{\mathfrak{Cact}^{k+l}_{\circlearrowright}},    
\end{equation}
which is a map of $\Lambda^{op}\times (\Delta^{op})^{k+l}\times\Lambda$-sets. \\
\begin{mydef}\label{thm:pictorial composition of cyclic cacti with spines}
\begin{itemize}
    \item Given $x_1\in \widehat{\mathfrak{Cact}^{l}}_{m_j;m_1',\cdots,m'_{l}}$ and $x_2\in (\widehat{\mathfrak{Cact}^k_{\circlearrowright}})_{n;m_1,\cdots,m_k;m_0}$, their composition $\widehat{\theta^{\circlearrowright,0}_j}(x_1,x_2)\in (\widehat{\mathfrak{Cact}^{k+l-1}_{\circlearrowright}})_{n;m_1,\cdots,m_{j-1},m_1',\cdots,m_l',m_j+1,\cdots,m_k;m_0}$ is the cyclic cactus with spines obtained by inserting $x_1$ into the $j$-th lobe of $x_2$ so that 1) the basepoint of $x_1$ matches the local input basepoint of the $j$-th lobe of $x_2$ and 2) the spines of $x_1$ precisely match up the marked points on the $j$-th lobe of $x_2$. 
    \item Given $x_1\in (\widehat{\mathfrak{Cact}^{l}_{\circlearrowright}})_{n';m_1',\cdots,m'_{l};m_0'}$ and $x_2\in (\widehat{\mathfrak{Cact}^k_{\circlearrowright}})_{m_0';m_1,\cdots,m_k;m_0}$, their composition $\widehat{\theta^{\circlearrowright,1}}(x_1,x_2)\in (\widehat{\mathfrak{Cact}^{k+l}_{\circlearrowright}})_{n';m_1,\cdots,m_k,m_1',\cdots,m_l';m_0}$ is the cyclic cactus with spines obtained by inserting $x_2$ into the base lobe of $x_1$ so that 1) the input basepoint of $x_2$ matches the output basepoint of $x_1$ and 2) the spines of $x_2$ precisely match up the marked points on the base lobe of $x_1$.
\end{itemize}
See Figure \ref{fig:composition_of_cyclic_cacti_with_spines} for an illustration.\par\indent
One easily checks that the structure maps $\widehat{\theta^{\circlearrowright,0}}, \widehat{\theta^{\circlearrowright,1}}$, together with $\widehat{{\theta}}$ from \eqref{eq:operadic structure on cacti with spine version 2}, satisfy (strict) associativity of composition. In particular, they endow the collection $\{\widehat{\mathfrak{Cact}^l},\widehat{\mathfrak{Cact}^k_{\circlearrowright}}\}_{l\geq 1,k\geq 0}$ with the structure of a two colored operad in the (co)cyclic/(co)simplicial multi-simplicial setting. Applying Lemma \autoref{thm:induced operadic structure on totalization-realization} and Lemma \autoref{thm:circle actions of realization of (co)cyclic sets} gives the following.\\
\begin{lemma}\label{thm:operadic structure on totalization-realization of cyclic cacti with spines}
Upon taking totalization-realization in the category of topological spaces (resp. $R$-chain complexes), $\{\widehat{\theta
},\widehat{\theta^{\circlearrowright,0}}, \widehat{\theta^{\circlearrowright,1}}\}$ induces the structure of a two colored topological (resp. dg) operad on the collection 
\begin{equation}
\{\widehat{\mathrm{Cact}^l},\widehat{\mathrm{Cact}^k_{\circlearrowright}}\}_{l\geq 1,k\geq 0} \quad(\mathrm{resp.}\;\{\widehat{\mathrm{Cact}^l_R},\widehat{\mathrm{Cact}^k_{\circlearrowright,R}}\}_{l\geq 1,k\geq 0}  ).
\end{equation}
Moreover, the structure maps induced from $\widehat{\theta^{\circlearrowright,1}}$ factors through the quotient $\widehat{\mathrm{Cact}^l}\times_{S^1}\widehat{\mathrm{Cact}^k_{\circlearrowright}}$ (resp. $\widehat{\mathrm{Cact}^l_R}\times_{C_*(S^1;R)}\widehat{\mathrm{Cact}^k_{\circlearrowright,R}}$)  and the structure maps induced from $\widehat{\theta^{\circlearrowright,0}}, \widehat{\theta^{\circlearrowright,1}}$ are $S^1\times S^1$-equivariant. \qed\\
\end{lemma}
\begin{mydef}\label{thm:the full KS operad}
Define $\{\mathbf{KS}(l,0),\mathbf{KS}(k,1)\}_{l\geq 1,k\geq 0}$ (resp. $\{\mathbf{KS}(l,0)_R,\mathbf{KS}(k,1)_R\}_{l\geq 1,k\geq 0}$) to be the topological (resp. dg) operad from Lemma \autoref{thm:operadic structure on totalization-realization of cyclic cacti with spines}.   
\end{mydef}
\begin{figure}[H]
 \centering
 \includegraphics[width=1.0\textwidth]{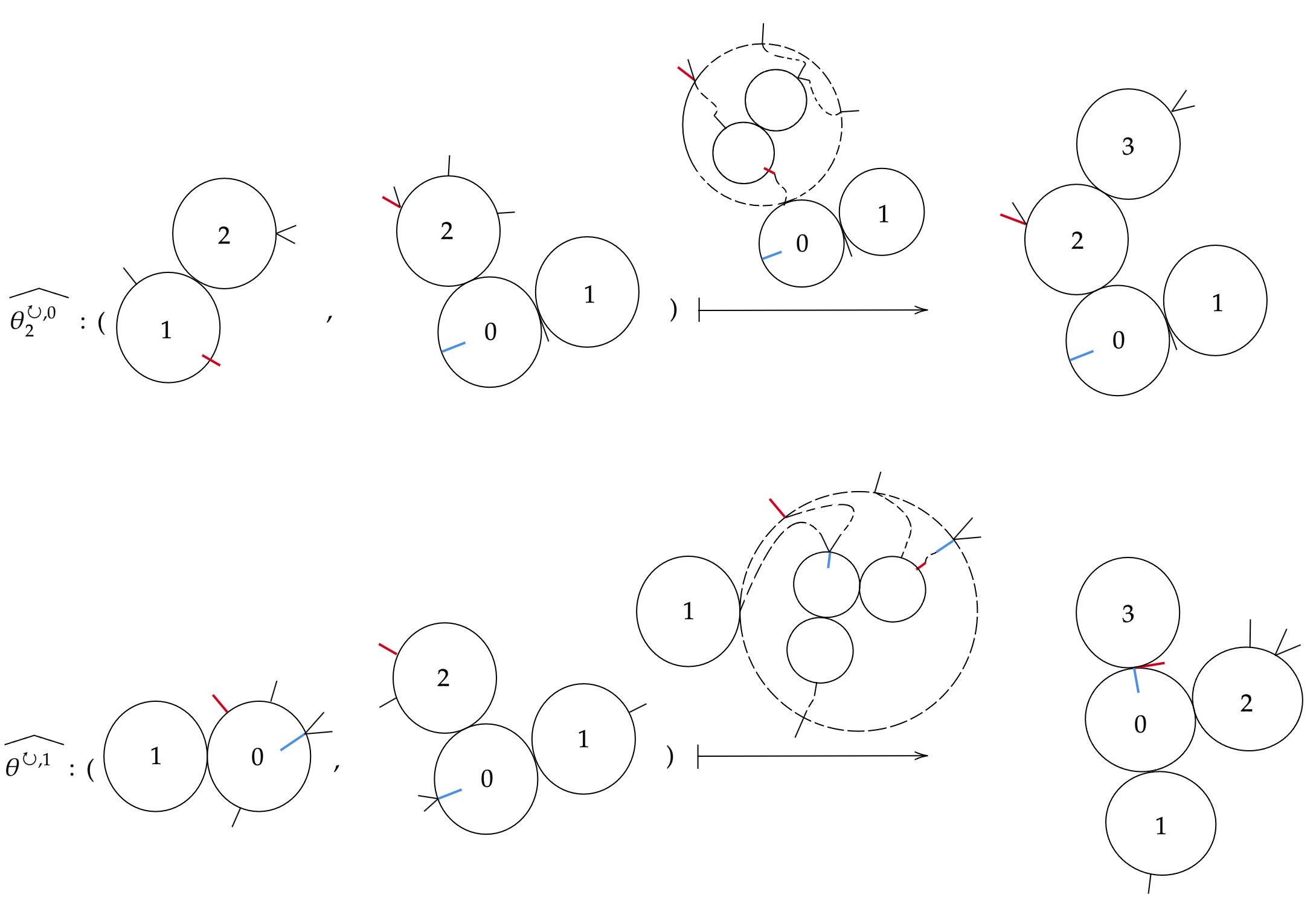}
 \caption{}
 \label{fig:composition_of_cyclic_cacti_with_spines}
\end{figure}
\end{mydef}
\begin{thm}\label{thm:KS equals weighted cyclic cacti}
There is a quasi-equivalence of two-colored topological operads 
\begin{equation}\label{eq:KS equals weighted cyclic cacti}
\{\mathbf{KS}(l,0),\mathbf{KS}(k,1)\}_{l\geq 1,k\geq 0}\simeq \{\widetilde{\mathrm{Cact}^l},\widetilde{\mathrm{Cact}^k_{\circlearrowright}}\}_{l\geq 1,k\geq 0}.
\end{equation}
extending the quasi-equivalence in Proposition \autoref{thm:KS equals weighted cacti}. Moreover, for each $k\geq 0$, the weak equivalence $\mathbf{KS}(k,1)\simeq \widetilde{\mathrm{Cact}^k_{\circlearrowright}}$ is $S^1\times (S^1)^{op}$-equivariant, i.e. equivariance in the second $S^1$-factor holds after composing with $S^1\xrightarrow{(\;\;)^{-1}} S^1$.
\end{thm}
\emph{Proof.} There is a homeomorphism 
\begin{equation}\label{eq:homeo from cyclic cacti with spines to cyclic cacti times Delta}
g=(g_1,g_2,g_3):\widehat{\mathrm{Cact}^{k}_{\circlearrowright,n}}\xrightarrow{\cong} \mathrm{Cact}^{k+1}\times \Delta^n\times S^1    
\end{equation}
where $g_1$ records the underlying (non-cyclic) cactus, $g_2$ records the consecutive distances between the spines (recall the convention is that the input basepoint is the $0$-th spine) and $g_3$ records the position of the output basepoint on the base lobe. We claim that \eqref{eq:homeo from cyclic cacti with spines to cyclic cacti times Delta} is $(S^1)^{op}$-equivariant: that is, it is $S^1$-equivariant where the $S^1$-action on the left hand side is induced from the cyclic structure in $\widehat{\mathfrak{Cact}^k_{\circlearrowright,n}}:(\Delta^{op})^k\times\Lambda\rightarrow \mathrm{Sets}$ via Lemma \autoref{thm:circle actions of realization of (co)cyclic sets}, and on the right hand side is induced from the \emph{inverse} regular representation on the third factor. \par
\begin{figure}[H]
 \centering
 \includegraphics[width=1.0\textwidth]{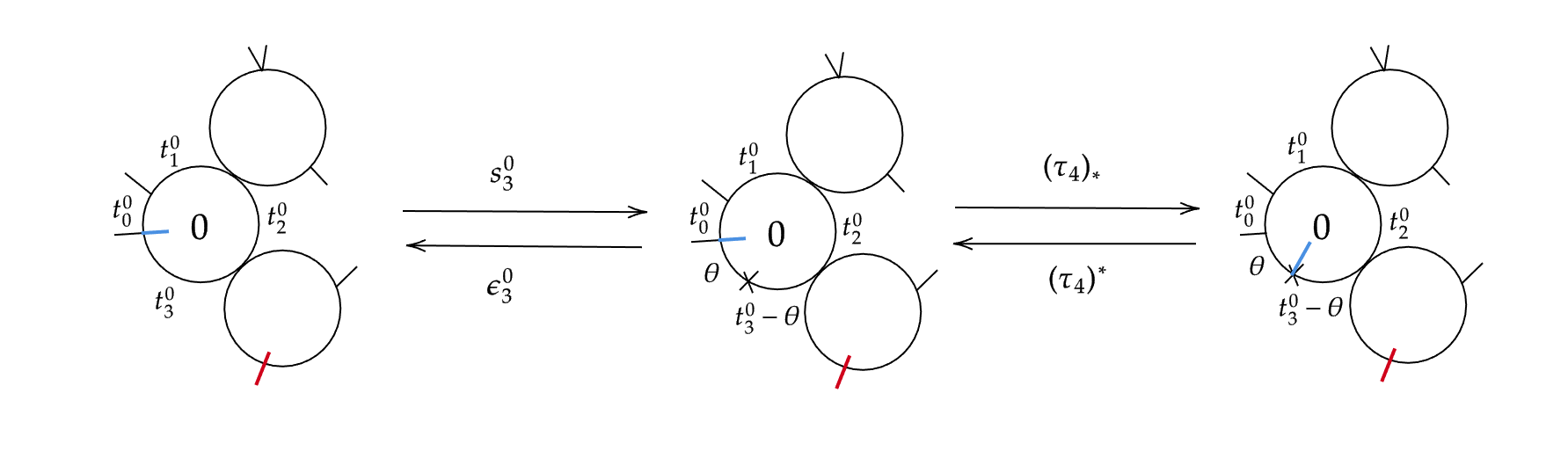}
 \caption{}
 \label{fig:circle_action_from_cyclic_structure_for_cyclic_cacti_with_spines}
\end{figure}

\indent
To see this, let 
\begin{equation}\
(x=(f,S,\iota),t^1,\cdots,t^k,t^0)\in (\widehat{\mathfrak{Cact}^k_{\circlearrowright}})_{n;m_1,\cdots,m_k;m_0}\times\Delta^{m_1}\times\cdots\times\Delta^{m_k}\times\Delta^{m_0},  
\end{equation}
viewed as a subspace of $(\widehat{\mathfrak{Cact}^k_{\circlearrowright}})_{n;m_1,\cdots,m_k;m_0}\times\Delta^{m_1}\times\cdots\times\Delta^{m_k}\times S^1\times\Delta^{m_0}$.
Let $f^{-1}(0)=\{a_0<\cdots<a_{m_0}\}$ and $t^0=(t^0_0,\cdots,t^0_{m_0})$. Geometrically, $t^0$ represents the consecutive distances between $a_{\iota},a_{\iota+1},\cdots,a_{m_0},\cdots,a_{\iota-1}$ (indices are considered cyclically mod $m_0$). Let $\theta\in S^1$ and without loss of generality assume $\theta\in[0,t_n]$ (the general case follows by induction). By the proof of Lemma \autoref{thm:circle actions of realization of (co)cyclic sets}, $\theta\cdot(x,t^1,\cdots,t^n,t^0)=$
$$\big(x,t^1,\cdots,t^k,(e^{-2\pi i\theta},(t^0_0,\cdots,t^0_{m_0}))\big)=\big(x,t^1,\cdots,t^k,(e^{-2\pi i\theta},\epsilon^0_{m_0}(t^0_0,\cdots,t^0_{m_0},t^0_{m_0}-\theta,\theta))\big)=$$
$$\big(s^0_{m_0}x,t^1,\cdots,t^k,(e^{-2\pi i\theta},(t^0_0,\cdots,t^0_{m_0},t^0_{m_0}-\theta,\theta))\big)=  \big(s^0_{m_0}x,t^1,\cdots,t^k,\tau_{m_0+1}^*(1,(\theta,t^0_0,\cdots,t^0_{m_0},t^0_{m_0}-\theta))\big)= $$
\begin{equation}
\big((\tau_{m_0+1})_*s^0_{m_0}x,t^1,\cdots,t^k,(1,(\theta,t^0_0,\cdots,t^0_{m_0},t^0_{m_0}-\theta))\big)\leftrightarrow\big((\tau_{m_0+1})_*s^0_{m_0}x,t^1,\cdots,t^k,(\theta,t^0_0,\cdots,t^0_{m_0},t^0_{m_0}-\theta)\big),
\end{equation}
which corresponds exactly to rotating the output basepoint \emph{counterclockwise} by angle $\theta$; see Figure \ref{fig:circle_action_from_cyclic_structure_for_cyclic_cacti_with_spines}. 
On the other hand, let $\mathrm{cact}^{k+1}:=\mathrm{Cact}^{k+1}/S^1$ be the space of \emph{unbased cacti}. The principal $S^1$-bundle $\mathrm{Cact}^{k+1}\rightarrow\mathrm{cact}^{k+1}$ admits a canonical section by letting the basepoint be the intersection of the $0$-th lobes with the connected component of its complement closure containing the $1$-st lobe (except for $k=0$, which requires a separate but easy discussion and is left to the reader); hence there is a homeomorphism $\mathrm{Cact}^{k+1}\cong \mathrm{cact}^{k+1}\times S^1$. \eqref{eq:homeo from cyclic cacti with spines to cyclic cacti times Delta} induces a homeomorphism 
\begin{equation} \label{eq:homeo from cyclic cacti with spines to unbased cacti times stuff}
\widehat{\mathrm{Cact}^k_{\circlearrowright,n}}\xrightarrow{\cong} \mathrm{cact}^{k+1}\times (S^1\times\Delta^n)  \times S^1.  
\end{equation}
From the geometric picture, one readily observes that \eqref{eq:homeo from cyclic cacti with spines to unbased cacti times stuff} gives rise to an isomorphism of cocyclic spaces $\widehat{\mathrm{Cact}^k_{\circlearrowright,\bullet}}\xrightarrow{\cong} \mathrm{cact}^{k+1}\times \Lambda^{\bullet}  \times S^1$, which inherits the previous $(S^1)^{op}$-equivariance. But then Lemma \autoref{thm:circle action on Mon(S^1)} implies that after totalization, the homeomorphism
\begin{equation}\label{eq:from KS to unbased cact times mon(S^1) times S^1}
\mathbf{KS}(k,1)=\mathrm{Tot}(\widehat{\mathrm{Cact}^{k}_{\circlearrowright,\bullet}})\cong  \mathrm{cact}^{k+1}\times \mathrm{Tot}(\Lambda^{\bullet})  \times S^1\cong \mathrm{cact}^{k+1}\times \mathrm{Mon}(S^1)  \times S^1   
\end{equation}
is $S^1\times (S^1)^{op}$-equivariant. \par\indent
Given an unbased cactus $x\in \mathrm{cact}^{k+1}$, by an abuse of notation let $c_x=(c_x^0,c_x^1,\cdots,c_x^k):S^1\rightarrow S^1\times(S^1)^{\times k}$ denote the cactus map of the image of $x$ under the canonical section of $\mathrm{Cact}^{k+1}\rightarrow\mathrm{cact}^{k+1}$. Then, the assignment
\begin{equation}
(x,f,\theta)\mapsto (c_x^0\circ f+\theta,c_x^1\circ f,\cdots, c_x^k\circ f)    
\end{equation}
defines an $S^1\times S^1$-equivariant embedding $\mathrm{cact}^{k+1}\times \mathrm{Mon}(S^1)  \times S^1\hookrightarrow \mathrm{CoEnd}_{S^1}(k,1)=\mathrm{Map}(S^1,S^1\times(S^1)^{\times k})$. Precomposing with \eqref{eq:from KS to unbased cact times mon(S^1) times S^1} one obtains an $S^1\times (S^1)^{op}$-equivariant embedding $\mathbf{KS}(k,1)\hookrightarrow \mathrm{CoEnd}_{S^1}(k,1)$ which (together with its counterpart from Proposition \autoref{thm:KS equals weighted cacti}), by a tedious but straightforward verification, induce an embedding of two-colored topological operads
\begin{equation}\label{eq:from KS to CoEnd two colored}
\{\mathbf{KS}(l,0),\mathbf{KS}(k,1)\}_{l\geq 1,k\geq 0}\hookrightarrow \{\mathrm{CoEnd}_{S^1}(l,0),\mathrm{CoEnd}_{S^1}(k,1)\}_{l\geq 1,k\geq 0}.   
\end{equation}
Proceeding entirely analogously to Proposition \autoref{thm:KS equals weighted cacti}, one shows that there is a further two-colored suboperad characterized by an additional piecewise linear condition (cf. Condition 2') of Proposition \autoref{thm:KS equals weighted cacti}) which is the image of an embedding $\{\widetilde{\mathrm{Cact}^l},\widetilde{\mathrm{Cact}^k_{\circlearrowright}}\}_{l\geq 1,k\geq 0}\hookrightarrow\{\mathrm{CoEnd}_{S^1}(l,0),\mathrm{CoEnd}_{S^1}(k,1)\}_{l\geq 1,k\geq 0}$ such that each $\widetilde{\mathrm{Cact}^k_{\circlearrowright}}\hookrightarrow\mathrm{CoEnd}_{S^1}(k,1)$ is $S^1\times S^1$-equivariant. This gives the desired embedding 
\begin{equation}
\{\widetilde{\mathrm{Cact}^l},\widetilde{\mathrm{Cact}^k_{\circlearrowright}}\}_{l\geq 1,k\geq 0}\hookrightarrow \{\mathbf{KS}(l,0),\mathbf{KS}(k,1)\}_{l\geq 1,k\geq 0}.  
\end{equation} 
Finally, that $\widetilde{\mathrm{Cact}^k_{\circlearrowright}}\hookrightarrow\mathbf{KS}(k,1)$ is a weak equivalence follows from a commutative diagram
\begin{equation}
\begin{tikzcd}[row sep=1.2cm, column sep=0.8cm]
\widetilde{\mathrm{Cact}^k_{\circlearrowright}}\arrow[hookrightarrow,rr]\arrow[dr]& &\mathbf{KS}(k,1)\arrow[dl]\\
 & \mathrm{Cact}^k_{\circlearrowright}  &
\end{tikzcd}
\end{equation}
whose downward arrows are both homotopy equivalences.\qed\par\indent
Finally, we remark that it is evident from the proof that this quasi-equivalence is equivariant with respect to the symmetric group actions permuting the labels of the lobes of a (cyclic) cactus.\par\indent
\textbf{Action on Hochschild chain and cochain}. Let $\mathcal{A}$ be a dg algebra over $R$ and $\mathcal{A}_{\sharp}$ (resp. $\mathcal{A}^{\sharp}$) be the (resp. co)simplicial chain complex from \eqref{eq:HH functor} (resp. \eqref{eq:HH cochain functor}) whose realization (resp. totalization) is the Hochschild (resp. co)chain complex of $\mathcal{A}$. \\
\begin{mydef}\label{thm: action of KS(n,1)_R on Hochschild (co)chains}
Fix $x=(f,S,\iota)\in (\widehat{\mathfrak{Cact}^k_{\circlearrowright}})_{n;m_1,\cdots,m_k;m_0}$, $\varphi_j\in\mathcal{A}^{\sharp}([m_j])=\mathrm{Hom}_R(\mathcal{A}^{\otimes m_j},\mathcal{A}), 1\leq j\leq k$ and $a_0\otimes a_1\otimes\cdots\otimes a_n\in \mathcal{A}_{\sharp}([n])=\mathcal{A}^{\otimes n+1}$, we define 
\begin{equation}
\mathrm{Act}_{x}(a_0\otimes\cdots\otimes a_n;\varphi_1,\cdots,\varphi_k)=b_0\otimes b_1\otimes\cdots\otimes b_{m_0} \in \mathcal{A}^{\otimes m_0+1}  
\end{equation}
as follows. Attach $a_0$ to the input basepoint, attach $a_i$ ($1\leq i\leq n$) to the $i$-th spine of $x$ (and attach the unit to an empty spine), and insert $\varphi_j$ into the $j$-th lobe of $x$. The cyclic cactus with spines $x$ specifies a way of successively applying $\varphi_j$ or multiplication in $\mathcal{A}$ to the tuple $(a_0,a_1,\cdots,a_n)$, and $b_0\otimes\cdots\otimes b_{m_0}$ is defined so that $b_0$ is the resulting element at the output basepoint and $b_i$ is the resulting element at the $i$-th marked point on the base lobe (counted clockwise from the output basepoint); again, the result is to be adjusted by the appropriate Koszul sign. See Figure \ref{fig:action_of_cyclic_cacti_with_spines_on_HH_co_chain} for an illustration.
\begin{figure}[H]
 \centering
 \includegraphics[width=1.0\textwidth]{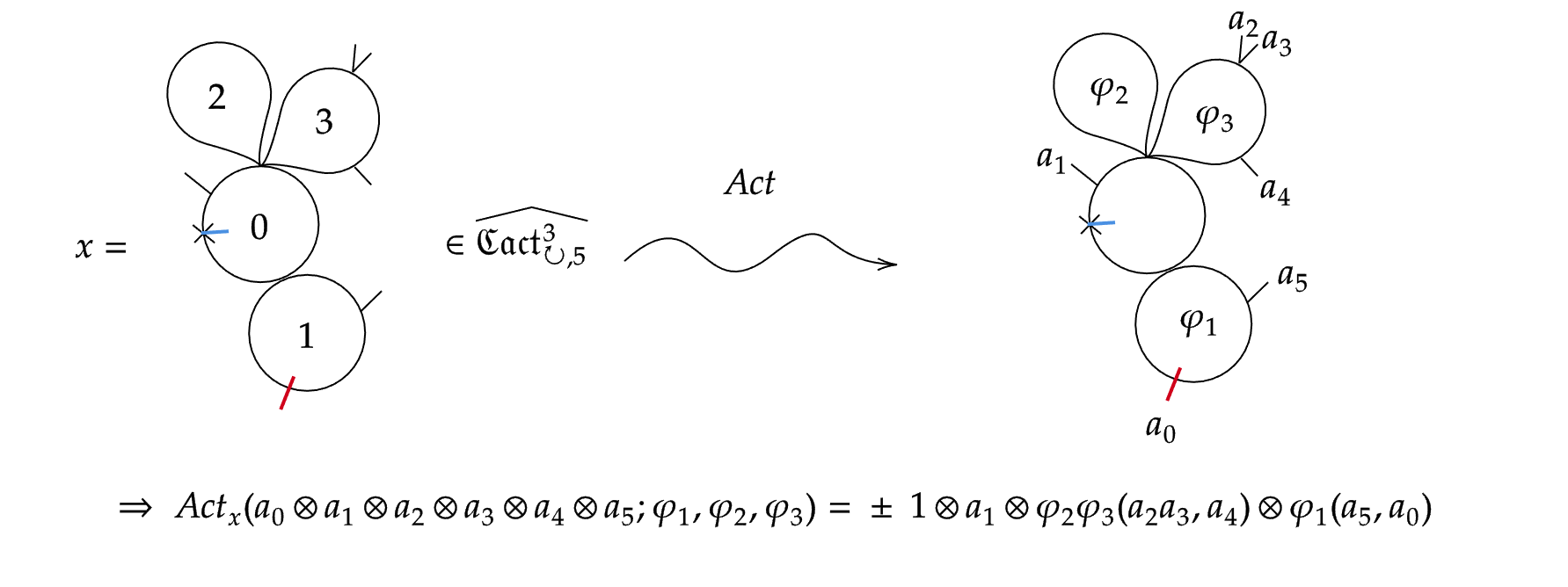}
 \caption{}
 \label{fig:action_of_cyclic_cacti_with_spines_on_HH_co_chain}
\end{figure}
\end{mydef}
The following lemma is immediate from Definition \autoref{thm: action of KS(n,1)_R on Hochschild (co)chains}. \\
\begin{lemma}
$\mathrm{Act}: \widehat{\mathfrak{Cact}^k_{\circlearrowright}}\rightarrow\mathrm{Hom}_R(\mathcal{A}_\sharp\otimes(\mathcal{A}^{\sharp})^{\otimes k},\mathcal{A}_{\sharp})$ defines a map of $\Lambda^{op}\times(\Delta^{op})^k\times\Lambda$-sets. Moreover, it intertwines the operadic structures $\widehat{\theta^{\circlearrowright,0}}$ and $ \widehat{\theta^{\circlearrowright,1}}$ with the operadic structures on $\{\mathrm{Hom}_R((\mathcal{A}^{\sharp})^{\otimes l},\mathcal{A}^{\sharp}),\mathrm{Hom}_R(\mathcal{A}_\sharp\otimes(\mathcal{A}^{\sharp})^{\otimes k},\mathcal{A}_{\sharp})\}_{l\geq 1, k\geq 0}$ given by insertion of multi-linear maps and composition.  \qed
\end{lemma}
Taking totalization-realization in the category of $R$-chain complexes and combining with \eqref{eq:KS(n,0)_R action on Hocschild cochains, after taking totalization} give a map of dg operads
\begin{equation}\label{eq:KS(n,0)_R and KS(n,1)_R action on Hocschild (co)chains, after taking totalization}    
\mathrm{Act}: (\{\mathbf{KS}(l,0)_R,\mathbf{KS}(k,1)_R\}_{l\geq 1,k\geq 0}\rightarrow \{\mathrm{Hom}_R(CC^*(\mathcal{A})^{\otimes l},CC^*(\mathcal{A})),\mathrm{Hom}_R(CC_*(\mathcal{A})\otimes CC^*(\mathcal{A})^{\otimes k},CC_*(\mathcal{A}))\}_{l\geq 1,k\geq 0}
\end{equation}
which is furthermore $S^1\times S^1$-equivariant when restricted to $\mathbf{KS}(k,1)_R$. We also note the evident fact that $\mathrm{Act}$ is equivariant with respect to the symmetric group actions, where on the left hand side it is permuting the labels of the lobes and on the right hand side it is permuting the tensor factors of the Hochschild cochains. \\
\begin{rmk}
In this paper we only consider the action of $\mathbf{KS}_R$ on the Hochschild cochain/chain of a $R$-linear dg algebra. On the other hand, by using a two-colored version of the argument in \cite[Section 5]{MS} one may readily check that the topological totalization/realization $\mathbf{KS}$ of Lemma \ref{thm:operadic structure on totalization-realization of cyclic cacti with spines} acts on the topological Hochschild cohomology/homology of any associative ring spectrum, and thus combined with Theorem \ref{thm:comparison of KS with Cyl} gives a `formulae-based proof' of a two-colored version of the topological Deligne conjecture. We remark that a more conceptual proof using factorization homology can be found in \cite{Hor}. 
\end{rmk}

\subsection{A little disks model for $\mathbf{KS}$}
In this subsection, we recall the construction of a two-colored topological operad $\mathrm{Cyl}$ based on configuration spaces \cite{KS1}\cite{DTT1}, which is more accessible to computations. Later in Section 6, we show that $\mathrm{Cyl}$ recovers the Kontsevich-Soibelman operad $\mathbf{KS}$ from Section 4.1 and 4.2.\par\indent
The values of the two-colored operad $\{\mathrm{Cyl}(l,0),\mathrm{Cyl(k,1)}\}_{l\geq 1,k\geq 0}$ are defined as follows. We set $\{\mathrm{Cyl}(l,0)\}_{l\geq 1}$ to be the little disk operad.\par
\begin{figure}[H]
 \centering
 \includegraphics[width=0.9\textwidth]{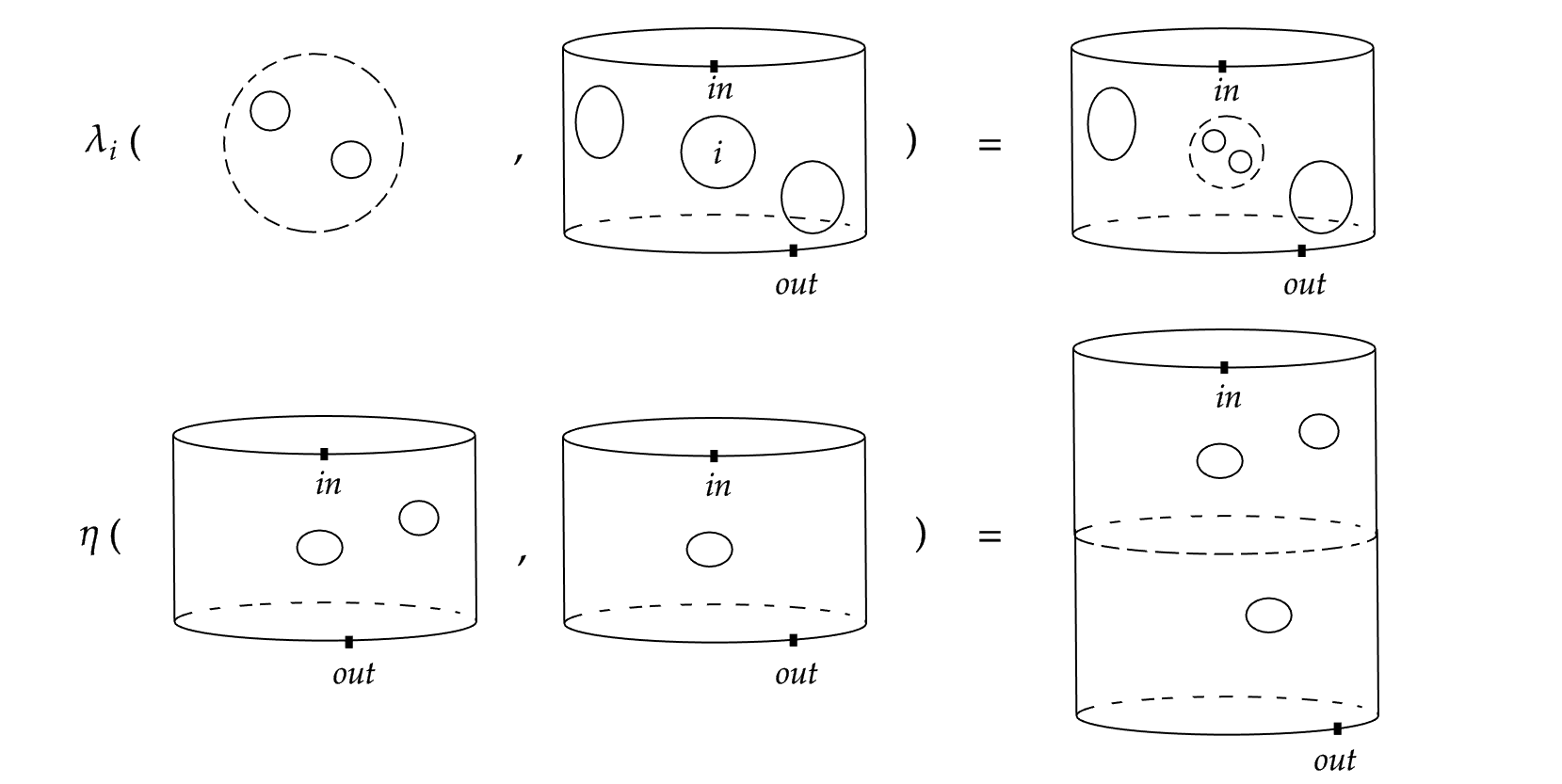}
 \caption{}
 \label{fig:figure3.5}
\end{figure}

\indent
On the other hand, $\mathrm{Cyl}(n,1), n\geq 0$ is the topological space of 
\begin{itemize}
    \item When $n\geq1$, cylinder $S^1\times [a,c],a<c$ together with a configuration of $n$ disks on $S^1\times(a,c)$ and two marked points $in\in S^1\times\{c\}, out\in S^1\times\{a\}$. When $n=0$, $\mathrm{Cyl}(0,1)$ is the space of two points $in,out\in S^1$ modulo simultaneous rotation. 
    \item These configurations are consider up to equivalence generated by parallel shifts $S^1\times[a,c]\rightarrow S^1\times [a+l,c+l]$ and overall $S^1$-rotation of the cylinder. 
\end{itemize}
The two kinds of operadic compositions involving $\mathrm{Cyl}(n,1)$
\begin{equation}\label{eq:cyl operadic 1}
\lambda_i:\mathrm{Cyl}(n,0)\times_i \mathrm{Cyl}(n',1)\rightarrow \mathrm{Cyl}(n+n'-1,1)  
\end{equation}
\begin{equation}\label{eq:cyl operadic 2}
\eta:\mathrm{Cyl}(n,1)\times \mathrm{Cyl}(n',1)\rightarrow \mathrm{Cyl}(n+n',1)    
\end{equation}
are given by insertion of configuration of little disks, and stacking two cylinders on top of each other while matching the $out$ of the first cylinder with the $in$ of the second cylinder (after a rotation), respectively; see Figure \ref{fig:figure3.5}. These should be viewed as the counterparts of the operadic compositions $\widehat{\theta^{\circlearrowright,0}_j},\widehat{\theta^{\circlearrowright,1}}$ of $\mathbf{KS}$, cf. Definition \autoref{thm:pictorial composition of cyclic cacti with spines}.\par\indent
Note that there is a homotopy equivalence 
\begin{equation}\label{eq:cyl as conf}
\mathrm{Cyl}(n,1)\simeq \mathrm{Conf}_n(\mathbb{C}^*)\times S^1,
\end{equation}
given as follows. Up to homotopy, we may assume everything lies on the standard cylinder $S^1\times [0,1]$ and up to $S^1$-rotation, we fix $in$ to be at $(0,1)\in S^1\times\{1\}$. Then, upon identifying $S^1\times (0,1)\cong \mathbb{C}^*$, we record the positions of the centers of the disks and the position of $out$, which gives the desired homotopy equivalence. 
The following comparison result will be crucial for our classification of Kontsevich-Soibelman operations.\\
\begin{thm}\label{thm:comparison of KS with Cyl}
$\mathbf{KS}$ is quasi-equivalent to $\mathrm{Cyl}$. Moreover, the quasi-equivalence is $\Sigma_n$-equivariant when restricted to $\mathbf{KS}(n,0)\simeq \mathrm{Cyl}(n,0)$ and $\Sigma_n\times S^1\times (S^1)^{op}$-equivariant when restricted to $\mathbf{KS}(n,1)\simeq \mathrm{Cyl}(n,1)$.
\end{thm}
In view of Theorem \autoref{thm:KS equals weighted cyclic cacti}, it suffices to prove a (suitably equivariant) quasi-equivalence of $\mathrm{Cyl}$ with the two-colored operad of weighted (cyclic) cacti, a task that we will defer to Section 6. The quasi-equivalence $\mathrm{Cyl}(l,0)\simeq \widetilde{Cact}^l$ is well-known. For instance, it was proved in \cite{MS} using Fiedorowicz's recognition principle and in \cite{Sal1} using flows of meromorphic differentials in the complex plane; in fact, it is the latter approach that we will adapt to give a proof of Theorem \autoref{thm:comparison of KS with Cyl}. 

\subsection{Kontsevich-Soibelman operations and the Getzler-Gauss-Manin connection}
In this subsection, we introduce a class of `natural'\footnote{In the aesthetic, not functorial, sense of the word.} operations on the periodic cyclic homology $HH_*^{per}(\mathcal{A})$ induced by the action of the Kontsevich-Soibelman operad on the pair $(CC^*(\mathcal{A}),CC_*(\mathcal{A}))$. We moreover show that they are automatically covariantly constant with respect to the Getzler-Gauss-Manin connection. Over an appropriate coefficient ring, these give rise to a large class of covariantly constant endomorphisms of $HH^{per}_*(\mathcal{A})$ which we will classify in Section 5.   \par\indent
Recall from Section 4.1 and 4.2 that there are, for $n\geq 0$, $\Sigma_n\times S^1\times S^1$-equivariant maps
\begin{equation}\label{eq:ksn1 action again}
\mathrm{Act}:\mathbf{KS}(n,1)_R\rightarrow\mathrm{Hom}_R(CC^*(\mathcal{A})^{\otimes n}\otimes CC_*(\mathcal{A}),CC_*(\mathcal{A}))
\end{equation}
that are moreover compatible with operadic structures. In particular, \eqref{eq:ksn1 action again} is $\Sigma_n\times S^1$-equivariant where $S^1$ acts diagonally via $S^1\rightarrow S^1\times S^1: x\mapsto (-x,x)$. We record the following elementary lemma which will be used freely later on. \\
\begin{lemma}\label{thm:lemma about fixed point and orbits}
Fix $G$ a compact Lie group, $X,Y$ $R$-linear chain complexes with $G$-actions, and $Z$ an $R$-linear chain complex with the trivial $G$-action. Then there are natural chain maps
\begin{equation}\label{eq:Fixpoint Hom}
\mathrm{Hom}_R(X,Y)^{hG}\rightarrow \mathrm{Hom}_{R^{hG}}(X^{hG},Y^{hG})  , 
\end{equation}
\begin{equation}\label{eq:Orbit Hom}
\mathrm{Hom}_R(X,Z)_{hG}\rightarrow \mathrm{Hom}_R(X^{hG},Z), 
\end{equation}
and 
\begin{equation}\label{eq:fixed point orbit to orbit}
X^{hG}\otimes Y_{hG}\rightarrow(X\otimes Y)_{hG}.    
\end{equation}
Moreover, \eqref{eq:Fixpoint Hom} is $R^{hG}$-equivariant and the following diagram commutes
\begin{equation}
\begin{tikzcd}[row sep=1.2cm, column sep=0.8cm]
\mathrm{Hom}_R(X,Y)^{hG}\otimes \mathrm{Hom}_R(Y,Z)_{hG}\arrow[d,"{\eqref{eq:fixed point orbit to orbit}}"]\arrow[rr,"{\eqref{eq:Fixpoint Hom}\otimes \eqref{eq:Orbit Hom}}"]& &\mathrm{Hom}_{R^{hG}}(X^{hG},Y^{hG})\otimes\mathrm{Hom}_R(Y^{hG},Z)\arrow[dd,"{\circ}"]\\
(\mathrm{Hom}_R(X,Y)\otimes\mathrm{Hom}_R(Y,Z))_{hG}\arrow[d,"{\circ_{hG}}"]&&\\
\mathrm{Hom}_R(X,Z)_{hG}\arrow[rr,"{\eqref{eq:Orbit Hom}}"]& & \mathrm{Hom}_R(X^{hG},Z)
\end{tikzcd}\label{eq:Hom orbit fixed point compatibility}
\end{equation}
\end{lemma}
\emph{Proof}. \eqref{eq:Fixpoint Hom} is just the induced map on homotopy fixed points of a homotopy $G$-equivariant chain map. For the other two maps, we fix a free resolution $Q\simeq R$ of right $R[G]$-modules (which we also view as left $R[G]$-modules via the inversion map $g\mapsto g^{-1}$) and a map of resolutions $\Delta:Q\rightarrow Q\otimes_R Q$ (which is unique up to homotopy).\par\indent
Then, \eqref{eq:Orbit Hom} is the map $\mathrm{Hom}_R(X,Z)_{hG}\otimes X^{hG}=$
\begin{equation}
Q\otimes_{R[G]}\mathrm{Hom}_R(X,Z)\otimes  \mathrm{Hom}_{R[G]}(Q,X)\rightarrow Z 
\end{equation}
given by
\begin{equation}
q\otimes F\otimes \varphi\mapsto (-1)^{|q|\cdot(|F|+|\varphi|))}F(\varphi(q)).    
\end{equation}
\eqref{eq:fixed point orbit to orbit} is the map $X^{hG}\otimes Y_{hG}=$
\begin{equation}
\mathrm{Hom}_{R[G]}(Q,X)\otimes Q\otimes_{R[G]}Y\rightarrow Q\otimes_{R[G]}(X\otimes Y)
\end{equation}
given by the composition 
\begin{equation}
\mathrm{Hom}_{R[G]}(Q,X)\otimes Q\otimes_{R[G]}Y\xrightarrow{\mathrm{id}\otimes \Delta\otimes\mathrm{id}} \mathrm{Hom}_{R[G]}(Q,X)\otimes (Q\otimes Q)\otimes_{R[G]}Y\rightarrow Q\otimes_{R[G]}(X\otimes Y),   
\end{equation}
where the second map above is given by
\begin{equation}
F\otimes q_1\otimes q_2 \otimes y\mapsto (-1)^{(|F|+|q_1|)\cdot|q_2|}q_2\otimes F(q_1) \otimes y.   
\end{equation}
One easily checks that these are well defined. The commutativity of \eqref{eq:Hom orbit fixed point compatibility} follows from straightforward diagram chasing. \qed

In what follows, we will save space by writing
\begin{equation}
(\mathbf{KS}(n,1)_R)^{hS^1}_{h\Sigma_n}:=((\mathbf{KS}(n,1)_R)_{h\Sigma_n})^{hS^1}.    
\end{equation}
Applying \eqref{eq:Orbit Hom} and then \eqref{eq:Fixpoint Hom} (note the $\Sigma_n$ and $S^1$ actions commute on both sides) to \eqref{eq:ksn1 action again}, one obtains an $R^{hS^1}$-equivariant map
$$\Xi:(\mathbf{KS}(n,1)_R)^{hS^1}_{h\Sigma_n}\rightarrow\mathrm{Hom}_R(CC^*(\mathcal{A})^{\otimes n}\otimes CC_*(\mathcal{A}),CC_*(\mathcal{A}))^{hS^1}_{h\Sigma_n}$$
$$\rightarrow \mathrm{Hom}_{R}((CC^*(\mathcal{A})^{\otimes n})^{h\Sigma_n}\otimes CC_*(\mathcal{A}),CC_*(\mathcal{A}))^{hS^1}$$
\begin{equation}\label{eq:ksn1 equiv}
\rightarrow \mathrm{Hom}_{R^{hS^1}}((CC^*(\mathcal{A})^{\otimes n})^{h\Sigma_n}\otimes CC_*(\mathcal{A})^{hS^1},CC_*(\mathcal{A})^{hS^1})
\end{equation}
Finally, taking cohomology of \eqref{eq:ksn1 equiv} gives an $H^*_{S^1}(pt;R)=R[[t]]$-equivariant map
\begin{equation}\label{eq:ksn1 equiv cohomology}
\Xi:H^*((\mathbf{KS}(n,1)_R)^{hS^1}_{h\Sigma_n})\rightarrow   \mathrm{Hom}_{H^*_{S^1}(pt;R)}(H^*_{\Sigma_n}(CC^*(\mathcal{A})^{\otimes n})\otimes HH_*^{-}(\mathcal{A}), HH_*^{-}(\mathcal{A})).
\end{equation}
At this stage, it is convenient to pass to the periodic cyclic homology. Let $t\in H^2_{S^1}(pt;R)$ denote the standard generator. By $H^*_{S^1}(pt;R)=R[[t]]$-equivariance, one may replace $HH_*^-$ by $HH_*^{per}=t^{-1}HH_*^-$ in \eqref{eq:ksn1 equiv cohomology}, in which case the resulting action factors through
\begin{equation}\label{eq:ksn1 equiv cohomology periodic}
\Xi:H^*((\mathbf{KS}(n,1)_R)^{tS^1}_{h\Sigma_n})=t^{-1}H^*((\mathbf{KS}(n,1)_R)^{hS^1}_{h\Sigma_n})\rightarrow \mathrm{Hom}_{H^*_{S^1}(pt;R)}(H^*_{\Sigma_n}(CC^*(\mathcal{A})^{\otimes n})\otimes HH_*^{per}(\mathcal{A}), HH_*^{per}(\mathcal{A})).
\end{equation}\par\indent
\begin{mydef}\label{thm:KS operations on HHper}
An endomorphism $L\in \mathrm{End}_{H^*_{S^1}(pt;R)}(HH^{per}_*(\mathcal{A}))$ is called a \emph{Kontsevich-Soibelman operation of arity $n$} (abbrev. \emph{$\mathrm{KS}^{n}$ operation}) if there exists $[\alpha]\in H^*((\mathbf{KS}(n,1)_R)^{tS^1}_{h\Sigma_n})$ and $[\phi]\in H^*_{\Sigma_n}(CC^*(\mathcal{A})^{\otimes n})$ such that
\begin{equation}
L=\Xi([\alpha])([\phi],-).    
\end{equation}
A \emph{Kontsevich-Soibelman operation} (abbrev. \emph{$\mathrm{KS}$ operation}) is an element of the sub $t^{-1}H^*_{S^1}(pt;R)$-algebra of $\mathrm{End}_{H^*_{S^1}(pt;R)}(HH^{per}_*(\mathcal{A}))$ generated by $\mathrm{KS}^{n}$, $n\geq 0$.\\ 
\end{mydef}
\begin{mydef}\label{thm:varphi KS operations on HHper}
Fix a cocycle $\varphi\in CC^{even}(\mathcal{A})$. An endomorphism $L\in \mathrm{End}_{H^*_{S^1}(pt;R)}(HH^{per}_*(\mathcal{A}))$ is called a \emph{$\varphi$-Kontsevich-Soibelman operation of arity $n$} (abbrev. \emph{$\mathrm{KS}^{n}_{\varphi}$ operation}) if there exists $[\alpha]\in H^*((\mathbf{KS}(n,1)_R)^{tS^1}_{h\Sigma_n})$ such that\footnote{It is well-known that $[\varphi^{\otimes n}]$ only depends on the cohomology class $[\varphi]$, though the assignment $[\varphi]\mapsto [\varphi^{\otimes n}]$ need not be additive.}
\begin{equation}
L=\Xi([\alpha])([\varphi^{\otimes n}],-).    
\end{equation}
\end{mydef}
\begin{rmk}
\begin{enumerate}[label=\arabic*)]
    \item It is possible that an endomorphism of $HH^{per}_*(\mathcal{A})$ is both a $\mathrm{KS}^{n}$ operation and a $\mathrm{KS}^m$  operation for $n\neq m$.
    \item The composition of a $\mathrm{KS}^{n}$ operation and a $\mathrm{KS}^{n'}$ operation may not be a $\mathrm{KS}^{n+n'}$ operation. However, there is a commutative diagram
\begin{equation}
\begin{tikzpicture}[commutative diagrams/every diagram]\label{eq:composition of KS operations diagram}
\node (P0) at (0,3) {$(\mathbf{KS}(n,1)_R)^{hS^1}_{h\Sigma_n}$};
\node (P1) at (0,2.5) {$\otimes\;(\mathbf{KS}(n',1)_R)^{hS^1}_{h\Sigma_n'}$} ;
\node (P2) at (10,3) {$\mathrm{Hom}_{R^{hS^1}}((CC^*(\mathcal{A})^{\otimes n})^{h\Sigma_n}\otimes CC_*(\mathcal{A})^{hS^1},CC_*(\mathcal{A})^{hS^1})$};
\node (P3) at (10,2.5) {$\otimes\;\mathrm{Hom}_{R^{hS^1}}((CC^*(\mathcal{A})^{\otimes n'})^{h\Sigma_{n'}}\otimes CC_*(\mathcal{A})^{hS^1},CC_*(\mathcal{A})^{hS^1})$};
\node (P4) at (0,-1) {$(\mathbf{KS}(n+n',1)_R)^{hS^1}_{h(\Sigma_n\times\Sigma_{n'})}$};
\node (P5) at (10,-2.5) {$\mathrm{Hom}_{R^{hS^1}}((CC^*(\mathcal{A})^{\otimes n+n'})^{h\Sigma_{n+n'}}\otimes CC_*(\mathcal{A})^{hS^1},CC_*(\mathcal{A})^{hS^1})$};
\node (P6) at (0,-2.5) {$(\mathbf{KS}(n+n',1)_R)^{hS^1}_{h\Sigma_{n+n'}}$};
\node (P7) at (10, -1) {$\mathrm{Hom}_{R^{hS^1}}((CC^*(\mathcal{A})^{\otimes n+n'})^{h(\Sigma_n\times\Sigma_{n'})}\otimes CC_*(\mathcal{A})^{hS^1},CC_*(\mathcal{A})^{hS^1})$};
\node (P8) at (10,1) {$\mathrm{Hom}_{R^{hS^1}}((CC^*(\mathcal{A})^{\otimes n})^{h\Sigma_{n}}\otimes (CC^*(\mathcal{A})^{\otimes n'})^{h\Sigma_{n'}}\otimes CC_*(\mathcal{A})^{hS^1},CC_*(\mathcal{A})^{hS^1})$};

\path[commutative diagrams/.cd, every arrow, every label]
(P1) edge node {$\Xi\,\otimes\,\Xi$} (P3)
(P6) edge node {$\Xi$} (P5)
(P1) edge node {$\eta$} (P4)
(P4) edge node {$\mathrm{Ind}_{\Sigma_n\times \Sigma_{n'}\subset \Sigma_{n+n'}}$} (P6)
(P7) edge node {$(\mathrm{Res}_{\Sigma_n\times \Sigma_{n'}\subset \Sigma_{n+n'}}\otimes\,\mathrm{id})\;\circ$} (P5)
(P4) edge node {$\Xi|_{\Sigma_n\times\Sigma_{n'}}$} (P7)
(P3) edge node {$\circ$} (P8)
(P7) edge node {($*$)} (P8);
\end{tikzpicture}
\end{equation}
where $\mathrm{Ind}, \mathrm{Res}$ denote taking further homotopy quotient/restricting homotopy fixed points along $\Sigma_n\times\Sigma_{n'}\subset \Sigma_{n+n'}$, respectively, and ($*$) is pre-composing with the natural map
   \begin{equation}\label{eq:tensor of homotopy fixed point to homotopy fixed point of tensor}
  (CC^*(\mathcal{A})^{\otimes n})^{h\Sigma_n}\otimes (CC^*(\mathcal{A})^{\otimes n'})^{h\Sigma_{n'}}\rightarrow(CC^*(\mathcal{A})^{\otimes n+n'})^{h(\Sigma_n\times\Sigma_{n'})}.
\end{equation}
In particular, when \eqref{eq:tensor of homotopy fixed point to homotopy fixed point of tensor} is a quasi-isomorphism, we may invert ($*$) in the derived category and the upper pentagon in \eqref{eq:composition of KS operations diagram} becomes a commutative square.
\item In contrast, fixing a $\varphi$, then the collection of $KS^n_{\varphi}$-operations ($n\geq 1$) is already closed under composition. This follows from diagram \eqref{eq:composition of KS operations diagram} and the fact that the restriction of $[\varphi^{\otimes n+n'}]$ to $(CC^*(\mathcal{A})^{\otimes n+n'})^{h(\Sigma_n\times\Sigma_{n'})}$ is exactly the image of $[\varphi^{\otimes n}]\otimes [\varphi^{\otimes n'}]$ under \eqref{eq:tensor of homotopy fixed point to homotopy fixed point of tensor}.
\end{enumerate} 
\end{rmk}

We now go back to the relative situation of a dg algebra $\mathcal{A}$ over a ring $R$ (considered as the base of a family) over a field $\mathbf{k}$ as in Section 3. Recall that a Kontsevich-Soibelman operation is said to be $\mathbf{k}$-linear if it is induced from a cohomology class in the image of  $H^*((\mathbf{KS}(n,1)_{\mathbf{k}})^{tS^1}_{h\Sigma_n})\rightarrow H^*((\mathbf{KS}(n,1)_{R})^{tS^1}_{h\Sigma_n})$ and an equivariant Hochschild cohomology class in the image of $H^*_{\Sigma_n}(CC^*(\mathcal{A}/R)^{\otimes_\mathbf{k} n})\rightarrow H^*_{\Sigma_n}(CC^*(\mathcal{A}/R)^{\otimes_R n})$. Below we use $CC_*(\mathcal{A})$ as shorthand for $CC_*(\mathcal{A}/R)$. \\
\begin{thm}\label{thm:automatic covariant constancy}
Any $\mathbf{k}$-linear Kontsevich-Soibelman operation $L\in \mathrm{End}_{H^*_{S^1}(pt;R)}(HH^{per}_*(\mathcal{A}))$ satisfies
\begin{equation}
L\circ \nabla^{GGM}=\nabla^{GGM}\circ L.    
\end{equation}
\end{thm}
\emph{Proof}. The proof is a straightforward consequence of Kaledin's formulation of the Getzler-Gauss-Manin connection \eqref{eq:GGM diagram}, which implies that, in the notation of loc.cit. a chain map $f: CC^{tS^1}_*(\mathcal{A})\rightarrow CC_*^{tS^1}(\mathcal{A})$ is covariantly constant with respect to $\nabla^{GGM}$ at the level of cohomology if there exists a chain map 
\begin{equation}f^{[2]}:I^0CC^{tS^1}_*({p_1}_*p_2^*\mathcal{A})/I^2CC^{tS^1}_*({p_1}_*p_2^*\mathcal{A})\rightarrow I^0CC^{tS^1}_*({p_1}_*p_2^*\mathcal{A})/I^2CC^{tS^1}_*({p_1}_*p_2^*\mathcal{A})
\end{equation}
making all the squares below commute (where all Hochschild chains are relative over $R$).
\begin{equation}
\begin{tikzcd}[row sep=1.2cm, column sep=0.8cm]
 I^0/I^2\arrow[dashed,dr,"{f^{[2]}}"]\arrow[r,"\pi"] \arrow[d,"m"]&CC^{tS^1}_*(\mathcal{A})\arrow[dr,"f"] &\\
{p_1}_*p_2^*CC^{tS^1}_*(\mathcal{A})\arrow[dr,"{{p_1}_*p_2^*f}"]  & I^0/I^2\arrow[r,"\pi"]\arrow[d,"m"]& CC^{tS^1}_*(\mathcal{A})\\
& {p_1}_*p_2^*CC^{tS^1}_*(\mathcal{A}) &
\end{tikzcd}\label{eq:GGM covariant constant diagram}
\end{equation}
Let $L=\Xi([\alpha])([\phi],-)$ be a $\mathbf{k}$-linear $\mathrm{KS}$ operation. Fix the explicit resolution 
\begin{equation}
\cdots\xrightarrow{\epsilon}\mathbf{k}[\epsilon]/\epsilon^2\xrightarrow{\epsilon}   \mathbf{k}[\epsilon]/\epsilon^2\rightarrow \mathbf{k}\rightarrow 0 
\end{equation}
of $\mathbf{k}$ as an $\mathbf{k}[\epsilon]/\epsilon^2\simeq \mathbf{k}[S^1]$-module, and some choice of free resolution $P_{\bullet}\rightarrow \mathbf{k}\rightarrow 0$ as right $\mathbf{k}[\Sigma_n]$-module. 
Then one can write down a chain level representative $\alpha$ of $[\alpha]$ of the form
\begin{equation}
\alpha=\sum_{i=-N}^{\infty} \alpha_i t^i,    
\end{equation}
where $\alpha_i\in P_{\bullet}\otimes_{\mathbf{k}[\Sigma_n]} \mathbf{KS}(n,1)_{\mathbf{k}}$ and a chain level representative $\phi\in \mathrm{Hom}_{\mathbf{k}[\Sigma_n]}(P_{\bullet},CC^*(\mathcal{A})^{\otimes_{\mathbf{k}} n})$ of $[\phi]$. $f=\Xi(\alpha)(\phi,-)$ gives a chain level endomorphism representing $L$. \par\indent
Recall that $R^{[2]}:=(R\otimes_{\mathbf{k}} R)/I^2$ is the first order neighborhood of the diagonal and $\mathcal{A}^{[2]}:={p_1}_*p_2^*\mathcal{A}=\mathcal{A}\otimes_RR^{[2]}$ is the base-changed dg algebra. Given a Hochschild cochain $\varphi \in CC^*(\mathcal{A}/R)$, let $\varphi^{[2]}\in CC^*(\mathcal{A}^{[2]}/R)$ denote the Hochschild cochain whose length $d$ part is the composition
\begin{equation}\label{eq:lift HH cochain}
(\mathcal{A}\otimes_{R}{R}^{[2]})^{\otimes d}\rightarrow (\mathcal{A}^{\otimes d})\otimes_{R}{R}^{[2]}\xrightarrow{\varphi\otimes_{R}{R}^{[2]}} \mathcal{A}\otimes_{R}{R}^{[2]}.     
\end{equation}
This induces a $\mathbf{k}$-linear map which sends a cocycle $\phi\in\mathrm{Hom}_{\mathbf{k}[\Sigma_n]}(P_{\bullet},CC^*(\mathcal{A})^{\otimes_{\mathbf{k}} n})$ to the cocycle $\phi^{[2]}\in\mathrm{Hom}_{\mathbf{k}[\Sigma_n]}(P_{\bullet},CC^*(\mathcal{A}^{[2]})^{\otimes_{\mathbf{k}} n})$. \par\indent
Set 
\begin{equation}
f^{[2]}:=\Xi(\alpha)(\phi^{[2]},-): I^0=CC^{tS^1}_*(\mathcal{A}^{[2]})\rightarrow CC^{tS^1}_*(\mathcal{A}^{[2]}).
\end{equation}
The fact that $f^{[2]}$ descends to a well-defined map $I^0/I^2\rightarrow I^0/I^2$ can be seen as follows. Since $\varphi^{[2]}$ is defined by extension of scalars \eqref{eq:lift HH cochain}, any Hochschild cochain of the form $\varphi^{[2]}$ satisfies $\varphi^{[2]}(I^1(\mathcal{A}^{[2]})^{\otimes_R d})\subset I^1\mathcal{A}^{[2]}$ and $\varphi^{[2]}(I^2(\mathcal{A}^{[2]})^{\otimes_R d})=0$ (since $I\subset R^{[2]}$ is square zero). As an immediate consequence, for any $x\in(\widehat{\mathfrak{Cact}^n_{\circlearrowright}})_{l;m_1,\cdots,m_n;m}$ and Hochschild cochains $\phi_i^{[2]}\;(1\leq i\leq n)$ of length $m_i$, the operation (cf. Figure \ref{fig:action_of_cyclic_cacti_with_spines_on_HH_co_chain})
\begin{equation}
\mathrm{Act}_x(-;\phi_1^{[2]},\cdots,\phi_n^{[2]}):(\mathcal{A}^{[2]})^{\otimes_R l+1}\rightarrow (\mathcal{A}^{[2]})^{\otimes_R m+1}   
\end{equation}
sends $I^i(\mathcal{A}^{[2]})^{\otimes_R l+1}$ to $I^i(\mathcal{A}^{[2]})^{\otimes_R m+1}$. This in particular implies that $f^{[2]}=\Xi(\alpha)(\phi^{[2]},-)$ also preserves the filtration $I^i$. \par\indent
By construction $f^{[2]}$ makes the diagram \eqref{eq:GGM covariant constant diagram} commute \footnote{Note however that if we replace $\mathbf{k}$-linear $\mathrm{KS}$ operations by $R$-linear ones, the left square in \eqref{eq:GGM covariant constant diagram} will not in general commute due to the switching of $R$-module structure in ${p_1}_*p_2^*$ (this should be intuitive as the property of being covariantly constant is not preserved by multiplication with an element in $R$).}. \qed

\section{Classification of Kontsevich-Soibelman operations}
The goal of this section is to identify topologically a set of generators for all KS operations on the periodic cyclic homology under operadic compositions. This is done by computing the (localized) equivariant cohomology of $\mathbf{KS}(n,1)$ and its interaction with operadic compositions. In the next section, we will write down explicit formulae for these generating operations. For the rest of the paper, we assume that the base ring $R$ is of pure characteristic, i.e. either $\mathbb{Q}\subset R$ or $\mathbb{F}_p\subset R$. 
\subsection{Equivariant homology of $\mathbf{KS}(n,1)$}
In light of Definition \autoref{thm:KS operations on HHper}, a first step towards understanding KS operations is a computation of the homology of
\begin{equation}\label{eq:equiv coh of Cyl}
(\mathbf{KS}(n,1)_R)^{tS^1}_{h\Sigma_n}\simeq C_{-*}(\mathrm{Cyl}(n,1);R)^{tS^1}_{h\Sigma_n},
\end{equation}
where the quasi-equivalence comes from Theorem \autoref{thm:comparison of KS with Cyl}. Note that the homotopy equivalence \eqref{eq:cyl as conf} is $\Sigma_n\times S^1$ equivariant, where $S^1$ acts on $\mathrm{Conf}_n(\mathbb{C}^*)\times S^1$ by rotating the ambient $\mathbb{C}^*$ clockwise in the first factor and trivially on the second factor. In particular, one further reduces \eqref{eq:equiv coh of Cyl} to computing the Tate equivariant homology
\begin{equation}
H^*(C_{-*}(\mathrm{UConf}_n(\mathbb{C}^*);R)^{tS^1})\otimes H_{-*}(S^1;R),
\end{equation}
where $\mathrm{UConf}_n$ (the \emph{unordered configuration space}) is defined as the free quotient $\mathrm{Conf}_n/\Sigma_n$. As will be clear from the proof of Theorem \autoref{thm:comparison of KS with Cyl}, the standard generator of $H^0(S^1;R)$ corresponds to the identity endomorphism, and the standard generator of $H^1(S^1;R)$ corresponds to postcomposition with the Connes operator (which is nullhomotopic as an endomorphism of the periodic cyclic complex). Thus it suffices to classify the operations coming from $H^*(C_{-*}(\mathrm{UConf}_n(\mathbb{C}^*);R)^{tS^1})$. \par\indent
\textbf{Orbit type classification}. We start by classifying all orbit types of the action 
\begin{equation}\label{eq:circle action on UConf}
S^1\acts \mathrm{UConf}_n(\mathbb{C}^*)
\end{equation}
induced by the ambient action $S^1\acts\mathbb{C}^*$ given by $e^{i\theta}\cdot z:=e^{-i\theta}z$. Suppose $\underline{z}=[(z_1,\cdots,z_n)]\in\mathrm{UConf}_n(\mathbb{C}^*)$ has nontrivial isotropy group $G_{\underline{z}}\subset S^1$, and let $e^{i\theta}\in G_{\underline{z}}$ be a non-identity element. Since $n$ is finite, there are coprime integers $(k,l),l\geq2$ such that $\theta=2\pi k/l, l\leq n$. We find an element $e^{i\theta_0}\in G_{\underline{z}}$ for which $l$ achieves its maximum $l_0$, and replacing $\theta_0$ by a multiple we may without loss of generality that $\theta_0=2\pi/l_0$. Hence for any $z\in \{z_1,\cdots,z_n\}$ we have $\{z,e^{-i\theta_0}z,\cdots,e^{-i(l_0-1)\theta_0}\}\subset \{z_1,\cdots,z_n\}$. By the maximality of $l_0$ and induction, we deduce that $G_{\underline{z}}=C_{l_0}\subset S^1, l_0|n$, and that $\underline{z}$ is a disjoint union of $n/l_0$ free $C_{l_0}$-orbits. \par\indent
An important tool in computing Tate equivariant homology is the localization theorem, which we recall below in both its topological and dg versions. \\
\begin{thm}\label{thm:localization}(Localization)
\begin{itemize}
    \item 
Let $G$ be a compact Lie group and $X$ be a $G$-space, paracompact
and with finite cohomological dimension. Then for any $s\in H^*(BG)$, the localized restriction morphism
\begin{equation}\label{eq:localized restriction}
s^{-1}H^*_G(X)\rightarrow s^{-1}H^*_G(X^s)    
\end{equation}
is an isomorphism, where $X^s=\{x\in X\,|\,\textrm{the image of $s$ under $H^*(BG)\rightarrow H^*(BG_x)$ is nonzero}\}$ and $G_x$ is the stabilizer of $x$.
\item Let $X$ be $C_p$-space such that the inclusion $X^{C_p}\hookrightarrow X$ is a good pair. Then there is a quasi-isomorphism
\begin{equation}
C_{-*}(X^{C_p};R)\otimes R^{tC_p}=C_{-*}(X^{C_p};R)^{tC_p}\xrightarrow{\simeq} C_{-*}(X;R)^{tC_p}.
\end{equation}
\end{itemize}
\end{thm}
\emph{Proof.} The first statement can be found in \cite[Theorem III.1']{Hsi}. For the second statement, consider the exact sequence 
\begin{equation}\label{eq:exact sequence of a pair}
C_{-*}(X^{C_p};R)\rightarrow  C_{-*}(X;R)\rightarrow \tilde{C}_{-*}(X/X^{C_p};R).
\end{equation}
The induced $C_p$-action on $X/X^{C_p}$ is free except for the basepoint; in particular, $\tilde{C}_{-*}(X/X^{C_p};R)$ is a complex of free $R[C_p]$-modules. Since the Tate construction is exact and vanishes on free $R[C_p]$-modules, the statement follows by applying $(-)^{tC_p}$ to \eqref{eq:exact sequence of a pair}.\qed \par\indent
For the computation of Tate equivariant homology of $\mathrm{UConf}_n(\mathbb{C}^*)$, we discuss two cases depending on the characteristic of the base ring $R$.
\subsubsection{Case I: $\mathbb{Q}\subset R$} 
In this case, by the orbit type classification, any isotropy group of the action \eqref{eq:circle action on UConf} is a finite cyclic group. In particular, in the context of Theorem \autoref{thm:localization}, $X^t=\emptyset$ and therefore $t^{-1}H^*_{S^1}(X;R)=0$. As an immediate consequence, we have the following.\\
\begin{lemma}\label{thm:rational KS operations are trivial}
Over a ground ring $R$ such that $\mathbb{Q}\subset R$, there are no nontrivial Kontsevich-Soibelman operations.    \qed
\end{lemma}

\subsubsection{Case II: $\mathbb{F}_p\subset R$}
\textbf{Step II.1.} It is convenient to first replace the $S^1$-action by the sub $C_p$-action (this will be justified below). Namely, restricting \eqref{eq:ksn1 equiv cohomology periodic} along $C_p\subset S^1$ gives an $H^*_{C_p}(pt;R)$-equivariant map
\begin{equation}\label{eq:ksn1 Cp equiv periodic}
\Xi^p:H^*((\mathbf{KS}(n,1)_R)^{tC_p}_{h\Sigma_n})\rightarrow \mathrm{Hom}_{H^*_{C_p}(pt;R)}(H^*_{\Sigma_n}(CC^*(\mathcal{A})^{\otimes n})\otimes HH_*^{C_p,per}(\mathcal{A}), HH_*^{C_p,per}(\mathcal{A})).
\end{equation}
\begin{mydef}\label{thm:C_p KS operations on HHper}
An endomorphism $L\in \mathrm{End}_{H^*_{C_p}(pt;R)}(HH^{C_p,per}_*(\mathcal{A}))$ is called a \emph{$C_p$-Kontsevich-Soibelman operation of arity $n$} (abbrev. \emph{$\mathrm{KS}^{n}_p$ operation}) if there exists $[\alpha]\in H^*((\mathbf{KS}(n,1)_R)^{tC_p}_{h\Sigma_n})$ and $[\phi]\in H^*_{\Sigma_n}(CC^*(\mathcal{A})^{\otimes n})$ such that
\begin{equation}
L=\Xi^p([\alpha])([\phi],-).    
\end{equation}
A \emph{$C_p$-Kontsevich-Soibelman operation} (abbrev. \emph{$\mathrm{KS}_p$ operation}) is an element of the sub $t^{-1}H^*_{C_p}(pt;R)$-algebra of $\mathrm{End}_{H^*_{C_p}(pt;R)}(HH^{C_p,per}_*(\mathcal{A}))$ generated by $\mathrm{KS}^{n}_p$ operations, $n\geq 0$. 
\end{mydef}
Fixing some $[\phi]\in H^*_{\Sigma_n}(CC^*(\mathcal{A})^{\otimes n})$, there is a commutative diagram whose vertical maps are restriction of homotopy/Tate fixed points along $C_p\subset S^1$.
\begin{equation}\label{eq:restricting ksn1 action}
\begin{tikzcd}[row sep=1.2cm, column sep=0.8cm]
H^*((\mathbf{KS}(n,1)_R)^{tS^1}_{h\Sigma_n})\otimes_{H^*_{S^1}(pt;R)} HH_*^{per}(\mathcal{A})\arrow[d,"{\mathrm{Res}_{C_p\subset S^1}}\;\otimes\; {\mathrm{Res}_{C_p\subset S^1}}"]\arrow[rrr,"{\Xi(-)([\phi],-)}"]& & & HH_*^{per}(\mathcal{A})\arrow[d,"{\mathrm{Res}_{C_p\subset S^1}}"]    \\
H^*((\mathbf{KS}(n,1)_R)^{tC_p}_{h\Sigma_n})\otimes_{H^*_{C_p}(pt;R)} HH_*^{C_p,per}(\mathcal{A})\arrow[rrr,"{\Xi^p(-)([\phi],-)}"]& & &HH_*^{C_p,per}(\mathcal{A})
\end{tikzcd}    
\end{equation}
We record the following lemma, which follows directly from the commutativity of diagram \eqref{eq:restricting ksn1 action}.\\
\begin{lemma}\label{thm:restriction of ksn1 action preserves circle fixed points}
If both $[\alpha]\in H^*((\mathbf{KS}(n,1)_R)^{tC_p}_{h\Sigma_n})$ and $x\in HH^{C_p,per}_*(\mathcal{A})$ lie in the image of ${\mathrm{Res}_{C_p\subset S^1}}$, then so does $\Xi^p([\alpha])([\phi],x)$, for any $[\phi]\in H^*_{\Sigma_n}(CC^*(\mathcal{A})^{\otimes n})$. \qed    
\end{lemma}
Similar as before, there is an analogue of diagram \eqref{eq:restricting ksn1 action} and Lemma \autoref{thm:restriction of ksn1 action preserves circle fixed points} for the twisted KS operations, which we omit.\par\indent
Finally, given an $R$-linear chain complex $C$ with $S^1$-action, the restriction $H^*_{S^1}(C)\rightarrow H^*_{C_p}(C)\cong H^*_{S^1}(C)\oplus H^*_{S^1}(C)\theta\;(|\theta|=1)$ is precisely the inclusion as the first factor (cf. the proof of Proposition \autoref{thm:pGysin for CC}). In particular, the vertical maps of \eqref{eq:restricting ksn1 action} are injections, and hence the KS operations are fully determined by the $\mathrm{KS}_p$ operations. \par\indent
\textbf{Step II.2.} By Theorem \autoref{thm:localization}, the inclusion of $C_p$-fixed points induces a quasi-isomorphism $C_*(\mathrm{UConf}_{pk}(\mathbb{C}^*)^{C_p};R)\otimes R^{tC_p}\cong C_*(\mathrm{UConf}_{pk}(\mathbb{C}^*);R)^{tC_p}$. We now give a computation of $H_*(\mathrm{UConf}_{pk}(\mathbb{C}^*)^{C_p};R)$. \par
\begin{figure}[H]
 \centering
 \includegraphics[width=0.9\textwidth]{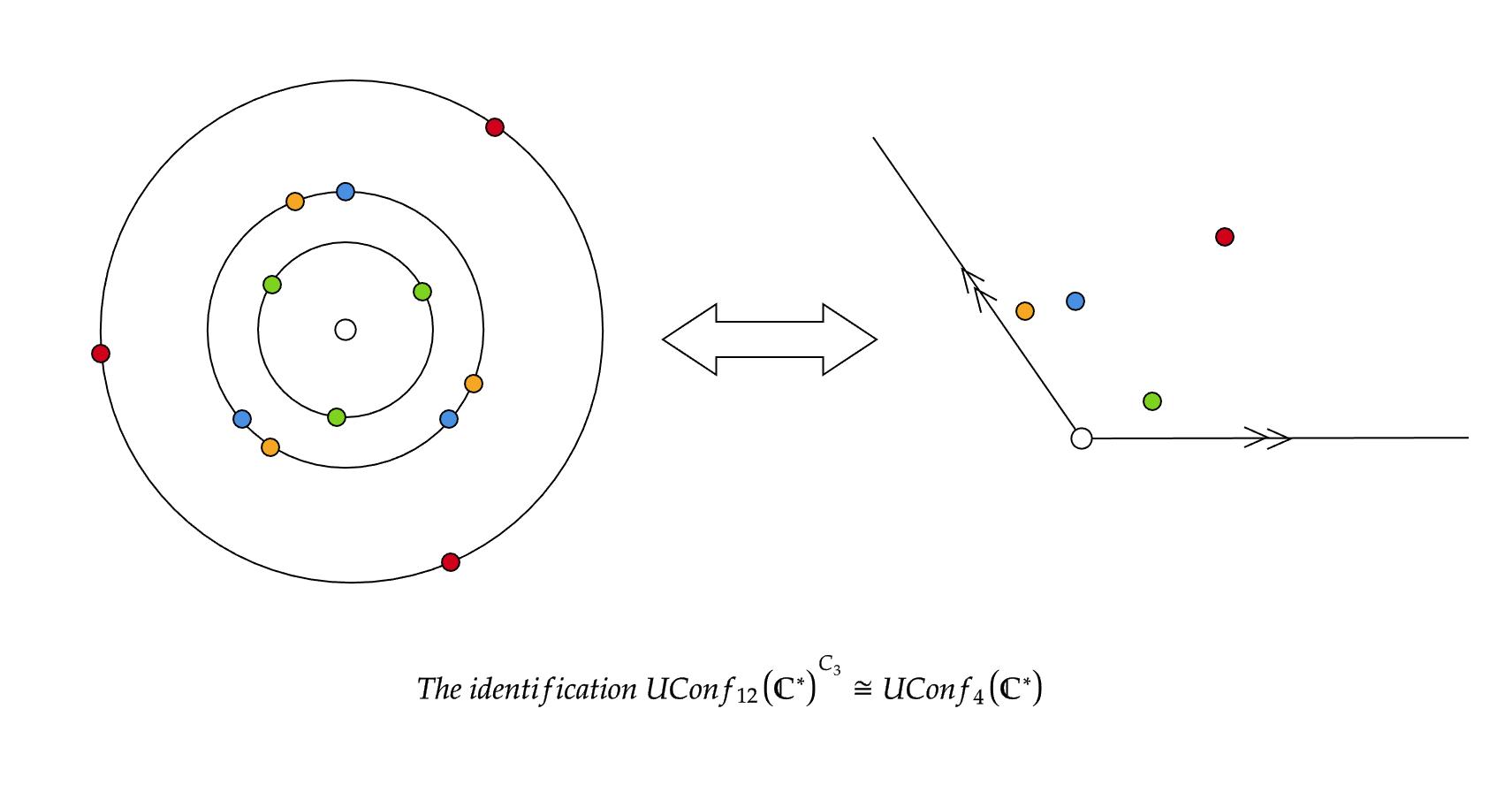}
 \caption{}
 \label{fig:figure4}
\end{figure}

\begin{figure}[H]
 \centering
 \includegraphics[width=0.9\textwidth]{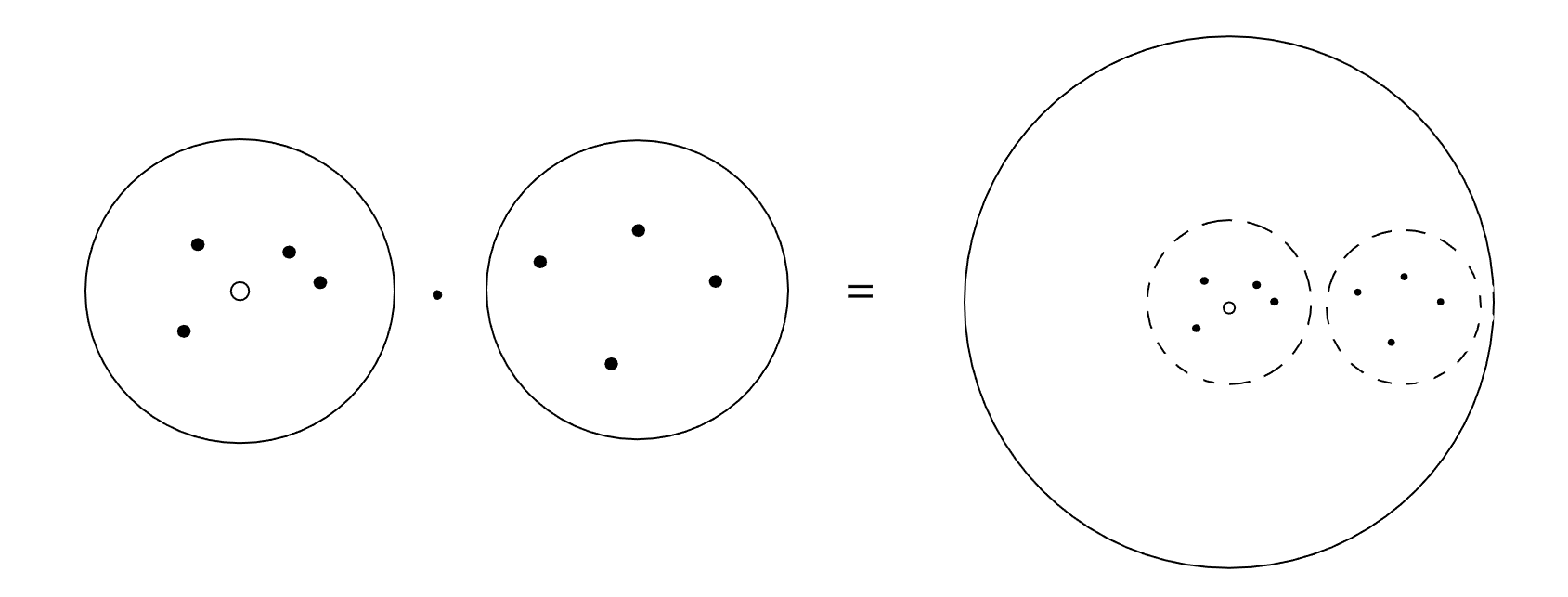}
 \caption{}
 \label{fig:figure5}
\end{figure}

\begin{figure}[H]
 \centering
 \includegraphics[width=0.9\textwidth]{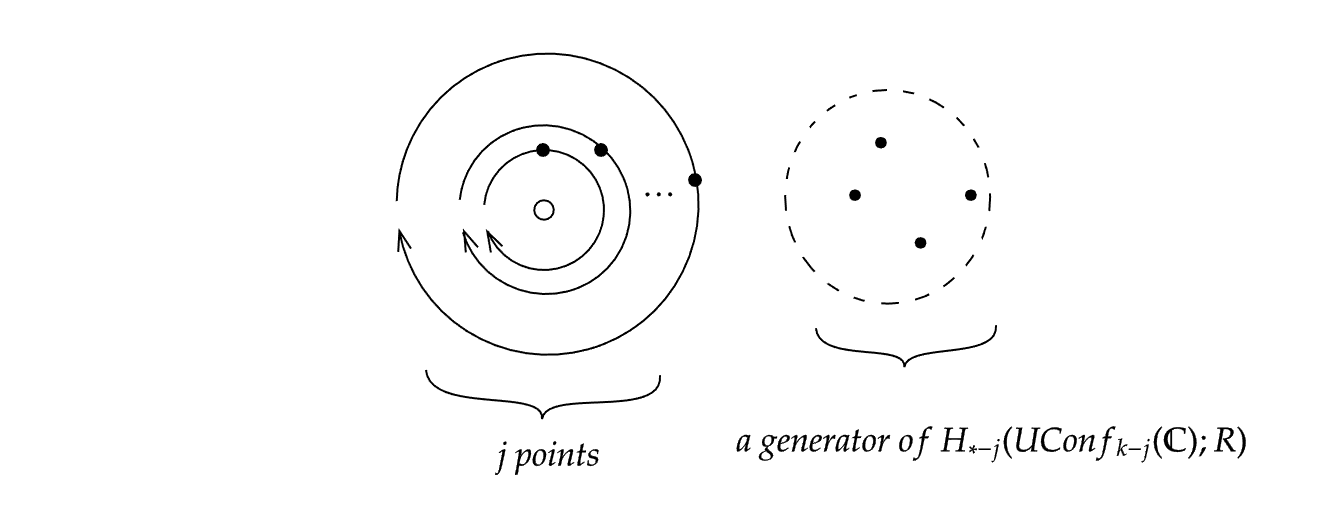}
 \caption{}
 \label{fig:figure6}
\end{figure}

\indent
Similar to \cite[Lemma 4.12]{Ros}, we have the following identification.\\
\begin{lemma}\label{thm:easy lemma}
There is a homeomorphism
\begin{equation}\label{eq:Cp fixed point is also UConf}
\mathrm{UConf}_{pk}(\mathbb{C}^*)^{C_p}\cong  \mathrm{UConf}_k(\mathbb{C}^*)   
\end{equation}    
\end{lemma}
\emph{Proof.} Let $\zeta=e^{2\pi i/p}$, and define
\begin{equation}
H=\{z\in\mathbb{C}^*:\arg(z)\in[0,2\pi/p]\}/(x\sim\zeta x,x\in\mathbb{R}_{>0})\end{equation}
The map $\mathrm{UConf}_k(H)\rightarrow \mathrm{UConf}_{pk}(\mathbb{C}^*)^{C_p}$ given by $\{z_1,\cdots,z_k\}\mapsto \{z_1,\zeta z_1,\cdots,\zeta^{p-1}z_1,\cdots,z_k,\zeta z_k,\cdots,\zeta^{p-1}z_k\}$, see Figure \ref{fig:figure4}, is a homeomorphism. Finally, observe that $H\cong \mathbb{C}^*$ via the $p$-fold map $z\mapsto z^p$.\qed
There is a product
\begin{equation}\label{E2 product}
\cdot:\mathrm{UConf}_k(\mathbb{C}^*)\times     \mathrm{UConf}_l(\mathbb{C})\rightarrow\mathrm{UConf}_{k+l}(\mathbb{C}^*)
\end{equation}
defined pictorially in Figure \ref{fig:figure5}, which comes from part of the $E_2$-structure on labeled configurations in $\mathbb{C}$, cf. Definition \autoref{thm:labeled configuration space}. 
Using Cohen's computation of the equivariant homology of configuration spaces \cite{Coh}, Rossi obtained the following description for $H_*(\mathrm{UConf}_k(\mathbb{C}^*);R)$. \\ 
\begin{lemma}(\cite[Corollary 4.15]{Ros})\label{thm:Rossi}
$H_*(\mathrm{UConf}_k(\mathbb{C}^*);R)$ is equal to the direct sum
\begin{equation}
H_*(\mathrm{UConf}_k(\mathbb{C});R)+ b^1\cdot H_{*-1}(\mathrm{UConf}_{k-1}(\mathbb{C});R)+\cdots +b^{k-1}\cdot H_1(\mathrm{UConf}_1(\mathbb{C});R)+b^k,
\end{equation}
where $b^j\in H_j(\mathrm{UConf}_j(\mathbb{C}^*);R)$ is the image of the fundamental class of $(S^1)^j$ under the map
\begin{equation}
(S^1)^j\rightarrow \mathrm{UConf}_j(\mathbb{C}^*): (\theta_1,\cdots,\theta_j)\mapsto \{e^{i\theta_1},2e^{i\theta_2},\cdots,je^{i\theta_j}\}.    
\end{equation}\qed
\end{lemma}
Figure \ref{fig:figure6} gives a pictorial description of the generators of $H_*(\mathrm{UConf}_k(\mathbb{C}^*);R)$.
Using the identification \eqref{eq:Cp fixed point is also UConf}, one can transfer Lemma \autoref{thm:Rossi} back to a statement about the $C_p$-equivariant cohomology of the Kontsevich-Soibelman operad (or equivalently, its topological partner $\mathrm{Cyl}$). We first fix some notations. \par\indent
Let $\mathrm{Cyl}^{\circ}(n,1)\subset\mathrm{Cyl}(n,1)$ denote the subspace with $in=(1,c)\in S^1\times\{c\}$ and $out=(1,a)\in S^1\times\{a\}$; in particular, $\mathrm{Cyl}^{\circ}(n,1)\times S^1\cong \mathrm{Cyl}(n,1)$. Let $\mathrm{UCyl}(n,0)=\mathrm{Cyl}(n,0)/\Sigma_n,\mathrm{UCyl}(n,1)=\mathrm{Cyl}(n,1)/\Sigma_n, \mathrm{UCyl}^{\circ}(n,1)=\mathrm{Cyl}^{\circ}(n,1)/\Sigma_n$ be the unordered versions. $(\mathrm{Cyl}(n,0),\mathrm{Cyl}^{\circ}(n,1))\subset (\mathrm{Cyl}(n,0),\mathrm{Cyl}(n,1))$ form a sub-operad, and there are homotopy equivalences
\begin{equation}
\mathrm{Cyl}(n,0)\simeq \mathrm{Conf}_n(\mathbb{C})\;,\quad \mathrm{UCyl}(n,0)\simeq \mathrm{UConf}_n(\mathbb{C})  
\end{equation} 
where the first equivalence is $\Sigma_n$-equivariant, and homotopy equivalences
\begin{equation}
\mathrm{Cyl}^{\circ}(n,1)\simeq \mathrm{Conf}_n(\mathbb{C}^*)\;,\quad \mathrm{UCyl}^{\circ}(n,1)\simeq \mathrm{UConf}_n(\mathbb{C}^*)  
\end{equation} 
which are $\Sigma_n\times S^1$ and  $S^1$-equivariant, respectively.\par\indent
Let $\lambda_{\Delta}$ denote the composition
\begin{equation}\label{eq:lambda delta}
\lambda_{\Delta}: \mathrm{Cyl}(l,0)\times \mathrm{Cyl}^{\circ}(p,1)\xrightarrow{\Delta\times \mathrm{id}} \mathrm{Cyl}(l,0)^{\times p}\times \mathrm{Cyl}^{\circ}(p,1)\xrightarrow{\prod_{i=1}^p\lambda_i}\mathrm{Cyl}^{\circ}(pl,1),
\end{equation}
where $\lambda_i$ is the first operadic composition in Figure \ref{fig:figure3.5}. Hence $\lambda_{\Delta}$ is a `simultaneous insertion of configurations of disks onto the cylinder'. $\lambda_{\Delta}$ is clearly $\Sigma_p\times \Sigma_l$-equivariant, where $\Sigma_p\times \Sigma_l$ acts on $\mathrm{Cyl}^{\circ}(pl,1)$ via the inclusion $\Sigma_p\times\Sigma_l\xrightarrow{\mathrm{id}\times \Delta}\Sigma_p\ltimes\Sigma_l^{\times p}\subset \Sigma_{pl}$. Hence $\lambda_{\Delta}$ induces a map
\begin{equation}\label{eq:lambda delta equiv}
\lambda_{\Delta}:\mathrm{UCyl}(l,0)\times \mathrm{UCyl}^{\circ}(p,1)\rightarrow \mathrm{Cyl}^{\circ}(pl,1)_{\Sigma_{p}\ltimes\Sigma_l^{\times p}}\rightarrow \mathrm{Cyl}^{\circ}(pl,1)_{\Sigma_{pl}}=\mathrm{UCyl}^{\circ}(pl,1),
\end{equation}
which by an abuse of notation is still denoted the same letter.\par\indent
On the other hand, the second operadic composition $\eta$ from Figure \ref{fig:figure3.5} is clearly $\Sigma_n\times \Sigma_{n'}\subset \Sigma_{n+n'}$-equivariant, and therefore induces a map 
\begin{equation}\label{eq:eta equiv}
\eta:\mathrm{UCyl}^{\circ}(n,1)\times\mathrm{UCyl}^{\circ}(n',1)\rightarrow \mathrm{Cyl}^{\circ}(n+n',1)_{\Sigma_{n}\times\Sigma_{n'}}\rightarrow \mathrm{Cyl}^{\circ}(n+n',1)_{\Sigma_{n+n'}}=\mathrm{UCyl}^{\circ}(n+n',1).   
\end{equation}

\begin{figure}[H]
 \centering
 \includegraphics[width=0.8\textwidth]{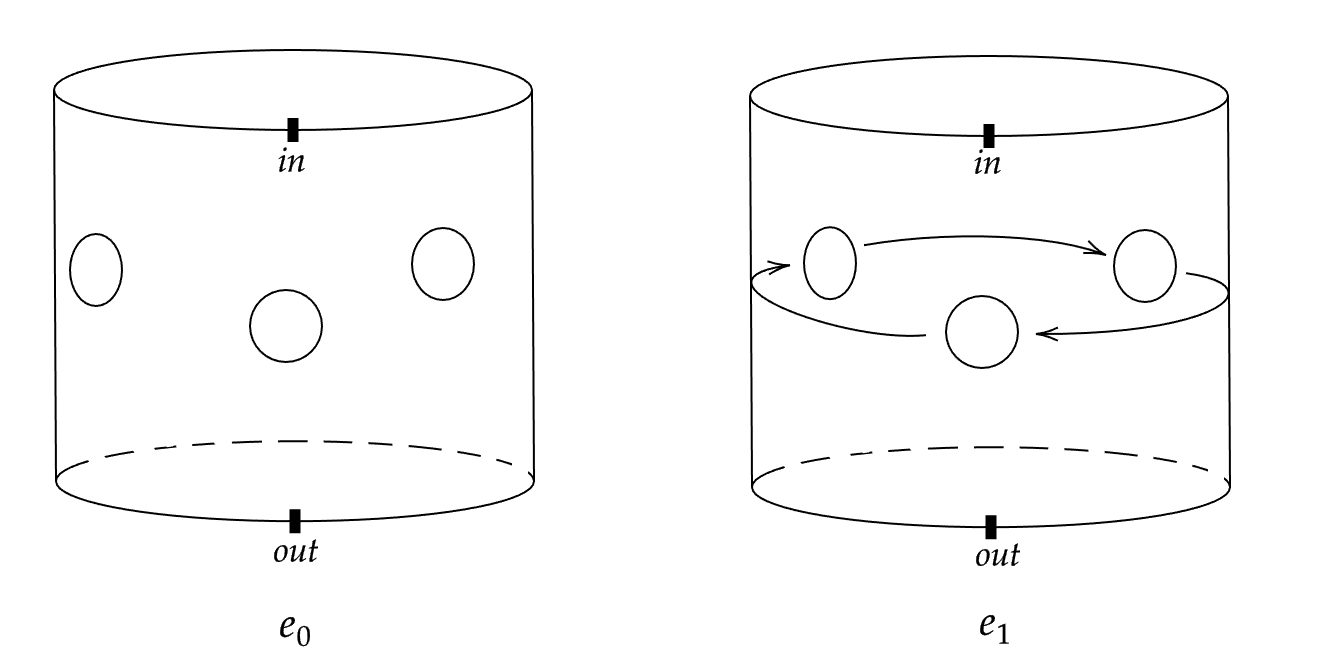}
 \caption{}
 \label{fig:figure7}
\end{figure}

\begin{lemma}\label{thm:generation of equiv cohomology of Cyl}
$\bigoplus_{k\geq 1}H^*(C_{-*}(\mathrm{UCyl}^{\circ}(pk,1);R)^{tC_p})$ is generated by $H^*(C_{-*}(\mathrm{UCyl}^{\circ}(p,1);R)^{tC_p})$ and $H_{-*}(\mathrm{UCyl}(l,0);R),l\geq 1$ via the (induced map on $C_p$-Tate cohomology of) following two types of operations:
\begin{equation}\label{eq:insert configuration}
\lambda_{\Delta}:\mathrm{UCyl}(l,0)\times \mathrm{UCyl}^{\circ}(p,1)\rightarrow \mathrm{UCyl}^{\circ}(pl,1)
\end{equation}
and
\begin{equation}\label{eq:stack cylinder}
\eta: \mathrm{UCyl}^{\circ}(pk,1)\times \mathrm{UCyl}^{\circ}(p,1)\rightarrow \mathrm{UCyl}^{\circ}(p(k+1),1).
\end{equation}
\end{lemma}
\emph{Proof}. The localization theorem implies that there is an isomorphism
\begin{equation}\label{eq:Cyl and UCyl, localization}
H^*(C_{-*}(\mathrm{Cyl}^{\circ}(pk,1);R)^{tC_p}_{h\Sigma_{pk}})\cong H_{-*}(\mathrm{UCyl}^{\circ}(pk,1)^{C_p};R)\otimes t^{-1}H_{C_p}^*(pt;R),    
\end{equation} 
which is furthermore compatible with operadic structures. From Lemma \autoref{thm:easy lemma}, one deduces that
\begin{equation}
H^*(C_{-*}(\mathrm{Cyl}^{\circ}(p,1);R)^{tC_p}_{h\Sigma_p})\cong H_{-*}(\mathrm{UCyl}^{\circ}(p,1)^{C_p})\otimes t^{-1}H^*_{C_p}(pt;R)=[e_0] t^{-1}H^*_{C_p}(pt;R)\oplus [e_1] t^{-1}H^*_{C_p}(pt;R),  
\end{equation}
where $[e_i]\in H_i(\mathrm{UCyl}^{\circ}(p,1)^{C_p})$ are the two generators drawn in Figure \ref{fig:figure7} (for $p=3$).
The proof of Lemma \autoref{thm:generation of equiv cohomology of Cyl} then follows from the picture below.
\begin{figure}[H]
 \centering
 \includegraphics[width=1.1\textwidth]{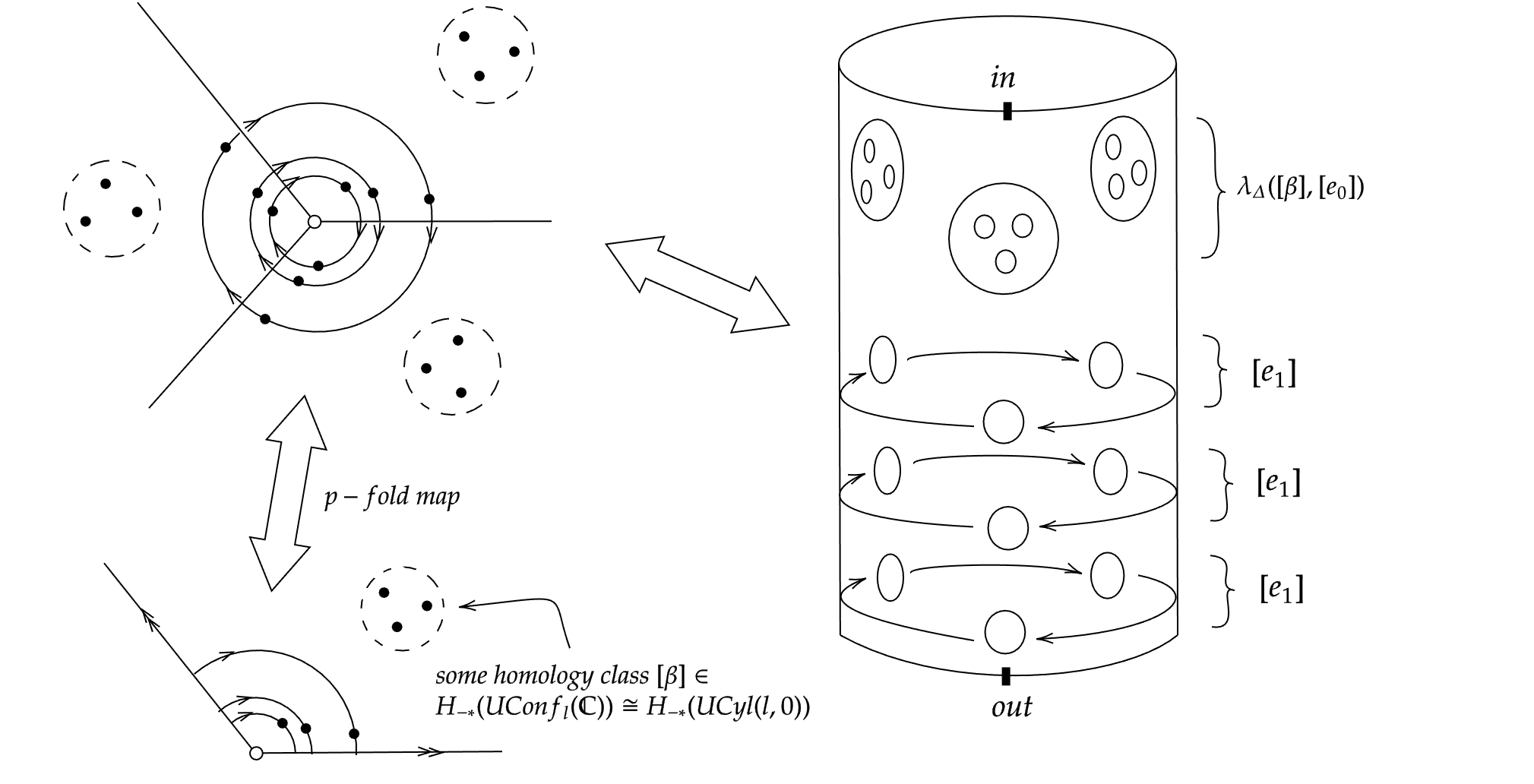}
 \caption{}
 \label{fig:figure8}
\end{figure}
We elaborate on Figure \ref{fig:figure8}. On top left we have a generator of $H_*(\mathrm{UConf}_{pk}(\mathbb{C}^*)^{C_p})$ obtained from a generator of $H_*(\mathrm{UConf}_{k}(\mathbb{C}^*))$ (bottom left; see also Figure \ref{fig:figure6}) under the `$p$-fold map' $H_*(\mathrm{UConf}_{k}(\mathbb{C}^*))\xrightarrow{\cong}H_*(\mathrm{UConf}_{pk}(\mathbb{C}^*)^{C_p})$ of Lemma \autoref{thm:easy lemma}. On the right we have the corresponding generator of $H_*(\mathrm{UCyl}^{\circ}(pk,1)^{C_p})$. From the right figure it is evident that one can first decompose such a generator via \eqref{eq:stack cylinder} (`stacking cylinders') into $[e_1]$'s followed by a $[e_0]$ at the top, where one may further insert a class $[\beta]\in H_*(\mathrm{UCyl}^{\circ}(l,0))$ diagonally into all the $p$ disks (cf. \eqref{eq:insert configuration}).\qed\par\indent
Before stating the main result of this section, we recall that the discussion at the beginning of Section 5.1 implies that all nontrivial ($C_p$-) Kontsevich Soibelman operations come from the subspace 
\begin{equation}
H^*(C_{-*}(\mathrm{Cyl}^{\circ}(n,1);R)^{hC_p}_{h\Sigma_n})\subset H^*(C_{-*}(\mathrm{Cyl}(n,1);R)^{hC_p}_{h\Sigma_n})\cong H^*((\mathbf{KS}(n,1)_R)^{hC_p}_{h\Sigma_n}).    
\end{equation}
We give the following notations for the operadic actions induced from the actions of $\mathbf{KS}$ in \eqref{eq:KS(n,0)_R and KS(n,1)_R action on Hocschild (co)chains, after taking totalization}:
\begin{equation}\label{eq:act from cyln}
\Phi_n: C_{-*}(\mathrm{Cyl}(n,0);R)\rightarrow \mathrm{Hom}_R(CC^*(\mathcal{A})^{\otimes n},CC^*(\mathcal{A})),  
\end{equation}
\begin{equation}\label{eq:act from cyln1}
\Phi_{n,1}: C_{-*}(\mathrm{Cyl}^{\circ}(n,1);R)\rightarrow \mathrm{Hom}_R(CC^*(\mathcal{A})^{\otimes n}\otimes CC_*(\mathcal{A}),CC_*(\mathcal{A})).
\end{equation}\par\indent
\begin{thm}\label{thm:classification of untwisted KSp operations}
\begin{enumerate}[label=\Roman*)]
    \item Assume that the natural map
\begin{equation}
(CC^*(\mathcal{A})^{\otimes k_1})^{h\Sigma_{k_1}}\otimes\cdots\otimes (CC^*(\mathcal{A})^{\otimes k_l})^{h\Sigma_{k_l}}\rightarrow (CC^*(\mathcal{A})^{\otimes k_1+\cdots+k_l})^{h(\Sigma_{k_1}\times \cdots\times\Sigma_{k_l})}     
\end{equation}
is a quasi-isomorphism for all $k_1,\cdots,k_l,l\geq 1$ and that the left hand side above satisfies the K\"{u}nneth isomorphism. Then all $C_p$-Kontsevich-Soibelman operations are generated under composition and linear combination by elements of the form
\begin{equation}\label{eq:e0 operation}
\Xi^p([e_0])([\phi],-)    
\end{equation}
and 
\begin{equation}\label{eq:e1 operation}
\Xi^p([e_1])([\phi],-),    
\end{equation}
where $[\phi]\in H^*_{\Sigma_p}(CC^*(\mathcal{A})^{\otimes p})$, and $[e_0],[e_1]$ are viewed as classes in $H^*(C_{-*}((\mathbf{KS}(p,1)_R)^{tC_p}_{h\Sigma_p})$ via the inclusions
\begin{equation}
H_{-*}(\mathrm{UCyl}^{\circ}(p,1)^{C_p};R)\subset H^*(C_{-*}(\mathrm{Cyl}^{\circ}(p,1);R)^{tC_p}_{h\Sigma_p})\subset H^*(C_{-*}(\mathrm{Cyl}(p,1);R)^{tC_p}_{h\Sigma_p})\cong H^*(C_{-*}(\mathbf{KS}(p,1)_R)^{tC_p}_{h\Sigma_p}).    
\end{equation}
\item  Fix cocycle $\varphi\in CC^{even}(\mathcal{A})$. All ($C_p$-)$\varphi$-Kontsevich-Soibelman operations are linear combination of elements of the form
\begin{equation}
\Xi^p([e_1])([\varphi^{\otimes p}],-)\circ\cdots\circ\Xi^p([e_1])([\varphi^{\otimes p}],-)\circ \Xi^p([e_0])([\phi],-),
\end{equation}
for some $[\phi]\in H^*_{\Sigma_p}(CC^*(\mathcal{A})^{\otimes p})$.
\end{enumerate}   
\end{thm}
\emph{Proof}. We discuss two separate scenarios. \par\indent
\textbf{Vertical decomposition}. Fix elements $[\alpha_1]\in H^*(C_{-*}(\mathrm{Cyl}^{\circ}(pk,1);R)^{tC_p}_{h\Sigma_{pk}})$, $[\alpha_2]\in H^*(C_{-*}(\mathrm{Cyl}^{\circ}(p,1);R)^{tC_p}_{h\Sigma_{p}})$ and $[\phi]\in H^*_{\Sigma_{p(k+1)}}(CC^*(\mathcal{A})^{\otimes p(k+1)})$.\par\indent
In the situation of I), the assumptions imply we can apply the Kunneth isomorphism and write 
\begin{equation}
\mathrm{Res}_{\Sigma_{pk}\times \Sigma_p\subset \Sigma_{p(k+1)}}[\phi]=  \sum_i[\phi_{1i}]\otimes[\phi_{2i}],\quad\mathrm{where}\,\;[\phi_{1i}]\in H^*_{\Sigma_{pk}}(CC^*(\mathcal{A})^{\otimes pk})\;\mathrm{and}\;[\phi_{2i}]\in H^*_{\Sigma_{p}}(CC^*(\mathcal{A})^{\otimes p}).
\end{equation}
Therefore, (the $C_p$-version of) diagram \eqref{eq:composition of KS operations diagram} implies that
\begin{equation}\label{eq:Xi decomposition vertical}
\Xi^p(\eta(\alpha_1,\alpha_2))(\phi,-)=\sum_i \pm \Xi^p(\alpha_2)(\phi_{2i})\circ \Xi^p(\alpha_1)(\phi_{1i})\end{equation}
at the level of homology.\par\indent
In the situation of II), since
\begin{equation}
\mathrm{Res}_{\Sigma_{pk}\times \Sigma_p\subset \Sigma_{p(k+1)}}[\varphi^{\otimes p(k+1)}]=[\varphi^{\otimes pk}]\otimes [\varphi^{\otimes p}],
\end{equation}
we have
\begin{equation}\label{eq:Xi decomposition vertical varphi}
\Xi^p(\eta(\alpha_1,\alpha_2))([\varphi^{\otimes p(k+1)}])=  \pm\Xi^p(\alpha_2)([\varphi^{\otimes p}])\circ \Xi^p(\alpha_1)([\varphi^{\otimes pk}]).    
\end{equation}
\textbf{Horizontal decomposition}. Fix elements $[\beta]\in H^{\Sigma_l}_{-*}(\mathrm{Cyl}(l,0);R)$, $[\alpha]\in H^*(C_{-*}(\mathrm{Cyl}^{\circ}(p,1);R)^{tC_p}_{h\Sigma_{p}})$, $[\phi]\in  H^*_{\Sigma_{pl}}(CC^*(\mathcal{A})^{\otimes pl})$. \par\indent
Since $\Phi_n, \Phi_{n,1}$ are compatible with operadic structures, we have the following commutative diagram (the ground ring $R$ will be omitted from the notations)
\begin{equation}\label{eq:Cyl big diagram}
\begin{tikzcd}[row sep=1.2cm, column sep=0.8cm]
C_{-*}(\mathrm{Cyl}(l,0))\otimes C_{-*}(\mathrm{Cyl}(p,1))\arrow[d,"{\lambda_{\Delta}}"]\arrow[rr,"{\Delta\otimes \mathrm{id}}"] &  & C_{-*}(\mathrm{Cyl}(l,0)^{\times p})\otimes C_{-*}(\mathrm{Cyl}(p,1))\arrow[d,"{\mathrm{EZ}^{-1}\otimes \mathrm{id}\footnotemark} "]\\
C_{-*}(\mathrm{Cyl}(pl,1))\arrow[d,"{\Phi_{pl,1}}"] & &C_{-*}(\mathrm{Cyl}(l,0))^{\otimes p}\otimes C_{-*}(\mathrm{Cyl}(p,1))\arrow[ll,"{\prod_{i=1}^p\lambda_i}"]\arrow[d,"{\Phi_l^{\otimes p}\otimes \Phi_{p,1}}"]\\
\mathrm{Hom}(CC^*(\mathcal{A})^{\otimes pl}\otimes  CC_*(\mathcal{A}), CC_*(\mathcal{A}))& & \mathrm{Hom}(CC^*(\mathcal{A})^{\otimes l},CC^*(\mathcal{A}))^{\otimes p}\otimes\mathrm{Hom}(CC^*(\mathcal{A})^{\otimes p}\otimes CC_*(\mathcal{A}), CC_*(\mathcal{A}))\arrow[ll,"{\circ}"]
\end{tikzcd}    
\end{equation}
\footnotetext{$\mathrm{EZ}^{-1}$ denotes a derived inverse to the Eilenberg-Zilber map and $(\mathrm{EZ}^{h\Sigma_p})^{-1}$ denotes a derived inverse to its induced map on homotopy fixed points. In particular, $(\mathrm{EZ}^{h\Sigma_p})^{-1}$ encodes the classical (dual) Steenrod operations.}
Upon applying suitable homotopy quotients/(Tate) fixed points to diagram \eqref{eq:Cyl big diagram}, we obtain diagram \eqref{eq:Cyl big diagram equiv}.\par\indent
The unlabeled arrows are induced by $\mathrm{Hom}(-,-)^{\otimes p}\rightarrow \mathrm{Hom}((-)^{\otimes p},(-)^{\otimes p})$; the arrows labeled $(\star\star)$ are induced by pre-composing with
\begin{equation}\label{eq:natural map involving fixed points of copies of CC^*}
((CC^*(\mathcal{A})^{\otimes l})^{h\Sigma_l)})^{\otimes p}\rightarrow     (CC^*(\mathcal{A})^{\otimes pl})^{h\Sigma_l^{\times p}}.
\end{equation}
The commutativity of the rightmost pentagon in \eqref{eq:Cyl big diagram equiv} follows from Lemma \autoref{thm:lemma about fixed point and orbits}, and the rest follows from either straightforward naturality or various compatibilities with operadic structures. \par\indent
By definition, the image of $\beta\otimes \alpha$ under the total composition of the left vertical maps in \eqref{eq:Cyl big diagram equiv} is $\Xi^p(\lambda_{\Delta}(\beta,\alpha))(-)$.\par\indent
In the situation of I), the assumptions imply that the bottom $(\star\star)$-arrow in \eqref{eq:Cyl big diagram equiv} has a derived inverse (indicated by a dashed arrow). Therefore, the commutativity of the diagram implies that 
\begin{equation}\label{eq:Xi decomposition horizontal}
 \Xi^p(\lambda_{\Delta}(\beta,\alpha))(\phi)=\Xi^p(\alpha)(\phi') 
\end{equation}
at the level of homology for 
\begin{equation}
[\phi']=\big(\Phi_l^{\otimes p}\circ(\mathrm{EZ}^{h\Sigma_p})^{-1}(\Delta(\beta))\big)([\phi]).   
\end{equation}
In the situation of II), even though the dashed arrow may not exist, the fact that  $\mathrm{Res}_{\Sigma_l^{\times p}\subset \Sigma_{pl}}[\varphi^{\otimes pl}]$ lies in the image of \eqref{eq:natural map involving fixed points of copies of CC^*} together with diagram \eqref{eq:Cyl big diagram equiv} still implies that 
\begin{equation}\label{eq:Xi decomposition horizontal varphi}
 \Xi^p(\lambda_{\Delta}(\beta,\alpha))([\varphi^{\otimes pl}])=\Xi^p(\alpha)([\phi]), 
\end{equation}
where 
\begin{equation}
[\phi]=\big(\Phi_l^{\otimes p}\circ(\mathrm{EZ}^{h\Sigma_p})^{-1}(\Delta(\beta))\big)([\varphi^{\otimes l}]^{\otimes p}).    
\end{equation}
Finally, combining \eqref{eq:Xi decomposition vertical} (resp. \eqref{eq:Xi decomposition vertical varphi}) and \eqref{eq:Xi decomposition horizontal} (resp. \eqref{eq:Xi decomposition horizontal varphi}) with Lemma \autoref{thm:generation of equiv cohomology of Cyl} gives the desired statements.\qed\par\indent
\clearpage
\newlength{\savedpw}\setlength{\savedpw}{\pdfpagewidth}
\newlength{\savedph}\setlength{\savedph}{\pdfpageheight}

\begin{figure}[H]
 \centering
 \includegraphics[width=0.9\textwidth]{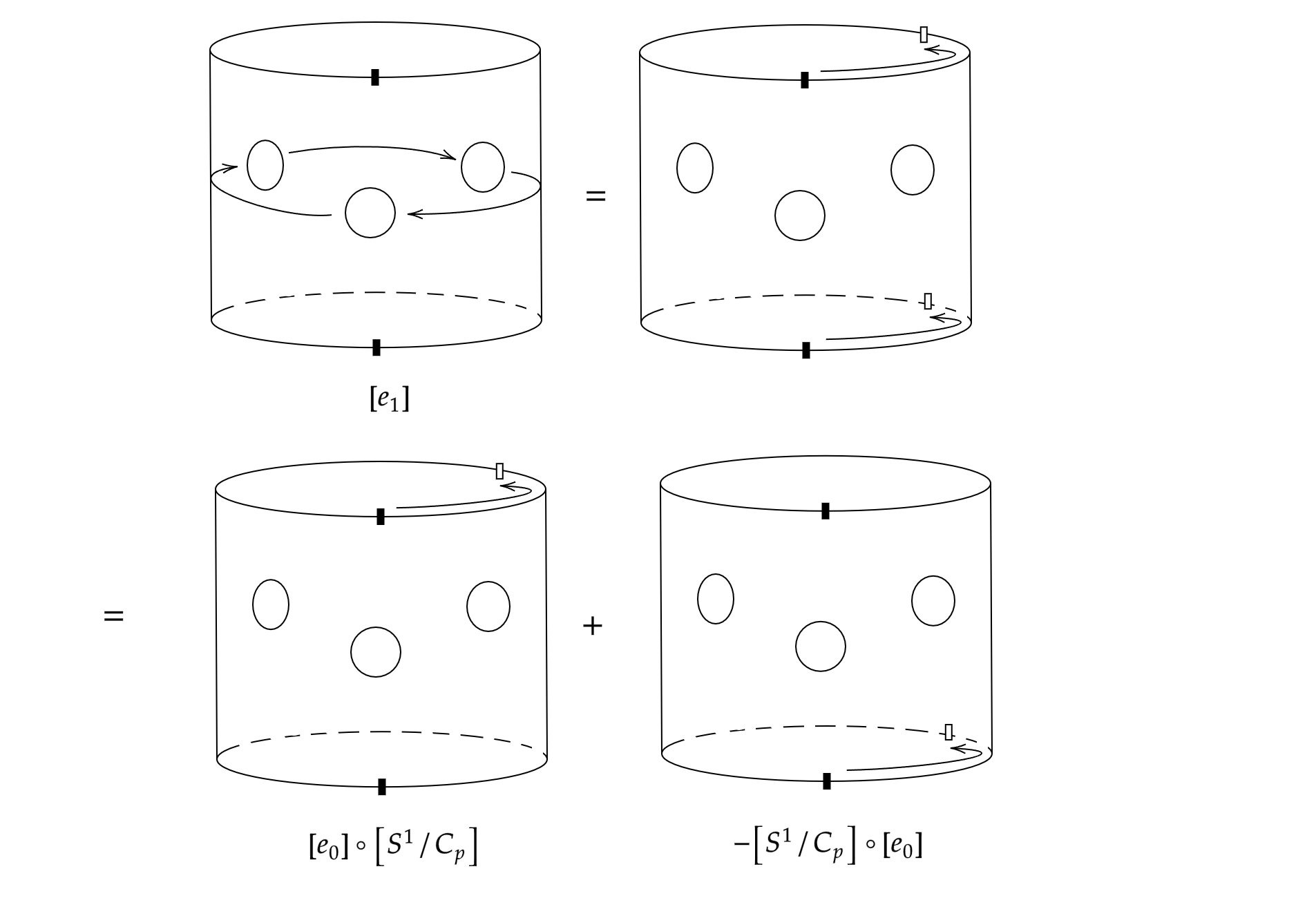}
 \caption{}
\label{fig:C_p_Cartan_formula}
\end{figure}

\pdfpageheight=16in      
\begin{landscape}
\newgeometry{margin=0.5in}
  \setlength{\textheight}{13in}   
  \thispagestyle{empty}
\begin{equation}\label{eq:Cyl big diagram equiv}
\hspace{-19cm}
\begin{tikzcd}[row sep=1.2cm, column sep=1.0cm]
C_{-*}(\mathrm{UCyl}(l,0))\otimes C_{-*}(\mathrm{UCyl}(p,1))^{tC_p}\arrow[rr,"{\Delta\otimes\mathrm{id}}"]\arrow[d,"{\lambda_{\Delta}}"] & & C_{-*}(\mathrm{UCyl}(l,0)^{\times p})^{h \Sigma_p}\otimes C_{-*}(\mathrm{UCyl}(p,1))^{tC_p}\arrow[d,"{(\mathrm{EZ}^{h\Sigma_p})^{-1}\otimes\mathrm{id}}"] & &\\
C_{-*}(\mathrm{Cyl}(pl,1))^{tC_p}_{h(\Sigma_p\ltimes\Sigma_l^{\times p})} \arrow[ddd,"{\Phi_{pl,1}}"]& &(C_{-*}(\mathrm{UCyl}(l,0))^{\otimes p})^{h\Sigma_p}\otimes C_{-*}(\mathrm{Cyl}(p,1))^{tC_p}_{h\Sigma_p}\arrow[d,"\eqref{eq:fixed point orbit to orbit}"]\arrow[r,"{\Phi_l^{\otimes p}\otimes \Phi_{p,1}}"]&\begin{array}{cc}
\big(\mathrm{Hom}((CC^*(\mathcal{A})^{\otimes l})^{h\Sigma_l},CC^*(\mathcal{A}))^{\otimes p}\big)^{h\Sigma_p}\otimes\\
\big(\mathrm{Hom}(CC^*(\mathcal{A})^{\otimes p}\otimes CC_*(\mathcal{A}),CC_*(\mathcal{A}))\big)^{tC_p}_{h\Sigma_p}    
\end{array}\arrow[d,"{\eqref{eq:fixed point orbit to orbit}}"]\arrow[r] & \begin{array}{cc}
\mathrm{Hom}\big(((CC^*(\mathcal{A})^{\otimes l})^{h\Sigma_l})^{\otimes p},CC^*(\mathcal{A})^{\otimes p}\big)^{h\Sigma_p}\otimes\\
\big(\mathrm{Hom}(CC^*(\mathcal{A})^{\otimes p}\otimes CC_*(\mathcal{A}),CC_*(\mathcal{A}))\big)^{tC_p}_{h\Sigma_p}    
\end{array}\arrow[ddl,"{\eqref{eq:fixed point orbit to orbit}}"]\arrow[dd,"{\eqref{eq:Fixpoint Hom}\;\textrm{and}\;\eqref{eq:Orbit Hom}}"]\\
& & (C_{-*}(\mathrm{UCyl}(l,0))^{\otimes p})\otimes C_{-*}(\mathrm{Cyl}(p,1))^{tC_p}_{h\Sigma_p}\arrow[d,"{\cong}"]]\arrow[r,"{\Phi_l^{\otimes p}\otimes \Phi_{p,1}}"]&\begin{array}{cc}
\big(\mathrm{Hom}((CC^*(\mathcal{A})^{\otimes l})^{h\Sigma_l},CC^*(\mathcal{A}))^{\otimes p}\otimes\\
\mathrm{Hom}(CC^*(\mathcal{A})^{\otimes p}\otimes CC_*(\mathcal{A}),CC_*(\mathcal{A}))\big)^{tC_p}_{h\Sigma_p}    
\end{array}\arrow[d] & \\
& &\big(C_{-*}(\mathrm{Cyl}(l,0))^{\otimes p}\otimes C_{-*}(\mathrm{Cyl}(p,1))\big)^{tC_p}_{h(\Sigma_p\ltimes \Sigma_l^{\times p})}\arrow[uull,"{\prod_{i=1}^p\lambda_i}"]\arrow[d,"{\Phi_l^{\otimes p}\otimes \Phi_{p,1}}"]& \begin{array}{cc}
\big(\mathrm{Hom}(((CC^*(\mathcal{A})^{\otimes l})^{h\Sigma_l})^{\otimes p},CC^*(\mathcal{A})^{\otimes p})\otimes\\
\mathrm{Hom}(CC^*(\mathcal{A})^{\otimes p}\otimes CC_*(\mathcal{A}),CC_*(\mathcal{A}))\big)^{tC_p}_{h\Sigma_p}    
\end{array}\arrow[d,"{\circ}"]&\begin{array}{cc}
\mathrm{Hom}\big((((CC^*(\mathcal{A})^{\otimes l})^{h\Sigma_l})^{\otimes p})^{h\Sigma_p},(CC^*(\mathcal{A})^{\otimes p})^{h\Sigma_p}\big)\otimes\\
\big(\mathrm{Hom}((CC^*(\mathcal{A})^{\otimes p})^{h\Sigma_p}\otimes CC_*(\mathcal{A}),CC_*(\mathcal{A}))\big)^{tC_p}
\end{array}\arrow[ddl,bend left=10, "{\circ}"] \\
\mathrm{Hom}(CC^*(\mathcal{A})^{\otimes pl}\otimes CC_*(\mathcal{A}), CC_*(\mathcal{A}))^{tC_p}_{h(\Sigma_p\ltimes \Sigma_l^{\times p})}\arrow[dd,"{\eqref{eq:Fixpoint Hom}\;\mathrm{and}\;\eqref{eq:Orbit Hom}}"]& &
{\begin{array}{l}
\big(\mathrm{Hom}(CC^*(\mathcal{A})^{\otimes l}, CC^*(\mathcal{A}))^{\otimes p}\otimes \\
\mathrm{Hom}(CC^*(\mathcal{A})^{\otimes p}\otimes CC_*(\mathcal{A}), CC_*(\mathcal{A}))\big)^{tC_p}_{h(\Sigma_p\ltimes \Sigma_l^{\times p})}\arrow[ll,"{\circ}"]
\end{array}}\arrow[d,"{\eqref{eq:Orbit Hom}}"]& \begin{array}{cc}
\mathrm{Hom}\big(((CC^*(\mathcal{A})^{\otimes l})^{h\Sigma_l})^{\otimes p}\otimes CC_*(\mathcal{A}),CC_*(\mathcal{A})\big)^{tC_p}_{h\Sigma_p}    
\end{array}\arrow[d,"{\eqref{eq:Orbit Hom}}"] &\\
 & & \begin{array}{cc}
\big(\mathrm{Hom}((CC^*(\mathcal{A})^{\otimes pl})^{h(\Sigma_l^{\times p})},CC^*(\mathcal{A})^{\otimes p})\otimes\\
\mathrm{Hom}(CC^*(\mathcal{A})^{\otimes p}\otimes CC_*(\mathcal{A}),CC_*(\mathcal{A}))\big)^{tC_p}_{h\Sigma_p}    
\end{array}\arrow[d,"{\circ}"]\arrow[uur, "(\star\star)"]& \begin{array}{cc}
\mathrm{Hom}\big((((CC^*(\mathcal{A})^{\otimes l})^{h\Sigma_l})^{\otimes p})^{h\Sigma_p}\otimes CC_*(\mathcal{A}),CC_*(\mathcal{A})\big)^{tC_p}\end{array}\arrow[dd,"{\eqref{eq:Fixpoint Hom}}"]&\\
\mathrm{Hom}((CC^*(\mathcal{A})^{\otimes pl})^{h(\Sigma_p\ltimes\Sigma_l^{\times p})}\otimes CC_*(\mathcal{A})^{tC_p}, CC_*(\mathcal{A})^{tC_p}) \arrow[drrr, bend right=5,  "{(\star\star)}"]\arrow[d,"{\mathrm{Res}_{\Sigma_p\ltimes \Sigma_l^{\times p}\subset \Sigma_{pl}}\otimes\mathrm{id}\;\circ}"] & & \mathrm{Hom}((CC^*(\mathcal{A})^{\otimes pl})^{h(\Sigma_l^{\times p})}\otimes CC_*(\mathcal{A}),CC_*(\mathcal{A}))^{tC_p}_{h\Sigma_p}\arrow[ll,"{\eqref{eq:Fixpoint Hom}\;\textrm{and}\;\eqref{eq:Orbit Hom}}"]\arrow[uur,"{(\star\star)}"] & \\
\mathrm{Hom}((CC^*(\mathcal{A})^{\otimes pl})^{h\Sigma_{pl}}\otimes CC_*(\mathcal{A})^{tC_p}, CC_*(\mathcal{A})^{tC_p}) & & &\begin{array}{cc}
\mathrm{Hom}\big((((CC^*(\mathcal{A})^{\otimes l})^{h\Sigma_l})^{\otimes p})^{h\Sigma_p}\otimes CC_*(\mathcal{A})^{tC_p},CC_*(\mathcal{A})^{tC_p}\big)\arrow[ulll, dashed]
\end{array}   
\end{tikzcd}
\end{equation}
\end{landscape}
\clearpage
\global\pdfpagewidth=\paperwidth
\global\pdfpageheight=\paperheight
\clearpage
On the other hand, there is a residual $S^1\cong S^1/C_p$-action on $CC^{C_p,per}_*(\mathcal{A})$; let 
\begin{equation}
B^p:HH_*^{C_p,per}(\mathcal{A})\rightarrow HH_{*-1}^{C_p,per}(\mathcal{A})    
\end{equation}
be the induced endomorphism corresponding to the fundamental class $[S^1/C_p]$. Consider the relation indicated in Figure \ref{fig:C_p_Cartan_formula}
Intuitively, $\Xi^p([e_1])([\phi],-)$ can be viewed as a `$C_p$-equivariant Lie action' by $[\phi]$, whereas Figure \ref{fig:C_p_Cartan_formula} can be viewed as a `$C_p$-equivariant (noncommutative) Cartan formula'. This observation, combined with Theorem \autoref{thm:classification of untwisted KSp operations}, immediately gives the following corollary.\\
\begin{cor}\label{thm:KSp operations are generated by e_0 operations, up to bracketing with S^1=S^1/C_p action}
Under the same assumption and notation as Theorem \autoref{thm:classification of untwisted KSp operations}, all $C_p$-Kontsevich-Soibelman operations are generated under composition and linear combination by elements of the form
\begin{equation}\label{eq:e0 operation}
\Xi^p([e_0])([\phi],-)\quad\mathrm{and}\quad [\Xi^p([e_0])([\phi],-),B^p],
\end{equation}
where $[\phi]\in H^*_{\Sigma_p}(CC^*(\mathcal{A}))$.\qed
\end{cor}
In particular, the understanding of ($C_p$-)Kontsevich-Soibelman operations reduces to understanding a single type of endomorphisms, namely, $\Xi^p([e_0])([\phi],-)$. We will complete this task in Section 6 by giving explicit formulae for these operations.  \par\indent
To conclude this section, we use the previous classification of $C_p$-Kontsevich-Soibelman operations to obtain a classification of ($S^1$-)Kontsevich-Soibelman operation on the periodic cyclic homology. First, note that by $C_p$-localization and Corollary \ref{thm:equivalence of S^1 vs C_p fixed points}, the natural map induced by inclusion of $C_p$-fixed points
\begin{equation}
C_{-*}(\mathrm{UCyl}^{\circ}(pk,1)^{C_p};R)^{tS^1} \rightarrow C_{-*}(\mathrm{UCyl}^{\circ}(pk,1);R)^{tS^1}    
\end{equation}
is also a quasi-isomorphism. Moreover, since the induced $S^1$-action on $\mathrm{UCyl}^{\circ}(pk,1)^{C_p}$ factors through the $p$-fold map $S^1\xrightarrow{\times p} S^1$ and $\mathbb{F}_p\subset R$, there is a quasi-isomorphism 
\begin{equation}\label{eq:tS^1 for C_p fixed points of UCyl^circ is trivial}
C_{-*}(\mathrm{UCyl}^{\circ}(pk,1)^{C_p};R)^{tS^1} \simeq (C_{-*}(\mathrm{UCyl}^{\circ}(pk,1)^{C_p};R)((t)),d),
\end{equation}
where $d$ denotes the differential of $C_{-*}(\mathrm{UCyl}^{\circ}(pk,1)^{C_p};R)$. Moreover, the quasi-isomorphism of \eqref{eq:tS^1 for C_p fixed points of UCyl^circ is trivial} is compatible with operadic compositions. Therefore, using the exact same idea as Lemma \ref{thm:generation of equiv cohomology of Cyl} and the decomposition of $H_{-*}(\mathrm{UCyl}^{\circ}(pk,1)^{C_p};R)$ as in Figure \ref{fig:figure8} gives the following lemma.\\
\begin{lemma}\label{thm:generation of S^1 equiv cohomology of Cyl}
$\bigoplus_{k\geq 1}H^*(C_{-*}(\mathrm{UCyl}^{\circ}(pk,1);R)^{tS^1})$ is generated by the two elements $[\epsilon_0],[\epsilon_1]\in H^*(C_{-*}(\mathrm{UCyl}^{\circ}(p,1);R)^{tS^1})$ and $H_{-*}(\mathrm{UCyl}(l,0);R),l\geq 1$ via the (induced map on $S^1$-Tate cohomology of) operadic compositions $\lambda_{\Delta}$ and $\eta$ (cf. \eqref{eq:lambda delta equiv}, \eqref{eq:eta equiv}); where $[\epsilon_0],[\epsilon_1]$ are the two generators $[e_0],[e_1]$ of $H_{-*}(\mathrm{UCyl}^{\circ}(p,1)^{C_p};R)$ indicated in Figure \ref{fig:figure7}, \emph{but viewed as elements of $H^*(C_{-*}(\mathrm{UCyl}^{\circ}(p,1);R)^{tS^1})$ via the inclusion}
\begin{equation}
H_{-*}(\mathrm{UCyl}^{\circ}(p,1)^{C_p};R)\xrightarrow{\textrm{include as constant term}} H_{-*}(\mathrm{UCyl}^{\circ}(p,1)^{C_p};R)((t))\cong    H^*(C_{-*}(\mathrm{UCyl}^{\circ}(p,1);R)^{tS^1}). 
\end{equation}\qed
\end{lemma}
Then, in exact parallel to Theorem \ref{thm:classification of untwisted KSp operations} we obtain\\
\begin{thm}\label{thm:classification of untwisted KS operations}
\begin{enumerate}[label=\Roman*)]
    \item Assume that the natural map
\begin{equation}
(CC^*(\mathcal{A})^{\otimes k_1})^{h\Sigma_{k_1}}\otimes\cdots\otimes (CC^*(\mathcal{A})^{\otimes k_l})^{h\Sigma_{k_l}}\rightarrow (CC^*(\mathcal{A})^{\otimes k_1+\cdots+k_l})^{h(\Sigma_{k_1}\times \cdots\times\Sigma_{k_l})}     
\end{equation}
is a quasi-isomorphism for all $k_1,\cdots,k_l,l\geq 1$ and that the left hand side above satisfies the K\"{u}nneth isomorphism. Then all Kontsevich-Soibelman operations on $HH^{per}_*(\mathcal{A})$ are generated under composition and linear combination by elements of the form
\begin{equation}\label{eq:e0 operation}
\Xi([\epsilon_0])([\phi],-)    
\end{equation}
and 
\begin{equation}\label{eq:e1 operation}
\Xi([\epsilon_1])([\phi],-),    
\end{equation}
where $[\phi]\in H^*_{\Sigma_p}(CC^*(\mathcal{A})^{\otimes p})$.
    \item Fix cocycle $\varphi\in CC^{even}(\mathcal{A})$. All $\varphi$-Kontsevich-Soibelman operations are linear combination of elements of the form
\begin{equation}
\Xi([e_1])([\varphi^{\otimes p}],-)\circ\cdots\circ\Xi([e_1])([\varphi^{\otimes p}],-)\circ \Xi([e_0])([\phi],-),
\end{equation}
for some $[\phi]\in H^*_{\Sigma_p}(CC^*(\mathcal{A})^{\otimes p})$.
\end{enumerate}\qed
\end{thm}
On the other hand, how does $\Xi([\epsilon_0])([\phi],-) $ and $\Xi([\epsilon_1])([\phi],-)$ act on $HH_*^{per}(\mathcal{A})$ explicitly?\par\indent  
As a preliminary step in obtaining such explicit formula, we pass to $(-)^{tC_p}$ again and study how the generators $[\epsilon_0],[\epsilon_1]\in H^0(C_{-*}(\mathrm{UCyl}^{\circ}(p,1);R)^{tS^1})$ restrict to $H^*(C_{-*}(\mathrm{UCyl}^{\circ}(p,1);R)^{tC_p})$. \par\indent
Recall that there is an $S^1$-equivariant homeomorphism $\mathrm{UCyl}^{\circ}(p,1)\cong \mathrm{UConf}_p(\mathbb{C}^*)$.\\
\begin{lemma}\label{thm:simplifying S^1 Tate fixed points of UConf_p}
The inclusion $\mathrm{UConf}_{p}(\mathbb{C}^*)^{C_p}\hookrightarrow \mathrm{UConf}_{p}(\mathbb{C}^*)$ induces a quasi-isomorphism
\begin{equation}
 C_{-*}(\mathrm{UConf}_{p}(\mathbb{C}^*)^{C_p};R)^{tS^1}\xrightarrow{\simeq} C_{-*}(\mathrm{UConf}_{p}(\mathbb{C}^*);R)^{tS^1}.   
\end{equation}
\end{lemma} 
\emph{Proof}. By $C_p$-localization (cf. Theorem \ref{thm:localization}) the induced map
\begin{equation}
  C_{-*}(\mathrm{UConf}_{p}(\mathbb{C}^*)^{C_p};R)^{tC_p}\xrightarrow{\simeq} C_{-*}(\mathrm{UConf}_{p}(\mathbb{C}^*);R)^{tC_p}  
\end{equation}
is a quasi-isomorphism. The statement then follows from Corollary \ref{thm:equivalence of S^1 vs C_p fixed points}.\qed\par\indent
By Lemma \ref{thm:easy lemma}, there is an $S^1$-equivariant homotopy equivalence $\mathrm{UConf}_p(\mathbb{C}^*)\simeq S^1/C_p$.\\
\begin{lemma}\label{thm:generator of S^1 fixed points of S^1 mod C_p}
There is a generator of $H^0(C_{-*}(S^1/C_p;R)^{tS^1})\cong R$ whose image under
\begin{equation}
H^0(C_{-*}(S^1/C_p;R)^{tS^1})\xrightarrow{\mathrm{Res}_{C_p\subset S^1}}  H^0(C_{-*}(S^1/C_p;R)^{tC_p})=H_0(S^1/C_p;R)\oplus    H_1(S^1/C_p;R)\theta
\end{equation}
is 
\begin{equation}
[pt]-[S^1/C_p]\theta,    
\end{equation} 
where we identified (since the $C_p$-action on $S^1/C_p$ is trivial) 
\begin{equation}
H^*(C_{-*}(S^1/C_p;R)^{tC_p})=H_{-*}(S^1/C_p;R)\otimes H^*(R^{tC_p})=H_{-*}(S^1/C_p;R)((t,\theta)).     
\end{equation}    
Moreover, there is a generator $H^{-1}(C_{-*}(S^1/C_p;R)^{tS^1})\cong R$ whose image under $\mathrm{Res}_{C_p\subset S^1}$ is 
\begin{equation}
 [S^1/C_p]\in H_{-*}(S^1/C_p;R)((t,\theta)).
\end{equation}
\end{lemma} 
\emph{Proof}. The action of $C_{-*}(S^1;R)$ on $C_{-*}(S^1/C_p;R)$ can be modeled by the dg algebra action $R[\tau,\sigma]\acts R[\Lambda]$ (cf. Appendix A.3) given by
\begin{equation}
\tau\cdot1=1\;,\;\tau\cdot\Lambda=\Lambda\;,\;\sigma\cdot 1=\Lambda\;,\;\sigma\cdot \Lambda=0. 
\end{equation}
In particular, $(1+\tau+\cdots+\tau^{p-1})\sigma\cdot 1=0$ and
\begin{equation}
C_{-*}(S^1/C_p;R)^{tS^1}\simeq (R[\Lambda]((t)), 0),   
\end{equation}
which implies that $H^0(C_{-*}(S^1/C_p;R)^{tS^1})\cong H^{-1}(C_{-*}(S^1/C_p;R)^{tS^1})\cong R$. On the other hand, by the explicit formula for $\mathrm{Res}_{C_p\subset S^1}$ in \eqref{eq:Gysin comparison chain map}, one sees that the image of $1\in R=H^0(C_{-*}(S^1/C_p;R)^{tS^1})$ under
\begin{equation}
(R[\Lambda]((t)),0)\xrightarrow{\mathrm{Res_{C_p\subset S^1}}}(R[\Lambda]((t,\theta)),0)\simeq \mathrm{RHom}_{R[\tau]}(R,R[\Lambda])    
\end{equation}
is $1-(\sigma\cdot 1)\theta=1-\Lambda\theta$, which exactly corresponds to the element $[pt]-[S^1/C_p]\theta$; and the image of $1\in R=H^{-1}(C_{-*}(S^1/C_p;R)^{tS^1})$ is given by $\Lambda-(-1)\sigma\cdot\Lambda\theta=\Lambda$, which corresponds to the fundamental class $[S^1/C_p]$. \qed\par\indent
Moreover, it is easy to see that the two generators of $H^i(C_{-*}(S^1/C_p;R)^{tS^1})$ ($i=0,-1$) in question in Lemma \ref{thm:generator of S^1 fixed points of S^1 mod C_p} correspond exactly to $[\epsilon_0]$ and $[\epsilon_1]$ under the identification $S^1/C_p\simeq \mathrm{UConf}_p(\mathbb{C}^*)\simeq \mathrm{UCyl}^{\circ}(p,1)$. In particular,\\
\begin{cor}\label{thm:restriction of epsilon classes}
$\mathrm{Res}_{C_p\subset S^1}([\epsilon_0])=[e_0]-[e_1]\theta$ and $\mathrm{Res}_{C_p\subset S^1}([\epsilon_1])=[e_1]$. \qed
\end{cor}
Recall from Corollary \ref{thm:residual circle action on C_p fixed point via Gysin} that $[S^1/C_p]$ acts on $HH^{C_p,per}_*(\mathcal{A})$ as 
\begin{equation}
([\alpha],[\beta]\theta)\mapsto([\beta],0) \quad\in \mathrm{End}(HH^{C_p,per}_*(\mathcal{A}))\cong \mathrm{End}\big(HH_*^{per}(\mathcal{A}) \oplus HH_*^{per}(\mathcal{A})\theta\big).
\end{equation}
Combining this with Corollary \ref{thm:restriction of epsilon classes} and Figure \ref{fig:C_p_Cartan_formula}, one obtains the following. \\
\begin{cor}\label{thm:formula for action of epsilon classes}
\begin{itemize}
    \item The endomorphism (when paired with an undetermined element of $H^*(CC^*(\mathcal{A})^{\otimes p})^{h\Sigma_p})$) 
\begin{equation}
\Xi([\epsilon_0])    
\end{equation}
is given by the composition
\begin{equation}\label{eq:action of e_0-e_1theta}
HH_*^{per}(\mathcal{A})\xrightarrow{\mathrm{Res}_{C_p\subset S^1}}    HH^{C_p,per}_*(\mathcal{A})\xrightarrow{\Xi^p([e_0])}  HH^{C_p,per}_*(\mathcal{A})=HH_*^{per}(\mathcal{A}) \oplus HH_*^{per}(\mathcal{A})\theta\xrightarrow{pr_1}HH_*^{per}(\mathcal{A}).
\end{equation}
      \item The endomorphism 
\begin{equation}
\Xi^p([e_1])
\end{equation}
 preserves the decomposition $HH^{C_p,per}_*(\mathcal{A})=HH_*^{per}(\mathcal{A}) \oplus HH_*^{per}(\mathcal{A})\theta$ and moreover 
\begin{equation}
 \Xi^p([e_1])|_{HH^{per}_*(\mathcal{A})}=\Xi([\epsilon_1])   
\end{equation} 
 is given by the composition
 $$HH_*^{per}(\mathcal{A})\xrightarrow{\mathrm{Res}_{C_p\subset S^1}}    HH^{C_p,per}_*(\mathcal{A})\xrightarrow{\Xi^p([e_0])}  HH^{C_p,per}_*(\mathcal{A})=HH_*^{per}(\mathcal{A}) \oplus HH_*^{per}(\mathcal{A})\theta$$
\begin{equation}\label{eq:action of e_1}
\xrightarrow{pr_2}HH_*^{per}(\mathcal{A})\theta\xrightarrow{\textrm{divide by $\theta$}} HH^{per}_*(\mathcal{A}).
\end{equation}
\end{itemize}  \qed
\end{cor}
In Section 6, we give precise formula for the operation $\Xi^p([e_0])$ as a `$p$-fold equivariant cap product' on an explicit chain model computing $HH^{C_p,per}_*(\mathcal{A})$, which completes our understanding of Kontsevich-Soibelman operations.
\subsubsection{Proof of Theorem \ref{thm:main theorem} and Theorem \ref{thm:main theorem without homological bounds assumption}}
Theorem \ref{thm:main theorem} (resp. Theorem \ref{thm:main theorem without homological bounds assumption}) follows from Theorem \ref{thm:classification of untwisted KS operations} I) (resp. II)) combined with Corollary \ref{thm:formula for action of epsilon classes}, modulo the fact that $\Xi^p([e_0])$ agrees with $\prod\nolimits^{C_p}$, which will be proved in Proposition \ref{thm:e0 operation is p fold equivariant cap product}.\qed

\subsection{Twisted operations}
Let $p$ be an odd prime. A priori, it seems natural to consider a parallel set of operations, which we tentatively call \emph{twisted $C_p$-Kontsevich-Soibelman operations}, arising from the map
\begin{equation}\label{eq:twisted ksn1 equiv cohomology periodic}
\Xi^{\pm}:H^*((\mathbf{KS}(n,1)_R\otimes R(n))^{tS^1}_{h\Sigma_n})\rightarrow   \mathrm{Hom}_{H^*_{S^1}(pt;R)}(H^*_{\Sigma_n}(CC^*(\mathcal{A})^{\otimes n}\otimes R(n))\otimes HH_*^{per}(\mathcal{A}), HH_*^{per}(\mathcal{A})),
\end{equation}
obtained from \eqref{eq:ksn1 equiv cohomology periodic} by twisting both sides by the sign representation $R(n)$ of $\Sigma_n$.\\
\begin{mydef}\label{thm:twisted KSp operation}
An endomorphism $L$ is called a \emph{twisted $C_p$-Kontsevich-Soibelman operation of arity $n$} (abbrev. \emph{$\mathrm{KS}^{\pm,n}_p$ operation}) if there exists $[\alpha]\in H^*((\mathbf{KS}(n,1)_R\otimes R(n))^{tC_p}_{h\Sigma_n})$ and $[\phi]\in H^*_{\Sigma_n}(CC^*(\mathcal{A})^{\otimes n}\otimes R(n))$ such that
\begin{equation}
L=\Xi^{\pm,p}([\alpha])([\phi],-).    
\end{equation}
\end{mydef}
The goal of the following discussion is to show that in fact:
\begin{center}
\emph{the twisted KS operations are fully subsumed by the untwisted ones}.    
\end{center}
As most of the arguments here are complete analogues of the untwisted case, we will often omit full proofs. The only difference here is that the operations in question are parametrized by the homology of
\begin{equation}\label{eq:twisted equiv cohomology of cyl and conf}
C_{-*}(\mathrm{Cyl}^{\circ}(pk,1);R)\otimes R(pk))^{tC_p}_{h\Sigma_{pk}}\cong   (C_{-*}(\mathrm{Conf}_{pk}(\mathbb{C}^*);R)\otimes R(pk))^{tC_p}_{h\Sigma_{pk}}.
\end{equation}
The computation of \eqref{eq:twisted equiv cohomology of cyl and conf} also follows from the classical work of Cohen \cite{Coh}, and is facilitated by the notion of labeled configuration space \cite{Bod}.\\
\begin{mydef}\label{thm:labeled configuration space}
Let $M$ be a topological space and $(X,*)$ be a based space. The \emph{space of (unordered) configurations of $M$ labeled by $X$} is defined as
\begin{equation}
\mathrm{UConf}(M;X):=\bigsqcup_{n\geq 0} \mathrm{Conf}_n(M)\times_{\Sigma_n} X^n/\sim,
\end{equation}
where the equivalence relation is
\begin{equation}
(m_1,\cdots,m_n;x_1,\cdots,x_n)\sim (m_1,\cdots,\hat{m}_i,\cdots,m_n;x_1,\cdots,\hat{x}_i,\cdots,x_n)\;\;\mathrm{if}\;\;x_i=*.    
\end{equation}
\end{mydef}
There is an obvious increasing filtration $\mathrm{UConf}_{\leq m}(M;X)\subset \mathrm{UConf}(M;X)$ given by
\begin{equation}
\mathrm{UConf}_{\leq m}(M;X):=\bigsqcup_{n\leq m} \mathrm{Conf}_n(M)\times_{\Sigma_n} X^n/\sim
\end{equation}
and we define $\mathrm{UConf}_n(M;X):=\mathrm{UConf}_{\leq n}(M;X)/\mathrm{UConf}_{\leq n-1}(M;X)$. \par\indent
Let $\tilde{C}_*(-)$ (resp. $\tilde{H}_*(-)$) denote the reduced singular chains (resp. homology) of a topological space. The relation of \eqref{eq:twisted equiv cohomology of cyl and conf} to labeled configuration is made by the following observation. Namely, as a $\Sigma_n$ representation, $R(n)\cong \tilde{H}_{-*}(S^1;R)^{\otimes n}$ up to a degree shift of $-n$. Therefore, there is a chain of quasi-isomorphisms
$$
(C_{-*}(\mathrm{Conf}_{n}(\mathbb{C}^*);R)\otimes R(n))_{h\Sigma_n}\simeq C_{-*}(\mathrm{Conf}_{n}(\mathbb{C}^*);R)\otimes_{\Sigma_n} \tilde{H}_{-*}(S^1;R)^{\otimes n}[-n]\simeq C_{-*}(\mathrm{Conf}_{n}(\mathbb{C}^*);R)\otimes_{\Sigma_n} \tilde{C}_{-*}(S^1;R)^{\otimes n}[-n]
$$
\begin{equation}\label{eq:twisted conf as labeled conf}
\simeq  \tilde{C}_{-*-n}(\mathrm{UConf}_n(\mathbb{C}^*;S^1);R).   
\end{equation}
\textbf{Observation}. Consider the based space $(S^0\vee S^1,*)$. Color the copy of $S^1$ black and color the point of $S^0$ which is not $*$ white. Then, up to homotopy, $\mathrm{UConf}(\mathbb{C}^*;S^1)$ may be identified as the subspace of $\mathrm{UConf}(\mathbb{C};S^0\vee S^1)$ consisting of exactly one white point (namely, up to homotopy, place the white point at the origin and think of it as a puncture).\par\indent
The advantage of re-interpreting $\mathrm{UConf}(\mathbb{C}^*;S^1)$ as labeled configurations in $\mathbb{C}$ is that the latter is naturally an $E_2$ space. Recall as part of the $E_2$ structure one has a (commutative up to homotopy) product $\cdot$ and a Lie bracket $[-,-]$. The equivariant homology of $E_2$ spaces (or more generally $E_n$ spaces) was extensively studied in the classical work of Cohen \cite{Coh}, and we recall the following result summarized in \cite[Theorem 4.8]{Ros}.    \\
\begin{thm}(\cite{Coh})\label{thm:Cohen theorem on labeled con}
Let $Q:H_q(\mathrm{UConf}(\mathbb{C};X);R)\rightarrow H_{pq+p-1}(\mathrm{UConf}(\mathbb{C};X);R)$ be the first Dyer-Lashof operation (when $p$ is odd, $Q$ only applies to odd degree classes).\par\indent
When $p=2$, $H_*(\mathrm{UConf}(\mathbb{C};X);R)$ is a free graded commutative algebra (under $\cdot$) by elements of the form $Q^i(x)$, where $x$ is a basic bracket.\par\indent
When $p>2$, $H_*(\mathrm{UConf}(\mathbb{C};X);R)$ is a free graded commutative algebra (under $\cdot$) by elements of the form $Q^i(x), \beta Q^i(x)$, where $\beta$ is the Bockstein and $x$ is a basic bracket of odd degree.\par\indent
Recall that a basic bracket is an element obtained from $\tilde{H}_*(X;R)\cong\tilde{H}_*(\mathrm{UConf}_1(\mathbb{C};X);R)$ by iterated Lie brackets. \qed
\end{thm}
Denote the white point in $S^0$ by $b$, and denote $a=[S^1]$ the fundamental class of $S^1$, both viewed as elements of $\tilde{H}_*(\mathrm{UConf}_1(\mathbb{C};S^0\vee S^1);R)$. Applying Theorem \autoref{thm:Cohen theorem on labeled con} to $X=S^0\vee S^1$, one concludes that\\
\begin{lemma}\label{thm:homology of conf labeled by S0 wedge S1}
$\tilde{H}_*(\mathrm{UConf}_n(\mathbb{C}^*;S^1);R)\subset H_*(\mathrm{UConf}(\mathbb{C};S^0\vee S^1);R)$ is equal to the direct sum
\begin{equation}
b\cdot \tilde{H}_*(\mathrm{UConf}_n(\mathbb{C};S^1);R)+ [a,b]\cdot\tilde{H}_*(\mathrm{UConf}_{n-1}(\mathbb{C};S^1);R)+\cdots+[\overbrace{a,[a,\cdots,[a}^{n\;\mathrm{times}},b]]]
\end{equation}    
\end{lemma}
\emph{Proof.} Note that applying a Dyer-Lashof operation scales the number of marked points (of each color) by $p$. In particular, since elements in $\mathrm{UConf}_n(\mathbb{C}^*;S^1)$ correspond to labeled configurations with exactly one white point, only generators of the form 1) $Q^i(x)$ (or $\beta Q^i(x)$), where $x$ is a basic Lie bracket built out of $a=[S^1]$ or 2) a basic Lie bracket with only one $b$ as input are relevant. \qed  \par\indent
Therefore, pictorially the generators of $\tilde{H}_*(\mathrm{UConf}_n(\mathbb{C}^*;S^1);R)$ are exactly shown as in Figure \ref{fig:figure6} except that all the black points are now labeled by the fundamental class $[S^1]$.\par\indent
Similar to the proof of Lemma \autoref{thm:easy lemma}, one can use the $p$-fold map to show that\\
\begin{lemma}\label{thm:easy lemma twisted case}
There is a homeomorphism 
\begin{equation}
\mathrm{UConf}_{pk}(\mathbb{C}^*;S^1)^{C_p}\cong    \mathrm{UConf}_{k}(\mathbb{C}^*;S^1).
\end{equation}
In particular, combining with \eqref{eq:twisted conf as labeled conf} and the localization theorem we obtain a chain of quasi-isomorphisms 
$$C_{-*}(\mathrm{Conf}_{pk}(\mathbb{C}^*;R)\otimes R(pk))^{tC_p}_{h\Sigma_{pk}}\simeq \tilde{C}_{pk-*}(\mathrm{UConf}_{pk}(\mathbb{C}^*;S^1);R))^{tC_p}\simeq \tilde{C}_{pk-*}(\mathrm{UConf}_{pk}(\mathbb{C}^*;S^1)^{C_p};R))$$
\begin{equation}
\simeq \tilde{C}_{-pk-*}(\mathrm{UConf}_{k}(\mathbb{C}^*;S^1);R)).
\end{equation}\qed
\end{lemma}
Let $\mathrm{UCyl}(n,0;S^1)$ (resp. $\mathrm{UCyl}(n,1;S^1)$) be the space of (unordered) embeddings of disks in a disk (resp. cylinder), with each embedded disk labeled by $S^1$ analogous to Definition \autoref{thm:labeled configuration space}; it is a homotopy equivalent model for $\mathrm{UConf}_n(\mathbb{C};S^1)$ (resp. $\mathrm{UConf}_n(\mathbb{C}^*;S^1)$). \par\indent 
Based on previous notions \eqref{eq:lambda delta equiv}, \eqref{eq:cyl operadic 2}, there are operadic compositions (which just further carry over the labelings)
\begin{equation}\label{eq:twisted cyl operadic 1}
\lambda_{\Delta}^{\pm}:\mathrm{UCyl}(l,0;S^1)\times \mathrm{UCyl}^{\circ}(p,1)\rightarrow \mathrm{UCyl}^{\circ}(pl,1;S^1)  
\end{equation}
\begin{equation}\label{eq:twisted cyl operadic 2}
\eta^{\pm}:\mathrm{UCyl}^{\circ}(n,1;S^1)\times \mathrm{UCyl}^{\circ}(n',1;S^1)\rightarrow \mathrm{UCyl}^{\circ}(n+n',1;S^1).
\end{equation}
Similar to the proof of Lemma \autoref{thm:generation of equiv cohomology of Cyl}, one can show that\\
\begin{lemma}\label{thm:generation of twisted equiv cohomology of Cyl}
$\bigoplus_{k\geq 1}H^*(\tilde{C}_{-*}(\mathrm{UCyl}^{\circ}(pk,1;S^1);R)^{tC_p})$ is generated by 
\begin{equation}
\tilde{H}_{-*}(\mathrm{UCyl}(l,0;S^1);R),    
\end{equation}
\begin{equation}
[e_0]\in H_0(\mathrm{UCyl}^{\circ}(p,1)^{C_p};R),   
\end{equation} and 
\begin{equation}
[e_1]\in H_1(\mathrm{UCyl}^{\circ}(p,1)^{C_p};R)   
\end{equation} 
under (induced map on $C_p$-Tate cohomology of) the two operations $\lambda_{\Delta}^{\pm}$ and $\eta^{\pm}$ (cf. Figure \ref{fig:figure7} for the definition of $e_0,e_1$).
\qed\\
\end{lemma}
\begin{thm}\label{thm:classification of twisted KSp operations}
Assume that the natural map
\begin{equation}
(CC^*(\mathcal{A})^{\otimes k_1}\otimes R(k_1))^{h\Sigma_{k_1}}\otimes\cdots\otimes (CC^*(\mathcal{A})^{\otimes k_l}\otimes R(k_l))^{h\Sigma_{k_l}}\rightarrow (CC^*(\mathcal{A})^{\otimes k_1+\cdots+k_l}\otimes R(k_1)\otimes\cdots\otimes R(k_n))^{h(\Sigma_{k_1}\times \cdots\times\Sigma_{k_l})}     
\end{equation}
is a quasi-isomorphism for all $k_1,\cdots,k_l,l\geq 1$ and that the left hand side above satisfies the K\"{u}nneth isomorphism. Then all twisted $C_p$-Kontsevich-Soibelman operations are generated under composition and linear combination by elements of the form
\begin{equation}\label{eq:e0 operation}
\Xi^p([e_0])([\phi],-)    
\end{equation}
and 
\begin{equation}\label{eq:e1 operation}
\Xi^p([e_1])([\phi],-),    
\end{equation}
where $[\phi]\in H^*_{\Sigma_p}(CC^*(\mathcal{A})^{\otimes p})$.  \qed
\end{thm}
Note that these are the same generators as in the untwisted case! We omit the proof of Theorem \autoref{thm:classification of twisted KSp operations}, as it follows from Lemma \autoref{thm:generation of twisted equiv cohomology of Cyl} using exactly the same argument as Theorem \autoref{thm:classification of untwisted KSp operations}. However, we give a simple example to illustrate the idea.\\
\begin{ex}\label{thm:example of twisted e0 operation}
 Consider the generator 
\begin{equation}
[e_0^{\pm}]:=\lambda^{\pm}_{\Delta}([S^1],[e_0])\in \tilde{H}_1(\mathrm{UCyl}^{\circ}(p,1;S^1)^{C_p};R)\subset H_1(\tilde{C}_*(\mathrm{UCyl}^{\circ}(p,1;S^1);R)^{tC_p}).
\end{equation}
where $[S^1]$ denotes the image of the fundamental class under $\tilde{H}_1(S^1;R)\cong \tilde{H}_1(\mathrm{Cyl}(1,0;S^1);R)$; see Figure \ref{fig:figure9}.
\begin{figure}[H]
 \centering
 \includegraphics[width=0.8\textwidth]{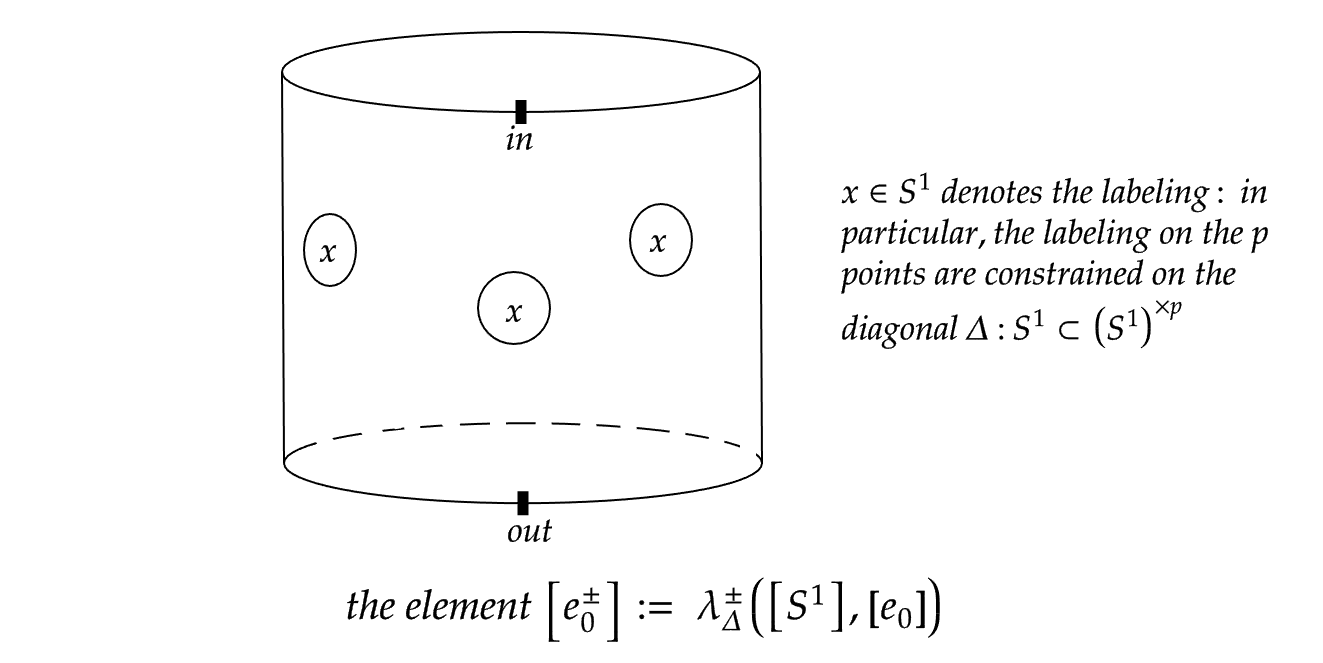}
 \caption{}
 \label{fig:figure9}
\end{figure}
Let $\delta_{[S^1]}$ be the image of $[S^1]$ under the map
\begin{equation}
\tilde{H}_1(S^1;R)\cong \tilde{H}_1(\mathrm{Cyl}(1,0;S^1);R)\xrightarrow{\Delta} H^{-1}(\tilde{C}_{-*}(\mathrm{Cyl}(1,0;S^1);R)^{\times p})^{h\Sigma_p})\cong H^{p-1}\big((C_{-*}(\mathrm{Cyl}(1,0);R)^{\otimes p}\otimes R(p))^{h\Sigma_p}\big).
\end{equation}
By the $E_2$ structure on Hochschild cochains, $\delta_{[S^1]}$ induces a map
\begin{equation}
\mathrm{Act}_{\delta_{[S^1]}}: H^*_{\Sigma_p}(CC^*(\mathcal{A})^{\otimes p}\otimes R(p))\rightarrow H^{*+p-1}_{\Sigma_p}(CC^*(\mathcal{A})^{\otimes p}).  
\end{equation}
In fact, since $\mathrm{Cyl(1,0)}$ is contractible, $\mathrm{Act}_{\delta_{[S^1]}}$ is just cup product with the diagonal $\Delta_{[S^1]}\in H^{-1}_{\Sigma_p}(\tilde{C}_{-*}((S^1)^{\times p};R))\cong H^{-1}_{\Sigma_p}(\tilde{C}_{-*}(S^1;R)^{\otimes p})\cong H^{p-1}_{\Sigma_p}(R(p))$.
From the operadic structure, it is easy to see that for $[\phi^{\pm}]\in H^*_{\Sigma_p}(CC^*(\mathcal{A})^{\otimes p}\otimes R(p))$, 
\begin{equation}
\Xi^{\pm,p}([e_0^{\pm}])([\phi^{\pm}])=\Xi^p([e_0])(\mathrm{Act}_{\delta_{[S^1]}}([\phi^{\pm}]))=\Xi^p([e_0])([\phi^{\pm}]\cup \Delta_{[S^1]}).
\end{equation}
\end{ex}\par\indent
\begin{rmk}
Another heuristic explanation for why there are no \emph{genuinely} twisted Kontsevich-Soibelman operations is the following. By the discussion on orbit types at the start of Section 5.1, one sees that for any representative $(z_1,\cdots,z_{pk})$ of a $C_p$-fixed point of $C_p\acts \mathrm{UConf}_{pk}(\mathbb{C}^*)$, $(e^{2\pi i/p} z_1, \cdots,e^{2\pi i/p}z_{pk})=\tau\cdot(z_1,\cdots,z_{pk})$, where $\tau\in \Sigma_{pk}$ is a permutation given by $k$ disjoint cycles of length $p$. If $p$ is odd, then $\tau$ must be an even permutation. In particular, for the sake of computing (Tate) fixed points, one may replace $\Sigma_{pk}$ by the alternating subgroup, which has trivial sign representation. 
\end{rmk}


\section{From configurations to cacti}
The goal of this section is to prove the comparison Theorem \autoref{thm:comparison of KS with Cyl} between weighted (cyclic) cacti and configurations on disks/cylinders. In Section 6.1, we review the key ideas in Salvatore's proof \cite{Sal1} that the operad of weighted cacti is quasi-equivalent to the little disk operad (or equivalently, the Fulton-MacPherson operad). In Section 6.2, we adapt the methods in loc.cit. to give a proof of the $S^1\times S^1$-equivariant equivalence $\mathrm{Cyl}(k,1)\simeq \widetilde{\mathrm{Cact}^k}$ that is furthermore compatible with (two-colored) operadic structures. In Section 6.3, we use the results of Section 6.1 and 6.2 to give explicit formulae computing the generators \eqref{eq:e0 operation} and \eqref{eq:e1 operation} of Kontsevich-Soibelman operations on the periodic cyclic homology.  
\subsection{Review of Salvatore's results}
To a configuration $(z_1,\cdots,z_n)\in \mathrm{Conf}_n(\mathbb{C})$ and a set of weights $(a_1,\cdots,a_n)\in \mathcal{P}(n)$, \cite{Sal1} first associates a combinatorial object called \emph{labeled tree}. Roughly speaking, take a vector field on the complex plane that is the superposition of radial vector fields centered at $z_i$ (with charges $a_i$) whose magnitudes are inversely proportional to their respective distances from $z_i$; the labeled tree is obtained as the critical graph (i.e. the union of flow lines connecting zeros and poles) of the vector field.\par
\begin{figure}[H]
 \centering
 \includegraphics[width=1.0\textwidth]{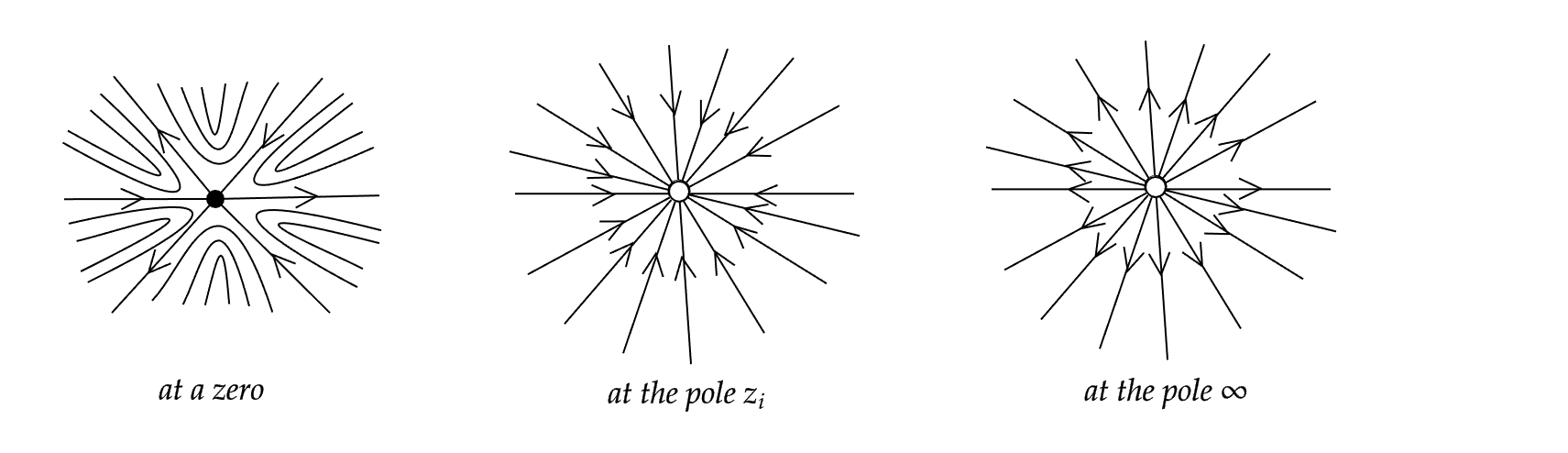}
 \caption{}
 \label{fig:local_behaviour_of_E}
\end{figure}

\indent
More precisely, let
\begin{equation}\label{eq:vector field E}
E(z):=\sum_{i=1}^n a_i\frac{z-z_i}{|z-z_i|^2}.
\end{equation}
$E$ has a potential of the form $E(z)=-\nabla \log|h(z)|$, $h(z)=\prod_{i=1}^n(z-z_i)^{a_i}$. In particular, the flow lines of $E$ correspond to level curves of $\mathrm{Im}(\log h(z))=\arg h(z)$. From this description, one can summarize the local behaviors of the flow lines of $E$ as follows, cf. \cite[Lemma 3.7]{Sal1}; see Figure \ref{fig:local_behaviour_of_E}.
\begin{itemize}
    \item Near a zero $w$ of $E$, if $\log h(z)-\log h(w)$ has a zero of order $m$ at $w$, then locally there are $m$ incoming and $m$ outgoing flow lines, in alternating fashion, at $w$.
    \item Near the pole $z_i, 1\leq i\leq n$, all flow lines are incoming.
    \item Near the pole $\infty$, all flow lines are outgoing. 
\end{itemize}
These local behaviors imply that the critical graph of $E$ satisfies certain properties which are summarized in the following definition.\\
\begin{mydef}\label{thm:admissible and labeled trees}
An \emph{admissible tree $T$} with $n$ leaves is an oriented tree such that
\begin{itemize}
    \item The set of vertices is equipped with a partition $V(T)=B\bigsqcup W$ into black and white vertices, such that $|W|=n$. The white vertices are also called \emph{leaves}.
    \item The set of edges $E_v$ incident to a fixed vertex $v$ is equipped with a cyclic ordering (i.e. $T$ is a \emph{ribbon graph}).
\end{itemize}
Furthermore, one requires the following conditions.
\begin{enumerate}[label=\arabic*)]
    \item A white vertex is not a source.
    \item If $(v,w)\in E(T)$, then $(w,v)\notin E(T)$ (i.e. an edge cannot be inverted).
    \item Any black vertex $v\in B$ is the source of at least two distinct edges. Moreover it cannot be the target of two edges next to each other in the cyclic ordering of $E_v$.
\end{enumerate}
Denote by $T_n$ the set of isomorphism classes of admissible trees with
$n$ leaves.  
A \emph{labeled tree} is an admissible tree $T$ with the following additional data:
\begin{enumerate}[label=\roman*)]
    \item A function $f: B\rightarrow (0,1]$ attaining the maximum $1$ and $f(b)\geq f(b')$ along an edge $(b,b')$ connecting two black vertices.
    \item To each white vertex $w\in W$, a function $g: E_w\rightarrow [0,1]$ such that $\sum_{e\in E_w} g(e)=1$. We view $2\pi g(e)$ as the angle between the edge $e$ and the next edge in the cyclic ordering.
\end{enumerate}
\end{mydef}
\textbf{From configurations to labeled trees}. Given a configuration $(z_1,\cdots,z_n)$, one associates to it a labeled tree with $n$ leaves as follows, cf. \cite[Proposition 3.9, Definition 3.10]{Sal1}. 
\begin{itemize}
    \item The underlying admissible tree $T$ has vertices the zeros (colored black) and poles (corresponding to $z_i$ and colored white) of $E$ and edges the flow lines of $E$, with the natural cyclic ordering of $E_v$ induced from its embedding into $\mathbb{C}$. 
    \item $f(b):=|h(b)|/M$, where $M:=\max_{b\in B}|h(b)|$.
    \item For each $z_i$, $E_{z_i}$ corresponds to the set of flow lines connecting a zero of $E$ to $z_i$. Then put $g(e):=\theta(e,e')/2\pi$, where $\theta(e,e')$ is the angle at $z_i$ between $e$ and the next flow line $e'$ in the cyclic ordering. Equivalently, $g(e)=(\arg h(e)-\arg h(e'))/2\pi a_i$. 
\end{itemize}
It turns out this assignment, denoted 
\begin{equation}\label{eq:Phi from configuration to labeled tree}
\Phi_{a_1,\cdots,a_n}:(z_1,\cdots,z_n)\mapsto (T,f,g),    
\end{equation} gives a homeomorphism from a suitable quotient of $\mathrm{Conf}_n(\mathbb{C})$ to a `space' of labeled trees, which we now make precise.\par\indent
As preliminary notations, given an admissible tree $T$, let $\sigma_B$ be the space of functions $B\rightarrow (0,1]$ satisfying Condition i) of Definition \autoref{thm:admissible and labeled trees}. Given a white vertex $w\in W$, let $\Delta_w$ denote the space of functions $g: E_w\rightarrow[0,1]$ satisfying Condition ii) of Definition \autoref{thm:admissible and labeled trees} (it is a simplex of dimension $|E_w|-1$).
\begin{figure}[H]
 \centering
 \includegraphics[width=0.8\textwidth]{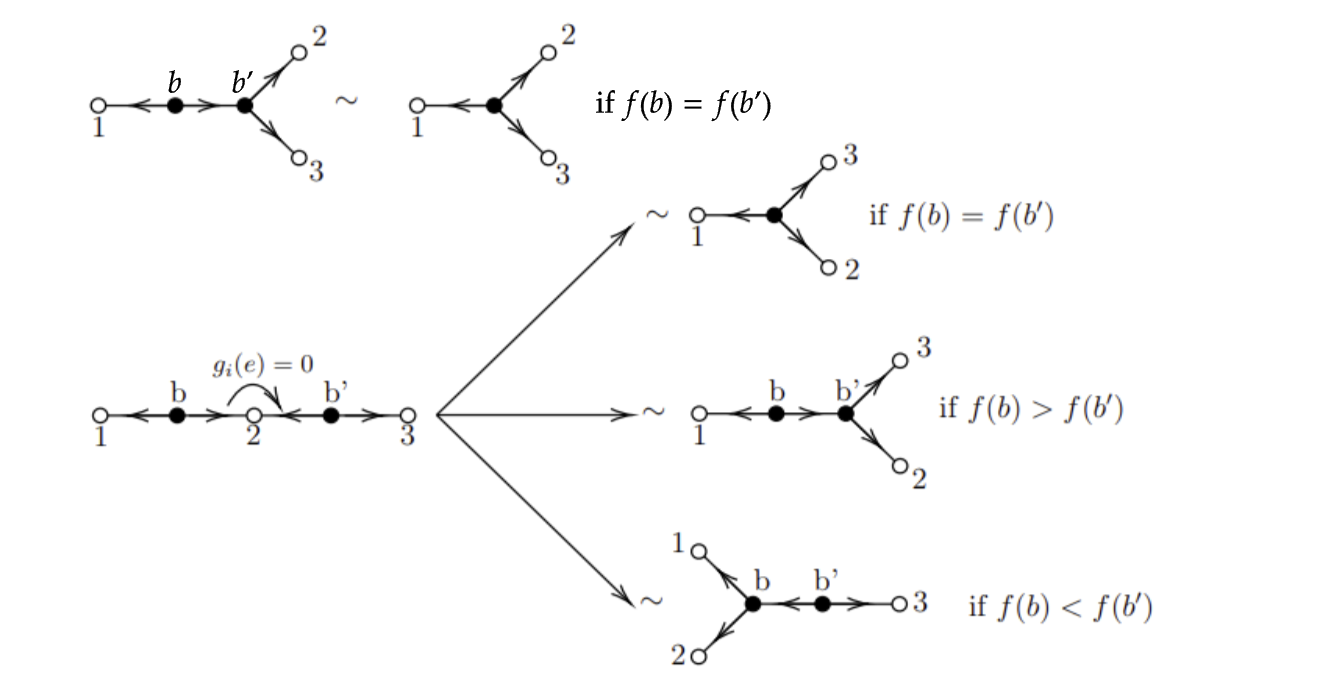}
 \caption{Equivalences of labeled trees}
 \label{fig:equivalence_in_spaces_of_labeled_trees}
\end{figure}
\begin{mydef}\label{thm:space of labeled trees}
The \emph{space of labeled trees with $n$ leaves} is 
\begin{equation}\label{eq:space of labeled trees}
Tr_n:=\coprod_{T\in T_n} \sigma_B\times \prod_{w\in W}\Delta_w\;\;/\sim,
\end{equation}
where $\sim$ is generated by two types of equivalences in Figure \ref{fig:equivalence_in_spaces_of_labeled_trees} (reproduced from \cite[Figure 5]{RS}). See \cite[Definition 2.8]{Sal1} for the precise definition.  \\
\end{mydef}
\begin{thm}(\cite[Theorem 3.13]{Sal1})\label{thm:from configuration to labeled trees, without S^1}
For each set of weights $(a_1,\cdots,a_n)\in \mathcal{P}(n)$, the map
\begin{equation}
\Phi_{a_1,\cdots,a_n}: \mathrm{Conf}_n(\mathbb{C})/\mathbb{C}\rtimes \mathbb{C}^*\cong Tr_n
\end{equation}
is a homeomorphism. \qed
\end{thm}
We briefly sketch the construction of the inverse $\Psi_{a_1,\cdots,a_n}:Tr_n\rightarrow \mathrm{Conf}_n(\mathbb{C})/\mathbb{C}\rtimes\mathbb{C}^*$, and refer the reader to \cite[Section 3]{Sal1} for more details. Starting with a labeled tree $(T,f,g)$, the construction goes in two steps; see Figure \ref{fig:from_labeled_tree_to_configuration} for an illustration. 
\begin{enumerate}[label=\arabic*)]
    \item Firstly, adjoin to $T$ an extra vertex called $\infty$. For each black vertex $b\in B$, adjoin edges from $\infty$ to $b$ in between any two consecutive outgoing edges from $b$. There is a natural cyclic ordering on adjoined edges from $\infty$. The resulting ribbon graph will be called $\overline{T}$. 
    \item Secondly, for each edge $e\in E_{v_i}$ going into a white vertex $v_i$, glue a strip $[-\infty,+\infty]\times[0,a_ig_i(e)]$ to $\overline{T}$, so that $\{-\infty\}\times [0,a_ig_i(e)]$ is glued to $v_i$ and $\{+\infty\}\times[0,a_ig_i(e)]$ is glued to $\infty$, with appropriate identifications on the edges. The resulting Riemann surface has genus $0$, and after identifying it with $\mathbb{C}P^1$ so that the vertex $\infty$ lies at $\infty$, one reads off the coordinates $(z_1,\cdots,z_n)$ of the white vertices.  
\end{enumerate}
\begin{figure}[H]
 \centering
 \includegraphics[width=0.8\textwidth]{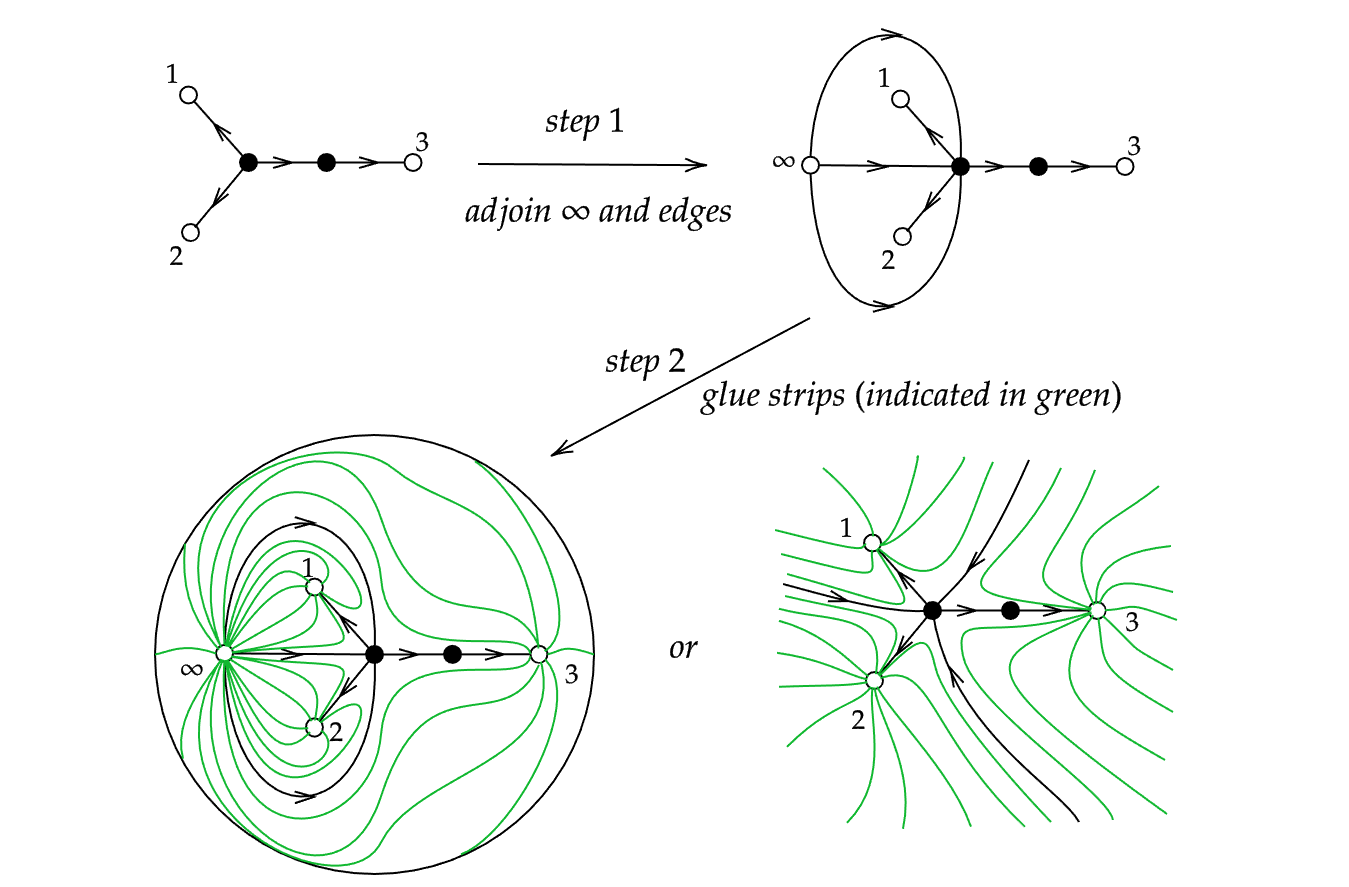}
 \caption{}
 \label{fig:from_labeled_tree_to_configuration}
\end{figure}
\textbf{From labeled trees to cacti.} Let $Tr_n^1\subset Tr_n$ denote the subspace where the function $f:B\rightarrow(0,1]$ is constant with value $1$. Note that $Tr_n$ deformation retracts onto $Tr_n^1$ by continuously increasing the value of $f$ on each $b\in B$ to $1$. By the first equivalence relation in Figure \ref{fig:equivalence_in_spaces_of_labeled_trees}, a labeled tree in $Tr^1_n$ (is equivalent to one that) has no edges between black vertices. In particular, there is a $\Sigma_n$-equivariant homeomorphism 
\begin{equation}\label{eq:homeo from Tr^1_n to cact^n}
Tr^1_n\cong \mathrm{cact}^n=\mathrm{Cact}^n/S^1  
\end{equation}
given as follows: each white vertex corresponds to a lobe, each black vertex corresponds to an intersection of (two or more, depending on its valency) lobes, and the values $g_i(e)$'s specify the lengths of arcs between intersection points of lobes; see Figure \ref{fig:from_labeled_tree_to_cactus} top.
\begin{figure}[H]
 \centering
 \includegraphics[width=1.0\textwidth]{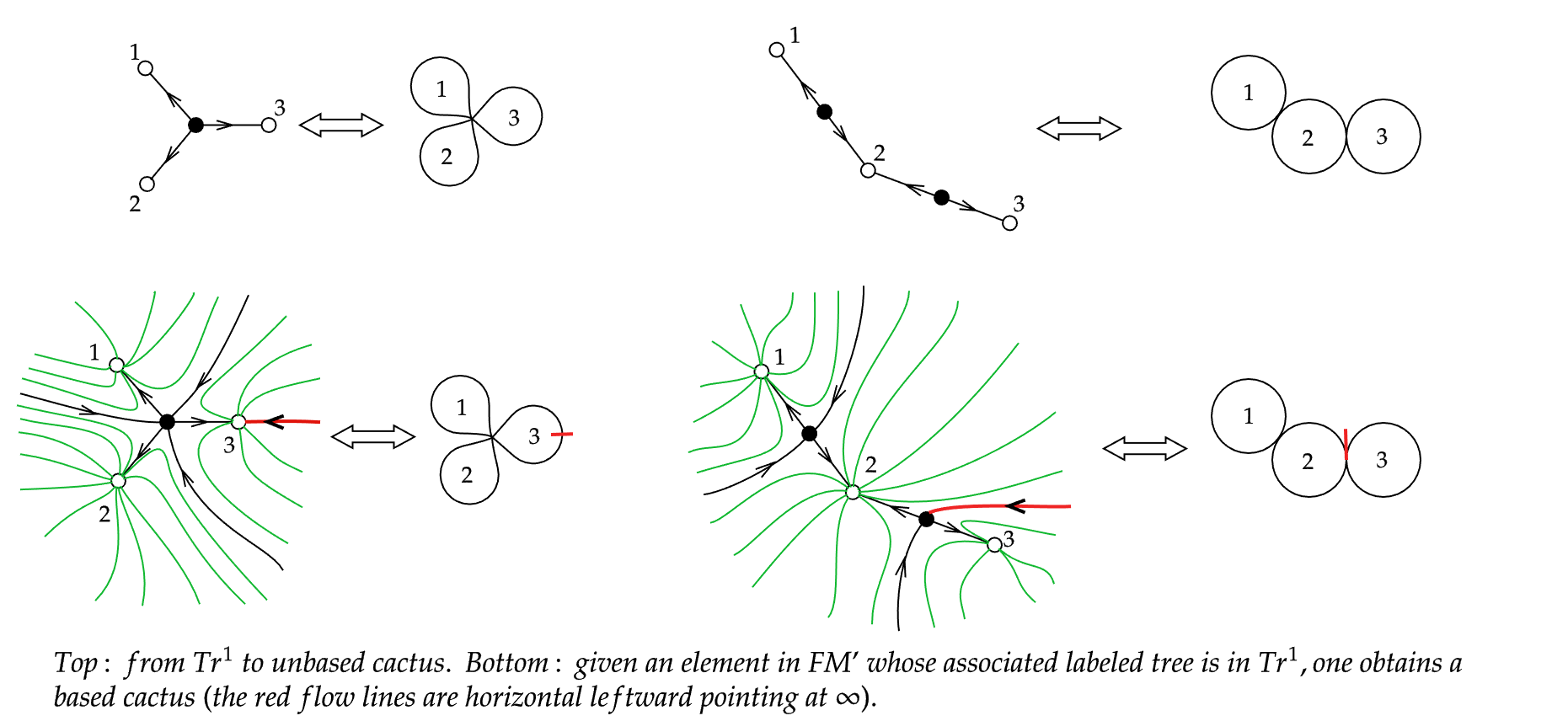}
 \caption{}
 \label{fig:from_labeled_tree_to_cactus}
\end{figure}
We now put the $S^1$ parameter in, i.e. compare the open stratum of the Fulton-MacPherson space 
\begin{equation}\label{eq:FM' open Fulton-MacPherson}
FM'(n):=\mathrm{Conf}_n(\mathbb{C})/\mathbb{C}\rtimes \mathbb{R}_+    
\end{equation} with the space of (based) cacti $\mathrm{Cact}^n$. We furthermore remove the dependence on weight parameters $(a_1,\cdots,a_n)$ by comparing instead $FM'(n)\times \mathcal{P}(n)$ with $\widetilde{\mathrm{Cact}^n}$. \par\indent
Namely, given $(z_1,\cdots,z_n)\in \mathrm{Conf}_n(\mathbb{C})$ such that $\Phi_{a_1,\cdots,a_n}[z_1,\cdots,z_n]=(T,f,g)$ lies in $Tr^1_n\subset Tr_n$, let $x'\in \widetilde{\mathrm{Cact}^n}/S^1$ be unbased cactus which is the image of $(T,f,g)$ under \eqref{eq:homeo from Tr^1_n to cact^n} whose lobes are weighted by $(a_1,\cdots,a_n)$. We assign a basepoint to $x'$ as follows. Consider the flow line $\gamma:(-\infty,\infty)\rightarrow \mathbb{C}$ of $E$ such that $\lim_{t\rightarrow-\infty}\gamma(t)=\infty$ and $\lim_{t\rightarrow-\infty}\gamma'(t)/|\gamma'(t)|=-1$. After potentially breaking at critical points, $\gamma$ eventually reaches some $z_i$; the basepoint on (the $i$-th lobe of) $x'$ records the infinitesimal tangent direction $\lim_{t\rightarrow+\infty}\gamma'(t)/|\gamma'(t)|$. The resulting (based weighted) cactus only depends on the image of $(z_1,\cdots,z_n)$ in $FM'(n)$. See Figure \ref{fig:from_labeled_tree_to_cactus} bottom. We summarize the result as follows.\\
\begin{prop}(compare \cite[Proposition 5.3]{Sal1})\label{thm:equivariant homotopy equivalence of FM' with Cact, weighted}
There is an $S^1\times \Sigma_n$-equivariant homotopy equivalence $FM'(n)\times \mathcal{P}(n)\simeq \widetilde{Cact^n}$.    \qed
\end{prop}
\textbf{Nested trees.} One can enhance the homotopy equivalence of Proposition \autoref{thm:equivariant homotopy equivalence of FM' with Cact, weighted} to a homeomorphism, at the cost of replacing $\widetilde{\mathrm{Cact}^n}$ by its cobar construction. Explicitly this is constructed in terms of nested trees, which we now recall.\\
\begin{mydef}\label{thm:nested trees}
A \emph{nested tree on $n$ leaves} ($n\geq 2$) is a collection $\mathcal{S}$ of subsets of $\{1,2,\cdots,n\}$ with cardinality at least $2$ such that
\begin{itemize}
    \item $\{1,\cdots,n\}\in \mathcal{S}$. This is the \emph{root} of the nested tree.
    \item For $S_1,S_2\in \mathcal{S}$, either $S_1\cap S_2=\emptyset$ or $S_1\subset S_2$ or $S_2\subset S_1$. 
\end{itemize}
The elements $S\in \mathcal{S}$ are the vertices of $\mathcal{S}$, and inclusions $S_1\subset S_2$ are the internal edges. Let $|S|_v$ denote the \emph{valency} of $S$, viewed as a vertex (which is not in general its cardinality $|S|$ as a set). Denote $\mathcal{S'}:=\mathcal{S}-\{1,2,\cdots,k\}$.\par\indent
Let $N_n$ denote the set of nested trees with $n$ leaves. See Figure \ref{fig:nested_trees} for examples of nested trees and their presentations (we caveat however that a nested tree should not be thought of as a planar tree, e.g. there are no cyclic ordering of edges incident to a vertex). 
\begin{figure}[H]
 \centering
 \includegraphics[width=0.8\textwidth]{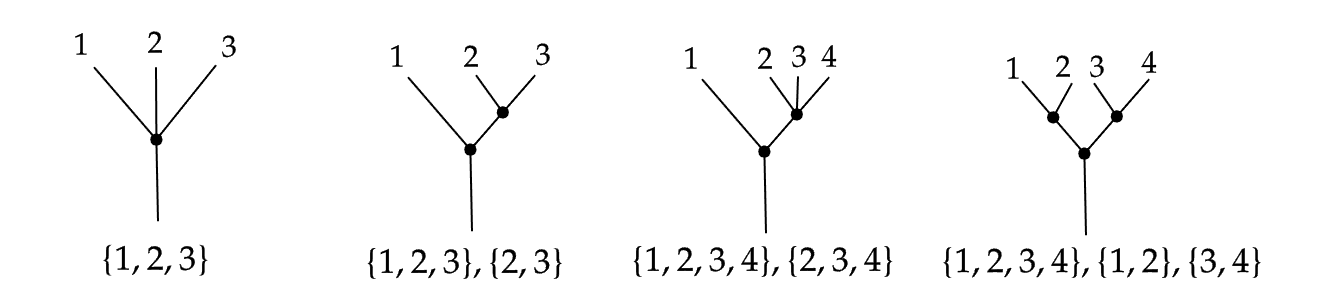}
 \caption{}
 \label{fig:nested_trees}
\end{figure}
\end{mydef}
Note that if $(O,\theta)$ is a topological operad, then a nested tree $\mathcal{S}\in N_k$ specifies a composition 
\begin{equation}\label{eq:nested tree specifies operadic composition}
\theta_{\mathcal{S}}: \prod_{S\in \mathcal{S}} O(|S|_v)\rightarrow O(k).
\end{equation}
We now combine this with the topological operad of weighted cacti, cf. \eqref{eq:operadic composition of weighted cacti}.\\
\begin{mydef}\label{thm:nested tree of weighted cacti}
The \emph{space of nested trees on $k$ leaves with vertices labeled by weighted cacti and internal edges labeled by $(0,1)$} is
\begin{equation}\label{eq:nested tree of weighted cacti}
N_k(\widetilde{\mathrm{Cact}}):=\coprod_{\mathcal{S}\in N_k}\prod_{S\in \mathcal{S}}\widetilde{\mathrm{Cact}^{|S|_v}}\times(0,1]^{\mathcal{S}'}\;/\sim ,   
\end{equation}
where $\sim$ is the equivalence relation generated by
\begin{equation}
((x_S)_{S\in\mathcal{S}},(\lambda_S)_{S\in\mathcal{S}'})\sim  ((x'_S)_{S\in\mathcal{S}-S_0},(\lambda_S)_{S\in\mathcal{S}'-S_0})\;\;,\;\;\;\lambda_{S_0}=1,
\end{equation}
where
\begin{equation}
\widetilde{\theta}_{\mathcal{S}}((x_S)_{S\in\mathcal{S}})=\widetilde{\theta}_{\mathcal{S}-S_0}((x'_S)_{S\in\mathcal{S}-S_0})\in\widetilde{\mathrm{Cact}^k}.
\end{equation}
\end{mydef}\par
\begin{figure}[H]
 \centering
 \includegraphics[width=0.9\textwidth]{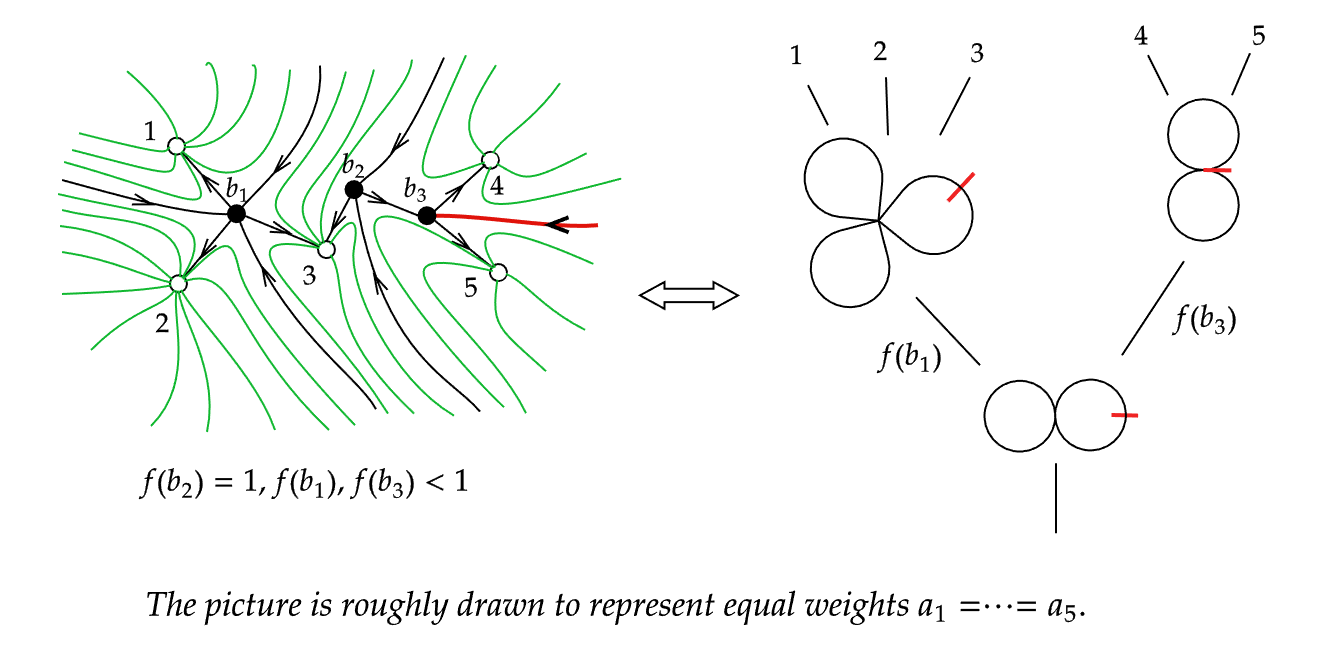}
 \caption{}
 \label{fig:from_FM__to_nested_trees_of_weighted_cacti}
\end{figure}

\begin{figure}[H]
 \centering
 \includegraphics[width=0.8\textwidth]{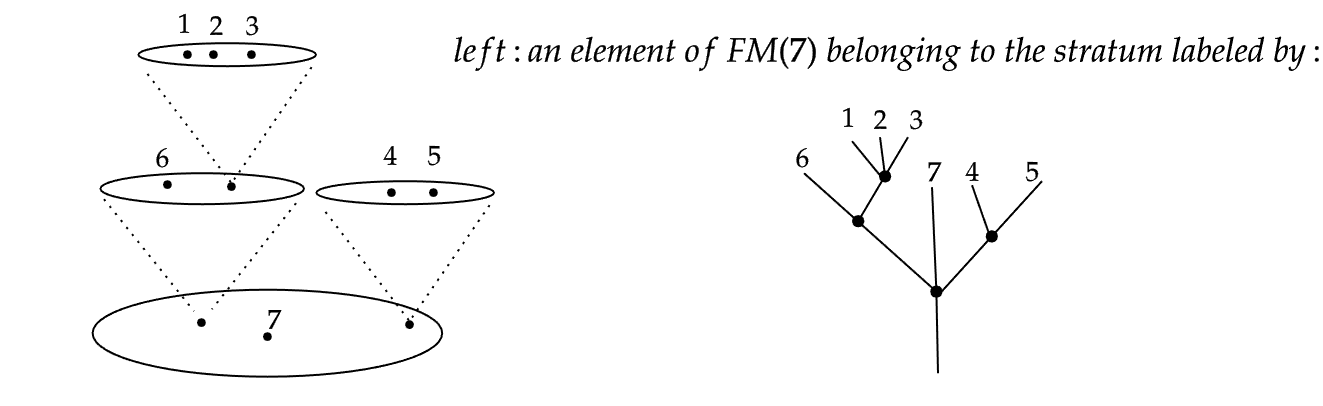}
 \caption{}
 \label{fig:Fulton_MacPherson_operad}
\end{figure}

\indent
\begin{lemma}(compare \cite[Lemma 5.8]{Sal1})\label{thm:homeo between FM' and nested trees of weighted cacti}
For $k\geq 1$, there is a $\Sigma_k$-equivariant homeomorphism 
\begin{equation}
FM'(k)\times \mathcal{P}(k)\cong N_k(\widetilde{\mathrm{Cact}}).    
\end{equation}  
\end{lemma}
\emph{Idea of Proof.} We indicate the map in the forward direction, and refer the readers to \cite[Lemma 5.8]{Sal1} for details. Start with $([z_1,\cdots,z_n],(a_1,\cdots,a_n))\in FM'(k)\times \mathcal{P}(k)$ and let $(T,f,g)$ be the associated labeled tree via Theorem \autoref{thm:from configuration to labeled trees, without S^1}. Remove from $T$ all black vertices such that $f(b)=1$ and edges incident to these vertices, the result of which is a disjoint union of subtrees. For each of the resulting subtree $T'$, let $S$ be its sets of white vertices; let $\lambda_S$ be the maximum value of a black vertex in $T'$. Finally, rescale $f|_{T'}$ by $1/\lambda_S$ (which produces black vertices with value $1$) and repeat the above procedure. This gives a nested tree $\mathcal{S}$ whose vertices correspond to the $T'$s appearing in the process, and internal edges labeled by the $\lambda_S$'s. Finally, there are unique (observe that for fixed $\mathcal{S}$, $\widetilde{\theta}_{\mathcal{S}}$ is injective) weighted cacti $(x_S)_{S\in\mathcal{S}}$ such that $\widetilde{\theta}_{\mathcal{S}}((x_S)_{\mathcal{S}})$ is the weighted cacti obtained by first retracting $(T,f,g)$ to $Tr^1_n$ (i.e. setting $f\equiv 1)$, then applying \eqref{eq:homeo from Tr^1_n to cact^n}, and finally attaching the basepoint according to the paragraph before Proposition \autoref{thm:equivariant homotopy equivalence of FM' with Cact, weighted}. See Figure \ref{fig:from_FM__to_nested_trees_of_weighted_cacti} for an illustration. \qed
\textbf{Extension to the full Fulton-MacPherson operad.} There exists a compactification of $FM'(k)=\mathrm{Conf}_k(\mathbb{C})/\mathbb{C}\rtimes\mathbb{R}_+, k\geq 1$, called the \emph{Fulton-MacPherson space} $FM(k)$ \cite{FM}, which is a manifold with corner whose strata are labeled by nested trees on $k$ leaves. Namely, 
\begin{equation}\label{eq:FM space}
FM(k):=\coprod_{\mathcal{S}\in N_k}\prod_{S\in \mathcal{S}} FM'(|S|_v).    
\end{equation}
See Figure \ref{fig:Fulton_MacPherson_operad} for an example. 
The collection $FM=\{FM(k)\}_{k\geq 1}$ has an operadic structure induced by grafting of nested trees \cite{GJ}, known as the \emph{Fulton-MacPherson operad}. For our purpose, it is more convenient to consider the \emph{weighted Fulton-MacPherson operad}  
$FM\times \mathcal{P}$ (see \eqref{eq:operadic composition of weights} for the definition of $\mathcal{P}$), which of course is quasi-equivalent to $FM$ via the projection since $\mathcal{P}$ is contractible. \par\indent
The smooth structure of $FM(k)$ near the corners is induced by `insertion of configurations', e.g. along a codimension one corner, there are gluing maps (for $\epsilon\ll 1$)
\begin{equation}
\Gamma:(0,\epsilon)\times FM'(l)\times_j FM'(k)\rightarrow FM'(k+l-1)\;,\;1\leq j\leq k
\end{equation}
given by 
\begin{equation}\label{eq:gluing of FM spaces}
\Gamma(t,[w_1,\cdots,w_l],[z_1,\cdots,z_k]):=\begin{cases}
 [z_1,\cdots,z_{j-1},z_j+t|z_j-z_{j-1}|w_1,\cdots,z_j+t|z_j-z_{j-1}|w_l,z_{j+1},\cdots,z_k],\quad\mathrm{if}\;k\geq 2\\
 [w_1,\cdots,w_l],\quad\mathrm{if}\;k=1.
\end{cases}  
\end{equation}
where we choose representative $(w_1,\cdots,w_l)$ such that $w_1=0$ and $\max|w_i|=1$ if $k\geq 2$. \par\indent
The homeomorphism of Lemma \autoref{thm:homeo between FM' and nested trees of weighted cacti} can be extended to the compactification $FM(k)$ in a way that is compatible with operadic structures.\\ 
\begin{mydef}\label{thm:compactified nested tree of weighted cacti}
Define 
 \begin{equation}\label{eq:nested tree of weighted cacti}
\overline{N}_k(\widetilde{\mathrm{Cact}}):=\coprod_{\mathcal{S}\in N_k}\prod_{S\in \mathcal{S}}\widetilde{\mathrm{Cact}^{|S|_v}}\times[0,1]^{\mathcal{S}'}\;/\sim ,   
\end{equation}
where $\sim$ is the same equivalence relation as in \autoref{thm:nested tree of weighted cacti}.
\end{mydef}
$\overline{N}_k(\widetilde{\mathrm{Cact}})$ is a compactification of $N_k(\widetilde{\mathrm{Cact}})$. Intuitively, whenever $\lambda_S=0$ for some internal edge, one should view the big nested tree as a `broken union' of two nested trees along this edge.\\
\begin{thm}(\cite{Sal1})\label{thm:FM homeomorphic to compactified nested tree of weighted cacti}
There is an isomorphism of topological operads
\begin{equation}
FM\times \mathcal{P}\cong \{\overline{N}_k(\widetilde{\mathrm{Cact}})\}_{k\geq 1},   
\end{equation}
where the operadic structure on the right hand side is given by grafting of nested trees with labels. \qed
\end{thm} 
\begin{figure}[H]
 \centering
 \includegraphics[width=0.9\textwidth]{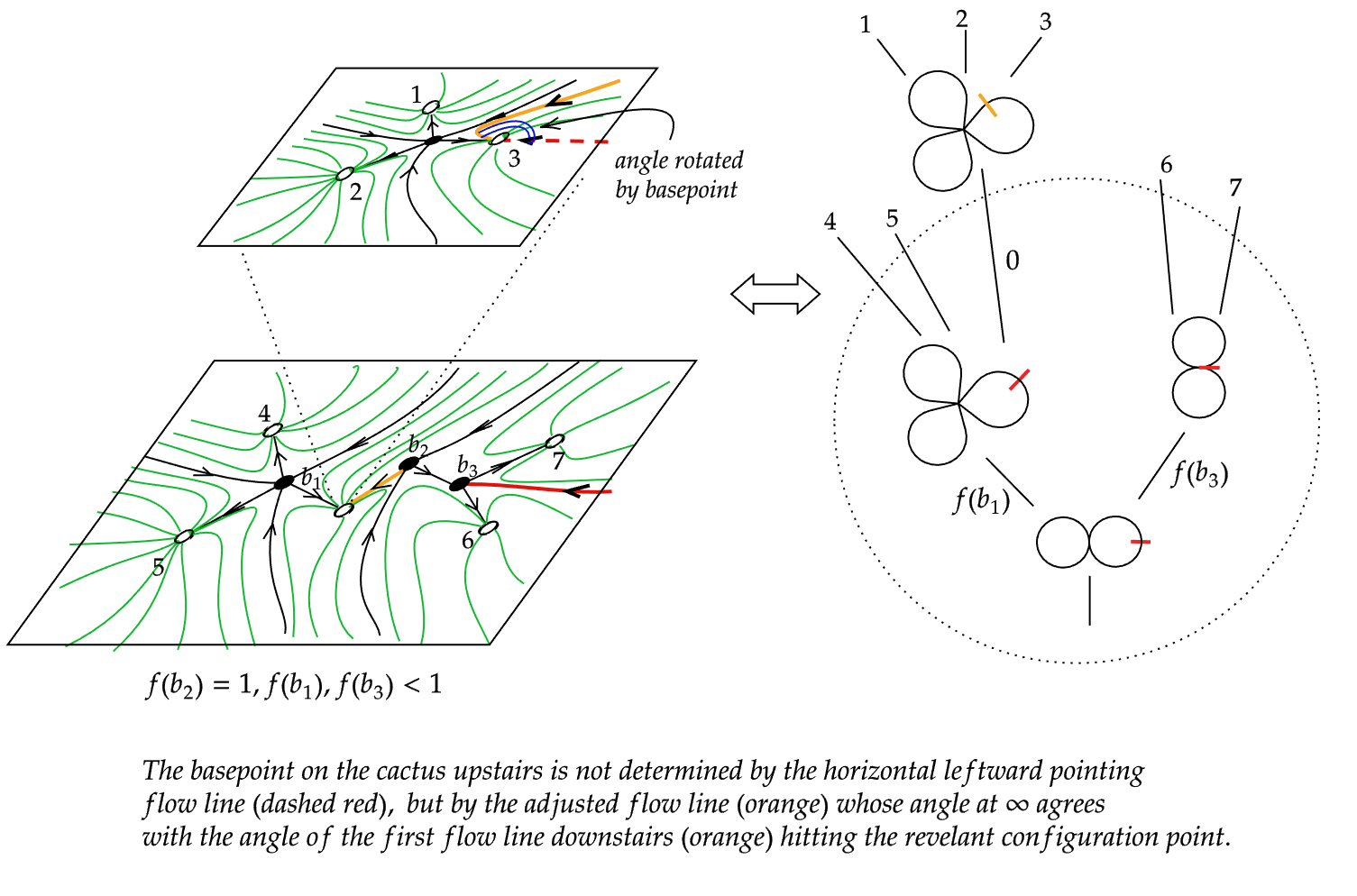}
 \caption{}
 \label{fig:from_FM_to_compactified_nested_trees}
\end{figure}
Another name for $\{\overline{N}_k(\widetilde{\mathrm{Cact}})\}_{k\geq 1}$ is the Boardman-Vogt resolution/W-construction \cite{BV} of the topological operad $\widetilde{\mathrm{Cact}}$. In other words, Theorem \autoref{thm:FM homeomorphic to compactified nested tree of weighted cacti} establishes the weighted Fulton-MacPherson operad as the W-construction of weighted cacti. \par\indent
Roughly speaking, the forward map in Theorem \autoref{thm:FM homeomorphic to compactified nested tree of weighted cacti} is obtained by `grafting together' instances of Lemma \autoref{thm:homeo between FM' and nested trees of weighted cacti} along a nested tree that represents a stratum of the Fulton-MacPherson space, whilst suitably rotating the basepoints of cacti on the nested tree (see Remark \autoref{thm:subtle rotations needed to prove FM = W(weighted cacti)} for a technical subtlety); see Figure \ref{fig:from_FM_to_compactified_nested_trees}. \\
\begin{rmk}\label{thm:subtle rotations needed to prove FM = W(weighted cacti)}
\begin{itemize}
    \item We remark that nested trees appeared in two different ways. On one hand, the strata of the Fulton-MacPherson space $FM(k)$ corresponds to $\mathcal{S}\in N_k$; on the other hand, the open stratum $FM'(k)$ further decomposes (up to taking product with $\mathcal{P}$) into a nested tree of cacti.
    \item There is a subtlety in constructing the isomorphism of Theorem \autoref{thm:FM homeomorphic to compactified nested tree of weighted cacti}: as the basepoints of cacti are rotated (this is to ensure continuity), the map depicted in Figure \ref{fig:from_FM_to_compactified_nested_trees} does not respect operadic structures (which recall, is given by grafting together nested trees). This is resolved in \cite[Section 7]{Sal1} by a technical argument that crucially relies on an isomorphism of topological operads $FM\cong WFM$ \cite{Sal3}.  
\end{itemize}
\end{rmk}
\subsection{Proof of Theorem \autoref{thm:comparison of KS with Cyl}}
In this subsection, we prove Theorem \autoref{thm:comparison of KS with Cyl} by adapting the methods of Section 6.1. Let 
\begin{equation}\label{eq:FM' circ}
FM_{\circlearrowright}'(n):=(\mathrm{Conf}_n(\mathbb{C}^*)/\mathbb{R}_+\times S^1\times S^1)/S^1, 
\end{equation}
where $S^1$ acts diagonally (it acts on $\mathrm{Conf}_n(\mathbb{C}^*)/\mathbb{R_+}$ by rotating configurations). In words, $FM'_{\circlearrowright}(n)$ is the configuration space of $n$ points on the infinite cylinder modulo translation, together with asymptotic markers (i.e. unit tangent directions) at the two ends, modulo overall rotation. There is a natural $S^1\times S^1$-action on $FM'_{\circlearrowright}(n)$ given by rotating the two asymptotic markers.  \par\indent
\textbf{From $FM'_{\circlearrowright}$ to cyclic cacti}. Given $(z_1,\cdots,z_n)\in \mathrm{Conf}_n(\mathbb{C}^*)$ and weights $(a_0,\cdots,a_n)\in \mathcal{P}(n+1)$, consider the vector field 
\begin{equation}\label{eq:vector field E circ}
E_{\circlearrowright}(z):=a_0\frac{z}{|z|^2}+\sum_{i=1}^n   a_i\frac{z-z_i}{|z-z_i|^2}  
\end{equation}
on the complex plane. By the local behavior analyzed in Section 6.1 (cf. Figure \ref{fig:local_behaviour_of_E}), we may view $E_{\circlearrowright}$ as a vector field on $\mathbb{C}^*$ with a sink at $0$ (as well as at the $z_i$'s) and a source at $\infty$. \par\indent
Let $Tr_{n+1}$ be the space of labeled trees with $n$ leaves, labeled by $\{0,1,\cdots,n\}$. Since $\mathrm{Conf}_n(\mathbb{C}^*)\cong\mathrm{Conf}_{n+1}(\mathbb{C})/\mathbb{C}$, Theorem \autoref{thm:from configuration to labeled trees, without S^1} immediately implies that\\
\begin{cor}\label{thm:from configuration on C^* to labeled trees, without S^1}
For each set of weights $(a_0,\cdots,a_n)\in\mathcal{P}(n+1)$, the map
\begin{equation}\label{eq:from configuration on C^* to labeled trees, without S^1}
\Phi_{a_0,\cdots,a_n}^\circlearrowright:\mathrm{Conf}_n(\mathbb{C}^*)/\mathbb{C}^*\cong\mathrm{Conf}_{n+1}(\mathbb{C})/\mathbb{C}\rtimes\mathbb{C}^*\xrightarrow[\cong]{\Phi_{a_0,\cdots,a_n}} Tr_{n+1}  
\end{equation}
is a homeomorphism. \qed
\end{cor}
We would now like to put in the $S^1\times S^1$-actions and upgrade \eqref{eq:from configuration on C^* to labeled trees, without S^1} to an equivariant comparison between $FM'_{\circlearrowright}(n)$ and the space of cyclic cacti $\mathrm{Cact}^n_{\circlearrowright}$.\\
\begin{prop}\label{thm:from FM' circ to weighted cyclic cacti, equivariantly}
There is an $S^1\times S^1\times \Sigma_n$-equivariant homotopy equivalence $FM'_{\circlearrowright}(n)\times \mathcal{P}(n+1)\simeq \widetilde{\mathrm{Cact}^n_{\circlearrowright}}$.    
\end{prop}
\emph{Proof}. Consider the $S^1\times S^1\times \Sigma_n$-invariant subspace $FM^{',1}_{\circlearrowright}(n)\subset FM'_{\circlearrowright}(n)$ consisting of equivalence classes of triples $([z_1,\cdots,z_n],\theta_0,\theta_{\infty})$ such that $\Phi_{a_0,\cdots,a_n}^{\circlearrowright}[z_1,\cdots,z_n]$ lies in $Tr^1_{n+1}\subset Tr_{n+1}$. The inclusion $FM^{',1}_{\circlearrowright}(n)\subset FM'_{\circlearrowright}(n)$ is a homotopy equivalence since $Tr_{n+1}$ deformation retracts onto $Tr^1_{n+1}$. \par\indent
Given $[([z_1,\cdots,z_n],\theta_0,\theta_{\infty})]\in FM^{',1}_{\circlearrowright}(n)$, denote $(T,f,g)=\Phi_{a_0,\cdots,a_n}^{\circlearrowright}[z_1,\cdots,z_n]$, and let $x'$ be the unbased cactus given by the image of $(T,f,g)$ under \eqref{eq:homeo from Tr^1_n to cact^n}, with its lobes weighted by $(a_0,\cdots,a_n)$. Recall that to upgrade $x'$ to a cyclic cactus (cf. Definition \autoref{thm:space of cyclic cacti}), we need to add an \emph{input basepoint} on its periphery (making it a based cactus) and an \emph{output basepoint} (lying on the $0$-th lobe). When $n=0$ (i.e. $x'$ only has a $0$-th lobe), the output and input basepoints are just given by $\theta_0+\frac{1}{2},\theta_{\infty}+\frac{1}{2}\in S^1=\mathbb{R}/\mathbb{Z}$, respectively.  \par\indent
Now consider the case $n\geq 1$. There is a unique zigzag of flow lines (as $T$ is a tree) $\gamma_0,\cdots,\gamma_m$ of $E_{\circlearrowright}$ connecting $0$ to $z_1$. Let $e^{-2\pi i\theta_0'}:=\lim_{t\rightarrow\infty}\frac{\gamma_0'(t)}{|\gamma_0'(t)|}$. Then the output basepoint is given by the difference $\theta_0-\theta_0'\in S^1$. \par
\begin{figure}[H]
 \centering
 \includegraphics[width=1.05\textwidth]{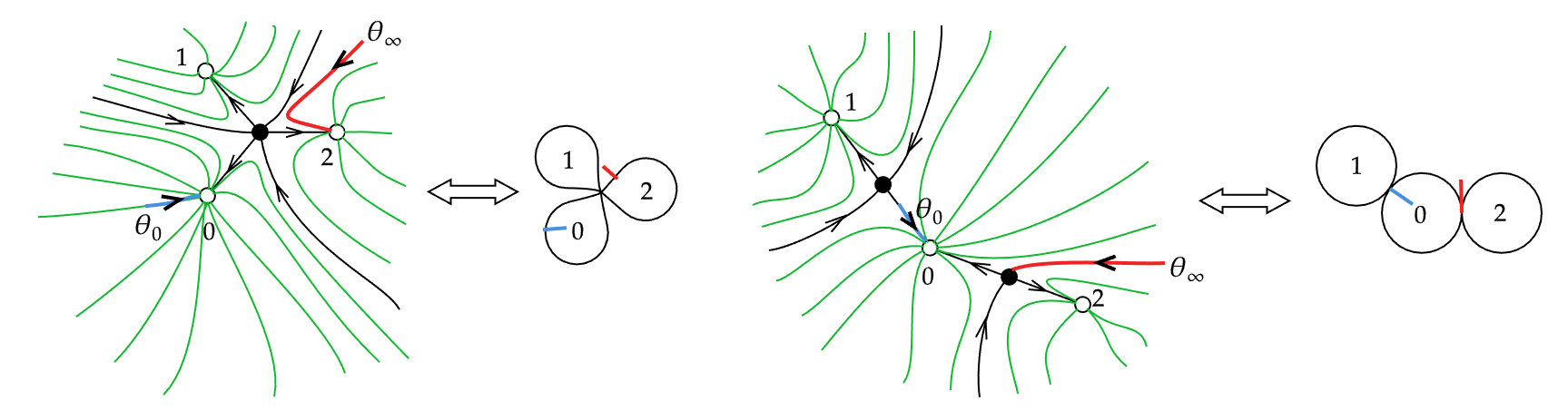}
 \caption{}
 \label{fig:from_FM__circ_to_cyclic_cacti}
\end{figure}

\indent
The input basepoint is determined as follows. Let $\gamma$ be the (broken) flow line of $E_{\circlearrowright}$ such that $\lim_{t\rightarrow-\infty}\gamma(t)=\infty$ and $\lim_{t\rightarrow-\infty}\frac{\gamma'(t)}{|\gamma'(t)|}=\theta_{\infty}$. Then after potentially breaking at critical points of $E_{\circlearrowright}$, $\gamma$ eventually reaches some $z_i$ (which might not be unique, but the resulting input basepoint is well-defined), and the input basepoint records minus the limiting tangent direction $-\lim_{t\rightarrow+\infty}\frac{\gamma'(t)}{|\gamma'(t)|}$; see Figure \ref{fig:from_FM__circ_to_cyclic_cacti}. It is easy to see that the resulting weighted cyclic cacti does not depend on an overall rotation of the initial data, and this assignment is $S^1\times S^1\times \Sigma_n$-equivariant. \qed
\textbf{Forest of nested trees}. Similar to Lemma \autoref{thm:homeo between FM' and nested trees of weighted cacti}, one can further upgrade Proposition \autoref{thm:from FM' circ to weighted cyclic cacti, equivariantly} to an equivariant \emph{homeomorphism}, at the cost of replacing $\widetilde{\mathrm{Cact}^n_{\circlearrowright}}$ by its cobar construction. Since $\widetilde{\mathrm{Cact}^n_{\circlearrowright}}$ is part of a two-colored topological operad, the appropriate combinatorial model for the cobar construction is based on the following variant of nested trees. \\
\begin{mydef}\label{thm:generalized nested tree}
The \emph{set of generalized nested trees with $n$ leaves $\hat{N}_n$} is
\begin{itemize}
    \item For $n\geq 2$: $\hat{N}_n=N_n^+\sqcup N_n^-$ where $N_n^-=N_n^+=N_n$ (cf. Definition \autoref{thm:nested trees}) are just two copies of nested trees with $n$ leaves. If $\mathcal{S}\in N_n$, we let $\mathcal{S}^{\pm}$ denote the same nested tree, but viewed as an element of $N_n^{\pm}$. 
    \item For $n=1$: $\hat{N}_1=N_1^-\sqcup N_1^+$, where $N_1^-=\emptyset$ and $N_1^+$ is a singleton, which we think of as the unique tree with one root, one internal vertex labeled by $\emptyset$, and one leaf.  
\end{itemize}
\end{mydef}\par\indent
\begin{mydef}\label{thm:forest of nested trees}
A \emph{forest of nested tree on $n$ leaves} is a partition $\{1,\cdots,n\}=\bigsqcup_{i=1}^m A_i$ together with a sequence $(\mathcal{S}_1^{\pm},\cdots,\mathcal{S}_m^{\pm})$ of generalized nested trees such that $\mathcal{S}_i^{\pm}\in \hat{N}_{n_i}$, where $|A_i|=n_i$ (i.e. the leaves of $\mathcal{S}_i^{\pm}$ are labeled by the set $A_i$). We often omit the $A_i$'s when referring to a forest of nested trees. Let $N_n^\circlearrowright$ denote the set of forests of nested trees on $n$ leaves.  
\end{mydef}
For $1\leq i\leq k$, one can \emph{graft} $\mathcal{S}\in N_l$ onto the $i$-th leaf of $(\mathcal{S}_1^{\pm},\cdots,\mathcal{S}_m^{\pm})\in N_k^{\circlearrowright}$ and obtain an element of $N^{k+l-1}_{\circlearrowright}$; one can also \emph{concatenate} two sequences $(\mathcal{T}_1^\pm,\cdots,\mathcal{T}^\pm_{m'})\in N^l_\circlearrowright$ and $(\mathcal{S}_1^{\pm},\cdots,\mathcal{S}_{m}^{\pm})\in N^k_\circlearrowright$ to get an element $(\mathcal{S}_1^{\pm},\cdots,\mathcal{S}_{m}^{\pm},\mathcal{T}_1^\pm,\cdots,\mathcal{T}^\pm_{m'})\in N^{k+l}_\circlearrowright$. Additionally, if $\{O,O_{\circlearrowright},\theta,\theta^{\circlearrowright,0},\theta^{\circlearrowright,1}\}$ (compare notations in \eqref{eq:cyclic cacti operadic structure type 0} and \eqref{eq:cyclic cacti operadic structure type 1}) is a two-colored topological operad, then a forest of nested trees on $n$ leaves $(\mathcal{S}_1^{\pm},\cdots,\mathcal{S}_m^{\pm})\in N^{\circlearrowright}_n, \mathcal{S}_i^{\pm}\in \hat{N}_{n_i}$ prescribes a well-defined composition
\begin{equation}\label{eq:operadic composition induced by forest of nested trees}
\theta^{\circlearrowright}_{(\mathcal{S}_1^{\pm},\cdots,\mathcal{S}_m^{\pm})}:\prod_{i=1}^m\begin{cases}
O_{\circlearrowright}(|A_i|_v)\times\prod_{S\in \mathcal{S}_i'}O(|S|_v)\quad,\;\mathrm{if}\;\mathcal{S}_i^{\pm}=\mathcal{S}_i^-\in N_{n_i}^-  \\
O_{\circlearrowright}(1)\times\prod_{S\in \mathcal{S}_i} O(|S|_v)\quad,\;\mathrm{if}\;\mathcal{S}_i^{\pm}=\mathcal{S}_i^+\in N_{n_i}^+
\end{cases}  \longrightarrow O_{\circlearrowright}(n).
\end{equation}\par\indent
\begin{mydef}\label{thm:forest of nested trees of weighted (cyclic) cacti}
The \emph{space of forests of nested trees on $n$ leaves with vertices labeled by weighted (cyclic) cacti and internal edges labeled by $(0,1)$} is
\begin{equation}\label{eq:forest of nested trees of weighted (cyclic) cacti}
N_n^{\circlearrowright}(\widetilde{\mathrm{Cact}},\widetilde{\mathrm{Cact}_{\circlearrowright}}):=\coprod_{(\mathcal{S}_1^\pm,\cdots,\mathcal{S}_m^\pm)\in N_n^{\circlearrowright}}(0,1]^{m-1}\times \prod_{i=1}^m\begin{cases}
\widetilde{\mathrm{Cact}_{\circlearrowright}}^{|A_i|_v}\times\prod_{S\in \mathcal{S}_i'}\widetilde{\mathrm{Cact}}^{|S|_v}\times(0,1]^{\mathcal{S}_i'}\;,\;\mathrm{if}\;\mathcal{S}_i^{\pm}=\mathcal{S}_i^-\in N_{n_i}^-  \\
\widetilde{\mathrm{Cact}^1_{\circlearrowright}}\times\prod_{S\in \mathcal{S}_i} \widetilde{\mathrm{Cact}}^{|S|_v}\times (0,1]^{\mathcal{S}_i}\;,\;\mathrm{if}\;\mathcal{S}_i^{\pm}=\mathcal{S}_i^+\in N_{n_i}^+
\end{cases}\Big/\sim
\end{equation}
where $\sim$ is the same equivalence relation as in Definition \autoref{thm:nested tree of weighted cacti} except one replaces $\widetilde{\theta}_{\mathcal{S}}$ by $\widetilde{\theta}^\circlearrowright_{(\mathcal{S}_1^\pm,\cdots,\mathcal{S}_m^\pm)}$. See Figure \ref{fig:forest_of_nested_trees_of_weighted_cyclic_cacti} for an example.
\begin{figure}[H]
 \centering
 \includegraphics[width=1.0\textwidth]{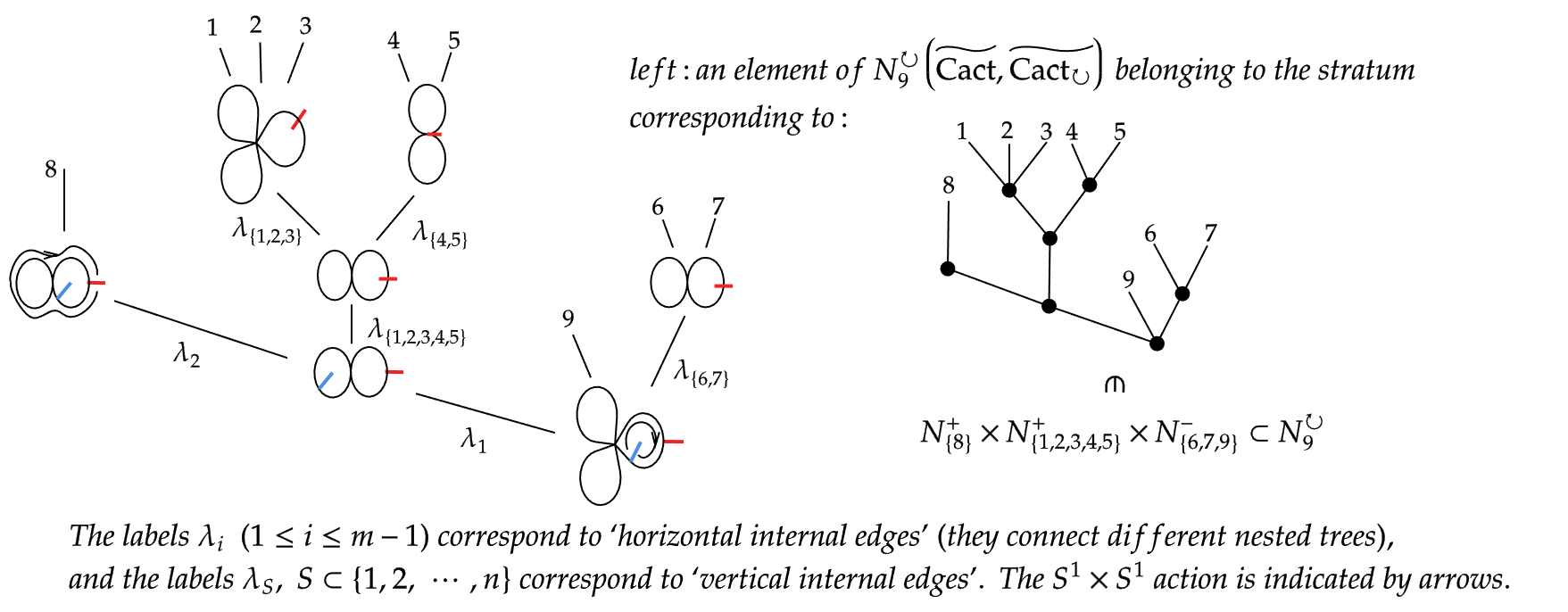}
 \caption{}
 \label{fig:forest_of_nested_trees_of_weighted_cyclic_cacti}
\end{figure}
In particular, an element 
\begin{equation}
[(\lambda_i)_{i=1}^{m-1} ,(x_i,(x^i_S)_{S\in \mathcal{S}_i'\;\mathrm{or}\;\mathcal{S}_i},(\lambda^i_S)_{S\in \mathcal{S}_i'\;\mathrm{or}\;\mathcal{S}_i})_{i=1}^m]  \in N_n^{\circlearrowright}(\widetilde{\mathrm{Cact}},\widetilde{\mathrm{Cact}_{\circlearrowright}}),
\end{equation}
contains a sequence of weighted cyclic cacti $(x_1,\cdots,x_m)$ sitting at the roots of the nested trees. In light of this, there is a natural action
\begin{equation}\label{S^1 times S^1 action on forest of nested trees of weighted (cyclic) cacti} 
S^1\times S^1\acts N_n^{\circlearrowright}(\widetilde{\mathrm{Cact}},\widetilde{\mathrm{Cact}_{\circlearrowright}})
\end{equation}
where the first $S^1$ factor rotates the input basepoint of $x_m$ and the second $S^1$ factor rotates the output basepoint of $x_1$. \\
\end{mydef}
\begin{lemma}\label{thm:homeo between FM' circ and forest of nested trees of weighted cyclic cacti}
There is an $S^1\times S^1\times \Sigma_k$-equivariant homeomorphism
\begin{equation}\label{eq:homeo between FM' circ and forest of nested trees of weighted cyclic cacti}
FM'_{\circlearrowright}(n)\times \mathcal{P}(n+1)\cong N_n^{\circlearrowright}(\widetilde{\mathrm{Cact}},\widetilde{\mathrm{Cact}_{\circlearrowright}}).
\end{equation}
\end{lemma}
\emph{Proof}. We produce the forward map in \eqref{eq:homeo between FM' circ and forest of nested trees of weighted cyclic cacti}. Given an element $([[z_1,\cdots,z_n],\theta_0,\theta_{\infty}],(a_0,\cdots,a_n))\in FM'_{\circlearrowright}(n)\times \mathcal{P}(n+1)$, let $(T,f,g)=\Phi_{a_0,\cdots,a_n}^\circlearrowright[z_1,\cdots,z_n]$ be the associated labeled tree of the configuration, whose $i$-th white point (which we abuse notation and also call $z_i$) has weight $a_i$. Additionally, the tree comes with two asymptotic markers $\theta_0$ and $\theta_{\infty}$, which are given by unit tangent directions induced from the flow of $E_\circlearrowright$, at $z_0$ and the outer periphery of the tree, respectively; see e.g. Figure \ref{fig:from_FM__circ_to_cyclic_cacti}.\par\indent 
Let $B_1=\{b\in B:f(b)=1\}$, which is nonempty by assumption. Pick a $b\in B_1$ which is farthest away from $z_0$ (where distance is measured by the length of the unique zigzag of edges connecting them); such $b$ might not be unique, but the end result will be independent of choices. We discuss two cases.  \par\indent
\textbf{I. Vertical decomposition}. Suppose the distance from $b$ to $z_0$ is $1$ (i.e. there is an edge from $b$ to $z_0$). This case is similar to Lemma \autoref{thm:homeo between FM' and nested trees of weighted cacti}. Namely, remove from $T$ the vertex $b$ and all edges incident to it, resulting in a disjoint union of trees. There is a unique tree containing $z_0$ (which naturally inherits asymptotic markers from this decomposition), and a few other trees $T_1,\cdots,T_l$. Let $T_0$ be the tree obtained from $T$ by collapsing each $T_i, 1\leq i\leq l$ to a single white vertex. For $1\leq i\leq l$, let $W_i$ be the set of white vertices of $T_i$, and let $\lambda_{W_i}$ be the maximum value of a black vertex in $T_i$; let $\lambda_0$ be the maximal value of a black vertex in $T_0$ (it is possible that $\lambda_{W_i}$ or $\lambda_0$ is $1$). Intuitively, we think of the result as a partial nested tree with $T_1,\cdots,T_l$ living over the root $T_0$. Finally, rescale $f|_{T_i}$ by $1/\lambda_{W_i}$ and $f|_{T_0}$ by $1/\lambda_0$. \par\indent
\textbf{II. Horizontal decomposition}. Suppose the distance from $b$ to $z_0$ is at least $2$. Again, we remove from $T$ the black vertex $b$ and edges incident to it, resulting a disjoint union of trees, of which $T_0$ contains $z_0$. Let $T_0'$ be the tree obtained from $T$ by collapsing $T_0$ to a single white vertex; this will be the new $0$-th white vertex of $T_0'$ and denoted $z_0'$. Intuitively, we think of $T$ being decomposed into a forest of nested trees $T_0'$ and $T_0$, which are connected by a horizontal internal edge; the length of this internal edge $\lambda_0$ is the maximal value of a black vertex in $T_0$. Furthermore, $T_0,T_0'$ naturally inherit asymptotic markers from this decomposition. Finally, we rescale $f|_{T_0}$ by $1/\lambda_0$ and $f|_{T_0'}$ by $1/\lambda_0'$. \par
\begin{figure}[H]
 \centering
 \includegraphics[width=1.0\textwidth]{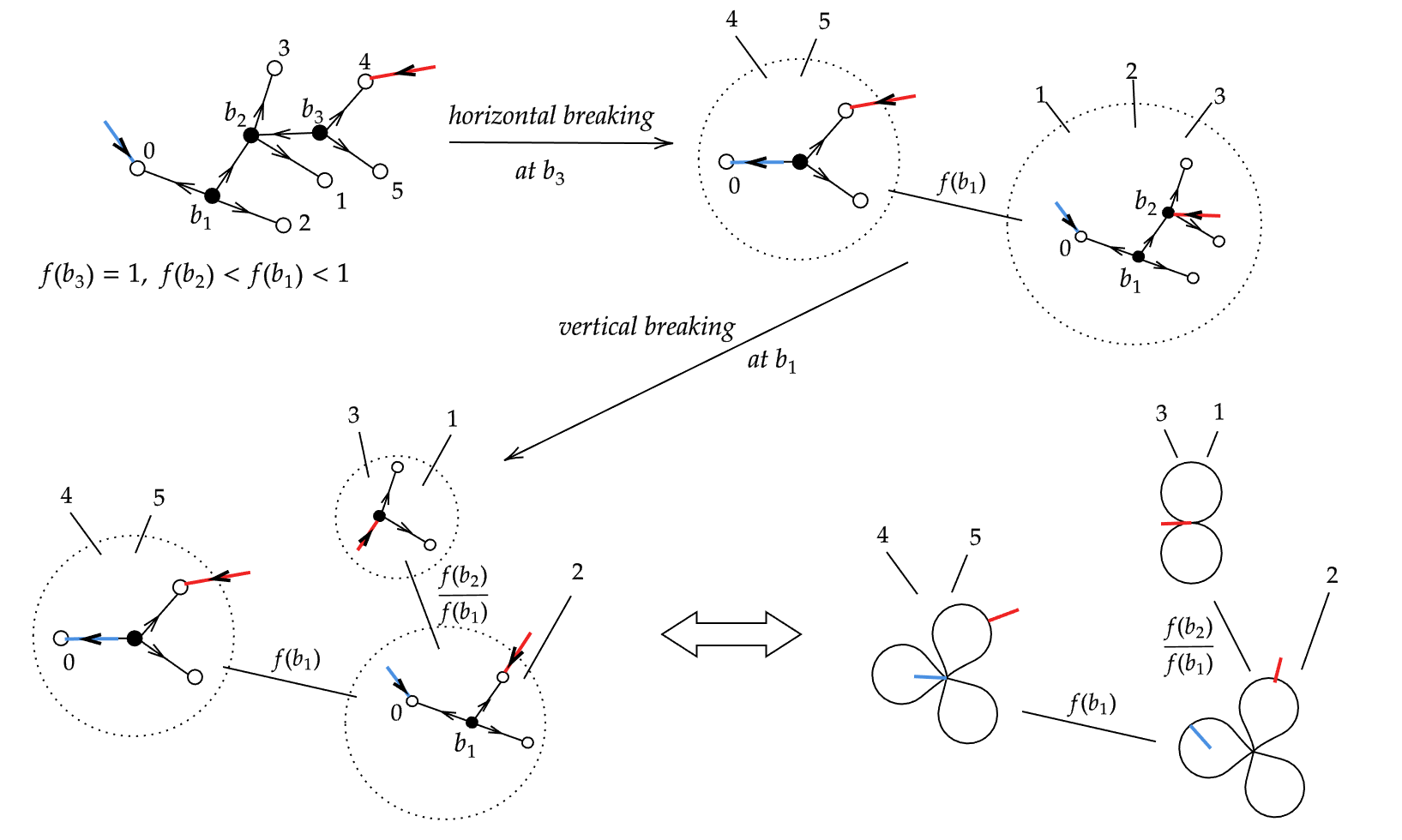}
 \caption{}
 \label{fig:from_FM__circ_to_forest_of_nested_trees_of_weighted_cyclic_cacti}
\end{figure}

\indent
Then we repeat. More precisely, for trees of type $T_i,i\geq 1$ that appeared in the process, we run the algorithm of Lemma \autoref{thm:homeo between FM' and nested trees of weighted cacti}; for trees of type $T_0,T_0'$ (i.e. has a distinguished $0$-th white vertex), we run the algorithm of Case I and II above. This eventually produces a forest $(\mathcal{S}_1,\cdots,\mathcal{S}_m)$ of nested trees, with internal edges labeled by $(0,1]$. The weighted (cyclic) cacti labeling the vertices of this forest are uniquely determined by the properties that 
\begin{itemize}
    \item Their total composition along this forest (cf. \eqref{eq:operadic composition induced by forest of nested trees}) agrees with the weighted cyclic cactus obtained from $(T,f,g)$ by first setting $f\equiv 1$ (in particular, collapsing all edges between black vertices) and then applying Proposition \autoref{thm:from FM' circ to weighted cyclic cacti, equivariantly}.
    \item The input and output basepoints of the weighted cyclic cacti lying at the roots of nested trees in the forest agree with the asymptotic markers of the relevant $T_0$ or $T_0'$ produced in the process. 
\end{itemize}
The property of being $S^1\times S^1\times\Sigma_n$-equivariant is evident from construction. See Figure \ref{fig:from_FM__circ_to_forest_of_nested_trees_of_weighted_cyclic_cacti} for an illustration. 
To argue that this map is a homeomorphism, one proceeds as in \cite[Lemma 5.8]{Sal1} to produce an explicit inverse; we leave the details to the readers.\qed  \par\indent
\textbf{Extension to the two-colored Fulton-MacPherson operad.} 
Similar to the usual Fulton-MacPherson space, there exists a compactification of $FM'_{\circlearrowright}(k)=(\mathrm{Conf}_k(\mathbb{C}^*)/\mathbb{R}_+\times S^1\times S^1)/S^1, k\geq 1$, which we denote by $FM_{\circlearrowright}(k)$.\par\indent
Intuitively, given a configuration in $FM'_{\circlearrowright}(k)$, two types of (codimension $1$) limiting behavior can occur. First, some marked points can get infinitesimally close to one another and `bubble off' at an interior point of the cylinder; second, some marked points can tend to one of the infinite ends with uniform speed, and then another cylinder `bubbles off' the original one.\par
\begin{figure}[H]
 \centering
 \includegraphics[width=1.0\textwidth]{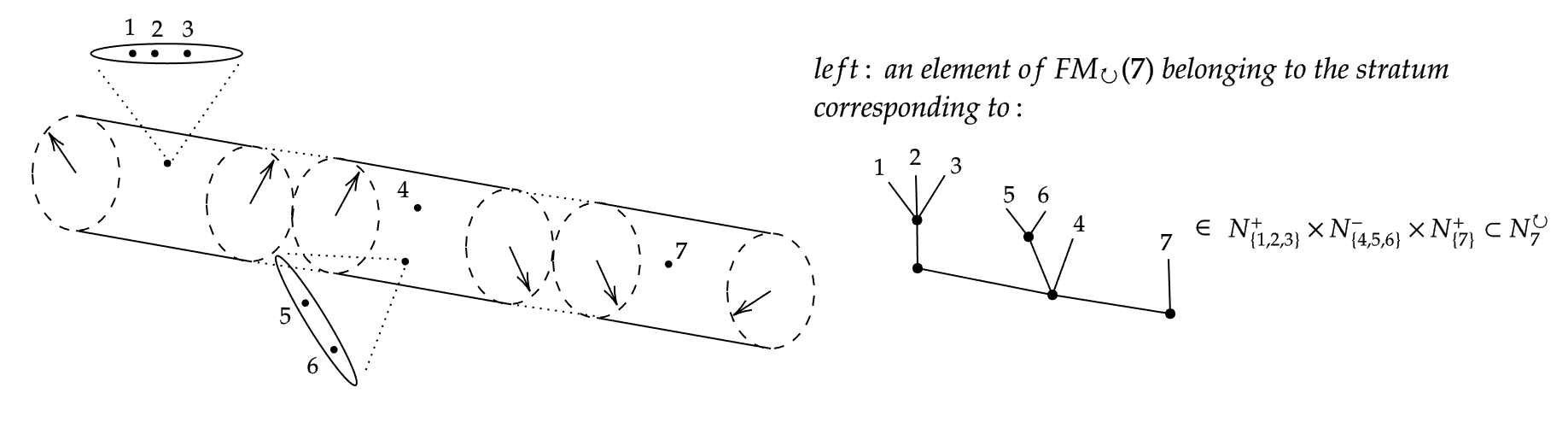}
 \caption{}
 \label{fig:FM_circ_space}
\end{figure}

\indent
Inductively, the above procedure indicates a manifold with corner structure on $FM_{\circlearrowright}(k)$ whose strata are labeled by forests of nested trees on $k$ leaves. Namely, using the notation of Definition \autoref{thm:forest of nested trees},
\begin{equation}\label{eq:FM_circ_space}
FM_{\circlearrowright}(k):=\coprod_{(\mathcal{S}_1^\pm,\cdots,\mathcal{S}_m^\pm)\in N_k^{\circlearrowright}}\prod_{i=1}^m\begin{cases}
FM'_{\circlearrowright}(|A_i|_v)\times\prod_{S\in \mathcal{S}_i'} FM'(|S|_v)\;\;,\;\mathrm{if}\;\mathcal{S}_i^\pm=\mathcal{S}_i^-\in N_{n_i}^-\\
FM'_{\circlearrowright}(1)\times\prod_{S\in \mathcal{S}_i} FM'(|S|_v)\;\;,\;\mathrm{if}\;\mathcal{S}_i^\pm=\mathcal{S}_i^+\in N_{n_i}^+.
\end{cases}    
\end{equation}
See Figure \ref{fig:FM_circ_space} for an example.
The usual Fulton-MacPherson operad can be extended to 
\begin{equation}\label{eq:FM two colored operad}
\{FM(l),FM_{\circlearrowright}(k)\}_{l\geq 1,k\geq 0}    
\end{equation}
form a two-colored topological operad, with two additional types of operadic structure maps induced by grafting of a nested tree onto a forest $FM(l)\times FM_{\circlearrowright}(k)\rightarrow FM_{\circlearrowright}(k+l-1)$ and concatenation of forests $FM(l)_{\circlearrowright}\times FM_{\circlearrowright}(k)\rightarrow FM_{\circlearrowright}(k+l)$. As in the case of Definition \autoref{thm:forest of nested trees of weighted (cyclic) cacti}, there is an $S^1\times S^1$-action on $FM_{\circlearrowright}(k)$ given by rotating the top asymptotic marker on the last cylinder in the forest, and the bottom asymptotic marker on the first cylinder. \\
\begin{mydef}\label{thm:compactified forest of nested trees of weighted (cyclic) cacti}
Define
\begin{equation}\label{eq:compactified forest of nested trees of weighted (cyclic) cacti}
\overline{N_n^{\circlearrowright}}(\widetilde{\mathrm{Cact}},\widetilde{\mathrm{Cact}_{\circlearrowright}}):=\coprod_{(\mathcal{S}_1^\pm,\cdots,\mathcal{S}_m^\pm)\in N_n^{\circlearrowright}}[0,1]^{m-1}\times \prod_{i=1}^m\begin{cases}
\widetilde{\mathrm{Cact}_{\circlearrowright}}^{|A_i|_v}\times\prod_{S\in \mathcal{S}_i'}\widetilde{\mathrm{Cact}}^{|S|_v}\times[0,1]^{\mathcal{S}_i'}\;,\;\mathrm{if}\;\mathcal{S}_i^{\pm}=\mathcal{S}_i^-\in N_{n_i}^-  \\
\widetilde{\mathrm{Cact}^1_{\circlearrowright}}\times\prod_{S\in \mathcal{S}_i} \widetilde{\mathrm{Cact}}^{|S|_v}\times [0,1]^{\mathcal{S}_i}\;,\;\mathrm{if}\;\mathcal{S}_i^{\pm}=\mathcal{S}_i^+\in N_{n_i}^+
\end{cases}\Big/\sim
\end{equation}
where $\sim$ is the same equivalence relation as in Definition \autoref{thm:forest of nested trees of weighted (cyclic) cacti}.
\end{mydef}
As before, if an internal edge (either horizontal or vertical) has length $0$, we think of it as `broken'. Therefore, an element of $\overline{N_n^{\circlearrowright}}(\widetilde{\mathrm{Cact}},\widetilde{\mathrm{Cact}_{\circlearrowright}})$ can be viewed as a forest of nested trees whose vertices are labeled by either elements of $N_k^{\circlearrowright}(\widetilde{\mathrm{Cact}},\widetilde{\mathrm{Cact}_{\circlearrowright}})$ (if it is a root) or $N_l(\widetilde{\mathrm{Cact}})$ (if it is not a root), i.e. these labels themselves are nested trees or forests of nested trees whose vertices are labeled by weighted (cyclic) cacti and edges labeled by $(0,1)$.  \par\indent
By grafting and concatenating together instances of Lemma \autoref{thm:homeo between FM' and nested trees of weighted cacti} and Lemma \autoref{thm:homeo between FM' circ and forest of nested trees of weighted cyclic cacti} along a forest of nested trees representing a stratum of $FM_\circlearrowright(k)$, one obtains the following theorem. Again, there is a technical issue involving rotation of basepoints of cacti that can be resolved similarly as in \cite[Section 7]{Sal1}. \\
\begin{thm}\label{thm:isomorphism of two colored FM and compactified (forests of) nested tree of
weighted (cyclic) cacti}
There is an isomorphism of two-colored topological operads
\begin{equation}\label{eq:from FM circ to compact (forest of) nested trees of (cyclic) cacti}
\{FM(l)\times\mathcal{P}(l),FM_{\circlearrowright}(k)\times \mathcal{P}(k+1)\}_{l\geq 1,k\geq 0}\cong \{\overline{N_l}(\widetilde{\mathrm{Cact}}),\overline{N_k^{\circlearrowright}}(\widetilde{\mathrm{Cact}},\widetilde{\mathrm{Cact}_{\circlearrowright}})\}_{l\geq 1,k\geq 0}
\end{equation}
which is equivariant with respect to symmetric group actions, and furthermore $S^1\times S^1$-equivariant when restricted to
\begin{equation}
FM_{\circlearrowright}(k)\times \mathcal{P}(k+1)\cong \overline{N_k^{\circlearrowright}}(\widetilde{\mathrm{Cact}},\widetilde{\mathrm{Cact}_{\circlearrowright}}).     
\end{equation}
\qed
\end{thm}
Since $\{\overline{N_l}(\widetilde{\mathrm{Cact}}),\overline{N_k^{\circlearrowright}}(\widetilde{\mathrm{Cact}},\widetilde{\mathrm{Cact}_{\circlearrowright}})\}_{l\geq 1,k\geq 0}$ is nothing other than the W-construction of the two-colored operad of weighted (cyclic) cacti $\{\widetilde{\mathrm{Cact}^l},\widetilde{\mathrm{Cact}^k_{\circlearrowright}})\}_{l\geq 1,k\geq 0}$, the total operadic composition induces a quasi-equivalence from the former to the latter that has the same equivariance properties as in Theorem \autoref{thm:isomorphism of two colored FM and compactified (forests of) nested tree of
weighted (cyclic) cacti}. \par\indent
On the other hand, the left hand side of \eqref{eq:from FM circ to compact (forest of) nested trees of (cyclic) cacti} is (equivariantly) quasi-equivalent to the two-colored operad $\{FM(l),FM_{\circlearrowright}(k)\}_{l\geq 1,k\geq 0}$ by projection. Finally, it is a standard argument to show that $\{FM(l),FM_{\circlearrowright}(k)\}_{l\geq 1,k\geq 0}$ and $\mathrm{Cyl}$ (cf. Section 4.3) are equivariantly quasi-equivalent (by a straightforward generalization of the classical quasi-equivalence between the little disks and Fulton-MacPherson operad \cite[Proposition 3.9]{Sal4}). This concludes the proof of Theorem \autoref{thm:comparison of KS with Cyl}.\qed

\subsection{$\Xi^p([e_0])([\phi],-)$ and equivariant cap products}
The goal of this subsection is to obtain explicit formulae for the endomorphisms (from which all $C_p$-Kontsevich-Soibelman operations can be derived, cf. Theorem \autoref{thm:classification of untwisted KSp operations})
\begin{equation}
\Xi^p([e_0])([\phi],-): CC^{C_p,per}_*(\mathcal{A})\rightarrow CC_*^{C_p,per}(\mathcal{A}),    
\end{equation}
where $CC^{C_p,per}_*(\mathcal{A})$ is the explicit chain complex (cf. Section 2.3) computing the $C_p$-Tate fixed points of $CC_*(\mathcal{A})$. This would require us to trace the generator $[e_0]\in H^0(C_{-*}(\mathrm{Cyl}(p,1)/\Sigma_p;R)^{tC_p})$ through (the induced map on equivariant homology of) the chain of quasi-equivalences proved in Section 4 and 6: 
\begin{equation}\label{eq:chain of quasi-equivalences from Cyl to KS}
\mathrm{Cyl}(p,1)\simeq FM_{\circlearrowright}(p)\simeq FM_{\circlearrowright}(p)\times \mathcal{P}(p+1)\simeq \widetilde{\mathrm{Cact}^p_{\circlearrowright}}\simeq \mathbf{KS}(p,1).  
\end{equation}
\textbf{Step I}. Under the first quasi-equivalence in \eqref{eq:chain of quasi-equivalences from Cyl to KS}, $[e_0]$ goes to
\begin{equation}\label{eq:equivariant class in FM' circ}
([1,\zeta,\cdots,\zeta^{p-1}],\theta_0,\theta_{\infty}=\theta_0)\in H^0(C_{-*}(FM'_{\circlearrowright}(p)/\Sigma_p;R)^{tC_p}),\;\;\mathrm{where}\;\zeta=e^{-2\pi i/p}.    
\end{equation}
\textbf{Step II}. Under the second quasi-equivalence in \eqref{eq:chain of quasi-equivalences from Cyl to KS}, \eqref{eq:equivariant class in FM' circ} corresponds to the class of any weighted configuration 
\begin{equation}\label{eq:equivariant class in weighted FM' circ}
 [([1,\zeta,\cdots,\zeta^{p-1}],\theta_0,\theta_{\infty}=\theta_0),(a_0,\cdots,a_p)]\in FM'_{\circlearrowright}(p)\times \mathcal{P}(p+1)\,/\Sigma_{p} ,\;\;\mathrm{where}\;a_1=\cdots=a_p.  
\end{equation}
\textbf{Step III}. To identify the image of \eqref{eq:equivariant class in weighted FM' circ} under the third equivalence in \eqref{eq:chain of quasi-equivalences from Cyl to KS}, we look back at the proof of Proposition \autoref{thm:from FM' circ to weighted cyclic cacti, equivariantly} and study the vector field $E_\circlearrowright(z)=a_0\frac{z}{|z|^2}+\sum_{j=1}^p a_j\frac{z-\zeta^{j-1}}{|z-\zeta^{j-1}|^2} $ with its associated function $h(z)=z^{a_0}\prod_{j=1}^p (z-\zeta^{j-1})^{a_j}=z^{a_0}(z^p-1)^a$, if we denote $a=a_1=\cdots=a_p$.\par
\begin{figure}[H]
 \centering
 \includegraphics[width=1.0\textwidth]{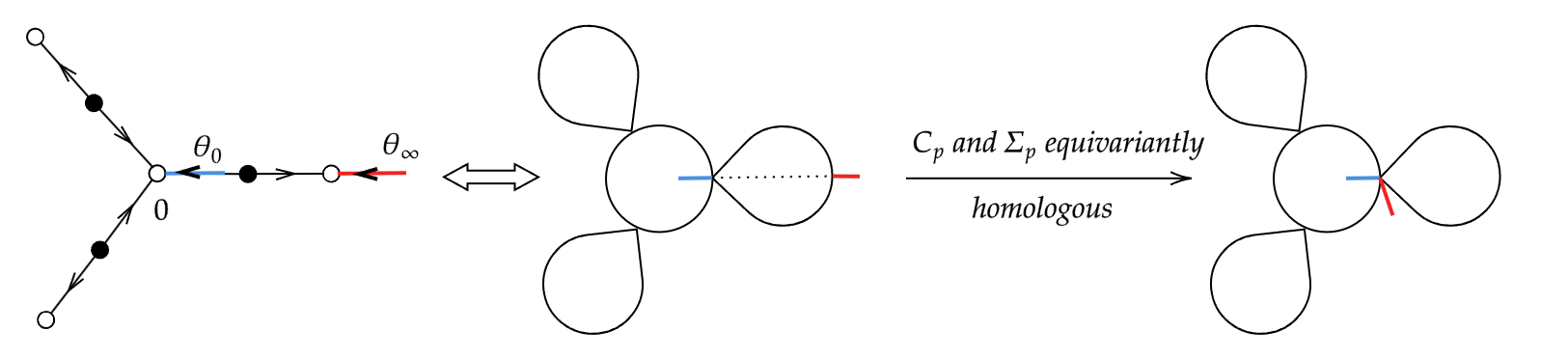}
 \caption{}
 \label{fig:weighted_cyclic_cactus_representing_the_class_e_0}
\end{figure}

\indent
 Then
\begin{equation}
\frac{h'(z)}{h(z)}= \frac{(a_0+ap)z^p-a_0}{z(z^p-1)},   
\end{equation}
which has $p$ distinct roots $p_j=\zeta^j\cdot\sqrt[\leftroot{0}\uproot{7}p]{\frac{a_0}{a_0+ap}}, 1\leq j\leq p$. Moreover, for each $1\leq j\leq p$, the flow lines of $E_{\circlearrowright}$ (which recall are curves on which $\arg h(z)$ is constant) connecting $0$ to $p_j$, $p_j$ to $\zeta^j$, and $\zeta^j$ to $\infty$ are straight lines; see Figure \ref{fig:weighted_cyclic_cactus_representing_the_class_e_0} left. Therefore, by the proof of Proposition \autoref{thm:from FM' circ to weighted cyclic cacti, equivariantly}, the corresponding (degree $0$) equivariant homology class of weighted cyclic cacti looks like Figure \ref{fig:weighted_cyclic_cactus_representing_the_class_e_0} middle, which is furthermore equivariantly homologous to Figure \ref{fig:weighted_cyclic_cactus_representing_the_class_e_0} right. We denote this equivariant homology class by $[\mathcal{E}_0]\in H^0(C_{-*}(\widetilde{\mathrm{Cact}^p_{\circlearrowright}}/\Sigma_p;R)^{tC_p})$.  
\textbf{Step IV}. We are now left with the last quasi-equivalence in \eqref{eq:chain of quasi-equivalences from Cyl to KS}, which is an instance of simplicial subdivision, cf. Theorem \autoref{thm:KS equals weighted cyclic cacti}. As the relevant action here is by finite cyclic subgroup $C_p\subset S^1$ instead of the full $S^1$, recall that is a map of combinatorial categories $j: {}_p\Lambda\rightarrow \Lambda$ which at the level of realization, recovers the sub $C_p\subset S^1$-action (cf. Section 2.4 b)). Consider the composition
\begin{equation}\label{eq:p-cyclic cacti with spine}
\widehat{{}_p\mathfrak{Cact}^p_{\circlearrowright}}:{}_p\Lambda^{op}\times (\Delta^{op})^p\times{}_p\Lambda\xrightarrow{j^{op}\times \mathrm{id}\times j} \Lambda^{op}\times (\Delta^{op})^p\times\Lambda \xrightarrow{\widehat{\mathfrak{Cact}^p_{\circlearrowright}}} \mathrm{Sets}.
\end{equation}
An element of $(\widehat{{}_p\mathfrak{Cact}^p_{\circlearrowright}})_{\bullet;\bullet,\cdots,\bullet;\bullet}$ is called a \emph{$p$-cyclic cactus with spines}. Pictorially, it is drawn like a cyclic cactus with spines except that it has $p$ input basepoints and $p$ output basepoints clockwisely ordered. \par\indent
By the Eilenberg-Zilber theorem, there is a weak equivalence 
\begin{equation}
 \widehat{{}_p\mathrm{Cact}^p_{\circlearrowright}}:=\mathrm{Tot}([k_1,\cdots,k_p]\mapsto |(\widehat{{}_p\mathfrak{Cact}^p_{\circlearrowright}})_{[k_1,\cdots,k_p];\bullet,\cdots,\bullet;\bullet} |)\simeq \mathrm{Tot}([n]\mapsto |(\widehat{\mathfrak{Cact}^p_{\circlearrowright}})_{n;\bullet,\cdots,\bullet;\bullet}|)=\mathbf{KS}(p,1),
\end{equation}
which is $\Sigma_p$-equivariant ($\Sigma_p$ acts by permuting labels of the lobes) and $C_p\times C_p$-equivariant (this follows from the homotopy exactness of \eqref{eq:two homotopy exact sq}). \par\indent
On the other hand, in the homotopy category of $R$-chain complexes, $ (\widehat{{}_p\mathrm{Cact}^p_{\circlearrowright}})_R:=\mathrm{Tot}_R([k_1,\cdots,k_p]\mapsto |(\widehat{{}_p\mathfrak{Cact}^p_{\circlearrowright}})_{[k_1,\cdots,k_p];\bullet,\cdots,\bullet;\bullet} |_R)$ can be explicitly computed as
\begin{equation}\label{eq:explicit chain complex computing Tot of a p-(co)cyclic object}
\prod_{k_1,\cdots,k_p} \bigoplus_{m_1,\cdots,m_p;k_1',\cdots,k_p'}R\langle(\widehat{{}_p\mathfrak{Cact}^p_{\circlearrowright}})_{[k_1,\cdots,k_p];m_1,\cdots,m_p;[k_1',\cdots,k_p']}\rangle [\sum m_i+\sum k_i'-\sum k_i], 
\end{equation}
where $R\langle-\rangle$ denotes the free $R$-module on a set of generators, and $[n]$ denotes degree shift down by $n$; it is equipped with the usual differential given by the alternating sum of multi-coface maps and multi-face maps.\par\indent
There is a degree $0$ cocycle $\mathfrak{e}_0$ of \eqref{eq:explicit chain complex computing Tot of a p-(co)cyclic object} defined as follows. Its entry corresponding to the multi-index $([k_1,\cdots,k_p],m_1,\cdots,m_p,[k_1',\cdots,k_p'])$ is given by:
\begin{itemize} 
    \item if $k_j=k_j'+m_j$ for all $j=1,\cdots,p$, it is the unique $p$-cyclic cactus with spines indicated in Figure \ref{fig:p_cyclic_cactus_with_spines_representing_e_0} (note the labeling of the lobes and the $m_i$'s differ by a cyclic permutation). 
    \begin{figure}[H]
 \centering
 \includegraphics[width=0.9\textwidth]{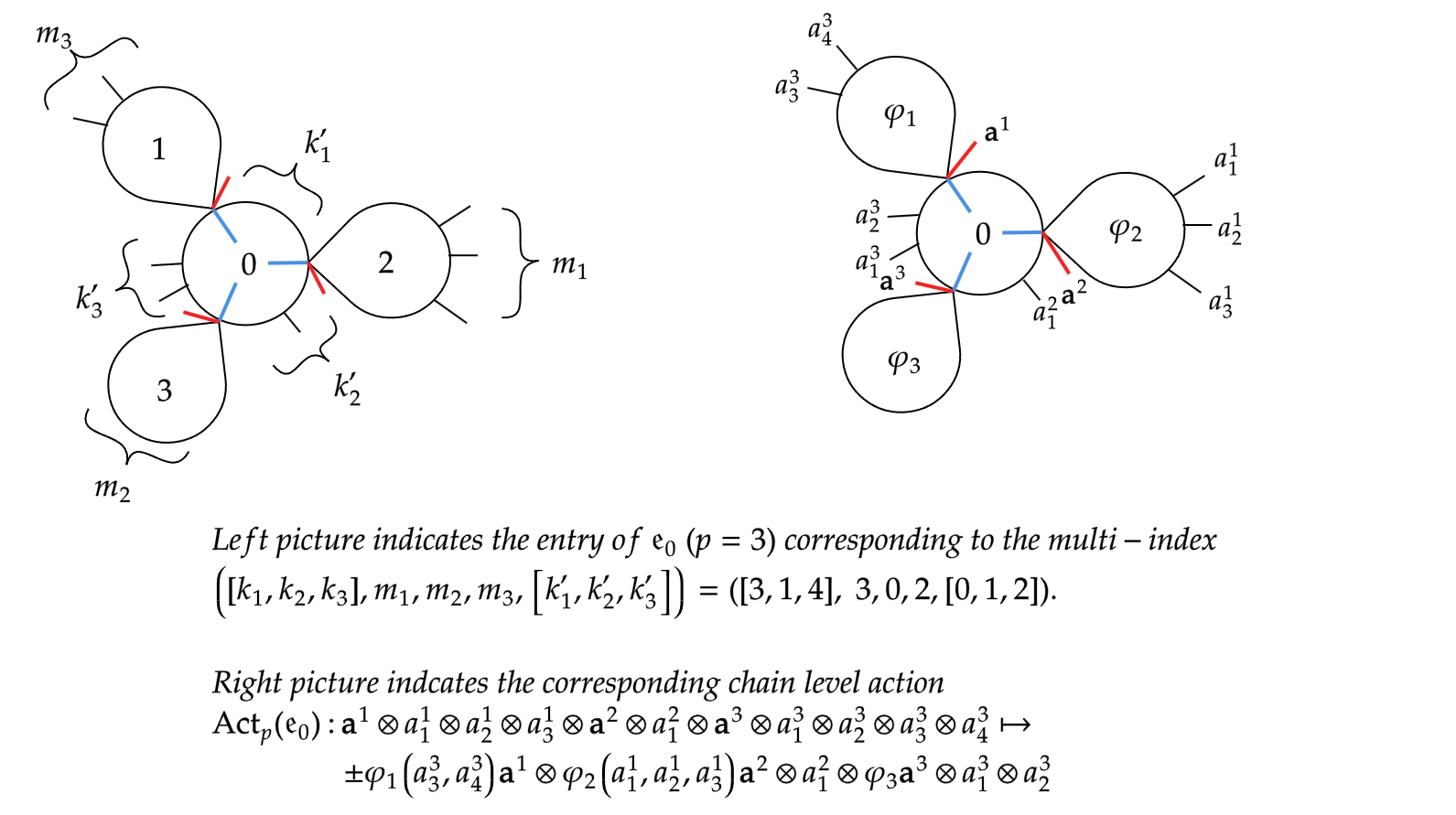}
 \caption{}
 \label{fig:p_cyclic_cactus_with_spines_representing_e_0}
\end{figure}
    To elaborate, the $p$ output basepoints on the base lobe are exactly its intersections with the lobes numbered $1,\cdots,p$; the $p$ input basepoints lie on the $p$ output basepoints and point leftward. The $j$-th lobe ($1\leq j\leq p$) has $m_{j-1}$ marked points, and on the base lobe, there are $k_j'$ marked points in between the $j$-th and $(j+1)$-th output basepoints. Finally, for each marked point that is not a (input or output) basepoint, there is exactly one spine on top of it. In particular, there are $k_j$ spines between the $j$-th and $(j+1)$-th input marked point.
    \item else, it is zero. 
\end{itemize}
One easily checks that $\mathfrak{e}_0$ defines a cocycle (cf. Figure \ref{fig:cactus_with_spine_co_face_and_co_degeneracy} for the definition of the (co)face maps of cacti with spines), and descends to a class in equivariant homology
\begin{equation}
[\mathfrak{e}_0]\in H^0(((\widehat{{}_p\mathrm{Cact}^p_{\circlearrowright}})_R)^{tC_p}_{h\Sigma_p})  .
\end{equation}
Moreover, under the quasi-equivalence $\widetilde{\mathrm{Cact}^p_{\circlearrowright}}\simeq \mathbf{KS}(p,1)$, $[\mathfrak{e}_0]$ exactly corresponds to $[\mathcal{E}_0]$. Hence we conclude that $[\mathfrak{e}_0]$ is the image of $[e_0]$ under the chain of quasi-equivalences in \eqref{eq:chain of quasi-equivalences from Cyl to KS}.  \par\indent
\textbf{The $C_p$-equivariant cap product}. Pull back the action map $\mathrm{Act}:\widehat{\mathfrak{Cact}^p_\circlearrowright}\rightarrow \mathrm{Hom}_R(\mathcal{A}_{\sharp}\otimes(\mathcal{A}^{\sharp})^{\otimes p},\mathcal{A}_{\sharp})$ of Definition \autoref{thm: action of KS(n,1)_R on Hochschild (co)chains} along ${}_p\Lambda^{op}\times (\Delta^{op})^p\times{}_p\Lambda\xrightarrow{j^{op}\times \mathrm{id}\times j} \Lambda^{op}\times (\Delta^{op})^p\times\Lambda$, one obtains 
\begin{equation}\label{eq:p fold action of of KS(n,1) on Hochschild (co)chains}
\mathrm{Act}_p:\widehat{{}_p\mathfrak{Cact}^p_\circlearrowright}\rightarrow \mathrm{Hom}_R(j^*\mathcal{A}_{\sharp}\otimes(\mathcal{A}^{\sharp})^{\otimes p},j^*\mathcal{A}_{\sharp}).    
\end{equation}
Taking realization/totalization induces an action map on equivariant homology (which by an abuse of notation is still denoted $\mathrm{Act}_p$)
\begin{equation}
\mathrm{Act}_p: ((\widehat{{}_p\mathfrak{Cact}^p_\circlearrowright})_R)^{tC_p}_{h\Sigma_p}\rightarrow\mathrm{Hom}_{R^{hC_p}}(CC_*^{C_p,per}(\mathcal{A})\otimes (CC^*(\mathcal{A})^{\otimes p})^{h\Sigma_p},CC_*^{C_p,per}(\mathcal{A}))
\end{equation}
and hence a map
\begin{equation}\label{eq:p fold action of e_0}
 \mathrm{Act}_p([\mathfrak{e}_0]):HH_*^{C_p,per}(\mathcal{A})\otimes H^*_{\Sigma_p}(CC^*(\mathcal{A})^{\otimes p})\longrightarrow HH_*^{C_p,per}(\mathcal{A}).   
\end{equation}
Finally, from the definition of $\mathrm{Act}_p$, one sees that $\mathrm{Act}_p([\mathfrak{e}_0])$ factors through the restriction (the restriction $H^*_{\Sigma_p}(CC^*(\mathcal{A})^{\otimes p})\hookrightarrow H^*_{C_p}(CC^*(\mathcal{A})^{\otimes p})$ is injective because $[\Sigma_p:C_p]$ is coprime to $p$)
\begin{equation}\label{eq:p fold action of e0 restricted from C_p to Sigma_p}
\begin{tikzcd}[row sep=1.2cm, column sep=0.8cm]
HH_*^{C_p,per}(\mathcal{A})\otimes H^*_{\Sigma_p}(CC^*(\mathcal{A})^{\otimes p})\arrow[rr,"{\mathrm{Act}_p([\mathfrak{e}_0])}"] & &HH_*^{C_p,per}(\mathcal{A}) \\
HH_*^{C_p,per}(\mathcal{A})\otimes H^*_{C_p}(CC^*(\mathcal{A})^{\otimes p}) \arrow[u,hookleftarrow,"{\mathrm{id}\otimes\mathrm{Res}_{C_p\subset\Sigma_p}}"]\arrow[urr,"{\prod^{C_p}}"]& &
\end{tikzcd}
\end{equation}
where $\prod^{C_p}$ is the explicit map defined below; compare Figure \ref{fig:p_cyclic_cactus_with_spines_representing_e_0} right.\\
\begin{mydef}(\cite{Che2})\label{thm:Cp equivariant cap product}
Given $\varphi_1\otimes\cdots\otimes \varphi_p\in CC^*(\mathcal{A})^{\otimes p}$, consider the formula 
\begin{equation*}
\mathbf{a}^1\otimes a_1^1\otimes \cdots\otimes a^1_{k_1}\otimes\mathbf{a}^2\otimes a_1^2\otimes \cdots\otimes a^2_{k_2}\otimes \cdots\otimes \mathbf{a}^p\otimes a_1^p\otimes \cdots\otimes a^p_{k_p}\mapsto
\end{equation*}
\begin{equation}\label{eq:p fold cap product}
\sum_{j_1,\cdots, j_p} (-1)^{\dagger}\varphi_1(a^p_{j_p+1}, \cdots,a^p_{k_p})\mathbf{a}^1\otimes a^1_1\otimes\cdots\otimes a^1_{j_1}\otimes\varphi_2(a^1_{j_1+1}, \cdots,a^1_{k_1})\mathbf{a}^2\otimes a^2_1\otimes\cdots \otimes\varphi_p(a^{p-1}_{j_{p-1}+1}, \cdots,a^{p-1}_{k_{p-1}})\mathbf{a}^p\otimes a^p_1\otimes\cdots\otimes a^p_{j_p},
\end{equation}
where 
$$\dagger=\sum_{i=j_p+1}^{k_p}\|a^p_{j_i}\|\cdot\big(\sum_{i=1}^{p-1}
(|\mathbf{a}^i|+\|a^i_1\|+\cdots+\|a^i_{k_i}\|)+|\mathbf{a}^p|+\|a^p_1\|+\cdots+\|a^p_{j_p}\|\big)$$
\begin{equation}
+\sum_{i=1}^p \|\varphi_i\|\cdot \big(\|a^{i-1}_{j_{i-1}+1}\|+\cdots+\|a^{i-1}_{k_{i-1}}\|+|\mathbf{a}^i|+\|a^i_1\|+\cdots+|\mathbf{a}^p|+\|a^p_1\|+\cdots+\|a^p_{j_p}\|\big)
\end{equation} 
gives the Koszul sign. \eqref{eq:p fold cap product} defines a chain map
\begin{equation}
{}_p\prod: CC^*(\mathcal{A})^{\otimes p}\otimes {}_pCC_*(\mathcal{A})\rightarrow {}_pCC_*(\mathcal{A}),
\end{equation}
which is called the \emph{$p$-fold cap product}. Moreover, it is $C_p$-equivariant (with respect to permuting the factors of $CC^*(\mathcal{A})^{\otimes p}$ and the action in \eqref{eq:tau}))
and hence descends to a map on $C_p$-homotopy fixed points
\begin{equation}
-\prod\nolimits^{C_p}-: H^*_{C_p}(CC^*(\mathcal{A})^{\otimes p})\otimes HH_*^{C_p,-}(\mathcal{A})\rightarrow HH_*^{C_p,-}(\mathcal{A}),  
\end{equation}
which we call \emph{the $p$-fold equivariant cap product} (one may also invert $t$ and replace $HH^{C_p,-}$ by $HH^{C_p,per}$). By precomposing with the Frobenius-linear map $HH^*(\mathcal{A})\rightarrow H^*_{C_p}(CC^*(\mathcal{A})^{\otimes p})$ given by $[\varphi]\mapsto [\varphi^{\otimes p}]$, one obtains another action denoted
\begin{equation}\label{eq:Cp equivariant cap product} 
\bigcap\nolimits^{C_p}:=HH^*(\mathcal{A})\otimes HH_*^{C_p,-}(\mathcal{A})\xrightarrow{([\varphi]\mapsto [\varphi^{\otimes p}])\otimes \mathrm{id}}    H^*_{C_p}(CC^*(\mathcal{A})^{\otimes p})\otimes HH_*^{C_p,-}(\mathcal{A})\xrightarrow{\prod^{C_p}} HH_*^{C_p,-}(\mathcal{A}),
\end{equation}
which we call the \emph{$C_p$-equivariant cap product}. 
\end{mydef}
Combining diagram \eqref{eq:p fold action of e0 restricted from C_p to Sigma_p} with the discussion preceding it, one obtains the following. \\
\begin{prop}\label{thm:e0 operation is p fold equivariant cap product}
Given $[\phi]\in H_*^{\Sigma_p}(CC^*(\mathcal{A})^{\otimes p})$, let $[\phi]_{C_p}\in H_*^{C_p}(CC^*(\mathcal{A})^{\otimes p})$ be its restriction to $C_p$ homotopy fixed points. Then
\begin{equation}\label{eq:e0 operation is Cp equivariant cap product}
\Xi^p([e_0])([\phi],-)=[\phi]_{C_p}\prod\nolimits^{C_p}-.
\end{equation}\qed
\end{prop}
We recall some basic properties of $\bigcap^{C_p}$ proved in \cite{Che2}:
\begin{enumerate}[label=C\arabic*)]
    \item  It is a graded multiplicative action: $(\phi\cup\varphi)\bigcap^{C_p} a=(-1)^{|\phi||\varphi|\frac{p(p-1)}{2}}\phi\bigcap^{C_p}(\varphi\bigcap^{C_p}a)$, where $\cup$ denotes the cup product on Hochschild cohomology. 
    \item It is additive in the second variable, and becomes additive in the first variable after multiplying by $t$ (the formal $S^1$-equivariant variable of degree $2$), the latter meaning $t(\phi+\varphi)\bigcap^{C_p}=t\phi\bigcap^{C_p}+t\varphi\bigcap^{C_p}$.
    \item $\bigcap^{C_p}$ is unital, i.e. $e\bigcap^{C_p}=\mathrm{id}$, where $e$ is the unit in the Hochschild cohomology algebra.
    \item It is Frobenius $p$-linear: for $a\in R$, $a\varphi\bigcap^{C_p}=a^p\varphi\bigcap^{C_p}$.
\end{enumerate}
The definition of the $C_p$-equivariant cap product easily extends to the case where $\mathcal{A}$ is an $A_{\infty}$-category, and it has several geometric incarnations.
\begin{itemize}
    \item When $\mathcal{A}$ is the Fukaya category of a closed monotone symplectic manifold, then $\bigcap^{C_p}$ corresponds to the \emph{quantum Steenrod operations} in equivariant Gromov-Witten theory, cf. \cite[Theorem 1.2]{Che2}. 
    \item When $\mathcal{A}$ is the dg category of matrix factorizations of an isolated singularity, $\bigcap^{C_p}$ is expected to correspond to a $p$-th power action of twisted functions on twisted differential forms, cf. \cite[Conjecture 6.1]{Che2}. 
\end{itemize}
We now state an analogue of Proposition \autoref{thm:e0 operation is p fold equivariant cap product} with $\bigcap^{C_p}$ in place of $\prod^{C_p}$. Before that, we recall the following basic computation. Suppose $p$ is an odd prime. By a standard computation (cf. e.g. \cite{Coh}), the diagonal class $\Delta_{[S^1]}\in H^1_{\Sigma_p}(\tilde{C}_*((S^1)^p;R))$ is sent to the element $t^{\frac{p-1}{2}}$ under 
\begin{equation}
H^1_{\Sigma_p}(\tilde{C}_*((S^1)^p;R))\cong H^{p-1}_{\Sigma_p}(R(p)) \xrightarrow{\mathrm{Res}_{C_p\subset\Sigma_p}} H^{p-1}_{C_p}(R)=R[[t,\theta]],   
\end{equation}
where $R(p)$ is the sign representation of $\Sigma_p$. Thus, in view of Example \autoref{thm:example of twisted e0 operation} and Proposition \autoref{thm:e0 operation is p fold equivariant cap product}, there is a commutative diagram
\begin{equation}\label{eq:e0 pm action big diagram}
\begin{tikzcd}[row sep=1.2cm, column sep=0.8cm]
HH_*^{C_p,per}(\mathcal{A})\otimes H^*_{\Sigma_p}(CC^*(\mathcal{A})^{\otimes p}\otimes R(p))\arrow[rrr,"{\Xi^{\pm,p}([e_0^\pm])}"] \arrow[dr,"{\mathrm{id}\otimes (-\cup \Delta_{[S^1]})}"] & & &HH_*^{C_p,per}(\mathcal{A})\arrow[dd,equal] \\
&HH_*^{C_p,per}(\mathcal{A})\otimes H^*_{\Sigma_p}(CC^*(\mathcal{A})^{\otimes p}) \arrow[urr,"{\Xi^{p}([e_0])}"]& &  \\
HH_*^{C_p,per}(\mathcal{A})\otimes H^*_{C_p}(CC^*(\mathcal{A})^{\otimes p}) \arrow[dr,"{\mathrm{id}\otimes t^{\frac{p-1}{2}}}"] \arrow[uu,hookleftarrow,"{\mathrm{id}\otimes\mathrm{Res}_{C_p\subset\Sigma_p}}"]& & & HH_*^{C_p,per}(\mathcal{A})\\
 &HH_*^{C_p,per}(\mathcal{A})\otimes H^*_{C_p}(CC^*(\mathcal{A})^{\otimes p})\arrow[uu,hookleftarrow,"{\mathrm{id}\otimes\mathrm{Res}_{C_p\subset\Sigma_p}}"]\arrow[urr,"{\prod^{C_p}}"] & &
\end{tikzcd}
\end{equation}\\
\begin{cor}\label{thm:Cp equivariant cap product are KS operations}
Fix $[\varphi]\in HH^l(\mathcal{A})$. If $l$ or $p$ is even, then $[\varphi^{\otimes p}]\in \mathrm{im}\big( H^*_{\Sigma_p}(CC^*(\mathcal{A})^{\otimes p})\subset H^*_{C_p}(CC^*(\mathcal{A})^{\otimes p})\big)$ and
\begin{equation}
[\varphi]\bigcap\nolimits^{C_p}-= \Xi^p([e_0])([\varphi^{\otimes p}],-).
\end{equation}
If $l$ and $p$ are both odd, then $[\varphi^{\otimes p}]\in \mathrm{im}\big( H^*_{\Sigma_p}(CC^*(\mathcal{A})^{\otimes p}\otimes R(p))\subset H^*_{C_p}(CC^*(\mathcal{A})^{\otimes p})\big)$ and
\begin{equation}
 t^{\frac{p-1}{2}}\cdot ([\varphi]\bigcap\nolimits^{C_p}-)=\Xi^p([e_0])([\varphi^{\otimes p}]\cup \Delta_{[S^1]},-).
\end{equation}
In particular, the $C_p$-equivariant cap product is a $C_p$-Kontsevich-Soibelman operation.
\end{cor}
\emph{Proof}. The case where $l$ or $p$ is even follows from Proposition \autoref{thm:e0 operation is p fold equivariant cap product}, and the case where $l$ and $p$ are odd follows from diagram \eqref{eq:e0 pm action big diagram} .\qed \\
Return to the relative setting $\mathcal{A}/R/\mathbf{k}$ as in Theorem \ref{thm:automatic covariant constancy}, where $\mathbf{k}$ is a field of characteristic $p>0$. \\
\begin{cor}\label{thm:covariant constancy of Cp equivariant cap product}
For any $[\varphi]\in HH^*(\mathcal{A})$, $[\varphi]\bigcap^{C_p}-$ is covariantly constant with respect to the ($C_p$-) Getzler-Gauss-Manin connection. 
\end{cor}
\emph{Proof.} In the context of Corollary \autoref{thm:Cp equivariant cap product are KS operations}, note that moreover $[\varphi^{\otimes p}]\in H^*_{\Sigma_p}(CC^*(\mathcal{A})^{\otimes_R p})$ lies in the image of $H^*_{\Sigma_p}(CC^*(\mathcal{A})^{\otimes_{\mathbf{k}} p})$; hence $[\varphi]\cap^{C_p}-$ is in fact a $\mathbf{k}$-linear Kontsevich-Soibelman operation. The statement then follows from (a $C_p$-counterpart of) Theorem \autoref{thm:automatic covariant constancy}.\qed\par\indent
Given the relationship between $C_p$-equivariant cap products (in the context of Fukaya categories) and quantum Steenrod operations proved in \cite{Che2}, Corollary \ref{thm:covariant constancy of Cp equivariant cap product} can be viewed as a non-commutative generalization of Seidel-Wilkins' covariant constancy of quantum Steenrod operations \cite{SW}.

\begin{appendices}
\renewcommand{\theequation}{A.\arabic{equation}}
\setcounter{equation}{0}
\section{Homotopy fixed points of circle actions}
In this appendix, we recall some basic facts about homotopy fixed points of chain complexes with circle and finite cyclic group actions. We work over a ring $\mathbf{k}$ of odd characteristic $p$ throughout, and all results work the same for either $\mathbb{Z}$ or $\mathbb{Z}/2$-graded chain complexes. 
\subsection{$S^1$-chain complexes}
Let $C_{-*}(S^1;\mathbf{k})$ be the dg algebra of singular chains on $S^1$ with coefficients in $\mathbf{k}$, with grading negated. An \emph{$S^1$-chain complex} is a dg module over $C_*(S^1;\mathbf{k})$. \par\indent
Given an $S^1$-chain complex $X$, the \emph{$S^1$-homotopy fixed point of $X$} is defined to be
\begin{equation}
X^{hS^1}:=\mathrm{RHom}_{C_{-*}(S^1;\mathbf{k})}(\mathbf{k}, X),   
\end{equation}
where $\mathbf{k}$ is equipped with the trivial $S^1$-action. \par\indent
To obtain an explicit chain model that computes $X^{hS^1}$, we note that by the homological perturbation lemma, there is an $A_{\infty}$-quasi-equivalence
\begin{equation}
C_{-*}(S^1;\mathbf{k})\simeq \mathbf{k}[\Lambda]/\Lambda^2,    
\end{equation}
where $\mathbf{k}[\Lambda]/\Lambda^2$ is viewed as a dg algebra with trivial differential and $|\Lambda|=-1$. We think of $\mathbf{k}[\Lambda]/\Lambda^2$ as coming from a cellular model for $S^1$ with one $0$-cell and one $1$-cell. \par\indent
Hence by transfer, one may equivalently view an $S^1$-chain complex $X$ as a dg module over $\mathbf{k}[\Lambda]/\Lambda^2$ (which by an abuse of notation still denoted $X$), and in particular
\begin{equation}\label{eq:S^1 fixed points of X}
X^{hS^1}\simeq \mathrm{RHom}_{\mathbf{k}[\Lambda]/\Lambda^2}(\mathbf{k},X).    
\end{equation}
On the other hand, there is a cofibrant resolution of the trivial $\mathbf{k}[\Lambda]/\Lambda^2$-module given by
\begin{equation}
\mathrm{Tot}(\cdots\xrightarrow{\Lambda} \mathbf{k}[\Lambda]/\Lambda^2\xrightarrow{\Lambda} \mathbf{k}[\Lambda]/\Lambda^2) \rightarrow\mathbf{k} 
\end{equation}
and therefore one can compute \eqref{eq:S^1 fixed points of X} by the chain complex
\begin{equation}\label{eq:explicit complex for S^1 fixed points}
\mathrm{Hom}_{\mathbf{k}[\Lambda]/\Lambda^2}(\mathrm{Tot}(\cdots\xrightarrow{\Lambda} \mathbf{k}[\Lambda]/\Lambda^2\xrightarrow{\Lambda} \mathbf{k}[\Lambda]/\Lambda^2) ,X) \cong  (X[[t]],d_X+t\Lambda)
\end{equation}
where $t$ is a formal variable of degree $2$ corresponding to the generator of $H^2(BS^1;\mathbf{k})$. Under this explicit chain model, the \emph{$S^1$-Tate fixed points} of $X$ can be computed as 
\begin{equation}
X^{tS^1}\simeq     (X((t)),d_X+t\Lambda).
\end{equation}

\subsection{$C_p$-chain complexes}
Let $C_p$ be the finite cyclic group of order $p$ (equal to the characteristic of $\mathbf{k}$). A \emph{$C_p$-chain complex} is a dg module over the group algebra $\mathbf{k}[C_p]$. \par\indent
Given a $C_p$-chain complex $X$, the \emph{$C_p$-homotopy fixed point of $X$} is defined to be
\begin{equation}\label{eq:C_p homotopy fixed point of X}
X^{hC_p}:=\mathrm{RHom}_{\mathbf{k}[C_p]}(\mathbf{k},X),    
\end{equation}
where $\mathbf{k}$ is equipped with the trivial $C_p$-action. \par\indent
To obtain an explicit chain model for $X^{hC_p}$, note that there is a cofibrant resolution of the trivial $\mathbf{k}[C_p]$-module given by
\begin{equation}
\mathrm{Tot}(\cdots\xrightarrow{\tau-1} \mathbf{k}[C_p]\xrightarrow{1+\tau+\cdots+\tau^{p-1}}\mathbf{k}[C_p]\xrightarrow{\tau-1} \mathbf{k}[C_p])\rightarrow \mathbf{k}.
\end{equation}
Therefore, 
\begin{equation}\label{eq:explicit complex for C_p fixed points}
X^{hC_p}\simeq \mathrm{Hom}_{\mathbf{k}[C_p]}( \mathrm{Tot}(\cdots\xrightarrow{\tau-1} \mathbf{k}[C_p]\xrightarrow{1+\tau+\cdots+\tau^{p-1}}\mathbf{k}[C_p]\xrightarrow{\tau-1} \mathbf{k}[C_p]),X)\cong(X[[t,\theta]],d_{eq}) ,  
\end{equation}
where $t,\theta$ are formal variables corresponding to the generators of $H_2(BC_p;\mathbf{k})$ and $H_1(BC_p;\mathbf{k})$, satisfying $|t|=2,|\theta|=1, t\theta=\theta t, \theta^2=0$, and the differential is given by
\begin{equation}
\begin{cases}
d_{eq}(x)=d_X(x)+(-1)^{|x|}(\tau x-x)\theta,\\
d_{eq}(x\theta)=d_X(x)\,\theta+(-1)^{|x|}(x+\tau x+\cdots+\tau^{p-1}x)t.     
\end{cases}      
\end{equation}
On $X^{hC_p}$, there is a `corrected action' by $\theta$ (cf. \cite[(2.38) and (2.39)]{SW}) given by
\begin{equation}\label{eq:corrected theta action}
\begin{cases}
xt^k\mapsto (-1)^{|x|}xt^k\theta \\
xt^k\theta\mapsto (-1)^{|x|+1}\tau(\tau-1)^{p-2}xt^{k+1}.
\end{cases}    
\end{equation}
One can check, using $\mathrm{char}(\mathbf{k})=p$, that \eqref{eq:corrected theta action} is a chain map whose square is nullhomotopic. This, combined with multiplication by $t$, descends to a $\mathbf{k}[[t,\theta]]$-module structure on $H^*(X^{hC_p})$. Finally, we remark that the \emph{$C_p$-Tate fixed points} of $X$ can be explicitly expressed as 
\begin{equation}
X^{tC_p}\simeq  (X((t,\theta)),d_{eq})
\end{equation}
with the same formula for the differential as $X^{hC_p}$.

\subsection{Comparing $S^1$ and $C_p$ homotopy fixed points}
Let $X$ be an $S^1$-chain complex. The goal of this subsection is to address the following simple question: how are $X^{hS^1}$ and $X^{hC_p}$ related, where the $C_p$-action on $X$ is induced from the inclusion $C_p\subset S^1$ as the $p$-th roots of unity?    \par\indent
To approach this problem, we use a suitable chain model for $C_{-*}(S^1)$ that is both computable and also fine enough to `see' the subgroup $C_p$ (for instance, the chain model $k[\Lambda]/\Lambda^2$ is too coarse for this purpose). Consider the following dg algebra (cf. \cite[Definition 2.3.1]{Sen}):
\begin{equation}
\mathbf{k}[\tau,\sigma]:=  \mathbf{k}[\tau,\sigma]/(\tau^p=1, d\sigma=\tau-1, \sigma^2=0)\;,\;|\tau|=0\;,\;|\sigma|=-1.
\end{equation}
There is a quasi-isomorphism of dg algebras (we write shorthand $\mathbf{k}[\Lambda]:=\mathbf{k}[\Lambda]/\Lambda^2$) given by
\begin{equation}
C_{-*}(S^1;\mathbf{k})\simeq \mathbf{k}[\Lambda]\xrightarrow[\simeq]{\Lambda\mapsto(1+\tau+\cdots+\tau^{p-1})\sigma}  \mathbf{k}[\tau,\sigma].  
\end{equation}
In particular, we can think of $\mathbf{k}[\tau,\sigma]$ as a chain model for $C_{-*}(S^1;\mathbf{k})$ coming from the cellular structure on $S^1$ with $p$ $0$-cells given by $1,\tau,\cdots,\tau^{p-1}$ and $p$ $1$-cells given by $\sigma,\tau\sigma,\cdots,\tau^{p-1}\sigma$. Indeed, there is an evident inclusion
\begin{equation}\label{eq:restriction from S^1 to C_p}
\mathbf{k}[C_p]\cong \mathbf{k}[\tau]\subset    \mathbf{k}[\tau,\sigma], 
\end{equation}
which corresponds to the standard inclusion $C_p\subset S^1$ of topological spaces. \par\indent
By transfer, we view an $S^1$-chain complex $X$ as a module over $\mathbf{k}[\tau,\sigma]$; by \eqref{eq:restriction from S^1 to C_p}, we may also view $X$ as a $\mathbf{k}[\tau]$-module. The question at the beginning of this subsection can then be reformulated as: what is the relation between $\mathrm{RHom}_{\mathbf{k}[\tau,\sigma]}(\mathbf{k},X)$ and $\mathrm{RHom}_{\mathbf{k}[\tau]}(\mathbf{k},X)$?
\par\indent
There is a natural restriction map $X^{hS^1}\rightarrow X^{hC_p}$ given by
\begin{equation}\label{eq:restriction}
 \mathrm{Res}_p: \mathrm{RHom}_{\mathbf{k}[\tau,\sigma]}(\mathbf{k},X)\cong \mathrm{Hom}_{\mathbf{k}[\tau,\sigma]}(T,X)\xrightarrow{\mathrm{for}}\mathrm{Hom}_{\mathbf{k}[\tau]}(T,X)\xrightarrow{(*)} \mathrm{Hom}_{\mathbf{k}[\tau]}(T',X)\simeq X^{hC_p},
\end{equation}
where $T$ (resp. $T'$) is a cofibrant resolution of $\mathbf{k}$ as a $\mathbf{k}[\tau,\sigma]$ (resp. $\mathbf{k}[\tau]$) module, `for' denotes the forgetful map, and $(*)$ denotes the map induced from a map of resolutions $T'\rightarrow T$ (which is unique up to homotopy).\\
\begin{lemma}\label{thm:C_p Gysin comparison}
$\mathrm{Res}_p$ is $\mathbf{k}[[t]]\simeq H^*(BS^1;\mathbf{k})$-linear, and its $(t,\theta)$-linear extension
\begin{equation}\label{eq:Gysin comparison}
\mathrm{Res}_p: X^{hS^1}\oplus X^{hS^1}\theta\rightarrow X^{hC_p}    
\end{equation}
is a quasi-isomorphism. 
\end{lemma}
\emph{Proof}. Recall from Appendix A.1 and A.2 that there are cofibrant resolutions
$$T= \mathrm{Tot}(\cdots \xrightarrow{(1+\tau+\cdots+\tau^{p-1})\sigma}\mathbf{k}[\tau,\sigma]\xrightarrow{(1+\tau+\cdots+\tau^{p-1})\sigma}\mathbf{k}[\tau,\sigma])$$
\begin{equation}
=\big(\mathbf{k}[\tau,\sigma][u,u^{-1}]/\sum_{k<0}\mathbf{k}[\tau,\sigma]u^k,d_T=d_{\mathbf{k}[\tau,\sigma]}+(1+\tau+\cdots+\tau^{p-1})\sigma u^{-1}\big),
\end{equation}
and
$$T'=\mathrm{Tot}(\cdots\xrightarrow{\tau-1} \mathbf{k}[\tau]\xrightarrow{1+\tau+\cdots+\tau^{p-1}}\mathbf{k}[\tau]\xrightarrow{\tau-1} \mathbf{k}[\tau])$$
\begin{equation}
=\big(\mathbf{k}[\tau][u,u^{-1},\epsilon]/(\sum_{k<0}\mathbf{k}[\tau]u^k+\mathbf{k}[\tau]u^k\epsilon), d_{T'}=\begin{cases}
u^{k}\mapsto (1+\tau+\cdots+\tau^{p-1})u^{k-1}\epsilon \\
u^{k}\epsilon\mapsto (\tau-1)u^{k}.
\end{cases}\big)
\end{equation}
of $\mathbf{k}$ as $\mathbf{k}[\tau,\sigma]$ (resp. $\mathbf{k}[\tau]$ modules), where $|u|=-2, |\epsilon|=-1, u\epsilon=\epsilon u, \epsilon^2=0$ are formal variables labeling the columns of the resolutions. There is a map from $T'$ to $T$ as $\mathbf{k}[\tau]$-resolutions of $\mathbf{k}$ given by
\begin{equation}\label{eq:comparison of k[tau] resolutions}
\begin{cases}
u^{k}\mapsto u^{k}\\
u^{k}\epsilon\mapsto \sigma u^{k}.
\end{cases}   
\end{equation}
Therefore, the composition $\mathrm{Res}_p:\mathrm{Hom}_{\mathbf{k}[\tau,\sigma]}(T,X)\xrightarrow{\mathrm{for}}\mathrm{Hom}_{\mathbf{k}[\tau]}(T,X)\xrightarrow{(*)} \mathrm{Hom}_{\mathbf{k}[\tau]}(T',X)$ in \eqref{eq:restriction} is given at the chain-level by
\begin{equation}\label{eq:Gysin comparison chain map}
\mathrm{Res}_p: xt^k\mapsto xt^k-(-1)^{|x|}\sigma xt^k\theta.  
\end{equation}
\eqref{eq:Gysin comparison chain map} is clearly $t$-linear. Moreover, the $\theta$-linear extension of $\mathrm{Res}_p$ satisfies the property that ${\mathrm{Res}_p}|_{t=\theta=0}=\mathrm{id}_X$. Since the $(t,\theta)$-filtration is complete, this shows that \eqref{eq:Gysin comparison} is a quasi-isomorphism. \qed\\
\begin{cor}\label{thm:equivalence of S^1 vs C_p fixed points}
Let $f:X\rightarrow Y$ be a map of $S^1$-chain complexes, then
\begin{equation}
f^{hS^1}:X^{hS^1}\rightarrow Y^{hS^1}    
\end{equation}
is a quasi-isomorphism if and only if 
\begin{equation}
f^{hC_p}:X^{hC_p}\rightarrow Y^{hC_p}    
\end{equation}
is a quasi-isomorphism.
\end{cor}
\emph{Proof}. This follows from the splitting of Lemma \ref{thm:C_p Gysin comparison} and the fact that $\mathrm{Res}_p$ is clearly natural with respect to $X$.\qed\par\indent
If $X$ is an $S^1$-chain complex, then $X^{hC_p}$ inherits a residual $S^1\simeq S^1/C_p$-action, which we think of as a `$p$-fold' circle action. This can be seen using the above explicit models following \cite[Lemma 2.3.4]{Sen}. Namely, the quasi-isomorphism $\mathbf{k}[\Lambda]\simeq \mathbf{k}[\tau,\sigma]$ induces a quasi-isomorphism of $C_p$-chain complexes
\begin{equation}\label{eq:rephrasing of induced C_p action}
X\simeq \mathrm{RHom}_{\mathbf{k}[\Lambda]}(\mathbf{k}[\tau,\sigma],X),    
\end{equation}
where the $C_p$-action on the right hand side is induced from the inclusion $\mathbf{k}[C_p]=\mathbf{k}[\tau]\subset \mathbf{k}[\tau,\sigma]$. \eqref{eq:rephrasing of induced C_p action} and tensor-hom adjunction gives rise to a chain of quasi-isomorphisms  
\begin{equation}\label{eq:rephrasing of X^hC_p}
X^{hC_p}\simeq\mathrm{RHom}_{\mathbf{k}[\tau]}(\mathbf{k}, \mathrm{RHom}_{\mathbf{k}[\Lambda]}(\mathbf{k}[\tau,\sigma],X))\simeq \mathrm{RHom}_{\mathbf{k}[\Lambda]}(\mathbf{k}[\tau,\sigma]\otimes^\mathbb{L}_{\mathbf{k}[\tau]}\mathbf{k},X)\simeq \mathrm{RHom}_{\mathbf{k}[\Lambda]}(\mathbf{k}[\Lambda^{(1)}],X),   
\end{equation}
where $\mathbf{k}[\Lambda^{(1)}]:=\mathbf{k}[\Lambda^{(1)}]/(\Lambda^{(1)})^2$ is just a copy of $\mathbf{k}[\Lambda]$, but equipped with the (unital) $\mathbf{k}[\Lambda]$-module structure where $\Lambda\cdot 1=p\Lambda^{(1)}=0$, corresponding to a $p$-fold circle action. From \eqref{eq:rephrasing of X^hC_p} there is an evident $C_{-*}(S^1;\mathbf{k})\simeq \mathbf{k}[\Lambda^{(1)}]$-action on $X^{hC_p}$ induced by multiplication by $\Lambda^{(1)}$ on the source of $\mathrm{RHom}_{\mathbf{k}[\Lambda]}(\mathbf{k}[\Lambda^{(1)}],X)$. From this description, it is also clear that the action of $\Lambda^{(1)}$ is $t$-linear (but not $\theta$-linear in general). \par\indent
Since the action of $\mathbf{k}[\Lambda]$ on $\mathbf{k}[\Lambda^{(1)}]$ is trivial, we may rewrite \eqref{eq:rephrasing of X^hC_p} as
$$X^{hC_p}\simeq \mathrm{RHom}_{\mathbf{k}[\Lambda]}(\mathbf{k}[\Lambda^{(1)}],X)=\mathrm{RHom}_{\mathbf{k}[\Lambda]}(\mathbf{k},X)\oplus \mathrm{RHom}_{\mathbf{k}[\Lambda]}(\mathbf{k}\Lambda^{(1)},X)$$
\begin{equation}\label{eq:second C_p Gysin interpretation}
=\mathrm{RHom}_{\mathbf{k}[\Lambda]}(\mathbf{k},X)\oplus \mathrm{RHom}_{\mathbf{k}[\Lambda]}(\mathbf{k},X)\theta\simeq X^{hS^1}\oplus X^{hS^1}\theta.
\end{equation}
Tracing through the quasi-isomorphisms in \eqref{eq:rephrasing of X^hC_p}, it is not hard to see that \eqref{eq:second C_p Gysin interpretation} recovers the splitting of Lemma \ref{thm:C_p Gysin comparison}. In particular,\\
\begin{cor}\label{thm:residual circle action on C_p fixed point via Gysin}
Under the (homological level) identification 
\begin{equation}
\mathrm{Res}_p: H^*(X^{hS^1})\oplus H^*(X^{hS^1})\theta\cong H^*(X^{hC_p})    
\end{equation}
of Lemma \ref{thm:C_p Gysin comparison}, the $H_{-*}(S^1/C_p;\mathbf{k})\simeq \mathbf{k}[\Lambda^{(1)}]$ action on $H^*(X^{hC_p})$ can be identified with the $\mathbf{k}[\Lambda^{(1)}]$-action on $H^*(X^{hS^1})\oplus H^*(X^{hS^1})\theta$ where $\Lambda^{(1)}$ acts by
\begin{equation}
([\alpha],[\beta]\theta)\mapsto ([\beta],0). 
\end{equation}\qed
\end{cor}

\end{appendices}

\end{document}